\documentclass[11pt,a4paper,twoside,final]{scrartcl}
\usepackage{a4wide}
\usepackage{amsfonts}
\usepackage{amsmath}
\usepackage{amssymb}
\usepackage{amsthm}
\usepackage{algorithm}
\usepackage{algorithmic}
\usepackage{pgfplots}
\pgfplotsset{compat=newest}
\usepackage{graphicx}
\usepackage{color}
\usepackage{colortbl}
\usepackage{showkeys}
\usepackage{enumerate}
\usepackage{subcaption}
\usepackage{verbatim}
\usepackage{bm}
\usepackage{hyperref}
\usepackage{tikz}
\usepackage{tikz-cd}
\pgfplotsset{colormap={whitered}{color(0cm)=(white);
	color(1cm)=(orange!75!red)}
}

\DeclareMathOperator{\esssup}{ess\,sup}

\newtheorem{theorem}{Theorem}[section]
\newtheorem{lemma}[theorem]{Lemma}
\newtheorem{remark}[theorem]{Remark}
\newtheorem{generalisation}[theorem]{Generalisation}
\newtheorem{definition}[theorem]{Definition}
\newtheorem{example}[theorem]{Example}
\newtheorem{corollary}[theorem]{Corollary}
\newtheorem{proposition}[theorem]{Proposition}

\makeatletter
\def\imod#1{\allowbreak\mkern10mu({\operator@font mod}\,\,#1)}
\makeatother

\makeatletter 
 
\@addtoreset{algorithm}{chapter} 
\makeatother

\numberwithin{equation}{section}
\numberwithin{table}{section}
\numberwithin{figure}{section}

\newcommand{\bend}{\hspace*{0ex} \hfill \hbox{\vrule height
    1.5ex\vbox{\hrule width 1.4ex \vskip 1.4ex\hrule  width 1.4ex}\vrule
    height 1.5ex}}

\long\def\symbolfootnote[#1]#2{\begingroup \def\thefootnote{\fnsymbol{footnote}}\footnote[#1]{#2}\endgroup}
\makeatletter
\let\@fnsymbol\@arabic
\makeatother

\title{Transformed {rank-1} lattices for high-dimensional approximation}

\date{}
\author{
	Robert Nasdala\thanks{
    	Faculty of Mathematics, Chemnitz University of Technology, D-09107 Chemnitz, Germany.\newline E-mail:                \href{mailto:robert.nasdala@math.tu-chemnitz.de}{robert.nasdala@math.tu-chemnitz.de}
	}
	\and
	Daniel Potts\thanks{
		Faculty of Mathematics, Chemnitz University of Technology, D-09107 Chemnitz, Germany.\newline E-mail:                \href{mailto:potts@math.tu-chemnitz.de}{potts@math.tu-chemnitz.de}
	}
}

\hypersetup{pdfauthor={Robert Nasdala, Daniel Potts},
		pdftitle={Transformed {rank-1} lattices for high-dimensional approximation},
		pdfsubject={},
		pdfcreator = {pdflatex},
		plainpages=false,
		pdfstartview=FitH,       
		pdfview=FitH,            
		pdfpagemode=UseOutlines, 
		bookmarksnumbered=true, 
		bookmarksopen=false,     
		bookmarksopenlevel=0,   
		colorlinks=true,       
		linkcolor=black,
		citecolor=black,
		urlcolor=black}
\begin{document}

\maketitle

\begin{abstract}
	This paper describes an extension of Fourier approximation methods for multivariate functions defined on the torus $\mathbb{T}^d$ to functions in a weighted Hilbert space $L_{2}(\mathbb{R}^d, \omega)$
	via a multivariate change of variables $\psi:\left(-\frac{1}{2},\frac{1}{2}\right)^d\to\mathbb{R}^d$.
	We establish sufficient conditions on $\psi$ and $\omega$ such that the composition of a function in such a weighted Hilbert space with $\psi$ yields a function in the Sobolev space $H_{\mathrm{mix}}^{m}(\mathbb{T}^d)$ of functions on the torus with mixed smoothness of natural order $m \in \mathbb{N}_{0}$.
	In this approach we adapt algorithms for the evaluation and reconstruction of multivariate trigonometric polynomials on the torus $\mathbb{T}^d$ 
	based on single and multiple reconstructing {rank-$1$} lattices.
	Since in applications it may be difficult to choose a related function space, we make use of dimension incremental construction methods for sparse frequency sets.
	Various numerical tests confirm obtained theoretical results for the transformed methods.
\medskip

\noindent {Keywords and phrases} : approximation on unbounded domains, change of variables,
sparse multivariate trigonometric polynomials, lattice rule,  multiple {rank-$1$} lattice, fast Fourier transform

\medskip

\noindent {2010 AMS Mathematics Subject Classification} : \text{
65T 42B05 }
\end{abstract}

\medskip

\section{Introduction}
The change of variables is a powerful tool in numerical analysis. Such transformations play an important role in spectral methods, numerical integration, and approximation of functions.
An excellent overview  can be found in \cite[Chapter 16 and 17]{boyd00} which contains many practical aspects of the mapped methods. 
In this paper we focus on change of variable mappings from multivariate bounded domains to unbounded ones in order to approximate functions defined on such unbounded domains. 
The main goal is to transfer the approximation error bounds of Fourier methods on the high-dimensional torus $\mathbb T^d\simeq[-\frac{1}{2},\frac{1}{2})^d$ to approximation methods on $\mathbb R^d$ with the help of an invertible transformation $\psi: (-\frac{1}{2},\frac{1}{2})^d\to\mathbb{R}^d$.

Regarding functions defined on the torus $\mathbb T^d$ there is well-developed theory, see \cite{Tem93, DuTeUl2018, PlPoStTa18}, concerned with the Wiener algebra $\mathcal{A}(\mathbb{T}^d)$, that contains all $L_{1}(\mathbb{T}^d)$-functions with absolutely summable Fourier coefficients
\begin{align*}
	{\hat f_{\mathbf k} := \int_{\mathbb T^d} f(\mathbf x)\,\mathrm{e}^{-2\pi\mathrm i \mathbf k \cdot \mathbf x} \,\mathrm d\mathbf x}
\end{align*}
with $\mathbf k = (k_1,\ldots,k_d)^{\top}\in\mathbb{Z}^d, \mathbf x = (x_1,\ldots,x_d)^{\top}\in\mathbb{R}^d$, and $\mathbf k \cdot \mathbf x := \sum_{j=1}^{d} k_j x_j$.
For $\beta\geq 0$ and the weight function 
\begin{align} \label{def:hyperbolic_cross_weight}
	\omega_{\mathrm{hc}}(\mathbf k) := \prod_{j=1}^{d} \max(1, |k_j|),
\end{align}
there are the subspaces of the Wiener algebra $\mathcal{A}(\mathbb{T}^d)$ in form of
\begin{align}
	\label{def:Aalphaspace}
	\mathcal{A}^{\beta}(\mathbb T^d) := \left\{ f \in L_{1}(\mathbb{T}^d) : \|f\|_{\mathcal{A}^{\beta}(\mathbb{T}^d)} := \sum_{\mathbf k\in\mathbb{Z}^d} \omega_{\mathrm{hc}}(\mathbf k)^{\beta} |\hat f_{\mathbf k}| < \infty \right\}
\end{align}
and the Hilbert spaces 
\begin{align}
	\label{def:HbetaRaum}
	\mathcal{H}^{\beta}(\mathbb T^d) := \left\{ f \in L_2(\mathbb{T}^d) : \|f\|_{\mathcal{H}^{\beta}(\mathbb{T}^d)} := \left( \sum_{\mathbf k\in\mathbb{Z}^d} \omega_{\mathrm{hc}}(\mathbf k)^{2\beta} |\hat f_{\mathbf k}|^2 \right)^{\frac{1}{2}} < \infty \right\},
\end{align}
whose norms contain information about the decay rate of the Fourier coefficients $\hat f_{\mathbf k}$ with respect to the weight function ${\omega_{\mathrm{hc}}}$.
For approximation purposes we consider non-empty frequency sets ${I \subset \mathbb{Z}^d}$ of finite cardinality ${|I|<\infty}$ and approximated Fourier partial sums 
\begin{align}\label{def:approx_FPS}
	S_{I}^{\Lambda}f(\mathbf x) := \sum_{\mathbf k\in I} \hat f_{\mathbf k}^{\Lambda}\,\mathrm{e}^{2\pi\mathrm i \mathbf k\cdot \mathbf x}
\end{align}
with approximated Fourier coefficients 
\begin{align}\label{def:approx_FC}
	\hat f_{\mathbf k}^{\Lambda} := \frac{1}{M} \sum_{j =0}^{M-1} f(\mathbf x_j) \,\mathrm{e}^{-2\pi\mathrm{i}\mathbf k\cdot\mathbf x_j} \approx\hat f_{\mathbf k}, 
\end{align}
that are sampled at the nodes $\mathbf x_j$ of a reconstructing rank-$1$ lattice $\Lambda(\mathbf z,M,I)$, whose definition is given in \eqref{eq:rank1latticecondition}.

For $N\in\mathbb{N}$ and hyperbolic crosses 
\begin{align} \label{def:hyperbolic_cross}
	I_{N}^{d}
	:= \left\{ \mathbf k \in\mathbb{Z}^d : \omega_{\mathrm{hc}}(\mathbf k) \leq N \right\},
\end{align}
it was shown in \cite[Theorem~3.3]{KaPoVo13} that when using single rank-$1$ lattices the error of approximating a continuous function $f\in\mathcal{A}^{\beta}(\mathbb{T}^d)$ by the approximated Fourier partial sum $S_{I_{N}^{d}}^{\Lambda}f$ measured in the $L_{\infty}(\mathbb{T}^d)$-norm is bounded above by $N^{-\beta} \|f\|_{\mathcal{A}^{\beta}(\mathbb{T}^d)}$.
The approximation of functions in the Hilbert spaces $\mathcal{H}^{\beta}(\mathbb{T}^d)$ was investigated by V. N. Temlyakov in e.g., \cite{Tem86, KaPoVo13}.
For certain $\beta \geq 0$ the error of approximating a continuous function $f\in\mathcal{H}^{\beta}(\mathbb{T}^d)$ by the approximated Fourier partial sum $S_{I_{N}^{d}}^{\Lambda}f$ measured in the $L_{2}(\mathbb{T}^d)$-norm is bounded above by $C_{d,\beta} N^{-\beta} (\log N)^{(d-1)/2} \|f\|_{\mathcal{H}^{\beta}(\mathbb{T}^d)}$ with some constant $C_{d,\beta} = C(d,\beta) > 0$ as shown in \cite[Theorem~2.30]{volkmerdiss}.

A major problem is that in general it's hard to calculate the Fourier coefficients $\hat f_{\mathbf k}$ in order to determine if they are absolutely or square summable.
Instead we utilize certain norm equivalences to get information about the decay rate of the Fourier coefficients $\hat f_{\mathbf k}$. 
Given a multi-index $\bm{\alpha} = (\alpha_1, \ldots, \alpha_d)^{\top}\in\mathbb{N}_{0}^d$ with $\|\bm\alpha\|_{\ell_{\infty}} := \max(|\alpha_1|,\ldots,|\alpha_d|)$ and the differential operator
\begin{align} \label{def:differential_def_mult}
	D^{\bm \alpha}[f](\mathbf x) = D^{(\alpha_1,\ldots,\alpha_d)}[f](x_1,\ldots,x_d) := \frac{\partial^{\alpha_1}}{\partial x_{1}^{\alpha_1}}\ldots\frac{\partial^{\alpha_d}}{\partial x_{d}^{\alpha_d}} [f](x_1,\ldots,x_d),
\end{align}
we define for ${\Omega \in \{ \mathbb{T}^d, \mathbb{R}^d \}}$ the norm
\begin{align} \label{def:Hmix_norm}
	\|f\|_{H_{\mathrm{mix}}^{m}(\Omega)} := \left( \sum_{\|\bm\alpha\|_{\ell_{\infty}}\leq m} \| D^{\bm\alpha} [f] \|_{L_{2}(\Omega)}^{2} \right)^{1/2}
\end{align}
of the Sobolev space $H_{\mathrm{mix}}^{m}(\Omega)$ of functions $f\in L_{2}(\Omega)$ with mixed natural smoothness ${m\in\mathbb{N}_{0}}$, that were discussed in \cite{MR891189, UllTDiss, VybiralDiss}.
As shown in \cite[Lemma~2.3]{KuSiUl15}, the norms $\|\cdot\|_{H_{\mathrm{mix}}^{m}(\mathbb{T}^d)}$
and $\|\cdot\|_{\mathcal{H}^{\beta}(\mathbb{T}^d)}$ are equivalent for ${\beta = m \in\mathbb{N}}$.
Furthermore, for all $\beta \geq 0$ and all $\lambda>\frac{1}{2}$ we have the continuous embedding ${\mathcal{H}^{\beta+\lambda}(\mathbb{T}^d) \hookrightarrow \mathcal{A}^{\beta}(\mathbb{T}^d)}$ as shown in \cite[Lemma~2.2]{KaPoVo13}.
Hence, for $m\in\mathbb{N}$ we can just check if $f$ is an element of a Sobolev space $H_{\mathrm{mix}}^{m}(\mathbb{T}^d)$ in order to determine if a function $f$ is in $\mathcal{A}^{m}(\mathbb{T}^d)$ or $\mathcal{H}^{m}(\mathbb{T}^d)$ instead of calculating all its Fourier coefficients $\hat f_{\mathbf k}$.

In order to utilize all these properties for functions defined on $\mathbb{R}^d$, we apply continuously differentiable and strictly increasing change of variables $\psi: (-\frac{1}{2},\frac{1}{2})^d\to\mathbb{R}^d$ componentwise to multivariate functions $h$ in a weighted Hilbert space $L_{2}(\mathbb{R}^d,\omega)$ as defined in \eqref{def:weighted_L2_space} with the weight function ${\omega:\mathbb{R}^d \to [0,\infty)}$.
As a result we consider transformed functions $f\in L_{2}(\mathbb{T}^d)$ of the form 
\begin{align*}
	f(\mathbf x) = h(\psi(\mathbf x))\,\sqrt{ \omega(\psi(\mathbf x))\,\psi'(\mathbf x) },
\end{align*}
so that we have the identity $\|h\|_{L_{2}(\mathbb{R}^d,\omega)} = \|f\|_{L_{2}(\mathbb{T}^d)}$.
Based on this connection we will later on observe that the inverse transformation $\psi^{-1}$ transforms the classical Fourier system $\{ \mathrm{e}^{2\pi\mathrm{i}\mathbf k \circ} \}$ into another orthonormal system of the form $\left\{ \sqrt{\frac{(\psi^{-1})'(\circ)}{\omega(\circ)}} \, \mathrm{e}^{2\pi\mathrm i \mathbf k\cdot\psi^{-1}(\circ)}\right\}$.
It's generally rather difficult to check if such a transformed function $f$ is in the Sobolev space $H_{\mathrm{mix}}^{m}(\mathbb{T}^d)$ by calculating its norm and checking the various $L_2$-integrability conditions.
Therefore we provide a set of sufficient $L_{\infty}$-conditions for $f$ being in $H_{\mathrm{mix}}^{m}(\mathbb{T}^d)$.
At first we prove these conditions for all possible transformations $\psi$ and weight functions $\omega$. Later on we consider families of parameterized transformations ${\psi(\circ) = \psi(\circ, \bm\eta)}$ and families of weight functions ${\omega(\circ) = \omega(\circ, \bm\mu)}$ with $\bm\eta,\bm\mu \in\mathbb{R}^d$.
Then we have parameterized transformed functions ${f(\circ) = f(\circ,\bm\eta,\bm\mu) \in L_{2}(\mathbb{T}^d)}$ and both parameters may impact the smoothness of these functions.
With the sufficient $L_{\infty}$-smoothness conditions we're then able to calculate lower bounds for $\bm\eta$ and $\bm\mu$ such that the smoothness degree $m$ of a function ${h\in L_{2}(\mathbb{R}^d,\omega(\circ,\bm\mu)) \cap H_{\mathrm{mix}}^{m}(\mathbb{R}^d)}$ does not change under composition with a family of transformations $\psi(\circ,\bm\eta)$ so that we end up with ${f\in H_{\mathrm{mix}}^{m}(\mathbb{T}^d)}$.
For two particular transformation families $\psi(\circ,\bm\eta)$ we explicitly calculate the resulting lower parameter bounds and observe a case in which the smoothness preservation under the transformation depends only on the parameter $\bm\mu\in\mathbb{R}^d$ appearing in the weight functions $\omega(\circ,\bm\mu)$, as far as the conditions are able to detect it. 
Furthermore, we present an example in which we compare the parameter bounds yielded by the $L_{\infty}$-conditions with the exact lower bounds resulting from calculating the Sobolev-norm $\|\cdot\|_{H_{\mathrm{mix}}^{m}(\mathbb{T}^d)}$. This will highlight that the easier to check $L_{\infty}$-conditions yield slightly coarser parameter bounds.
These conditions as a tool to determine when a transformed function $f$ is at least an $L_{2}(\mathbb{T}^d)$-function, enable us to prove upper bounds for the approximation error ${\|h - S_{I}^{\Lambda}h\|}$ measured in weighted $L_{2}$- and $L_{\infty}$-norms on $\mathbb{R}^d$.
These are based on the already established error bounds for ${\|f - S_{I}^{\Lambda}f\|}$ on the torus with respect to the $L_{2}(\mathbb{T}^d)$- and $L_{\infty}(\mathbb{T}^d)$-norms.

One advantage of the proposed method is the availability of fast algorithms for high dimensional approximation, see e.g., \cite{boyd00}, in contrast to function approximation based on for instance multivariate Hermite functions or Sinc methods.
To this end, there are lattice rules that in recent years became an important tool in numerical analysis for high dimensional integration and approximation of multivariate functions. 
An introduction to lattice rules can be found in \cite{Nie78, SlJo94, JoKuSl13}. 
These rules are used for the approximation of functions on the torus, see \cite{Tem93}. 
Recently, efficient algorithms based on component-by-component methods \cite{CoKuNu10, NuCo04} were presented in order to compute high-dimensional integrals.
For the approximation of high-dimensional functions there are efficient algorithms using sampling schemes based on {rank-$1$} lattices \cite{KaPoVo13, kaemmererdiss}, and furthermore these schemes provide good approximation properties, see also \cite{ByKaUlVo16}. 
We adapt these algorithms and incorporate the outlined use of transformations. Furthermore, we present numerical examples.

We note that it was recently suggested in \cite{Kae16, Kae17} to use multiple {rank-$1$} lattices which are obtained by taking a union of several single {rank-$1$} lattices.
This method overcomes the limitations of the single {rank-$1$} lattice approach. 
That is, for the reconstruction of multivariate trigonometric polynomials supported on an arbitrary frequency set $I$ of finite cardinality ${|I|<\infty}$ with a single reconstructing {rank-$1$} lattice, the lattice size $M$ is bounded by ${|I| \leq M \leq |I|^2}$ under certain mild assumptions, see \cite[Lemma~2.1]{KaPoVo13} and \cite[Corollary~1]{Kae2013}. 
Multiple {rank-$1$} lattices improve the upper bound to $M\leq C |I|\log|I|$ with high probability \cite{Kae16, KaPoVo17}.
Remarkably, in both cases the upper bound is independent of the dimension $d$.
Furthermore there are methods where the support of the Fourier coefficients $\hat f_{\mathbf k}$ is unknown. 
We adapt the methods presented in \cite{PoVo14} that describe a dimension incremental construction of a frequency set $I\subset\mathbb{Z}^d$ containing only non-zero or the approximately largest Fourier coefficients $\hat h_{\mathbf k}$, based on component-by-component construction of {rank-$1$} lattices. 
This is done with respect to a specific search space in form of a full integer grid $[-N,N]^d\cap\mathbb{Z}^d$ with refinement $N\in\mathbb{N}$ and a sparsity that bounds the cardinality of the support.
We incorporate the change of variables method into both the multiple {rank-$1$} lattice methods as well as the component-by-component construction method.
Let us note that instead of rank-$1$ lattice points one can use a dimensional incremental support identification technique based on randomly chosen sampling points, that was recently developed in \cite{ChIwKr18}.

The outline of the paper is as follows:
In Section~2 we establish the basic notions from classical Fourier approximation theory on the torus $\mathbb T^d$, the corresponding function spaces and important convergence properties. 
We introduce the Sobolev spaces $H_{\mathrm{mix}}^{m}(\mathbb T^d)$ of mixed natural smoothness order $m\in\mathbb{N}_0$ and the Wiener Algebra $\mathcal{A}(\mathbb T^d)$ of functions with absolutely summable Fourier coefficients.
Furthermore, we discuss certain properties of the subspaces $\mathcal{A}^{\beta}(\mathbb T^d)$ and $\mathcal{H}^{\beta}(\mathbb T^d)$ of the Wiener Algebra, in particular we highlight the norm equivalence of $\|\cdot\|_{\mathcal{H}^{m}(\mathbb{T}^d)}$ and $\|\cdot\|_{H_{\mathrm{mix}}^{m}(\mathbb{T}^d)}$ for all $m\in\mathbb{N}$, see \cite{KuSiUl15} .
Then we define {rank-$1$} lattices as introduced in \cite{Ko59},
discuss their importance in the context of Fourier approximation 
and recall two important approximation error bounds on the torus in Theorems~\ref{thm:L_infty_approx_error_torus} and \ref{eq:L_2_approximation_error_bound}.
In Section~3 we define the notion of a transformation $\psi:(-\frac{1}{2},\frac{1}{2})^d \to \mathbb{R}^d$ and provide a couple of examples that we will use later on.
Then we introduce weight functions ${\omega:\mathbb{R}^d\to[0,\infty)}$ and
describe the structure of the weighted Hilbert spaces $L_2(\mathbb{R}^d,\omega)$, the corresponding weighted scalar product $(\cdot, \cdot)_{L_2(\mathbb{R}^d,\omega)}$ and the resulting Fourier coefficients $\hat h_{\mathbf k}$.
Afterwards we prove sufficient $L_{\infty}$-conditions on the transformation $\psi$ and weight function $\omega$, such that a function $h\in L_2(\mathbb{R}^d,\omega)\cap H_{\mathrm{mix}}^{m}(\mathbb{R}^d)$ is transformed under composition with $\psi$ into a smooth function $f \in H_{\mathrm{mix}}^{m}(\mathbb{T}^d)$.
Then we are able to prove approximation error bounds on $\mathbb{R}^d$ in Theorems~\ref{thm:L_infty_approx_error_multivar} and \ref{thm:Hm_approx_error_decay_multivar} based on the theorems on the torus that were recalled in Section~2.
In Section~4 we incorporate the usage of transformations $\psi$ into the algorithms \cite[Algorithm~3.1 and 3.2]{kaemmererdiss} for the evaluation and the reconstruction of multivariate functions in Algorithms~\ref{alg:LFFT_eval} and \ref{alg:LFFT_recon} based on transformed {rank-$1$} lattices.
In Section~5 we discuss examples for the algebraic transformation \eqref{eq:algebraic_trafo} and the error function transformation \eqref{eq:error_function_trafo} that were introduced in Section~3.
In these examples we use a parameterized transformation ${\psi(\circ) = \psi(\circ,\bm\eta)}$ with $\bm\eta\in\mathbb{R}^d$ and a parameterized weight function ${\omega(\circ) = \omega(\circ,\bm\mu)}$ with $\bm\mu\in\mathbb{R}^d$ that fit their original definitions in Sections~2 and 3.
With the sufficient $L_{\infty}$-conditions from Section~3 we then calculate explicit lower bounds for $\bm\eta$ and $\bm\mu$ that determine the degree of smoothness $m\in\mathbb{N}$ of ${h\in L_{2}(\mathbb{R}^d,\omega(\circ,\bm\mu)) \cap H_{\mathrm{mix}}^{m}(\mathbb{R}^d)}$ that is preserved under composition with the family of transformations ${\psi(\circ) = \psi(\circ,\bm\eta)}$.
Then we use the algorithms of the previous section to illustrate the theoretical upper approximation error bounds.
For some special cases in which the Fourier coefficients $\hat h_{\mathbf k}$ are explicitly given, we compare those to the theoretically proposed rate of decay of their absolutely values.
In Section~6 we add some remarks on how the tool of change of variables is incorporated into the ideas of multiple {rank-$1$} lattices and sparse fast Fourier algorithms. 
Furthermore we present examples with various test functions and different transformation maps in up to dimension $d=12$.

\section{Fourier approximation on the torus}
At first we introduce weighted $L_{p}$-function spaces and Sobolev spaces of mixed smoothness,
recall some definitions of classical Fourier approximation theory
and define a space of functions that have absolute square-summable Fourier coefficients.
Finally, we reflect the ideas of {rank-$1$} lattices from \cite{SLKA87, CoKuNu10,kaemmererdiss}, the corresponding Fourier approximation methods, and approximation error bounds that were discussed in e.g., \cite{Tem86, KaPoVo13, ByKaUlVo16}.
\subsection{Preliminaries}
Let $\Omega \in \{ \mathbb{T}^d, \mathbb{R}^d \}$ with $\mathbb T^d \simeq [-\frac{1}{2},\frac{1}{2})^d$ being the $d$-dimensional torus. 
The space ${(\mathcal{C}(\Omega), \|\cdot\|_{L_{\infty}(\Omega)})}$ denotes the collection of all continuous functions ${f:\Omega \to \mathbb{C}}$, and $(\mathcal{C}_{0}(\mathbb{R}^d), \|\cdot\|_{L_{\infty}(\mathbb{R}^d)})$ denotes the space of all continuous functions ${f:\mathbb{R}^d \to \mathbb{C}}$ vanishing at infinity in every direction.
We define weighted function spaces $L_p(\mathbb{R}^d, \omega)$ for $1\leq p < \infty$
with the weight function ${\omega:\mathbb{R}^d\to[0,\infty)}$ as
\begin{align}
	\label{def:weighted_L2_space}
	L_p(\mathbb{R}^d, \omega) := \left\{ h\in L_p(\mathbb{R}^d) : \|h\|_{L_p(\mathbb{R}^d, \omega)} := \left( \int_{\mathbb{R}^d} |h(\mathbf x)|^p \, \omega(\mathbf x) \,\mathrm{d}\mathbf x \right)^{\frac{1}{p}} < \infty \right\}
\end{align}
with the usual adjustments for $p=\infty$.
We have $L_{p}(\mathbb{R}^d) \subset L_p(\mathbb{R}^d, \omega)$ if $\omega$ is bounded 
and $L_p(\mathbb{R}^d, \omega) \subset L_{p}(\mathbb{R}^d)$ if $\omega$ is unbounded.
For the constant weight function $\omega(\mathbf x) \equiv 1$ we have $L_p(\mathbb{R}^d, \omega) = L_{p}(\mathbb{R}^d)$ and the $L_{p}(\mathbb{T}^d)$-spaces are defined analogously.

For functions $f$ and $g$ in the Hilbert space $L_2(\mathbb T^d)$ we have the scalar product
\begin{align*}
	(f,g)_{L_2(\mathbb T^d)} := \int_{\mathbb T^d} f(\mathbf x)\, \overline{g(\mathbf x)}\,\mathrm d\mathbf x .
\end{align*}
The functions $\mathrm{e}^{2\pi\mathrm i \mathbf k \cdot\mathbf x} := \prod_{j=1}^{d}\mathrm{e}^{2\pi\mathrm i k_j x_j}$ with ${\mathbf{k}\in \mathbb{Z}^d}$ and ${\mathbf x \in\mathbb T^d}$ 
are orthogonal with respect to the $L_2(\mathbb T^d)$-scalar product.
For any frequency set $I \subset \mathbb{Z}^d$ of finite cardinality $|I|<\infty$ we denote the space of all multivariate trigonometric polynomials supported on $I$ by
\begin{align*} 
	\Pi_{I} := \mathrm{span}\{ \mathrm{e}^{2\pi\mathrm i \mathbf k \cdot \circ} : \mathbf k\in I\}.
\end{align*}
For all $\mathbf k\in \mathbb{Z}^d$ we denote the \textsl{Fourier coefficients} $\hat f_{\mathbf k}$ by
\begin{align*}
	\hat f_{\mathbf k} 
	= (f, \mathrm{e}^{2\pi\mathrm i \mathbf k \cdot \circ})_{L_2(\mathbb T^d)} 
	= \int_{\mathbb T^d} f(\mathbf x)\,\mathrm{e}^{-2\pi\mathrm i \mathbf k \cdot \mathbf x} \,\mathrm d\mathbf x,
\end{align*}
and the corresponding \textsl{Fourier partial sum} by $S_{I}f(\mathbf x) := \sum_{\mathbf k\in I} \hat f_{\mathbf k}\,\mathrm{e}^{2\pi\mathrm i \mathbf k\mathbf \cdot \mathbf x}$.

For multi-indices $\bm\alpha \in\mathbb{N}_{0}^d$ and the differential operator $D^{\bm \alpha}[f](\mathbf x)$ as defined in \eqref{def:differential_def_mult} we define the \emph{Sobolev spaces of mixed natural smoothness} of $L_2(\Omega)$-functions with smoothness order $m\in\mathbb{N}_{0}$, see \cite{MR891189,UllTDiss, VybiralDiss}, as
\begin{align*} 
H_{\mathrm{mix}}^{m}(\Omega)
	= \left\{ f\in L_2(\Omega) : \|f\|_{H_{\mathrm{mix}}^{m}(\Omega)} < \infty \right\}
\end{align*}
with $\|\cdot\|_{H_{\mathrm{mix}}^{m}(\Omega)}$ as given in \eqref{def:Hmix_norm}.
The univariate spaces are denoted as $H^{m}(\mathbb{T})$ and $H^{m}(\mathbb{R})$, respectively.
For $\Omega = \mathbb{T}^d$ we recall some notation introduced in \cite{KuSiUl15}. 
The $H_{\mathrm{mix}}^{m}(\mathbb{T}^d)$-norm is expressible in terms of the Fourier coefficients $\hat f_{\mathbf k}$, which leads to the equivalent norm
\begin{align*}
	\|f\|_{H^{m,+}(\mathbb{T}^{d})}
	&:= \left( \sum_{\mathbf k\in\mathbb{Z}^{d}} |\hat{f}_{\mathbf k}|^{2} \prod_{j=1}^{d}(1+|k_j|^{2})^{m} \right)^{\frac{1}{2}}.
\end{align*}
In \cite[Lemma~2.3]{KuSiUl15} it is specified that for $m\in\mathbb{N}$ and all $f\in H_{\mathrm{mix}}^{m}(\mathbb{T}^{d})$ we have
\begin{align*}
\|f\|_{H_{\mathrm{mix}}^{m}(\mathbb{T}^{d})}
	\leq \|f\|_{H^{m,+}(\mathbb{T}^{d})}
	\leq \left(\frac{2^m}{m+1}\right)^{\frac{d}{2}} \|f\|_{H_{\mathrm{mix}}^{m}(\mathbb{T}^{d})}.
\end{align*}
The norm $\|\cdot \|_{H^{m,+}(\mathbb{T}^{d})}$-norm and the $\|\cdot\|_{\mathcal{H}^{\beta}(\mathbb{T}^d)}$-norm, that is given in \eqref{def:HbetaRaum}, are also equivalent for $m = \beta$ because of the observation that
\begin{align*}
	\max(1, |k_j|)^2 \leq 1+|k_j|^{2} \leq 2 \max(1, |k_j|)^2
\end{align*}
for all $k_j\in\mathbb{Z}$.
In total, for $m\in\mathbb{N}$ we have the norm equivalences
\begin{align} \label{eq:Hs_norm_equivalence}
	\|\cdot\|_{\mathcal{H}^{m}(\mathbb{T}^d)} \sim \|\cdot \|_{H^{m,+}(\mathbb{T}^{d})} \sim \|\cdot\|_{H_{\mathrm{mix}}^{m}(\mathbb{T}^d)},
\end{align}
but we distinguish the related function spaces anyway, because $\mathcal{H}^{\beta}(\mathbb{T}^d)$ appears in results concerned with approximation error bounds, whereas $H_{\mathrm{mix}}^{m}(\mathbb{T}^{d})$ is considered later on when we discuss smoothness-preserving transformation mappings.
Considering furthermore the function spaces $\mathcal{A}^{\beta}(\mathbb T^d)$ as defined in \eqref{def:Aalphaspace} it was shown in \cite[Lemma~2.2]{KaPoVo13} that for ${\beta\geq 0, \lambda>\frac{1}{2}}$ and fixed $d\in\mathbb{N}$ there are the continuous embeddings
\begin{align}
	\label{eq:Wiener_algebra_inclusion}
	\mathcal{H}^{\beta+\lambda}(\mathbb{T}^d) \hookrightarrow
	\mathcal{A}^{\beta}(\mathbb T^d) \hookrightarrow
	\mathcal{A}(\mathbb T^d)
\end{align}
and for $f\in\mathcal{A}^{\beta}(\mathbb T^d)$ we have 
\begin{align}
	\label{eq:Wiener_algebra_inclusion2}
	\|f\|_{\mathcal{A}^{\beta}(\mathbb T^d)} \leq C_{d,\lambda} \|f\|_{\mathcal{H}^{\beta+\lambda}(\mathbb{T}^d)}
\end{align}
with a constant $C_{d,\lambda} := C(d,\lambda) > 1$.
Additionally, for each function in $\mathcal{A}(\mathbb T^d)$ there exists a continuous representative, as proven in \cite[Lemma~2.1]{kaemmererdiss}.
Later on, when we sample functions $f\in\mathcal{H}^{\beta+\lambda}(\mathbb T^d)$ we identify them with their continuous representatives given by their Fourier series $\sum_{\mathbf k\in \mathbb{Z}^d} \hat f_{\mathbf k}\,\mathrm{e}^{2\pi\mathrm i \mathbf k\mathbf \cdot \circ}$
and this identification will be denoted by ${f\in\mathcal{H}^{\beta+\lambda}(\mathbb{T}^d)\cap\mathcal{C}(\mathbb{T}^d)}$.

\subsection{{Rank-1} lattices and reconstructing {rank-1} lattices}
Before discussing the approximation of functions ${f\in\mathcal{H}^{\beta}(\mathbb T^d)\cap\mathcal{C}(\mathbb{T}^d)}$ we recollect some related objects and observations from \cite{SLKA87, CoKuNu10, kaemmererdiss, PlPoStTa18}.
For each frequency set $I \subset \mathbb{Z}^d$ there is the \textsl{difference set}
\begin{align*}
	\mathcal{D}(I) &:= \{\mathbf{k} \in \mathbb{Z}^d : \mathbf{k}=\mathbf{k}_1-\mathbf{k}_2 \text{ with } \mathbf{k}_1,\mathbf{k}_2 \in I \}.
\end{align*}
Furthermore, the set
\begin{align}
	\label{def:rank_one_lattice}
	\Lambda(\mathbf{z}, M) &:= \left\{ \mathbf{x}_j := \left(\frac{j}{M}\,\mathbf z \bmod \mathbf 1\right) \in \mathbb{T}^d: j = 0,1,\ldots M-1 \right\}
\end{align}
is called \textsl{{rank-$1$} lattice} with the \textsl{generating vector} $\mathbf{z}\in\mathbb{Z}^d$ and the \textsl{lattice size} $M \in \mathbb{N}$, where $\mathbf 1 := \left(1,\ldots,1\right)^{\top}\in\mathbb{Z}^d$.
To ensure that $\Lambda(\mathbf z,M)$ has exactly $M$ distinct elements it is pointed out in \cite[p.428]{PlPoStTa18} that we need to assume that $M$ is coprime with at least one component of the generating vector $\mathbf z$.
A \textsl{reconstructing {rank-$1$} lattice} $\Lambda(\mathbf{z}, M, I)$ is a {rank-$1$} lattice $\Lambda(\mathbf{z}, M)$ for which the condition
\begin{align} \label{eq:rank1latticecondition}
\mathbf{t} \cdot \mathbf{z} \not\equiv 0 \,(\bmod{M}) \quad \text{for all } \mathbf{t}\in\mathcal{D}(I)\setminus\{\mathbf 0\} 
\end{align}
holds.
Given a reconstructing {rank-$1$} lattice $\Lambda(\mathbf z,M,I)$, we have exact integration for all multivariate trigonometric polynomials $g \in \Pi_{\mathcal{D}(I)}$, see \cite{SLKA87}, so that
\begin{align*}
	\int_{\mathbb{T}^d} g(\mathbf x) \,\mathrm{d}\mathbf x 
	= \frac{1}{M} \sum_{j =0}^{M-1} g(\mathbf x_j), 
	\quad \mathbf x_j\in\Lambda(\mathbf z,M,I).
\end{align*}
In particular, for $f \in \Pi_{I}$ and $\mathbf k\in I$ we have $f(\circ)\,\mathrm{e}^{-2\pi\mathrm{i}\mathbf k\cdot\circ} \in \Pi_{\mathcal{D}(I)}$ and
\begin{align} 
	\label{eq:exact_integration_formula}
	\hat f_{\mathbf k}
	= \int_{\mathbb{T}^d} f(\mathbf x) \,\mathrm{e}^{-2\pi\mathrm{i}\mathbf k\cdot\mathbf x} \,\mathrm{d}\mathbf x 
	= \frac{1}{M} \sum_{j =0}^{M-1} f(\mathbf x_j) \,\mathrm{e}^{-2\pi\mathrm{i}\mathbf k\cdot\mathbf x_j}, 
	\quad \mathbf x_j\in\Lambda(\mathbf z,M,I).
\end{align}
For an arbitrary function ${f \in \mathcal{H}^{\beta}(\mathbb T^d)\cap\mathcal{C}(\mathbb{T}^d)}$ and lattice points ${\mathbf x_j \in \Lambda(\mathbf z,M,I)}$ we lose the former mentioned exactness and get \textsl{approximated Fourier coefficients} $\hat f_{\mathbf k}^\Lambda \approx \hat f_{\mathbf k}$ of the form \eqref{def:approx_FC}
leading to the \textsl{approximated Fourier partial sum} $S_{I}^{\Lambda} f(\mathbf x) \approx	S_{I} f(\mathbf x)$ as given in \eqref{def:approx_FPS}.
\subsection{Lattice based approximation on the torus}
\begin{figure}[t]
	\centering
\begin{tikzpicture}[scale=1]
	\begin{axis}[scatter/classes = { a = {mark=o, draw=black} },
	font=\footnotesize,
	grid = both,
	xmax = 20, xmin = -20, ymax = 20, ymin = -20,
	height=0.4\textwidth, width=0.4\textwidth,
unit vector ratio*=1 1 1
	]
	\addplot[scatter ,only marks, mark size=1, scatter src = explicit symbolic]coordinates{
		(-16,-1)(-16,0)(-16,1)(-15,-1)(-15,0)(-15,1)(-14,-1)(-14,0)(-14,1)(-13,-1)(-13,0)(-13,1)(-12,-1)(-12,0)(-12,1)(-11,-1)(-11,0)(-11,1)(-10,-1)(-10,0)(-10,1)(-9,-1)(-9,0)(-9,1)(-8,-2)(-8,-1)(-8,0)(-8,1)(-8,2)(-7,-2)(-7,-1)(-7,0)(-7,1)(-7,2)(-6,-2)(-6,-1)(-6,0)(-6,1)(-6,2)(-5,-3)(-5,-2)(-5,-1)(-5,0)(-5,1)(-5,2)(-5,3)(-4,-4)(-4,-3)(-4,-2)(-4,-1)(-4,0)(-4,1)(-4,2)(-4,3)(-4,4)(-3,-5)(-3,-4)(-3,-3)(-3,-2)(-3,-1)(-3,0)(-3,1)(-3,2)(-3,3)(-3,4)(-3,5)(-2,-8)(-2,-7)(-2,-6)(-2,-5)(-2,-4)(-2,-3)(-2,-2)(-2,-1)(-2,0)(-2,1)(-2,2)(-2,3)(-2,4)(-2,5)(-2,6)(-2,7)(-2,8)(-1,-16)(-1,-15)(-1,-14)(-1,-13)(-1,-12)(-1,-11)(-1,-10)(-1,-9)(-1,-8)(-1,-7)(-1,-6)(-1,-5)(-1,-4)(-1,-3)(-1,-2)(-1,-1)(-1,0)(-1,1)(-1,2)(-1,3)(-1,4)(-1,5)(-1,6)(-1,7)(-1,8)(-1,9)(-1,10)(-1,11)(-1,12)(-1,13)(-1,14)(-1,15)(-1,16)(0,-16)(0,-15)(0,-14)(0,-13)(0,-12)(0,-11)(0,-10)(0,-9)(0,-8)(0,-7)(0,-6)(0,-5)(0,-4)(0,-3)(0,-2)(0,-1)(0,0)(0,1)(0,2)(0,3)(0,4)(0,5)(0,6)(0,7)(0,8)(0,9)(0,10)(0,11)(0,12)(0,13)(0,14)(0,15)(0,16)(1,-16)(1,-15)(1,-14)(1,-13)(1,-12)(1,-11)(1,-10)(1,-9)(1,-8)(1,-7)(1,-6)(1,-5)(1,-4)(1,-3)(1,-2)(1,-1)(1,0)(1,1)(1,2)(1,3)(1,4)(1,5)(1,6)(1,7)(1,8)(1,9)(1,10)(1,11)(1,12)(1,13)(1,14)(1,15)(1,16)(2,-8)(2,-7)(2,-6)(2,-5)(2,-4)(2,-3)(2,-2)(2,-1)(2,0)(2,1)(2,2)(2,3)(2,4)(2,5)(2,6)(2,7)(2,8)(3,-5)(3,-4)(3,-3)(3,-2)(3,-1)(3,0)(3,1)(3,2)(3,3)(3,4)(3,5)(4,-4)(4,-3)(4,-2)(4,-1)(4,0)(4,1)(4,2)(4,3)(4,4)(5,-3)(5,-2)(5,-1)(5,0)(5,1)(5,2)(5,3)(6,-2)(6,-1)(6,0)(6,1)(6,2)(7,-2)(7,-1)(7,0)(7,1)(7,2)(8,-2)(8,-1)(8,0)(8,1)(8,2)(9,-1)(9,0)(9,1)(10,-1)(10,0)(10,1)(11,-1)(11,0)(11,1)(12,-1)(12,0)(12,1)(13,-1)(13,0)(13,1)(14,-1)(14,0)(14,1)(15,-1)(15,0)(15,1)(16,-1)(16,0)(16,1) 
	};
	\end{axis}
	\end{tikzpicture}
	\caption{The hyperbolic cross $I_{N}^{d}$ for $N=16$ and $d=2$.}
	\label{fig:HyperbolicCrosses}
\end{figure}
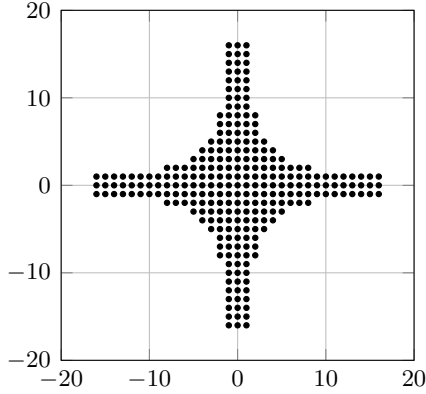

We discuss upper bounds for certain approximation errors $\|f - S_{I_N^d}^{\Lambda} f\|$ of functions $f$ in ${\mathcal{A}^{\beta}(\mathbb{T}^d)\cap\mathcal{C}(\mathbb{T}^d)}$ and ${\mathcal{H}^{\beta}(\mathbb{T}^d)\cap\mathcal{C}(\mathbb{T}^d)}$. 
For this matter the frequency sets are hyperbolic crosses $I_N^d$ as defined in \eqref{def:hyperbolic_cross} and are illustrated for $N=16$ in two dimensions in Figure~\ref{fig:HyperbolicCrosses}. 
For approximation purposes the existence of reconstructing {rank-$1$} lattices is secured by the arguments provided in \cite[Corollary~1]{Kae2013} and \cite[Lemma~2.1]{KaPoVo13}:
\begin{lemma}
	\label{lem:existence_of_M}
	Let $I\subset\mathbb{Z}^d$ be a frequency set of finite cardinality $4 \leq |I|<\infty$ and with ${I\subset \mathbb{Z}^d\cap(-M/2,M/2)^d}$.
	For all multivariate trigonometric polynomials $f\in\Pi_{I}$ there exists a reconstructing {rank-$1$} lattice $\Lambda(\mathbf z,M,I)$ 
	where the lattice size $M\in\mathbb{N}$ is bounded by ${|I| \leq M \leq |\mathcal{D}(I)| \leq |I|^2}$, 
	such that ${\hat f_{\mathbf k} = \hat f_{\mathbf k}^\Lambda}$.
	The generating vector $\mathbf z$ can be constructed using a component-by-component approach.
\end{lemma}

Then it's possible to prove an upper error bound for the $L_{\infty}$-approximation of functions in the subspace $\mathcal{A}^{\beta}(\mathbb{T}^d)$ of the Wiener Algebra, as seen in \cite[Theorem~3.3]{KaPoVo13}:
\begin{theorem} \label{thm:L_infty_approx_error_torus}
	Let $f \in \mathcal{A}^{\beta}(\mathbb{T}^d)\cap\mathcal{C}(\mathbb{T}^d)$ with $\beta\geq 0$ and $d\in\mathbb{N}$, 
	a hyperbolic cross $I_N^d$ with ${|I_N^d|<\infty}$ and $N\in\mathbb{N}$,
	and a reconstructing {rank-$1$} lattice ${\Lambda(\mathbf{z}, M, I_N^d)}$ be given.
	The approximation of $f$ by the approximated Fourier partial sum $S_{I_N^d}^{\Lambda} f$
	leads to an approximation error that is estimated by
	\begin{align} \label{eq:torus_infty_approx_bound}
		\|f - S_{I_N^d}^{\Lambda} f\|_{L_{\infty}(\mathbb{T}^d)}
		\leq 2 N^{-\beta} \|f\|_{\mathcal{A}^{\beta}(\mathbb{T}^d)}.
	\end{align}
\end{theorem}

The approximation of functions in the Hilbert spaces $\mathcal{H}^{\beta}(\mathbb{T}^d)$ was investigated by V. N. Temlyakov, see \cite{Tem86, KaPoVo13}. 
He showed that for $\beta > 1$ there exists a reconstructing {rank-$1$} lattice generated by a vector of Korobov form $\mathbf z := (1,z,z^2,\ldots,z^{d-1})^{\top}\in\mathbb{Z}^d$ such that the $L_2$-truncation error is bounded above by
\begin{align*}
	\|f - S_{I_{N}^{d}}^{\Lambda}f\|_{L_{2}(\mathbb{T}^d)} 
	\leq N^{-\beta} (\log N)^{(d-1)/2} \|f\|_{\mathcal{H}^{\beta}(\mathbb{T}^d)}.
\end{align*}
A generalization of this estimate as well as an upper bound for the corresponding aliasing error is stated in \cite[Theorem~2]{ByKaUlVo16} using dyadic hyperbolic cross frequency sets and a component-by-component approach to construct the generating vector $\mathbf z\in\mathbb{Z}^d$, which generally isn't of Korobov form anymore.
However, every dyadic hyperbolic cross is embedded in a non-dyadic one, see \cite[Lemma~2.29]{volkmerdiss}. 
Thus, the error estimates are easily translated in terms of non-dyadic hyperbolic crosses $I_{N}^{d}$, see \cite[Theorem~2.30]{volkmerdiss},
and we are particularly interested in the following special case:
\begin{theorem}\label{eq:L_2_approximation_error_bound}
	Let $\beta > \frac{1}{2}$,
	the dimension $d\in\mathbb{N}$,
	a function $f\in\mathcal{H}^{\beta}(\mathbb{T}^d)\cap\mathcal{C}(\mathbb{T}^d)$,
	a hyperbolic cross $I_{N}^{d}$ with $N\geq 2^{d+1}$,
	and a reconstructing {rank-$1$} lattice $\Lambda(\mathbf z, M, I_{N}^{d})$ be given.
	Then we have
	\begin{align} \label{eq:H_beta_error_bound}
		\|f - S_{I_{N}^{d}}^{\Lambda}f\|_{L_{2}(\mathbb{T}^d)} 
		\leq C_{d,\beta} N^{-\beta} (\log N)^{(d-1)/2} \|f\|_{\mathcal{H}^{\beta}(\mathbb{T}^d)}
	\end{align}
	with some constant $C_{d,\beta} := C(d,\beta)>0$.
\end{theorem}

As highlighted earlier in \eqref{eq:Hs_norm_equivalence}, for $\beta=m\in\mathbb{N}$ the norms $\|\cdot\|_{\mathcal{H}^{\beta}(\mathbb{T}^d)}$ and $\|\cdot\|_{H_{\mathrm{mix}}^{m}(\mathbb{T}^d)}$ are equivalent.
Eventually we utilize this norm equivalence in order to apply the above approximation error bounds for functions $f$ in the Sobolev space $H_{\mathrm{mix}}^{m}(\mathbb{T}^d)$ that are characterized by their derivatives.

\section{Torus-to-$\mathbb{R}$ transformation mappings}
Change of variables were discussed for example in \cite{boyd00, ShTaWa11} and were used for high dimensional integration in e.g.,\cite{KuWaWa06, KPPW18}.
In this chapter we define transformations $\psi:(-\frac{1}{2},\frac{1}{2})^d\to\mathbb{R}^d$ and provide examples that will reappear later in this paper. 
Afterwards we describe the weighted Hilbert spaces $L_{2}(\mathbb{R}^d,\omega)$ with weight functions $\omega:\mathbb{R}^d\to[0,\infty)$ and investigate their structure.
Then we prove sufficient conditions on $\psi$ and $\omega$ such that an initially chosen $h\in L_{2}(\mathbb{R}^d,\omega) \cap H_{\mathrm{mix}}^{m}(\mathbb{R}^d)$ is transformed by the change of variables $\psi$ into a function that is lying in a Sobolev space $H_{\mathrm{mix}}^{m}(\mathbb{T}^d)$ of mixed natural smoothness order $m\in\mathbb{N}_{0}$. 
Eventually we show that with an incorporated transformation $\psi$, we still have upper bounds for certain approximation errors on $\mathbb{R}^d$, which are based on the already established error bounds with respect to the $L_{\infty}(\mathbb{T}^d)$- and $L_{2}(\mathbb{T}^d)$-norms recalled in Theorems~\ref{thm:L_infty_approx_error_torus} and \ref{eq:L_2_approximation_error_bound} respectively.

\subsection{Transformations to $\mathbb{R}^d$}
We call a map $\psi:(-\frac{1}{2},\frac{1}{2}) \to \mathbb{R}$ a \textsl{transformation} or \textsl{change of variables} if it is continuously differentiable, strictly increasing, odd, and we have
\begin{align}\label{def:Trafo_def}
	\lim\limits_{x \to -\frac{1}{2}}\psi(x) = -\infty, 
	\quad \lim\limits_{x \to \frac{1}{2}}\psi(x) = \infty.
\end{align}
We denote its first derivative by $\psi'(x) := \frac{\mathrm{d}}{\mathrm{d}x}[\psi](x) $.
The respective inverse transformation is also continuously differentiable, increasing and is denoted by $\psi^{-1}:\mathbb{R}\to (-\frac{1}{2},\frac{1}{2})$ in the sense of $y=\psi(x) \Leftrightarrow x=\psi^{-1}(y)$.
We call the derivative of the inverse transformation the \emph{density function of $\psi$}, which we define as
\begin{align}\label{def:Varrho_def}
	\varrho(y) := (\psi^{-1})'(y) = \frac{1}{\psi'(\psi^{-1}(y))}
\end{align}
and for which we have $\varrho(y)\geq 0$ for all $y\in\mathbb{R}$.
Furthermore, we have $\lim\limits_{y\to -\infty}\psi^{-1}(y)= -\frac{1}{2}$ and $\lim\limits_{y\to \infty}\psi^{-1}(y)= \frac{1}{2}$.
We note that $\varrho$ is a bounded function with 
\begin{align*} 
	\left\|\varrho\right\|_{L_{1}(\mathbb{R})}
	= \int_{-\infty}^{\infty} \varrho(y) \,\mathrm{d}y 
	= 1 . \end{align*}
For multivariate transformations we put 
\begin{align*}
	\psi(\mathbf x) := ( \psi_1(x_1),\ldots,\psi_d(x_d) )^{\top}
	\quad \text{and} \quad
	\psi'(\mathbf x) := \prod_{j=1}^d \psi_j'(x_j)
\end{align*}
with ${\mathbf x = (x_1,\ldots,x_d)^{\top}\in(-\frac{1}{2},\frac{1}{2})^d}$, 
where we may use different transformations $\psi_j$ in each direction.
Similarly, we put ${\psi^{-1}(\mathbf y) := ( \psi_1^{-1}(y_1),\ldots,\psi_d^{-1}(y_d) )^{\top}}$ and 
\begin{align} \label{def:varrho_mult}
	{\varrho(\mathbf y) := \prod_{j=1}^d \varrho_j(y_j)}
\end{align}
with ${\mathbf y = (y_1,\ldots,y_d)^{\top}\in\mathbb{R}^d}$.

Later on we consider families of parameterized transformations 
\begin{align}
	\label{eq:psi_weighted_mult}
	\psi(\mathbf x,\bm \eta) := ( \psi_1(x_1,\eta_1),\ldots,\psi_d(x_d,\eta_d) )^{\top}
\end{align}
with $\bm \eta = (\eta_1,\ldots,\eta_d)^{\top}\in\mathbb{R}^d$.
We only consider parametrizations for which the transformation $\psi$, its inverse $\psi^{-1}$ and the density function $\varrho$ fit into the given definitions above despite being impacted by the parameter $\bm\eta$. 
On several occasions throughout this paper we will replace transformations $\psi(x)$ by 
\begin{align} \label{eq:param_trafo_explicit_example}
	\psi(\mathbf x,\bm\eta) := \bm\eta\cdot\psi(\mathbf x)
\end{align}
with $\bm\eta\in(0,\infty)^d$.
As the transformations are going to be composed with functions defined on $\mathbb{R}^d$ the parameter $\bm \eta$ may impact the smoothness of the resulting transformed functions, which we will discuss in depth later on.
For now, we omit the parameter in the notation for simplicity and proceed to just write $\psi(\circ)$ until we actually consider particular parameterized families of the form \eqref{eq:psi_weighted_mult} or \eqref{eq:param_trafo_explicit_example}.

\subsection{Exemplary transformations}
We list some feasible univariate transformations $\psi$ with either an algebraic or an exponential density function $\varrho$, some of which were suggested in the literature, see e.g.,\cite[Section~17.6]{boyd00} and \cite[Section~7.5]{ShTaWa11}.
With the remark on \eqref{eq:psi_weighted_mult} in mind, we list these transformations here in their univariate non-parameterized form with $\eta=1$ and $\psi(x)=\psi(x,1)$ for simplicity.
Later on when we fix a particular family of parameterized transformations $\psi(\circ,\eta), \eta\in\mathbb{R}$ we recall these definitions accordingly.

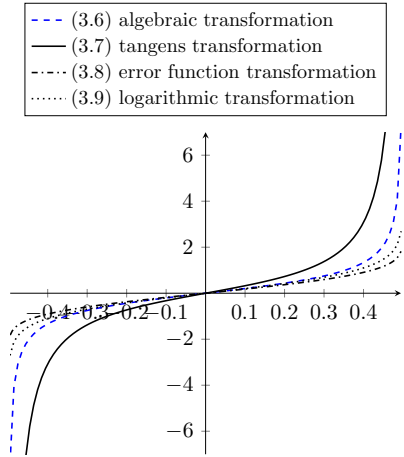
\begin{figure}[t]
\centering
\begin{tikzpicture}[scale=0.75]
\begin{axis}[samples=1000,  
xtick={-0.4,-0.3,-0.2,-0.1,0,0.1,0.2,0.3,0.4},
xmin=-0.495, 
xmax=0.495,
ytick={-6,-4,-2,0,2,4,6},
ymin=-7, 
ymax=7,
axis x line=center, 
axis y line=center,
every axis plot/.append style={thick},
legend style={at={(0.5,1.05)}, anchor=south,legend columns=1,legend cell align=left, font=\small}
]
\addplot[blue, dashed] { 2*x/sqrt(1-4*x^2) }; \addlegendentry{\eqref{eq:algebraic_trafo} algebraic transformation};
\addplot[black] {tan(pi*deg(x))}; \addlegendentry{\eqref{eq:tan_trafo} tangens transformation};
\addplot[dashdotted, mark size = 2] coordinates { 
	(-0.495, -1.820160e+00) (-0.4899, -1.643646e+00) (-0.4849, -1.533081e+00) (-0.4799, -1.450755e+00) (-0.4749, -1.384387e+00) (-0.4698, -1.328359e+00) (-0.4648, -1.279608e+00) (-0.4598, -1.236276e+00) (-0.4548, -1.197144e+00) (-0.4497, -1.161368e+00) (-0.4447, -1.128340e+00) (-0.4397, -1.097605e+00) (-0.4347, -1.068814e+00) (-0.4296, -1.041693e+00) (-0.4246, -1.016023e+00) (-0.4196, -9.916272e-01) (-0.4146, -9.683573e-01) (-0.4095, -9.460913e-01) (-0.4045, -9.247259e-01) (-0.3995, -9.041728e-01) (-0.3945, -8.843565e-01) (-0.3894, -8.652114e-01) (-0.3844, -8.466804e-01) (-0.3794, -8.287133e-01) (-0.3744, -8.112660e-01) (-0.3693, -7.942990e-01) (-0.3643, -7.777774e-01) (-0.3593, -7.616699e-01) (-0.3543, -7.459482e-01) (-0.3492, -7.305869e-01) (-0.3442, -7.155629e-01) (-0.3392, -7.008551e-01) (-0.3342, -6.864445e-01) (-0.3291, -6.723135e-01) (-0.3241, -6.584460e-01) (-0.3191, -6.448272e-01) (-0.3141, -6.314436e-01) (-0.309, -6.182824e-01) (-0.304, -6.053320e-01) (-0.299, -5.925815e-01) (-0.294, -5.800208e-01) (-0.2889, -5.676406e-01) (-0.2839, -5.554320e-01) (-0.2789, -5.433867e-01) (-0.2739, -5.314971e-01) (-0.2688, -5.197559e-01) (-0.2638, -5.081563e-01) (-0.2588, -4.966919e-01) (-0.2538, -4.853565e-01) (-0.2487, -4.741446e-01) (-0.2437, -4.630505e-01) (-0.2387, -4.520694e-01) (-0.2337, -4.411962e-01) (-0.2286, -4.304263e-01) (-0.2236, -4.197554e-01) (-0.2186, -4.091792e-01) (-0.2136, -3.986937e-01) (-0.2085, -3.882953e-01) (-0.2035, -3.779801e-01) (-0.1985, -3.677447e-01) (-0.1935, -3.575859e-01) (-0.1884, -3.475003e-01) (-0.1834, -3.374849e-01) (-0.1784, -3.275367e-01) (-0.1734, -3.176530e-01) (-0.1683, -3.078310e-01) (-0.1633, -2.980680e-01) (-0.1583, -2.883614e-01) (-0.1533, -2.787090e-01) (-0.1482, -2.691082e-01) (-0.1432, -2.595567e-01) (-0.1382, -2.500524e-01) (-0.1332, -2.405930e-01) (-0.1281, -2.311765e-01) (-0.1231, -2.218008e-01) (-0.1181, -2.124640e-01) (-0.1131, -2.031640e-01) (-0.108, -1.938991e-01) (-0.103, -1.846673e-01) (-0.09799, -1.754669e-01) (-0.09296, -1.662962e-01) (-0.08794, -1.571533e-01) (-0.08291, -1.480366e-01) (-0.07789, -1.389444e-01) (-0.07286, -1.298752e-01) (-0.06784, -1.208273e-01) (-0.06281, -1.117991e-01) (-0.05779, -1.027891e-01) (-0.05276, -9.379577e-02) (-0.04774, -8.481759e-02) (-0.04271, -7.585306e-02) (-0.03769, -6.690070e-02) (-0.03266, -5.795906e-02) (-0.02764, -4.902667e-02) (-0.02261, -4.010210e-02) (-0.01759, -3.118392e-02) (-0.01256, -2.227069e-02) (-0.007538, -1.336100e-02) (-0.002513, -4.453431e-03) (0.002513, 4.453431e-03) (0.007538, 1.336100e-02) (0.01256, 2.227069e-02) (0.01759, 3.118392e-02) (0.02261, 4.010210e-02) (0.02764, 4.902667e-02) (0.03266, 5.795906e-02) (0.03769, 6.690070e-02) (0.04271, 7.585306e-02) (0.04774, 8.481759e-02) (0.05276, 9.379577e-02) (0.05779, 1.027891e-01) (0.06281, 1.117991e-01) (0.06784, 1.208273e-01) (0.07286, 1.298752e-01) (0.07789, 1.389444e-01) (0.08291, 1.480366e-01) (0.08794, 1.571533e-01) (0.09296, 1.662962e-01) (0.09799, 1.754669e-01) (0.103, 1.846673e-01) (0.108, 1.938991e-01) (0.1131, 2.031640e-01) (0.1181, 2.124640e-01) (0.1231, 2.218008e-01) (0.1281, 2.311765e-01) (0.1332, 2.405930e-01) (0.1382, 2.500524e-01) (0.1432, 2.595567e-01) (0.1482, 2.691082e-01) (0.1533, 2.787090e-01) (0.1583, 2.883614e-01) (0.1633, 2.980680e-01) (0.1683, 3.078310e-01) (0.1734, 3.176530e-01) (0.1784, 3.275367e-01) (0.1834, 3.374849e-01) (0.1884, 3.475003e-01) (0.1935, 3.575859e-01) (0.1985, 3.677447e-01) (0.2035, 3.779801e-01) (0.2085, 3.882953e-01) (0.2136, 3.986937e-01) (0.2186, 4.091792e-01) (0.2236, 4.197554e-01) (0.2286, 4.304263e-01) (0.2337, 4.411962e-01) (0.2387, 4.520694e-01) (0.2437, 4.630505e-01) (0.2487, 4.741446e-01) (0.2538, 4.853565e-01) (0.2588, 4.966919e-01) (0.2638, 5.081563e-01) (0.2688, 5.197559e-01) (0.2739, 5.314971e-01) (0.2789, 5.433867e-01) (0.2839, 5.554320e-01) (0.2889, 5.676406e-01) (0.294, 5.800208e-01) (0.299, 5.925815e-01) (0.304, 6.053320e-01) (0.309, 6.182824e-01) (0.3141, 6.314436e-01) (0.3191, 6.448272e-01) (0.3241, 6.584460e-01) (0.3291, 6.723135e-01) (0.3342, 6.864445e-01) (0.3392, 7.008551e-01) (0.3442, 7.155629e-01) (0.3492, 7.305869e-01) (0.3543, 7.459482e-01) (0.3593, 7.616699e-01) (0.3643, 7.777774e-01) (0.3693, 7.942990e-01) (0.3744, 8.112660e-01) (0.3794, 8.287133e-01) (0.3844, 8.466804e-01) (0.3894, 8.652114e-01) (0.3945, 8.843565e-01) (0.3995, 9.041728e-01) (0.4045, 9.247259e-01) (0.4095, 9.460913e-01) (0.4146, 9.683573e-01) (0.4196, 9.916272e-01) (0.4246, 1.016023e+00) (0.4296, 1.041693e+00) (0.4347, 1.068814e+00) (0.4397, 1.097605e+00) (0.4447, 1.128340e+00) (0.4497, 1.161368e+00) (0.4548, 1.197144e+00) (0.4598, 1.236276e+00) (0.4648, 1.279608e+00) (0.4698, 1.328359e+00) (0.4749, 1.384387e+00) (0.4799, 1.450755e+00) (0.4849, 1.533081e+00) (0.4899, 1.643646e+00) (0.495, 1.820160e+00) 
}; \addlegendentry{\eqref{eq:error_function_trafo} error function transformation};
\addplot[dotted] {(1/2)*ln( (1/2+x)/(1/2-x )}; \addlegendentry{\eqref{eq:logarithmic_trafo} logarithmic transformation}  
\end{axis}
\end{tikzpicture}
\caption{Plots of exemplary transformations \eqref{eq:algebraic_trafo}-\eqref{eq:logarithmic_trafo}.}
\label{fig:Plot of trafo examples}
\end{figure}

Let $x\in(-\frac{1}{2}, \frac{1}{2})$ and $y\in\mathbb{R}$.
We are particularly interested in the following transformations:
\begin{itemize}
	\item 
	\textsl{algebraic transformation:}
	\begin{align} \label{eq:algebraic_trafo}
		\psi(x) &= \frac{2x}{(1-4x^2)^{\frac{1}{2}}} ,
		\quad \psi'(x) = \frac{2}{(1-4x^2)^{\frac{3}{2}}},\\
		\psi^{-1}(y) &= \frac{y}{2 (1+y^2)^{\frac{1}{2}}},
		\quad \varrho(y) = \frac{1}{2 (1+y^2)^{\frac{3}{2}}} \nonumber
		\end{align}
	\item
	\textsl{tangent transformation}:
	\begin{align} \label{eq:tan_trafo}
		\psi(x) &= \tan\left( \pi x \right),
		\quad \psi'(x) = \frac{\pi}{\cos^2(\pi x)}\\
		\psi^{-1}(y) &= \frac{1}{\pi} \arctan\left(y\right) ,
		\quad \varrho(y) = \frac{1}{\pi}\left( \frac{1}{1 + y^2} \right) \nonumber
	\end{align}
	\item 
	\textsl{error function transformation}:
	\begin{align} \label{eq:error_function_trafo}
		\psi(x) &= \mathrm{erf}^{-1}(2x) ,
		\quad \psi'(x) = \sqrt{\pi} \, \mathrm{e}^{(\mathrm{erf}^{-1}(2x))^2} \\
		\psi^{-1}(y) &= \frac{1}{2}\,\mathrm{erf}\left(y\right),
		\quad \varrho(y) = \frac{1}{\sqrt{\pi}}\,\mathrm{e}^{-y^2} \nonumber
	\end{align}
	with the error function 
	\begin{align*}
		\mathrm{erf}(x) = \frac{1}{\sqrt{\pi}} \int_{-x}^{x} \mathrm{e}^{-t^2} \,\mathrm{d}t, \quad x\in\mathbb{R},
	\end{align*}
	and $\mathrm{erf}^{-1}(\circ)$ denoting the inverse error function
	\item
	\textsl{logarithmic transformation}:
	\begin{align} \label{eq:logarithmic_trafo}
		\psi(x) 
		&= \frac{1}{2}\log\left(\frac{1+2x}{1-2x}\right) 
		= \tanh^{-1}(2x),
		\quad \psi'(x) = \frac{2}{1-4x^2}, \\
		\psi^{-1}(y) &= \frac{1}{2}\left(\frac{\mathrm{e}^{2y} - 1}{\mathrm{e}^{2y} + 1}\right) = \frac{1}{2}\tanh\left(y\right) ,
		\quad \varrho(y) = \frac{2\mathrm{e}^{2y}}{(\mathrm{e}^{2y}+1)^2}
		\nonumber
	\end{align}
\end{itemize}
For a side-by-side comparison of their individual slope see Figure~\ref{fig:Plot of trafo examples}.

\subsection{Weighted Hilbert spaces on $\mathbb{R}$}
We describe the structure of the weighted $L_2(\mathbb{R},\omega)$-function spaces as defined in $\eqref{def:weighted_L2_space}$.
In this section the weight function $\omega:\mathbb{R}\to[0,\infty)$ remains unspecified. 
However, similarly to the generalization \eqref{eq:psi_weighted_mult} of transformations $\psi$ defined in \eqref{def:Trafo_def}, we will later on consider families of non-negative parameterized weight functions $\omega(\circ,\mu)$ with $\mu\in\mathbb{R}$ for the purpose of controlling the smoothness of functions in $L_2(\mathbb{R},\omega(\circ,\mu))\cap H^{m}(\mathbb{R})$ and of the corresponding transformed functions on the torus $\mathbb{T}$.
Analogously, families of multivariate parameterized weight functions are defined as 
\begin{align}
	\label{eq:omega_weighted_mult}
	\omega(\mathbf y, \bm \mu) := \prod_{j=1}^{d}\omega_j(y_j,\mu_j), \quad \mathbf y, \bm \mu \in \mathbb{R}^d
\end{align}
with univariate weight functions $\omega_j(\circ,\mu_j):\mathbb{R}\to[0,\infty)$.

For now we remain in the univariate setting.
The system $\left\{\varphi_{k}\right\}_{k\in\mathbb{Z}}$ of weighted exponential functions
\begin{align}\label{eq:transformed_basis_functions}
	\varphi_{k}(y)
	:= \sqrt{\frac{\varrho(y)}{\omega(y)}} \, \mathrm{e}^{2\pi\mathrm i k\psi^{-1}(y)}, \quad y\in\mathbb{R}
\end{align}
forms an orthogonal system with respect to the scalar product
\begin{align} \label{def:weighted_scalar_product}
	(h_1, h_2)_{L_2\left(\mathbb{R}, \omega \right)}
	:= \int_{\mathbb{R}} \omega(y) \, h_1(y) \, \overline{h_2(y)} \, \mathrm dy
\end{align}
and for $k_1,k_2\in\mathbb{Z}$ we have
\begin{align*}
	(\varphi_{k_1}, \varphi_{k_2})_{L_2\left(\mathbb{R}, \omega\right)}
	= \delta_{k_1, k_2}.
\end{align*}
The weighted scalar product \eqref{def:weighted_scalar_product} induces the norm
\begin{align*}
	\|h\|_{L_2\left(\mathbb{R}, \omega \right)}
	:= \sqrt{ (h, h)_{L_2\left(\mathbb{R}, \omega \right)} }
\end{align*}
and in a natural way we have Fourier coefficients of the form
\begin{align} \label{def:FouCoeff_of_h}
	\hat h_{k}
	:= \left(h, \varphi_{k} \right)_{L_2\left(\mathbb{R}, \omega \right)}
	= \int_{\mathbb{R}} h(y) \, \sqrt{\varrho(y)\,\omega(y)} \, \mathrm{e}^{-2\pi\mathrm i k\psi^{-1}(y)} \, \mathrm dy,
\end{align}
as well as the respective Fourier partial sum for $I\subset\mathbb{Z}$ given by
\begin{align} \label{def:Fourier_part_sum_of_h}
	S_{I}h(y) 
	:= \sum_{k\in I} \hat h_{k} \, \varphi_{k}(y).
\end{align}
\begin{example}
\begin{itemize}
\item 
	For the algebraic transformation \eqref{eq:algebraic_trafo} with the density $\varrho(y) = \frac{1}{2 (1+y^2)^{\frac{3}{2}}}$ 
	and the parameterized weight function 
	\begin{align}\label{eq:alebraic_weight_function_parametrized_univar}
		\omega(y,\mu) = \left( \frac{1}{1+y^2} \right)^{\mu}, \quad \mu\in\mathbb{R},
	\end{align}
	the orthogonal system functions $\varphi_{k}$ as in $\eqref{eq:transformed_basis_functions}$ are of the form
	\begin{align*}
		\varphi_{k}(y)
		= \sqrt{\frac{1}{2}\left(\frac{1}{1+y^2}\right)^{\frac{3}{2}-\mu}} \, \mathrm{e}^{\pi\mathrm i k \frac{y}{\sqrt{1+y^2}}}.
	\end{align*}
	The graphs of their real and imaginary parts of these $\varphi_{k}$ are shown for $\mu=2$ and $k= 0,1,2,3$ in Figure~\ref{fig:phi_k_example_with_alg_trafo}.
\item 
	For the error function transformation \eqref{eq:error_function_trafo} with the density $\varrho(y) = \frac{1}{\sqrt{\pi}}\,\mathrm{e}^{-y^2}$
	and the Gaussian weight function 
	\begin{align}\label{eq:gaussian_weight_param_univar}
		\omega(y,\mu) = \frac{1}{\sqrt{\pi}}\mathrm{e}^{-\mu^2 y^2}, \quad \mu\in\mathbb{R}, 
	\end{align}
	the orthogonal system functions $\varphi_{k}$ as in $\eqref{eq:transformed_basis_functions}$ are of the form
	\begin{align*}
		\varphi_{k}(y)
		= \mathrm{e}^{\frac{1}{2}(\mu^2-1)y^2 + \pi\mathrm i k\,\mathrm{erf}\left(y\right)},
	\end{align*}
	with graphs of their real and imaginary parts for $\mu=\sqrt{2}$ and $k=0,1,2,3$ shown in Figure~\eqref{fig:phi_k_example_with_erf_trafo}
	and the corresponding weighted scalar product \eqref{def:weighted_scalar_product} reads as
	\begin{align*}
		(h_1, h_2)_{L_2\left(\mathbb{R}, \omega(\circ, \mu) \right)}
		= \frac{1}{\sqrt{\pi}}\,\int_{\mathbb{R}} \mathrm{e}^{-\mu^2 y^2} \, h_1(y) \, \overline{h_2(y)} \, \mathrm d y.
	\end{align*}
\end{itemize}
\end{example}

\begin{figure}[t]
	\begin{minipage}{.5\linewidth}
		\centering
\begin{tikzpicture}[scale=0.75]
		\pgfmathsetmacro\MU{2}
		\begin{axis}[samples=500,  
		xmin=-5, xmax=5, ymin=-1.5, ymax=1.5,
		title = {$\mathrm{Re}(\varphi_{k}(y)) = \sqrt{\frac{1}{2} \sqrt{1+y^2} } \, \cos\left(\pi k \frac{y}{\sqrt{1+y^2}}\right)$},
		axis x line=center, axis y line=center,
		every axis plot/.append style={thick},
		legend style={at={(0.5,1.25)}, anchor=south,legend columns=4,legend cell align=left, font=\small}
		]
		\addplot[dashdotted] { sqrt((1/2)*(1/(1+x^2))^(3/2-\MU))*cos(pi*0*(deg(x/(sqrt(1+x^2))))) };
		\addlegendentry{$k=0$};
		\addplot[black] { sqrt((1/2)*(1/(1+x^2))^(3/2-\MU))*cos(pi*1*(deg(x/(sqrt(1+x^2))))) };
		\addlegendentry{$k=1$};
		\addplot[dashed] { sqrt((1/2)*(1/(1+x^2))^(3/2-\MU))*cos(pi*2*(deg(x/(sqrt(1+x^2))))) }; 
		\addlegendentry{$k=2$};
		\addplot[dotted] { sqrt((1/2)*(1/(1+x^2))^(3/2-\MU))*cos(pi*3*(deg(x/(sqrt(1+x^2)))))  };
		\addlegendentry{$k=3$};   
		\end{axis}
		\end{tikzpicture}
	\end{minipage}
	\begin{minipage}{.5\linewidth}
		\centering
\begin{tikzpicture}[scale=0.75]
		\pgfmathsetmacro\MU{2}
		\begin{axis}[samples=500,  
		xmin=-5, xmax=5, ymin=-1.5, ymax=1.5,
		title = {$\mathrm{Im}(\varphi_{k}(y)) = \sqrt{\frac{1}{2} \sqrt{1+y^2} } \, \sin\left(\pi k \frac{y}{\sqrt{1+y^2}}\right)$},
		axis x line=center, axis y line=center,
		every axis plot/.append style={thick},
		legend style={at={(0.5,1.25)}, anchor=south,legend columns=4,legend cell align=left, font=\small}
		]
\addplot[black] { sqrt((1/2)*(1/(1+x^2))^(3/2-\MU))*sin(pi*1*(deg(x/(sqrt(1+x^2))))) };
		\addlegendentry{$k=1$};
		\addplot[dashed] { sqrt((1/2)*(1/(1+x^2))^(3/2-\MU))*sin(pi*2*(deg(x/(sqrt(1+x^2))))) }; 
		\addlegendentry{$k=2$};
		\addplot[dotted] { sqrt((1/2)*(1/(1+x^2))^(3/2-\MU))*sin(pi*3*(deg(x/(sqrt(1+x^2)))))  };
		\addlegendentry{$k=3$};  
\end{axis}
		\end{tikzpicture}
	\end{minipage}
	\caption{Real and imaginary part of the weighted exponential functions $\varphi_{k}, k=0,1,2,3$ in \eqref{eq:transformed_basis_functions} with the density function $\varrho$ of the algebraic transformation \eqref{eq:algebraic_trafo} and the algebraic parameterized weight function $\omega(y,\mu)$ as given in \eqref{eq:alebraic_weight_function_parametrized} for fixed $\mu = 2$.}
	\label{fig:phi_k_example_with_alg_trafo}
\end{figure}
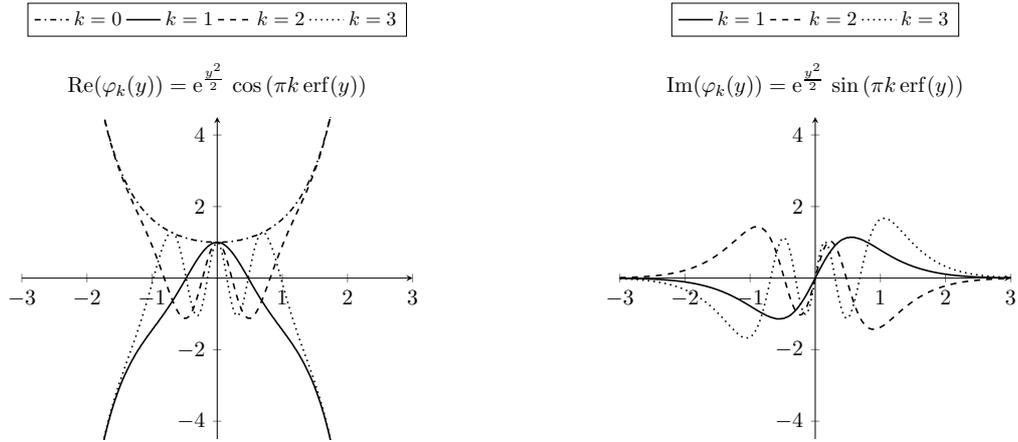
\begin{figure}[t]
	\begin{minipage}{.5\linewidth}
	\centering
\begin{tikzpicture}[scale=0.75]
	\pgfmathsetmacro\MU{2}
	\begin{axis}[samples=500,  
	xmin=-3, xmax=3, ymin=-4.5, ymax=4.5,
title = {$\mathrm{Re}(\varphi_{k}(y)) = \mathrm{e}^{\frac{y^2}{2}} \, \cos\left(\pi k \,\mathrm{erf}(y)\right)$},
	axis x line=center, axis y line=center,
	every axis plot/.append style={thick},
	legend style={at={(0.5,1.25)}, anchor=south,legend columns=4,legend cell align=left, font=\small}
	]
	\addplot[dashdotted] coordinates { 
		(-2, 7.39) (-1.987, 7.19) (-1.973, 7.01) (-1.96, 6.82) (-1.946, 6.65) (-1.933, 6.48) (-1.92, 6.31) (-1.906, 6.15) (-1.893, 6) (-1.88, 5.85) (-1.866, 5.71) (-1.853, 5.57) (-1.839, 5.43) (-1.826, 5.3) (-1.813, 5.17) (-1.799, 5.05) (-1.786, 4.93) (-1.773, 4.81) (-1.759, 4.7) (-1.746, 4.59) (-1.732, 4.48) (-1.719, 4.38) (-1.706, 4.28) (-1.692, 4.19) (-1.679, 4.09) (-1.666, 4) (-1.652, 3.92) (-1.639, 3.83) (-1.625, 3.75) (-1.612, 3.67) (-1.599, 3.59) (-1.585, 3.51) (-1.572, 3.44) (-1.559, 3.37) (-1.545, 3.3) (-1.532, 3.23) (-1.518, 3.17) (-1.505, 3.1) (-1.492, 3.04) (-1.478, 2.98) (-1.465, 2.92) (-1.452, 2.87) (-1.438, 2.81) (-1.425, 2.76) (-1.411, 2.71) (-1.398, 2.66) (-1.385, 2.61) (-1.371, 2.56) (-1.358, 2.51) (-1.344, 2.47) (-1.331, 2.43) (-1.318, 2.38) (-1.304, 2.34) (-1.291, 2.3) (-1.278, 2.26) (-1.264, 2.22) (-1.251, 2.19) (-1.237, 2.15) (-1.224, 2.12) (-1.211, 2.08) (-1.197, 2.05) (-1.184, 2.02) (-1.171, 1.98) (-1.157, 1.95) (-1.144, 1.92) (-1.13, 1.89) (-1.117, 1.87) (-1.104, 1.84) (-1.09, 1.81) (-1.077, 1.79) (-1.064, 1.76) (-1.05, 1.74) (-1.037, 1.71) (-1.023, 1.69) (-1.01, 1.67) (-0.9967, 1.64) (-0.9833, 1.62) (-0.9699, 1.6) (-0.9565, 1.58) (-0.9431, 1.56) (-0.9298, 1.54) (-0.9164, 1.52) (-0.903, 1.5) (-0.8896, 1.49) (-0.8763, 1.47) (-0.8629, 1.45) (-0.8495, 1.43) (-0.8361, 1.42) (-0.8227, 1.4) (-0.8094, 1.39) (-0.796, 1.37) (-0.7826, 1.36) (-0.7692, 1.34) (-0.7559, 1.33) (-0.7425, 1.32) (-0.7291, 1.3) (-0.7157, 1.29) (-0.7023, 1.28) (-0.689, 1.27) (-0.6756, 1.26) (-0.6622, 1.25) (-0.6488, 1.23) (-0.6355, 1.22) (-0.6221, 1.21) (-0.6087, 1.2) (-0.5953, 1.19) (-0.5819, 1.18) (-0.5686, 1.18) (-0.5552, 1.17) (-0.5418, 1.16) (-0.5284, 1.15) (-0.5151, 1.14) (-0.5017, 1.13) (-0.4883, 1.13) (-0.4749, 1.12) (-0.4615, 1.11) (-0.4482, 1.11) (-0.4348, 1.1) (-0.4214, 1.09) (-0.408, 1.09) (-0.3946, 1.08) (-0.3813, 1.08) (-0.3679, 1.07) (-0.3545, 1.06) (-0.3411, 1.06) (-0.3278, 1.06) (-0.3144, 1.05) (-0.301, 1.05) (-0.2876, 1.04) (-0.2742, 1.04) (-0.2609, 1.03) (-0.2475, 1.03) (-0.2341, 1.03) (-0.2207, 1.02) (-0.2074, 1.02) (-0.194, 1.02) (-0.1806, 1.02) (-0.1672, 1.01) (-0.1538, 1.01) (-0.1405, 1.01) (-0.1271, 1.01) (-0.1137, 1.01) (-0.1003, 1.01) (-0.08696, 1) (-0.07358, 1) (-0.0602, 1) (-0.04682, 1) (-0.03344, 1) (-0.02007, 1) (-0.006689, 1) (0.006689, 1) (0.02007, 1) (0.03344, 1) (0.04682, 1) (0.0602, 1) (0.07358, 1) (0.08696, 1) (0.1003, 1.01) (0.1137, 1.01) (0.1271, 1.01) (0.1405, 1.01) (0.1538, 1.01) (0.1672, 1.01) (0.1806, 1.02) (0.194, 1.02) (0.2074, 1.02) (0.2207, 1.02) (0.2341, 1.03) (0.2475, 1.03) (0.2609, 1.03) (0.2742, 1.04) (0.2876, 1.04) (0.301, 1.05) (0.3144, 1.05) (0.3278, 1.06) (0.3411, 1.06) (0.3545, 1.06) (0.3679, 1.07) (0.3813, 1.08) (0.3946, 1.08) (0.408, 1.09) (0.4214, 1.09) (0.4348, 1.1) (0.4482, 1.11) (0.4615, 1.11) (0.4749, 1.12) (0.4883, 1.13) (0.5017, 1.13) (0.5151, 1.14) (0.5284, 1.15) (0.5418, 1.16) (0.5552, 1.17) (0.5686, 1.18) (0.5819, 1.18) (0.5953, 1.19) (0.6087, 1.2) (0.6221, 1.21) (0.6355, 1.22) (0.6488, 1.23) (0.6622, 1.25) (0.6756, 1.26) (0.689, 1.27) (0.7023, 1.28) (0.7157, 1.29) (0.7291, 1.3) (0.7425, 1.32) (0.7559, 1.33) (0.7692, 1.34) (0.7826, 1.36) (0.796, 1.37) (0.8094, 1.39) (0.8227, 1.4) (0.8361, 1.42) (0.8495, 1.43) (0.8629, 1.45) (0.8763, 1.47) (0.8896, 1.49) (0.903, 1.5) (0.9164, 1.52) (0.9298, 1.54) (0.9431, 1.56) (0.9565, 1.58) (0.9699, 1.6) (0.9833, 1.62) (0.9967, 1.64) (1.01, 1.67) (1.023, 1.69) (1.037, 1.71) (1.05, 1.74) (1.064, 1.76) (1.077, 1.79) (1.09, 1.81) (1.104, 1.84) (1.117, 1.87) (1.13, 1.89) (1.144, 1.92) (1.157, 1.95) (1.171, 1.98) (1.184, 2.02) (1.197, 2.05) (1.211, 2.08) (1.224, 2.12) (1.237, 2.15) (1.251, 2.19) (1.264, 2.22) (1.278, 2.26) (1.291, 2.3) (1.304, 2.34) (1.318, 2.38) (1.331, 2.43) (1.344, 2.47) (1.358, 2.51) (1.371, 2.56) (1.385, 2.61) (1.398, 2.66) (1.411, 2.71) (1.425, 2.76) (1.438, 2.81) (1.452, 2.87) (1.465, 2.92) (1.478, 2.98) (1.492, 3.04) (1.505, 3.1) (1.518, 3.17) (1.532, 3.23) (1.545, 3.3) (1.559, 3.37) (1.572, 3.44) (1.585, 3.51) (1.599, 3.59) (1.612, 3.67) (1.625, 3.75) (1.639, 3.83) (1.652, 3.92) (1.666, 4) (1.679, 4.09) (1.692, 4.19) (1.706, 4.28) (1.719, 4.38) (1.732, 4.48) (1.746, 4.59) (1.759, 4.7) (1.773, 4.81) (1.786, 4.93) (1.799, 5.05) (1.813, 5.17) (1.826, 5.3) (1.839, 5.43) (1.853, 5.57) (1.866, 5.71) (1.88, 5.85) (1.893, 6) (1.906, 6.15) (1.92, 6.31) (1.933, 6.48) (1.946, 6.65) (1.96, 6.82) (1.973, 7.01) (1.987, 7.19) (2, 7.39)
	}; \addlegendentry{$k=0$};
	\addplot[black] coordinates{ 
		(-2, -7.39) (-1.987, -7.19) (-1.973, -7.01) (-1.96, -6.82) (-1.946, -6.65) (-1.933, -6.48) (-1.92, -6.31) (-1.906, -6.15) (-1.893, -6) (-1.88, -5.85) (-1.866, -5.7) (-1.853, -5.56) (-1.839, -5.43) (-1.826, -5.3) (-1.813, -5.17) (-1.799, -5.04) (-1.786, -4.92) (-1.773, -4.81) (-1.759, -4.7) (-1.746, -4.59) (-1.732, -4.48) (-1.719, -4.38) (-1.706, -4.28) (-1.692, -4.18) (-1.679, -4.09) (-1.666, -4) (-1.652, -3.91) (-1.639, -3.82) (-1.625, -3.74) (-1.612, -3.66) (-1.599, -3.58) (-1.585, -3.5) (-1.572, -3.43) (-1.559, -3.36) (-1.545, -3.29) (-1.532, -3.22) (-1.518, -3.15) (-1.505, -3.09) (-1.492, -3.02) (-1.478, -2.96) (-1.465, -2.9) (-1.452, -2.84) (-1.438, -2.79) (-1.425, -2.73) (-1.411, -2.68) (-1.398, -2.63) (-1.385, -2.58) (-1.371, -2.53) (-1.358, -2.48) (-1.344, -2.43) (-1.331, -2.38) (-1.318, -2.34) (-1.304, -2.29) (-1.291, -2.25) (-1.278, -2.21) (-1.264, -2.16) (-1.251, -2.12) (-1.237, -2.08) (-1.224, -2.04) (-1.211, -2) (-1.197, -1.97) (-1.184, -1.93) (-1.171, -1.89) (-1.157, -1.85) (-1.144, -1.82) (-1.13, -1.78) (-1.117, -1.75) (-1.104, -1.71) (-1.09, -1.68) (-1.077, -1.64) (-1.064, -1.61) (-1.05, -1.58) (-1.037, -1.54) (-1.023, -1.51) (-1.01, -1.48) (-0.9967, -1.44) (-0.9833, -1.41) (-0.9699, -1.38) (-0.9565, -1.34) (-0.9431, -1.31) (-0.9298, -1.28) (-0.9164, -1.25) (-0.903, -1.21) (-0.8896, -1.18) (-0.8763, -1.14) (-0.8629, -1.11) (-0.8495, -1.08) (-0.8361, -1.04) (-0.8227, -1.01) (-0.8094, -0.974) (-0.796, -0.939) (-0.7826, -0.903) (-0.7692, -0.868) (-0.7559, -0.832) (-0.7425, -0.795) (-0.7291, -0.758) (-0.7157, -0.721) (-0.7023, -0.684) (-0.689, -0.646) (-0.6756, -0.607) (-0.6622, -0.569) (-0.6488, -0.53) (-0.6355, -0.49) (-0.6221, -0.45) (-0.6087, -0.41) (-0.5953, -0.369) (-0.5819, -0.329) (-0.5686, -0.287) (-0.5552, -0.246) (-0.5418, -0.204) (-0.5284, -0.162) (-0.5151, -0.12) (-0.5017, -0.0782) (-0.4883, -0.0359) (-0.4749, 0.00639) (-0.4615, 0.0487) (-0.4482, 0.091) (-0.4348, 0.133) (-0.4214, 0.175) (-0.408, 0.217) (-0.3946, 0.258) (-0.3813, 0.299) (-0.3679, 0.34) (-0.3545, 0.38) (-0.3411, 0.419) (-0.3278, 0.458) (-0.3144, 0.496) (-0.301, 0.534) (-0.2876, 0.57) (-0.2742, 0.605) (-0.2609, 0.64) (-0.2475, 0.673) (-0.2341, 0.705) (-0.2207, 0.736) (-0.2074, 0.765) (-0.194, 0.793) (-0.1806, 0.819) (-0.1672, 0.844) (-0.1538, 0.867) (-0.1405, 0.889) (-0.1271, 0.909) (-0.1137, 0.926) (-0.1003, 0.943) (-0.08696, 0.957) (-0.07358, 0.969) (-0.0602, 0.979) (-0.04682, 0.987) (-0.03344, 0.994) (-0.02007, 0.998) (-0.006689, 1) (0.006689, 1) (0.02007, 0.998) (0.03344, 0.994) (0.04682, 0.987) (0.0602, 0.979) (0.07358, 0.969) (0.08696, 0.957) (0.1003, 0.943) (0.1137, 0.926) (0.1271, 0.909) (0.1405, 0.889) (0.1538, 0.867) (0.1672, 0.844) (0.1806, 0.819) (0.194, 0.793) (0.2074, 0.765) (0.2207, 0.736) (0.2341, 0.705) (0.2475, 0.673) (0.2609, 0.64) (0.2742, 0.605) (0.2876, 0.57) (0.301, 0.534) (0.3144, 0.496) (0.3278, 0.458) (0.3411, 0.419) (0.3545, 0.38) (0.3679, 0.34) (0.3813, 0.299) (0.3946, 0.258) (0.408, 0.217) (0.4214, 0.175) (0.4348, 0.133) (0.4482, 0.091) (0.4615, 0.0487) (0.4749, 0.00639) (0.4883, -0.0359) (0.5017, -0.0782) (0.5151, -0.12) (0.5284, -0.162) (0.5418, -0.204) (0.5552, -0.246) (0.5686, -0.287) (0.5819, -0.329) (0.5953, -0.369) (0.6087, -0.41) (0.6221, -0.45) (0.6355, -0.49) (0.6488, -0.53) (0.6622, -0.569) (0.6756, -0.607) (0.689, -0.646) (0.7023, -0.684) (0.7157, -0.721) (0.7291, -0.758) (0.7425, -0.795) (0.7559, -0.832) (0.7692, -0.868) (0.7826, -0.903) (0.796, -0.939) (0.8094, -0.974) (0.8227, -1.01) (0.8361, -1.04) (0.8495, -1.08) (0.8629, -1.11) (0.8763, -1.14) (0.8896, -1.18) (0.903, -1.21) (0.9164, -1.25) (0.9298, -1.28) (0.9431, -1.31) (0.9565, -1.34) (0.9699, -1.38) (0.9833, -1.41) (0.9967, -1.44) (1.01, -1.48) (1.023, -1.51) (1.037, -1.54) (1.05, -1.58) (1.064, -1.61) (1.077, -1.64) (1.09, -1.68) (1.104, -1.71) (1.117, -1.75) (1.13, -1.78) (1.144, -1.82) (1.157, -1.85) (1.171, -1.89) (1.184, -1.93) (1.197, -1.97) (1.211, -2) (1.224, -2.04) (1.237, -2.08) (1.251, -2.12) (1.264, -2.16) (1.278, -2.21) (1.291, -2.25) (1.304, -2.29) (1.318, -2.34) (1.331, -2.38) (1.344, -2.43) (1.358, -2.48) (1.371, -2.53) (1.385, -2.58) (1.398, -2.63) (1.411, -2.68) (1.425, -2.73) (1.438, -2.79) (1.452, -2.84) (1.465, -2.9) (1.478, -2.96) (1.492, -3.02) (1.505, -3.09) (1.518, -3.15) (1.532, -3.22) (1.545, -3.29) (1.559, -3.36) (1.572, -3.43) (1.585, -3.5) (1.599, -3.58) (1.612, -3.66) (1.625, -3.74) (1.639, -3.82) (1.652, -3.91) (1.666, -4) (1.679, -4.09) (1.692, -4.18) (1.706, -4.28) (1.719, -4.38) (1.732, -4.48) (1.746, -4.59) (1.759, -4.7) (1.773, -4.81) (1.786, -4.92) (1.799, -5.04) (1.813, -5.17) (1.826, -5.3) (1.839, -5.43) (1.853, -5.56) (1.866, -5.7) (1.88, -5.85) (1.893, -6) (1.906, -6.15) (1.92, -6.31) (1.933, -6.48) (1.946, -6.65) (1.96, -6.82) (1.973, -7.01) (1.987, -7.19) (2, -7.39)
	};
	\addlegendentry{$k=1$};
	\addplot[dashed] coordinates{ 
		(-2, 7.39) (-1.987, 7.19) (-1.973, 7) (-1.96, 6.82) (-1.946, 6.64) (-1.933, 6.47) (-1.92, 6.31) (-1.906, 6.15) (-1.893, 5.99) (-1.88, 5.84) (-1.866, 5.7) (-1.853, 5.56) (-1.839, 5.42) (-1.826, 5.29) (-1.813, 5.16) (-1.799, 5.04) (-1.786, 4.91) (-1.773, 4.8) (-1.759, 4.68) (-1.746, 4.57) (-1.732, 4.47) (-1.719, 4.36) (-1.706, 4.26) (-1.692, 4.16) (-1.679, 4.07) (-1.666, 3.98) (-1.652, 3.89) (-1.639, 3.8) (-1.625, 3.71) (-1.612, 3.63) (-1.599, 3.55) (-1.585, 3.47) (-1.572, 3.39) (-1.559, 3.32) (-1.545, 3.25) (-1.532, 3.17) (-1.518, 3.1) (-1.505, 3.04) (-1.492, 2.97) (-1.478, 2.9) (-1.465, 2.84) (-1.452, 2.78) (-1.438, 2.72) (-1.425, 2.65) (-1.411, 2.6) (-1.398, 2.54) (-1.385, 2.48) (-1.371, 2.42) (-1.358, 2.37) (-1.344, 2.31) (-1.331, 2.26) (-1.318, 2.2) (-1.304, 2.15) (-1.291, 2.09) (-1.278, 2.04) (-1.264, 1.99) (-1.251, 1.94) (-1.237, 1.88) (-1.224, 1.83) (-1.211, 1.78) (-1.197, 1.73) (-1.184, 1.67) (-1.171, 1.62) (-1.157, 1.57) (-1.144, 1.51) (-1.13, 1.46) (-1.117, 1.41) (-1.104, 1.35) (-1.09, 1.3) (-1.077, 1.24) (-1.064, 1.18) (-1.05, 1.13) (-1.037, 1.07) (-1.023, 1.01) (-1.01, 0.952) (-0.9967, 0.892) (-0.9833, 0.831) (-0.9699, 0.77) (-0.9565, 0.707) (-0.9431, 0.644) (-0.9298, 0.58) (-0.9164, 0.516) (-0.903, 0.45) (-0.8896, 0.384) (-0.8763, 0.318) (-0.8629, 0.251) (-0.8495, 0.183) (-0.8361, 0.115) (-0.8227, 0.0475) (-0.8094, -0.0207) (-0.796, -0.0887) (-0.7826, -0.157) (-0.7692, -0.224) (-0.7559, -0.291) (-0.7425, -0.357) (-0.7291, -0.422) (-0.7157, -0.487) (-0.7023, -0.549) (-0.689, -0.61) (-0.6756, -0.669) (-0.6622, -0.726) (-0.6488, -0.78) (-0.6355, -0.831) (-0.6221, -0.879) (-0.6087, -0.924) (-0.5953, -0.965) (-0.5819, -1) (-0.5686, -1.03) (-0.5552, -1.06) (-0.5418, -1.09) (-0.5284, -1.1) (-0.5151, -1.12) (-0.5017, -1.12) (-0.4883, -1.12) (-0.4749, -1.12) (-0.4615, -1.11) (-0.4482, -1.09) (-0.4348, -1.07) (-0.4214, -1.04) (-0.408, -1) (-0.3946, -0.958) (-0.3813, -0.909) (-0.3679, -0.854) (-0.3545, -0.794) (-0.3411, -0.728) (-0.3278, -0.657) (-0.3144, -0.582) (-0.301, -0.502) (-0.2876, -0.419) (-0.2742, -0.332) (-0.2609, -0.244) (-0.2475, -0.153) (-0.2341, -0.0608) (-0.2207, 0.0316) (-0.2074, 0.124) (-0.194, 0.215) (-0.1806, 0.304) (-0.1672, 0.391) (-0.1538, 0.475) (-0.1405, 0.555) (-0.1271, 0.63) (-0.1137, 0.699) (-0.1003, 0.763) (-0.08696, 0.82) (-0.07358, 0.87) (-0.0602, 0.912) (-0.04682, 0.947) (-0.03344, 0.973) (-0.02007, 0.99) (-0.006689, 0.999) (0.006689, 0.999) (0.02007, 0.99) (0.03344, 0.973) (0.04682, 0.947) (0.0602, 0.912) (0.07358, 0.87) (0.08696, 0.82) (0.1003, 0.763) (0.1137, 0.699) (0.1271, 0.63) (0.1405, 0.555) (0.1538, 0.475) (0.1672, 0.391) (0.1806, 0.304) (0.194, 0.215) (0.2074, 0.124) (0.2207, 0.0316) (0.2341, -0.0608) (0.2475, -0.153) (0.2609, -0.244) (0.2742, -0.332) (0.2876, -0.419) (0.301, -0.502) (0.3144, -0.582) (0.3278, -0.657) (0.3411, -0.728) (0.3545, -0.794) (0.3679, -0.854) (0.3813, -0.909) (0.3946, -0.958) (0.408, -1) (0.4214, -1.04) (0.4348, -1.07) (0.4482, -1.09) (0.4615, -1.11) (0.4749, -1.12) (0.4883, -1.12) (0.5017, -1.12) (0.5151, -1.12) (0.5284, -1.1) (0.5418, -1.09) (0.5552, -1.06) (0.5686, -1.03) (0.5819, -1) (0.5953, -0.965) (0.6087, -0.924) (0.6221, -0.879) (0.6355, -0.831) (0.6488, -0.78) (0.6622, -0.726) (0.6756, -0.669) (0.689, -0.61) (0.7023, -0.549) (0.7157, -0.487) (0.7291, -0.422) (0.7425, -0.357) (0.7559, -0.291) (0.7692, -0.224) (0.7826, -0.157) (0.796, -0.0887) (0.8094, -0.0207) (0.8227, 0.0475) (0.8361, 0.115) (0.8495, 0.183) (0.8629, 0.251) (0.8763, 0.318) (0.8896, 0.384) (0.903, 0.45) (0.9164, 0.516) (0.9298, 0.58) (0.9431, 0.644) (0.9565, 0.707) (0.9699, 0.77) (0.9833, 0.831) (0.9967, 0.892) (1.01, 0.952) (1.023, 1.01) (1.037, 1.07) (1.05, 1.13) (1.064, 1.18) (1.077, 1.24) (1.09, 1.3) (1.104, 1.35) (1.117, 1.41) (1.13, 1.46) (1.144, 1.51) (1.157, 1.57) (1.171, 1.62) (1.184, 1.67) (1.197, 1.73) (1.211, 1.78) (1.224, 1.83) (1.237, 1.88) (1.251, 1.94) (1.264, 1.99) (1.278, 2.04) (1.291, 2.09) (1.304, 2.15) (1.318, 2.2) (1.331, 2.26) (1.344, 2.31) (1.358, 2.37) (1.371, 2.42) (1.385, 2.48) (1.398, 2.54) (1.411, 2.6) (1.425, 2.65) (1.438, 2.72) (1.452, 2.78) (1.465, 2.84) (1.478, 2.9) (1.492, 2.97) (1.505, 3.04) (1.518, 3.1) (1.532, 3.17) (1.545, 3.25) (1.559, 3.32) (1.572, 3.39) (1.585, 3.47) (1.599, 3.55) (1.612, 3.63) (1.625, 3.71) (1.639, 3.8) (1.652, 3.89) (1.666, 3.98) (1.679, 4.07) (1.692, 4.16) (1.706, 4.26) (1.719, 4.36) (1.732, 4.47) (1.746, 4.57) (1.759, 4.68) (1.773, 4.8) (1.786, 4.91) (1.799, 5.04) (1.813, 5.16) (1.826, 5.29) (1.839, 5.42) (1.853, 5.56) (1.866, 5.7) (1.88, 5.84) (1.893, 5.99) (1.906, 6.15) (1.92, 6.31) (1.933, 6.47) (1.946, 6.64) (1.96, 6.82) (1.973, 7) (1.987, 7.19) (2, 7.39)
	}; 
	\addlegendentry{$k=2$};
	\addplot[dotted] coordinates{
	 	(-2, -7.38) (-1.987, -7.19) (-1.973, -7) (-1.96, -6.82) (-1.946, -6.64) (-1.933, -6.47) (-1.92, -6.3) (-1.906, -6.14) (-1.893, -5.98) (-1.88, -5.83) (-1.866, -5.69) (-1.853, -5.55) (-1.839, -5.41) (-1.826, -5.28) (-1.813, -5.15) (-1.799, -5.02) (-1.786, -4.9) (-1.773, -4.78) (-1.759, -4.66) (-1.746, -4.55) (-1.732, -4.44) (-1.719, -4.34) (-1.706, -4.24) (-1.692, -4.14) (-1.679, -4.04) (-1.666, -3.94) (-1.652, -3.85) (-1.639, -3.76) (-1.625, -3.67) (-1.612, -3.58) (-1.599, -3.5) (-1.585, -3.42) (-1.572, -3.34) (-1.559, -3.26) (-1.545, -3.18) (-1.532, -3.1) (-1.518, -3.03) (-1.505, -2.95) (-1.492, -2.88) (-1.478, -2.81) (-1.465, -2.74) (-1.452, -2.67) (-1.438, -2.6) (-1.425, -2.53) (-1.411, -2.46) (-1.398, -2.39) (-1.385, -2.32) (-1.371, -2.25) (-1.358, -2.19) (-1.344, -2.12) (-1.331, -2.05) (-1.318, -1.98) (-1.304, -1.91) (-1.291, -1.85) (-1.278, -1.78) (-1.264, -1.71) (-1.251, -1.64) (-1.237, -1.57) (-1.224, -1.49) (-1.211, -1.42) (-1.197, -1.35) (-1.184, -1.27) (-1.171, -1.2) (-1.157, -1.12) (-1.144, -1.04) (-1.13, -0.966) (-1.117, -0.886) (-1.104, -0.805) (-1.09, -0.723) (-1.077, -0.64) (-1.064, -0.556) (-1.05, -0.471) (-1.037, -0.385) (-1.023, -0.299) (-1.01, -0.211) (-0.9967, -0.123) (-0.9833, -0.0353) (-0.9699, 0.0529) (-0.9565, 0.141) (-0.9431, 0.229) (-0.9298, 0.315) (-0.9164, 0.401) (-0.903, 0.486) (-0.8896, 0.569) (-0.8763, 0.649) (-0.8629, 0.727) (-0.8495, 0.802) (-0.8361, 0.873) (-0.8227, 0.94) (-0.8094, 1) (-0.796, 1.06) (-0.7826, 1.11) (-0.7692, 1.16) (-0.7559, 1.2) (-0.7425, 1.23) (-0.7291, 1.25) (-0.7157, 1.26) (-0.7023, 1.27) (-0.689, 1.27) (-0.6756, 1.25) (-0.6622, 1.23) (-0.6488, 1.2) (-0.6355, 1.16) (-0.6221, 1.1) (-0.6087, 1.04) (-0.5953, 0.967) (-0.5819, 0.885) (-0.5686, 0.794) (-0.5552, 0.694) (-0.5418, 0.588) (-0.5284, 0.474) (-0.5151, 0.356) (-0.5017, 0.233) (-0.4883, 0.108) (-0.4749, -0.0192) (-0.4615, -0.146) (-0.4482, -0.27) (-0.4348, -0.391) (-0.4214, -0.507) (-0.408, -0.616) (-0.3946, -0.716) (-0.3813, -0.805) (-0.3679, -0.882) (-0.3545, -0.946) (-0.3411, -0.996) (-0.3278, -1.03) (-0.3144, -1.05) (-0.301, -1.05) (-0.2876, -1.03) (-0.2742, -0.993) (-0.2609, -0.941) (-0.2475, -0.872) (-0.2341, -0.788) (-0.2207, -0.69) (-0.2074, -0.58) (-0.194, -0.458) (-0.1806, -0.329) (-0.1672, -0.193) (-0.1538, -0.0531) (-0.1405, 0.0875) (-0.1271, 0.226) (-0.1137, 0.361) (-0.1003, 0.488) (-0.08696, 0.606) (-0.07358, 0.712) (-0.0602, 0.804) (-0.04682, 0.88) (-0.03344, 0.938) (-0.02007, 0.978) (-0.006689, 0.997) (0.006689, 0.997) (0.02007, 0.978) (0.03344, 0.938) (0.04682, 0.88) (0.0602, 0.804) (0.07358, 0.712) (0.08696, 0.606) (0.1003, 0.488) (0.1137, 0.361) (0.1271, 0.226) (0.1405, 0.0875) (0.1538, -0.0531) (0.1672, -0.193) (0.1806, -0.329) (0.194, -0.458) (0.2074, -0.58) (0.2207, -0.69) (0.2341, -0.788) (0.2475, -0.872) (0.2609, -0.941) (0.2742, -0.993) (0.2876, -1.03) (0.301, -1.05) (0.3144, -1.05) (0.3278, -1.03) (0.3411, -0.996) (0.3545, -0.946) (0.3679, -0.882) (0.3813, -0.805) (0.3946, -0.716) (0.408, -0.616) (0.4214, -0.507) (0.4348, -0.391) (0.4482, -0.27) (0.4615, -0.146) (0.4749, -0.0192) (0.4883, 0.108) (0.5017, 0.233) (0.5151, 0.356) (0.5284, 0.474) (0.5418, 0.588) (0.5552, 0.694) (0.5686, 0.794) (0.5819, 0.885) (0.5953, 0.967) (0.6087, 1.04) (0.6221, 1.1) (0.6355, 1.16) (0.6488, 1.2) (0.6622, 1.23) (0.6756, 1.25) (0.689, 1.27) (0.7023, 1.27) (0.7157, 1.26) (0.7291, 1.25) (0.7425, 1.23) (0.7559, 1.2) (0.7692, 1.16) (0.7826, 1.11) (0.796, 1.06) (0.8094, 1) (0.8227, 0.94) (0.8361, 0.873) (0.8495, 0.802) (0.8629, 0.727) (0.8763, 0.649) (0.8896, 0.569) (0.903, 0.486) (0.9164, 0.401) (0.9298, 0.315) (0.9431, 0.229) (0.9565, 0.141) (0.9699, 0.0529) (0.9833, -0.0353) (0.9967, -0.123) (1.01, -0.211) (1.023, -0.299) (1.037, -0.385) (1.05, -0.471) (1.064, -0.556) (1.077, -0.64) (1.09, -0.723) (1.104, -0.805) (1.117, -0.886) (1.13, -0.966) (1.144, -1.04) (1.157, -1.12) (1.171, -1.2) (1.184, -1.27) (1.197, -1.35) (1.211, -1.42) (1.224, -1.49) (1.237, -1.57) (1.251, -1.64) (1.264, -1.71) (1.278, -1.78) (1.291, -1.85) (1.304, -1.91) (1.318, -1.98) (1.331, -2.05) (1.344, -2.12) (1.358, -2.19) (1.371, -2.25) (1.385, -2.32) (1.398, -2.39) (1.411, -2.46) (1.425, -2.53) (1.438, -2.6) (1.452, -2.67) (1.465, -2.74) (1.478, -2.81) (1.492, -2.88) (1.505, -2.95) (1.518, -3.03) (1.532, -3.1) (1.545, -3.18) (1.559, -3.26) (1.572, -3.34) (1.585, -3.42) (1.599, -3.5) (1.612, -3.58) (1.625, -3.67) (1.639, -3.76) (1.652, -3.85) (1.666, -3.94) (1.679, -4.04) (1.692, -4.14) (1.706, -4.24) (1.719, -4.34) (1.732, -4.44) (1.746, -4.55) (1.759, -4.66) (1.773, -4.78) (1.786, -4.9) (1.799, -5.02) (1.813, -5.15) (1.826, -5.28) (1.839, -5.41) (1.853, -5.55) (1.866, -5.69) (1.88, -5.83) (1.893, -5.98) (1.906, -6.14) (1.92, -6.3) (1.933, -6.47) (1.946, -6.64) (1.96, -6.82) (1.973, -7) (1.987, -7.19) (2, -7.38) 
	};
	\addlegendentry{$k=3$};   
	\end{axis}
	\end{tikzpicture}
\end{minipage}
\begin{minipage}{.5\linewidth}
	\centering
\begin{tikzpicture}[scale=0.75]
	\pgfmathsetmacro\MU{2}
	\begin{axis}[samples=500,  
	xmin=-3, xmax=3, ymin=-4.5, ymax=4.5,
	title = {$\mathrm{Im}(\varphi_{k}(y)) = \mathrm{e}^{\frac{y^2}{2}} \, \sin\left(\pi k \,\mathrm{erf}(y)\right)$},
	axis x line=center, axis y line=center,
	every axis plot/.append style={thick},
	legend style={at={(0.5,1.25)}, anchor=south,legend columns=4,legend cell align=left, font=\small}
	]
\addplot[black] coordinates { 
		(-3, -0.00625) (-2.98, -0.00667) (-2.96, -0.00713) (-2.94, -0.00761) (-2.92, -0.00812) (-2.9, -0.00866) (-2.88, -0.00924) (-2.86, -0.00985) (-2.839, -0.0105) (-2.819, -0.0112) (-2.799, -0.0119) (-2.779, -0.0127) (-2.759, -0.0135) (-2.739, -0.0143) (-2.719, -0.0152) (-2.699, -0.0162) (-2.679, -0.0172) (-2.659, -0.0183) (-2.639, -0.0194) (-2.619, -0.0206) (-2.599, -0.0219) (-2.579, -0.0232) (-2.559, -0.0246) (-2.538, -0.0261) (-2.518, -0.0276) (-2.498, -0.0292) (-2.478, -0.031) (-2.458, -0.0328) (-2.438, -0.0347) (-2.418, -0.0367) (-2.398, -0.0387) (-2.378, -0.0409) (-2.358, -0.0433) (-2.338, -0.0457) (-2.318, -0.0482) (-2.298, -0.0509) (-2.278, -0.0537) (-2.258, -0.0566) (-2.237, -0.0597) (-2.217, -0.0629) (-2.197, -0.0663) (-2.177, -0.0698) (-2.157, -0.0735) (-2.137, -0.0773) (-2.117, -0.0813) (-2.097, -0.0855) (-2.077, -0.0899) (-2.057, -0.0945) (-2.037, -0.0993) (-2.017, -0.104) (-1.997, -0.109) (-1.977, -0.115) (-1.957, -0.121) (-1.936, -0.126) (-1.916, -0.133) (-1.896, -0.139) (-1.876, -0.146) (-1.856, -0.152) (-1.836, -0.16) (-1.816, -0.167) (-1.796, -0.175) (-1.776, -0.183) (-1.756, -0.191) (-1.736, -0.2) (-1.716, -0.209) (-1.696, -0.218) (-1.676, -0.228) (-1.656, -0.238) (-1.635, -0.248) (-1.615, -0.259) (-1.595, -0.27) (-1.575, -0.281) (-1.555, -0.293) (-1.535, -0.305) (-1.515, -0.318) (-1.495, -0.331) (-1.475, -0.344) (-1.455, -0.358) (-1.435, -0.372) (-1.415, -0.387) (-1.395, -0.402) (-1.375, -0.418) (-1.355, -0.434) (-1.334, -0.45) (-1.314, -0.467) (-1.294, -0.484) (-1.274, -0.502) (-1.254, -0.52) (-1.234, -0.539) (-1.214, -0.558) (-1.194, -0.577) (-1.174, -0.597) (-1.154, -0.617) (-1.134, -0.638) (-1.114, -0.659) (-1.094, -0.68) (-1.074, -0.701) (-1.054, -0.723) (-1.033, -0.745) (-1.013, -0.767) (-0.9933, -0.789) (-0.9732, -0.812) (-0.9532, -0.834) (-0.9331, -0.856) (-0.913, -0.879) (-0.893, -0.901) (-0.8729, -0.922) (-0.8528, -0.944) (-0.8328, -0.965) (-0.8127, -0.985) (-0.7926, -1) (-0.7726, -1.02) (-0.7525, -1.04) (-0.7324, -1.06) (-0.7124, -1.07) (-0.6923, -1.09) (-0.6722, -1.1) (-0.6522, -1.11) (-0.6321, -1.12) (-0.612, -1.13) (-0.592, -1.14) (-0.5719, -1.14) (-0.5518, -1.14) (-0.5318, -1.14) (-0.5117, -1.13) (-0.4916, -1.13) (-0.4716, -1.12) (-0.4515, -1.1) (-0.4314, -1.09) (-0.4114, -1.07) (-0.3913, -1.05) (-0.3712, -1.02) (-0.3512, -0.99) (-0.3311, -0.956) (-0.311, -0.92) (-0.291, -0.88) (-0.2709, -0.836) (-0.2508, -0.789) (-0.2308, -0.739) (-0.2107, -0.686) (-0.1906, -0.631) (-0.1706, -0.572) (-0.1505, -0.511) (-0.1304, -0.448) (-0.1104, -0.382) (-0.0903, -0.315) (-0.07023, -0.247) (-0.05017, -0.177) (-0.0301, -0.107) (-0.01003, -0.0356) (0.01003, 0.0356) (0.0301, 0.107) (0.05017, 0.177) (0.07023, 0.247) (0.0903, 0.315) (0.1104, 0.382) (0.1304, 0.448) (0.1505, 0.511) (0.1706, 0.572) (0.1906, 0.631) (0.2107, 0.686) (0.2308, 0.739) (0.2508, 0.789) (0.2709, 0.836) (0.291, 0.88) (0.311, 0.92) (0.3311, 0.956) (0.3512, 0.99) (0.3712, 1.02) (0.3913, 1.05) (0.4114, 1.07) (0.4314, 1.09) (0.4515, 1.1) (0.4716, 1.12) (0.4916, 1.13) (0.5117, 1.13) (0.5318, 1.14) (0.5518, 1.14) (0.5719, 1.14) (0.592, 1.14) (0.612, 1.13) (0.6321, 1.12) (0.6522, 1.11) (0.6722, 1.1) (0.6923, 1.09) (0.7124, 1.07) (0.7324, 1.06) (0.7525, 1.04) (0.7726, 1.02) (0.7926, 1) (0.8127, 0.985) (0.8328, 0.965) (0.8528, 0.944) (0.8729, 0.922) (0.893, 0.901) (0.913, 0.879) (0.9331, 0.856) (0.9532, 0.834) (0.9732, 0.812) (0.9933, 0.789) (1.013, 0.767) (1.033, 0.745) (1.054, 0.723) (1.074, 0.701) (1.094, 0.68) (1.114, 0.659) (1.134, 0.638) (1.154, 0.617) (1.174, 0.597) (1.194, 0.577) (1.214, 0.558) (1.234, 0.539) (1.254, 0.52) (1.274, 0.502) (1.294, 0.484) (1.314, 0.467) (1.334, 0.45) (1.355, 0.434) (1.375, 0.418) (1.395, 0.402) (1.415, 0.387) (1.435, 0.372) (1.455, 0.358) (1.475, 0.344) (1.495, 0.331) (1.515, 0.318) (1.535, 0.305) (1.555, 0.293) (1.575, 0.281) (1.595, 0.27) (1.615, 0.259) (1.635, 0.248) (1.656, 0.238) (1.676, 0.228) (1.696, 0.218) (1.716, 0.209) (1.736, 0.2) (1.756, 0.191) (1.776, 0.183) (1.796, 0.175) (1.816, 0.167) (1.836, 0.16) (1.856, 0.152) (1.876, 0.146) (1.896, 0.139) (1.916, 0.133) (1.936, 0.126) (1.957, 0.121) (1.977, 0.115) (1.997, 0.109) (2.017, 0.104) (2.037, 0.0993) (2.057, 0.0945) (2.077, 0.0899) (2.097, 0.0855) (2.117, 0.0813) (2.137, 0.0773) (2.157, 0.0735) (2.177, 0.0698) (2.197, 0.0663) (2.217, 0.0629) (2.237, 0.0597) (2.258, 0.0566) (2.278, 0.0537) (2.298, 0.0509) (2.318, 0.0482) (2.338, 0.0457) (2.358, 0.0433) (2.378, 0.0409) (2.398, 0.0387) (2.418, 0.0367) (2.438, 0.0347) (2.458, 0.0328) (2.478, 0.031) (2.498, 0.0292) (2.518, 0.0276) (2.538, 0.0261) (2.559, 0.0246) (2.579, 0.0232) (2.599, 0.0219) (2.619, 0.0206) (2.639, 0.0194) (2.659, 0.0183) (2.679, 0.0172) (2.699, 0.0162) (2.719, 0.0152) (2.739, 0.0143) (2.759, 0.0135) (2.779, 0.0127) (2.799, 0.0119) (2.819, 0.0112) (2.839, 0.0105) (2.86, 0.00985) (2.88, 0.00924) (2.9, 0.00866) (2.92, 0.00812) (2.94, 0.00761) (2.96, 0.00713) (2.98, 0.00667) (3, 0.00625)
	};
	\addlegendentry{$k=1$};
	\addplot[dashed] coordinates { 
		(-3, 0.0125) (-2.98, 0.0133) (-2.96, 0.0143) (-2.94, 0.0152) (-2.92, 0.0162) (-2.9, 0.0173) (-2.88, 0.0185) (-2.86, 0.0197) (-2.839, 0.021) (-2.819, 0.0224) (-2.799, 0.0238) (-2.779, 0.0253) (-2.759, 0.027) (-2.739, 0.0287) (-2.719, 0.0305) (-2.699, 0.0324) (-2.679, 0.0344) (-2.659, 0.0366) (-2.639, 0.0388) (-2.619, 0.0412) (-2.599, 0.0437) (-2.579, 0.0464) (-2.559, 0.0492) (-2.538, 0.0521) (-2.518, 0.0552) (-2.498, 0.0585) (-2.478, 0.0619) (-2.458, 0.0655) (-2.438, 0.0693) (-2.418, 0.0733) (-2.398, 0.0775) (-2.378, 0.0819) (-2.358, 0.0865) (-2.338, 0.0914) (-2.318, 0.0965) (-2.298, 0.102) (-2.278, 0.107) (-2.258, 0.113) (-2.237, 0.119) (-2.217, 0.126) (-2.197, 0.133) (-2.177, 0.14) (-2.157, 0.147) (-2.137, 0.155) (-2.117, 0.163) (-2.097, 0.171) (-2.077, 0.18) (-2.057, 0.189) (-2.037, 0.199) (-2.017, 0.209) (-1.997, 0.219) (-1.977, 0.23) (-1.957, 0.241) (-1.936, 0.253) (-1.916, 0.265) (-1.896, 0.278) (-1.876, 0.291) (-1.856, 0.305) (-1.836, 0.319) (-1.816, 0.334) (-1.796, 0.349) (-1.776, 0.365) (-1.756, 0.382) (-1.736, 0.399) (-1.716, 0.417) (-1.696, 0.435) (-1.676, 0.454) (-1.656, 0.474) (-1.635, 0.495) (-1.615, 0.516) (-1.595, 0.538) (-1.575, 0.56) (-1.555, 0.583) (-1.535, 0.607) (-1.515, 0.632) (-1.495, 0.657) (-1.475, 0.684) (-1.455, 0.71) (-1.435, 0.738) (-1.415, 0.766) (-1.395, 0.795) (-1.375, 0.824) (-1.355, 0.854) (-1.334, 0.884) (-1.314, 0.915) (-1.294, 0.947) (-1.274, 0.979) (-1.254, 1.01) (-1.234, 1.04) (-1.214, 1.07) (-1.194, 1.11) (-1.174, 1.14) (-1.154, 1.17) (-1.134, 1.2) (-1.114, 1.23) (-1.094, 1.26) (-1.074, 1.29) (-1.054, 1.32) (-1.033, 1.34) (-1.013, 1.36) (-0.9933, 1.38) (-0.9732, 1.4) (-0.9532, 1.42) (-0.9331, 1.43) (-0.913, 1.43) (-0.893, 1.43) (-0.8729, 1.43) (-0.8528, 1.42) (-0.8328, 1.41) (-0.8127, 1.39) (-0.7926, 1.36) (-0.7726, 1.33) (-0.7525, 1.29) (-0.7324, 1.24) (-0.7124, 1.19) (-0.6923, 1.12) (-0.6722, 1.05) (-0.6522, 0.971) (-0.6321, 0.883) (-0.612, 0.788) (-0.592, 0.685) (-0.5719, 0.576) (-0.5518, 0.461) (-0.5318, 0.342) (-0.5117, 0.219) (-0.4916, 0.0929) (-0.4716, -0.0339) (-0.4515, -0.16) (-0.4314, -0.285) (-0.4114, -0.405) (-0.3913, -0.52) (-0.3712, -0.627) (-0.3512, -0.725) (-0.3311, -0.812) (-0.311, -0.886) (-0.291, -0.946) (-0.2709, -0.99) (-0.2508, -1.02) (-0.2308, -1.03) (-0.2107, -1.02) (-0.1906, -0.99) (-0.1706, -0.945) (-0.1505, -0.882) (-0.1304, -0.802) (-0.1104, -0.707) (-0.0903, -0.598) (-0.07023, -0.478) (-0.05017, -0.348) (-0.0301, -0.212) (-0.01003, -0.0711) (0.01003, 0.0711) (0.0301, 0.212) (0.05017, 0.348) (0.07023, 0.478) (0.0903, 0.598) (0.1104, 0.707) (0.1304, 0.802) (0.1505, 0.882) (0.1706, 0.945) (0.1906, 0.99) (0.2107, 1.02) (0.2308, 1.03) (0.2508, 1.02) (0.2709, 0.99) (0.291, 0.946) (0.311, 0.886) (0.3311, 0.812) (0.3512, 0.725) (0.3712, 0.627) (0.3913, 0.52) (0.4114, 0.405) (0.4314, 0.285) (0.4515, 0.16) (0.4716, 0.0339) (0.4916, -0.0929) (0.5117, -0.219) (0.5318, -0.342) (0.5518, -0.461) (0.5719, -0.576) (0.592, -0.685) (0.612, -0.788) (0.6321, -0.883) (0.6522, -0.971) (0.6722, -1.05) (0.6923, -1.12) (0.7124, -1.19) (0.7324, -1.24) (0.7525, -1.29) (0.7726, -1.33) (0.7926, -1.36) (0.8127, -1.39) (0.8328, -1.41) (0.8528, -1.42) (0.8729, -1.43) (0.893, -1.43) (0.913, -1.43) (0.9331, -1.43) (0.9532, -1.42) (0.9732, -1.4) (0.9933, -1.38) (1.013, -1.36) (1.033, -1.34) (1.054, -1.32) (1.074, -1.29) (1.094, -1.26) (1.114, -1.23) (1.134, -1.2) (1.154, -1.17) (1.174, -1.14) (1.194, -1.11) (1.214, -1.07) (1.234, -1.04) (1.254, -1.01) (1.274, -0.979) (1.294, -0.947) (1.314, -0.915) (1.334, -0.884) (1.355, -0.854) (1.375, -0.824) (1.395, -0.795) (1.415, -0.766) (1.435, -0.738) (1.455, -0.71) (1.475, -0.684) (1.495, -0.657) (1.515, -0.632) (1.535, -0.607) (1.555, -0.583) (1.575, -0.56) (1.595, -0.538) (1.615, -0.516) (1.635, -0.495) (1.656, -0.474) (1.676, -0.454) (1.696, -0.435) (1.716, -0.417) (1.736, -0.399) (1.756, -0.382) (1.776, -0.365) (1.796, -0.349) (1.816, -0.334) (1.836, -0.319) (1.856, -0.305) (1.876, -0.291) (1.896, -0.278) (1.916, -0.265) (1.936, -0.253) (1.957, -0.241) (1.977, -0.23) (1.997, -0.219) (2.017, -0.209) (2.037, -0.199) (2.057, -0.189) (2.077, -0.18) (2.097, -0.171) (2.117, -0.163) (2.137, -0.155) (2.157, -0.147) (2.177, -0.14) (2.197, -0.133) (2.217, -0.126) (2.237, -0.119) (2.258, -0.113) (2.278, -0.107) (2.298, -0.102) (2.318, -0.0965) (2.338, -0.0914) (2.358, -0.0865) (2.378, -0.0819) (2.398, -0.0775) (2.418, -0.0733) (2.438, -0.0693) (2.458, -0.0655) (2.478, -0.0619) (2.498, -0.0585) (2.518, -0.0552) (2.538, -0.0521) (2.559, -0.0492) (2.579, -0.0464) (2.599, -0.0437) (2.619, -0.0412) (2.639, -0.0388) (2.659, -0.0366) (2.679, -0.0344) (2.699, -0.0324) (2.719, -0.0305) (2.739, -0.0287) (2.759, -0.027) (2.779, -0.0253) (2.799, -0.0238) (2.819, -0.0224) (2.839, -0.021) (2.86, -0.0197) (2.88, -0.0185) (2.9, -0.0173) (2.92, -0.0162) (2.94, -0.0152) (2.96, -0.0143) (2.98, -0.0133) (3, -0.0125) 
	}; 
	\addlegendentry{$k=2$};
	\addplot[dotted] coordinates { 
		(-3, -0.0187) (-2.98, -0.02) (-2.96, -0.0214) (-2.94, -0.0228) (-2.92, -0.0244) (-2.9, -0.026) (-2.88, -0.0277) (-2.86, -0.0295) (-2.839, -0.0315) (-2.819, -0.0335) (-2.799, -0.0357) (-2.779, -0.038) (-2.759, -0.0404) (-2.739, -0.043) (-2.719, -0.0457) (-2.699, -0.0486) (-2.679, -0.0517) (-2.659, -0.0549) (-2.639, -0.0583) (-2.619, -0.0618) (-2.599, -0.0656) (-2.579, -0.0696) (-2.559, -0.0738) (-2.538, -0.0782) (-2.518, -0.0828) (-2.498, -0.0877) (-2.478, -0.0929) (-2.458, -0.0983) (-2.438, -0.104) (-2.418, -0.11) (-2.398, -0.116) (-2.378, -0.123) (-2.358, -0.13) (-2.338, -0.137) (-2.318, -0.145) (-2.298, -0.153) (-2.278, -0.161) (-2.258, -0.17) (-2.237, -0.179) (-2.217, -0.189) (-2.197, -0.199) (-2.177, -0.209) (-2.157, -0.22) (-2.137, -0.232) (-2.117, -0.244) (-2.097, -0.257) (-2.077, -0.27) (-2.057, -0.283) (-2.037, -0.298) (-2.017, -0.313) (-1.997, -0.328) (-1.977, -0.345) (-1.957, -0.361) (-1.936, -0.379) (-1.916, -0.397) (-1.896, -0.416) (-1.876, -0.436) (-1.856, -0.457) (-1.836, -0.478) (-1.816, -0.5) (-1.796, -0.523) (-1.776, -0.547) (-1.756, -0.572) (-1.736, -0.598) (-1.716, -0.624) (-1.696, -0.652) (-1.676, -0.68) (-1.656, -0.709) (-1.635, -0.739) (-1.615, -0.771) (-1.595, -0.803) (-1.575, -0.836) (-1.555, -0.87) (-1.535, -0.904) (-1.515, -0.94) (-1.495, -0.977) (-1.475, -1.01) (-1.455, -1.05) (-1.435, -1.09) (-1.415, -1.13) (-1.395, -1.17) (-1.375, -1.21) (-1.355, -1.25) (-1.334, -1.29) (-1.314, -1.33) (-1.294, -1.37) (-1.274, -1.41) (-1.254, -1.44) (-1.234, -1.48) (-1.214, -1.51) (-1.194, -1.55) (-1.174, -1.58) (-1.154, -1.6) (-1.134, -1.63) (-1.114, -1.65) (-1.094, -1.66) (-1.074, -1.67) (-1.054, -1.67) (-1.033, -1.67) (-1.013, -1.65) (-0.9933, -1.63) (-0.9732, -1.61) (-0.9532, -1.57) (-0.9331, -1.52) (-0.913, -1.46) (-0.893, -1.39) (-0.8729, -1.3) (-0.8528, -1.21) (-0.8328, -1.1) (-0.8127, -0.98) (-0.7926, -0.85) (-0.7726, -0.709) (-0.7525, -0.559) (-0.7324, -0.401) (-0.7124, -0.237) (-0.6923, -0.0693) (-0.6722, 0.0995) (-0.6522, 0.266) (-0.6321, 0.428) (-0.612, 0.582) (-0.592, 0.723) (-0.5719, 0.848) (-0.5518, 0.954) (-0.5318, 1.04) (-0.5117, 1.09) (-0.4916, 1.12) (-0.4716, 1.12) (-0.4515, 1.08) (-0.4314, 1.01) (-0.4114, 0.915) (-0.3913, 0.787) (-0.3712, 0.633) (-0.3512, 0.458) (-0.3311, 0.267) (-0.311, 0.0656) (-0.291, -0.138) (-0.2709, -0.336) (-0.2508, -0.521) (-0.2308, -0.685) (-0.2107, -0.822) (-0.1906, -0.925) (-0.1706, -0.989) (-0.1505, -1.01) (-0.1304, -0.99) (-0.1104, -0.926) (-0.0903, -0.821) (-0.07023, -0.68) (-0.05017, -0.509) (-0.0301, -0.315) (-0.01003, -0.107) (0.01003, 0.107) (0.0301, 0.315) (0.05017, 0.509) (0.07023, 0.68) (0.0903, 0.821) (0.1104, 0.926) (0.1304, 0.99) (0.1505, 1.01) (0.1706, 0.989) (0.1906, 0.925) (0.2107, 0.822) (0.2308, 0.685) (0.2508, 0.521) (0.2709, 0.336) (0.291, 0.138) (0.311, -0.0656) (0.3311, -0.267) (0.3512, -0.458) (0.3712, -0.633) (0.3913, -0.787) (0.4114, -0.915) (0.4314, -1.01) (0.4515, -1.08) (0.4716, -1.12) (0.4916, -1.12) (0.5117, -1.09) (0.5318, -1.04) (0.5518, -0.954) (0.5719, -0.848) (0.592, -0.723) (0.612, -0.582) (0.6321, -0.428) (0.6522, -0.266) (0.6722, -0.0995) (0.6923, 0.0693) (0.7124, 0.237) (0.7324, 0.401) (0.7525, 0.559) (0.7726, 0.709) (0.7926, 0.85) (0.8127, 0.98) (0.8328, 1.1) (0.8528, 1.21) (0.8729, 1.3) (0.893, 1.39) (0.913, 1.46) (0.9331, 1.52) (0.9532, 1.57) (0.9732, 1.61) (0.9933, 1.63) (1.013, 1.65) (1.033, 1.67) (1.054, 1.67) (1.074, 1.67) (1.094, 1.66) (1.114, 1.65) (1.134, 1.63) (1.154, 1.6) (1.174, 1.58) (1.194, 1.55) (1.214, 1.51) (1.234, 1.48) (1.254, 1.44) (1.274, 1.41) (1.294, 1.37) (1.314, 1.33) (1.334, 1.29) (1.355, 1.25) (1.375, 1.21) (1.395, 1.17) (1.415, 1.13) (1.435, 1.09) (1.455, 1.05) (1.475, 1.01) (1.495, 0.977) (1.515, 0.94) (1.535, 0.904) (1.555, 0.87) (1.575, 0.836) (1.595, 0.803) (1.615, 0.771) (1.635, 0.739) (1.656, 0.709) (1.676, 0.68) (1.696, 0.652) (1.716, 0.624) (1.736, 0.598) (1.756, 0.572) (1.776, 0.547) (1.796, 0.523) (1.816, 0.5) (1.836, 0.478) (1.856, 0.457) (1.876, 0.436) (1.896, 0.416) (1.916, 0.397) (1.936, 0.379) (1.957, 0.361) (1.977, 0.345) (1.997, 0.328) (2.017, 0.313) (2.037, 0.298) (2.057, 0.283) (2.077, 0.27) (2.097, 0.257) (2.117, 0.244) (2.137, 0.232) (2.157, 0.22) (2.177, 0.209) (2.197, 0.199) (2.217, 0.189) (2.237, 0.179) (2.258, 0.17) (2.278, 0.161) (2.298, 0.153) (2.318, 0.145) (2.338, 0.137) (2.358, 0.13) (2.378, 0.123) (2.398, 0.116) (2.418, 0.11) (2.438, 0.104) (2.458, 0.0983) (2.478, 0.0929) (2.498, 0.0877) (2.518, 0.0828) (2.538, 0.0782) (2.559, 0.0738) (2.579, 0.0696) (2.599, 0.0656) (2.619, 0.0618) (2.639, 0.0583) (2.659, 0.0549) (2.679, 0.0517) (2.699, 0.0486) (2.719, 0.0457) (2.739, 0.043) (2.759, 0.0404) (2.779, 0.038) (2.799, 0.0357) (2.819, 0.0335) (2.839, 0.0315) (2.86, 0.0295) (2.88, 0.0277) (2.9, 0.026) (2.92, 0.0244) (2.94, 0.0228) (2.96, 0.0214) (2.98, 0.02) (3, 0.0187)
	};
	\addlegendentry{$k=3$};  
	\end{axis}
	\end{tikzpicture}
	\end{minipage}
	\caption{Real and imaginary part of the weighted exponential functions $\varphi_{k}, k=0,1,2,3$ in \eqref{eq:transformed_basis_functions} with the density function $\varrho$ of the error function transformation \eqref{eq:error_function_trafo} and the parameterized Gaussian weight function $\omega(y,\mu)$ as given in \eqref{eq:gaussian_weight_param} for fixed $\mu = \sqrt{2}$.}
	\label{fig:phi_k_example_with_erf_trafo}
\end{figure}

\subsection{Smoothness properties of composed functions in Sobolev spaces}
In this section we discuss the smoothness of univariate functions $h$ defined on $\mathbb{R}$ and of their resulting transformed versions $f$ on the torus $\mathbb{T}$.
In \cite{KPPW18} the authors used change of variables for integration problems with respect to a family of integrands with bounded $L_{p}$-norm of mixed first order partial derivatives with $1\leq p\leq \infty$ and provided sufficient conditions such that the transformed integrand belongs to a Sobolev space of mixed smootheness order one.
We will propose a specific set of sufficient conditions for $\psi$ and $\omega$ such that $f\in H_{\mathrm{mix}}^{m}(\mathbb{T}^d)$ with $m\in\mathbb{N}_{0}$.
These conditions are stated for both univariate and multivariate functions.
Afterwards we utilize the norm equivalence of the Sobolev space $H_{\mathrm{mix}}^{m}(\mathbb{T}^d)$ and the subspace $\mathcal{H}^{\beta}(\mathbb T^d)$ of the Wiener Algebra $\mathcal{A}(\mathbb{T}^d)$ for $m = \beta$ as described in \eqref{eq:Hs_norm_equivalence}
and combine it with the embedding $\mathcal{H}^{\beta+\lambda}(\mathbb{T}^d) \hookrightarrow \mathcal{A}^{\beta}(\mathbb T^d)$ in \eqref{eq:Wiener_algebra_inclusion} for all $\lambda > \frac{1}{2}$
in order to discuss high dimensional approximation problems in which we apply {rank-$1$} lattice based fast Fourier approximation methods. Throughout this section we still omit the parameters $\bm\eta,\bm\mu\in\mathbb{R}^d$ in the notation of the transformations $\psi$ and the weight functions $\omega$, as outlined in \eqref{eq:psi_weighted_mult} and \eqref{eq:omega_weighted_mult}.

For now we consider univariate transformed functions $f\in L_{2}(\mathbb{T})$ of the form 
\begin{align}\label{eq:f_is_transformed_h}
	f(x) 
	:= h(\psi(x)) \, \sqrt{ \omega(\psi(x)) \, \psi'(x) },
	\quad x\in\mathbb{T},
\end{align}
that are the result of applying the change of variables $y = \psi(x)$ as defined in \eqref{def:Trafo_def} to a function $h\in L_2(\mathbb{R},\omega)$ and for which we have the identity
\begin{align} \label{eq:transformed_function_natural}
	\|h\|_{L_{2}(\mathbb{R},\omega)}^2
	= \int_{\mathbb{R}} |h(y)|^2 \, \omega(y) \,\mathrm{d}y
	= \int_{\mathbb{T}} \left| h(\psi(x)) \right|^2 \, \omega(\psi(x)) \, \psi'(x) \, \mathrm{d}x
	= \|f\|_{L_{2}(\mathbb{T})}^2,
\end{align}
schematically shown in Figure~\ref{fig:DiagrammTrafo2}.
\begin{figure}[]
	\begin{center}
\begin{tikzpicture}[baseline= (a).base]
		\node[scale=1] (a) at (0,0)
		{
			\begin{tikzcd}[row sep=huge]
			\mathbb{T} \simeq [-\frac{1}{2},\frac{1}{2})  \arrow[r, phantom, "{\supset}"] &  (-\frac{1}{2},\frac{1}{2}) \arrow[shift left]{rr}{\psi(x)} \arrow{dr}[swap]{ L_{2}(\mathbb{T}) \,\ni\, h(\psi(x))\,\sqrt{\omega(\psi(x)) \psi'(x)} =: f(x) }
			& & \mathbb{R} \arrow{dl}{h(y) \,\in\, L_{2}\left(\mathbb{R}, \omega \right)} \arrow[shift left]{ll}{\psi^{-1}(y)} & 
			\phantom{\mathbb{T} \simeq [-\frac{1}{2},\frac{1}{2})  \arrow[l, phantom]} \\
			& & \mathbb{C} & &
			\end{tikzcd}
		};
		\end{tikzpicture}
	\end{center}
	\caption{Scheme of the relation between $f$ and $h$ caused by a transformation $\psi$.}
	\label{fig:DiagrammTrafo2}
\end{figure}
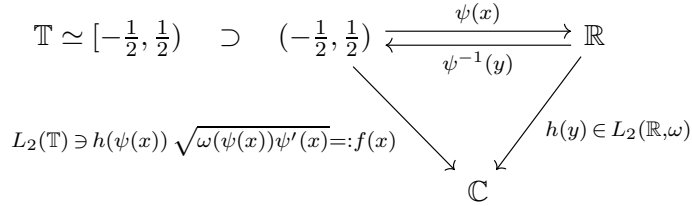
\begin{remark}
	Transformed functions $f$ as given in \eqref{eq:f_is_transformed_h} are generally not in $L_{2}(\mathbb{T})$ for all transformations $\psi$.
 	We will consider families of transformations $\psi(\circ,\eta), \eta\in\mathbb{R}$ as in \eqref{eq:psi_weighted_mult} and families of weight functions $\omega(\circ,\mu), \mu\in\mathbb{R}$ as in \eqref{eq:omega_weighted_mult}.
 	Generally, there are restrictions to the range of feasible parameters $\eta,\mu\in\mathbb{R}$ for which the transformed functions $f(\circ,\eta,\mu)$ as in \eqref{eq:f_is_transformed_h} are in $L_{2}(\mathbb{T})$.
 	Later on we present examples with multivariate functions ${h \in L_2(\mathbb{R}^d,\omega(\circ,\bm\mu))\cap H_{\mathrm{mix}}^{m}(\mathbb{R}^d)}$ and a fixed family of transformations $\psi(\circ,\bm\eta), \bm\eta\in\mathbb{R}$ and calculate parameter ranges of $\bm\eta$ and $\bm\mu$ for which the transformed functions $f$ are in $H_{\mathrm{mix}}^{m}(\mathbb{T}^d)$ for $m\in\mathbb{N}_0$.
\end{remark}

It is generally rather difficult to check if such transformed functions $f$ are in $H^{m}(\mathbb{T})$ for some fixed $m\in\mathbb{N}_0$ by calculating the individual $L_{2}(\mathbb{T})$-norms within the Sobolev norm $\|f\|_{H^{m}(\mathbb{T})}$.
Therefore we propose two different sets of sufficient conditions such that ${f\in H^{m}(\mathbb{T})}$ with $m\in\mathbb{N}_0$ by utilizing the product structure of the functions $f$ in \eqref{eq:f_is_transformed_h}.
At first we state conditions on $h\in L_{2}(\mathbb{R},\omega)$, the weight function $\omega$ and the transformation $\psi$ to preserve a certain degree of smoothness $m$ of $h$ under transformation with $\psi$, that slightly simplify the problem of the difficult evaluation of $L_{2}(\mathbb{T})$-integrals.

\begin{theorem} \label{thm:Hm_composition_criteria_simple}
	Let $m\in\mathbb{N}_0$, a transformation $\psi:(-\frac{1}{2},\frac{1}{2})\to\mathbb{R}$ as defined in \eqref{def:Trafo_def},
	and a function ${h\in L_2(\mathbb{R},\omega)}$ with a weight function $\omega:\mathbb{R}\to[0,\infty)$, 
	and the corresponding transformed functions $f$ of the form \eqref{eq:f_is_transformed_h} be given.
	We have $f \in H^{m}(\mathbb{T})$ 
	if either for all $k=0,1,\ldots,m$
	\begin{align*} \frac{\mathrm{d}^{k}}{\mathrm{d}x^{k}}\left[ h\circ\psi \right](x) \in L_{\infty}(\mathbb{T}) 
		\quad \text{and} \quad \frac{\mathrm{d}^{k}}{\mathrm{d}x^{k}}\left[ \sqrt{ (\omega\circ\psi) \, \psi' } \right](x) \in L_{2}(\mathbb{T})
	\end{align*}
	or if for all $k=0,1,\ldots,m$
	\begin{align*} \frac{\mathrm{d}^{k}}{\mathrm{d}x^{k}}\left[ h\circ\psi \right](x) \in L_{2}(\mathbb{T}) 
		\quad \text{and} \quad \frac{\mathrm{d}^{k}}{\mathrm{d}x^{k}}\left[ \sqrt{ (\omega\circ\psi) \, \psi' } \right](x) \in L_{\infty}(\mathbb{T}).
	\end{align*}
\end{theorem}
\begin{proof} 
	Let $f$ be of the form \eqref{eq:f_is_transformed_h} and $k=0,1,\ldots,m$.
	Using the well-known \emph{generalized Leibniz rule} for the $k$-th derivative of a product of two functions leads to
	\begin{align}
		\label{eq:Leibniz_estimate}
		\left\| \frac{\mathrm{d}^{k}}{\mathrm{d}x^{k}}[f](x) \right\|_{L_2(\mathbb{T})}
		&\leq \sum_{\ell=0}^{k} \binom{k}{\ell}
		\left\|
		\frac{\mathrm{d}^{\ell}}{\mathrm{d}x^{\ell}}\left[ h\circ\psi \right](x) 
		\frac{\mathrm{d}^{k-\ell}}{\mathrm{d}x^{k-\ell}}\left[ \sqrt{ (\omega\circ\psi) \, \psi' } \right](x) 
		\right\|_{L_{2}(\mathbb{T})}.
	\end{align}
	Now we either estimate
	\begin{align*}
		\left\| \frac{\mathrm{d}^{k}}{\mathrm{d}x^{k}}[f](x) \right\|_{L_2(\mathbb{T})}
		\leq \sum_{\ell=0}^{k} \binom{k}{\ell}
		\left\|\frac{\mathrm{d}^{\ell}}{\mathrm{d}x^{\ell}}\left[ h\circ\psi \right](x) \right\|_{L_{\infty}(\mathbb{T})}	
		\left\|\frac{\mathrm{d}^{k-\ell}}{\mathrm{d}x^{k-\ell}}\left[ \sqrt{ (\omega\circ\psi) \, \psi' } \right](x) \right\|_{L_{2}(\mathbb{T})}
	\end{align*}
	or
	\begin{align*}
		\left\| \frac{\mathrm{d}^{k}}{\mathrm{d}x^{k}}[f](x) \right\|_{L_2(\mathbb{T})}
		\leq \sum_{\ell=0}^{k} \binom{k}{\ell}
		\left\|\frac{\mathrm{d}^{\ell}}{\mathrm{d}x^{\ell}}\left[ h\circ\psi \right](x) \right\|_{L_{2}(\mathbb{T})}	
		\left\|\frac{\mathrm{d}^{k-\ell}}{\mathrm{d}x^{k-\ell}}\left[ \sqrt{ (\omega\circ\psi) \, \psi' } \right](x) \right\|_{L_{\infty}(\mathbb{T})} .
	\end{align*}
	If the $\left\| \cdot \right\|_{L_2(\mathbb{T})}$-norms of $k$-th derivatives of $f$ are finite for all $k=0,1,\ldots,m$, then their sum, i.e., the $H^{m}(\mathbb{T})$-norm is finite, too.
\end{proof}

Now, we derive a set of sufficient $L_{\infty}$-conditions on $\psi$ and $\omega$ that ensure that a function $h\in L_{2}(\mathbb{R},\omega) \cap H^{m}(\mathbb{R})$ can be transformed by $\psi$ into an $f \in H^{m}(\mathbb{T})$ of form \eqref{eq:f_is_transformed_h}.
This eliminates the necessity to evaluate $L_{2}$-integrals of various derivatives of $f$.
Furthermore, once we consider particular parameterized families of transformations $\psi(\circ,\eta)$ and families of weight functions $\omega(\circ,\mu)$,
these conditions enable us for each smoothness order $m\in\mathbb{N}$ to explicitly calculate how large the parameters $\eta,\mu\in\mathbb{R}$ have to be in order to preserve the fixed degree of smoothness $m$ when transforming ${h\in L_{2}(\mathbb{R},\omega(\circ,\mu)) \cap H^{m}(\mathbb{R})}$ into $f\in H^{m}(\mathbb{T})$ via $\psi(\circ,\eta)$.

For simplified notation we alternate between equivalent expressions for derivatives of the appearing functions and for improved readability we write explicit arguments within certain norms. 
We denote the $k$-th derivative of a function $f(x)$ with respect to $x$ by either $\frac{\mathrm{d}^{k}}{\mathrm{d}x^{k}}[f](x)$ or $f^{(k)}(x)$, and for $k=1,2,3$ we sometimes use the notation $f'(x), f''(x)$, and $f'''(x)$.
\begin{theorem} \label{thm:Hm_composition_criteria}
	Let $m\in\mathbb{N}_0$, a transformation $\psi:(-\frac{1}{2},\frac{1}{2})\to\mathbb{R}$ as defined in \eqref{def:Trafo_def} with the density function $\varrho$ of $\psi$ as in \eqref{def:Varrho_def}, 
	a function ${h\in L_2(\mathbb{R},\omega)\cap H^{m}(\mathbb{R})}$ with a non-negative weight function $\omega:\mathbb{R}\to[0,\infty)$ 
	and the corresponding transformed functions $f$ of the form \eqref{eq:f_is_transformed_h} be given.

	If for all $\ell=0,1,\ldots,m$ we have 
	\begin{align} &{\frac{\mathrm{d}^{\ell}}{\mathrm{d}y^{\ell}}\left[ \varrho \right](y)\in\mathcal{C}_{0}(\mathbb{R})}, \qquad
		{\frac{\mathrm{d}^{\ell}}{\mathrm{d}x^{\ell}}\left[ \psi \right](x)\in\mathcal{C}((-1/2,1/2))} \nonumber \\
		\text{and} \quad 
		&\max_{k=0,\ldots,\ell} \left\| \frac{\mathrm{d}^{\ell-k}}{\mathrm{d}x^{\ell-k}}\left[ \sqrt{ (\omega\circ\psi) \, \psi' } \right](x)\, \psi'(x)^{\max(-\frac{1}{2},2k-\frac{3}{2})} \right\|_{L_{\infty}(\mathbb{T})} < \infty,
	\end{align}
	then we have $f \in H^{m}(\mathbb{T})$.
\end{theorem}
\begin{proof}
	For $h\in L_{2}(\mathbb{R},\omega) \cap H^{m}(\mathbb{R})$ with $m\in\mathbb{N}_{0}$ and a transformation $\psi$ as defined in $\eqref{def:Trafo_def}$
	we consider the function $f$ as given in \eqref{eq:f_is_transformed_h}.
	In order to prove that $f\in H^{m}(\mathbb{T})$ we have to show that $\left\| \frac{\mathrm{d}^{n}}{\mathrm{d}x^{n}}[f](x) \right\|_{L_2(\mathbb{T})} < \infty$ for all $n=0,1,\ldots,m$. 
	We present the arguments for $n=m$ and they are applicable in the same way for $n=0,1,\ldots,m-1$, too.
	We consider $\left\| \frac{\mathrm{d}^{m}}{\mathrm{d}x^{m}}[f](x) \right\|_{L_2(\mathbb{T})}$ and apply the \emph{generalized Leibniz rule} as in \eqref{eq:Leibniz_estimate}, so that we now have to ensure that
	\begin{align} \label{eq:Leibniz_single_term}
		\left\|
		\frac{\mathrm{d}^{k}}{\mathrm{d}x^{k}}\left[ h\circ\psi \right](x) 
		\frac{\mathrm{d}^{m-k}}{\mathrm{d}x^{m-k}}\left[ \sqrt{ (\omega\circ\psi) \, \psi' } \right](x) 
		\right\|_{L_{2}(\mathbb{T})} < \infty
	\end{align}
	for all $k=0,\ldots,m$.
	We leave $h\circ\psi$ in the term corresponding to $k=0$ untouched for now. 
	For $k = 1,\ldots,m$ we use the \emph{Fa\'{a} di Bruno} formula to write the $k$-th derivative of the composition of functions $h$ and $\psi$ as
	\begin{align} \label{eq:FdiB_formula}
		\frac{\mathrm{d}^{k}}{\mathrm{d}x^{k}}\left[ h\circ\psi \right](x)
		= \sum_{\ell = 1}^{k} h^{(\ell)}(\psi(x)) \cdot B_{k,\ell}(\psi'(x), \psi''(x),\ldots, \psi^{(k-\ell+1)}(x))
	\end{align}
	and the well-known Bell polynomials $B_{k,\ell}$ for $k,\ell\in\mathbb{N}_{0}$ are given by
	\begin{align}
		\label{def:Bell_polynomial}
		B_{k,\ell}(\mathbf z) := \sum_{ \substack{j_1+j_2+\ldots+j_{k-\ell+1}=\ell,\\ j_1+2j_2+\ldots+(k-\ell+1)j_{k-\ell+1} = k} }
		\frac{\ell!}{j_1! \cdot \ldots \cdot j_{k-\ell+1}!}
		\prod_{r=1}^{k-\ell+1} \left(\frac{z_{r}}{r!}\right)^{j_{r}}
	\end{align}
	with $\mathbf z = (z_1,\ldots,z_{k-\ell+1})^{\top}$. 
	By differentiating both sides of ${\psi^{-1}(\psi(x)) = x}$ we obtain
	\begin{align*}
		\psi'(x) = \frac{1}{\varrho(\psi(x))}
		,\qquad 
		\psi''(x) = - \frac{\varrho'(\psi(x)) \psi'(x)}{\varrho(\psi(x))^2} = - \varrho'(\psi(x)) \psi'(x)^3.
	\end{align*}
	Based on this we also observe that for $k\in\mathbb{N}$
	\begin{align} \label{eq:powers_of_psi_prime}
		\frac{\mathrm{d}}{\mathrm{d}x} \left[ (\psi')^{k} \right](x) 
		= k \psi'(x)^{k-1} \psi''(x) = - k \psi'(x)^{k+2} \varrho'(\psi(x)).
	\end{align}
	Hence, the $k$-th derivative of $\psi$ can be expressed soley in terms of powers of $\psi'$ and the first $(k-1)$ derivatives of $\varrho$ by repeated insertion of the expression of $\psi''$.
	Formula \eqref{eq:powers_of_psi_prime} implies that the highest appearing power of $\psi'$ increases by $2$ with each differentiation.
	For example, 
	\begin{align*}
		\psi'''(x) 
		&= \psi'(x)^5 \left( - \frac{\varrho''(\psi(x))}{\psi'(x)} + 3 \varrho'(\psi(x)) \right), \\
\psi^{(4)}(x) 
		&= \psi'(x)^7 \left( -\frac{\varrho'''(\psi(x))}{\psi'(x)^2} + \frac{4 \varrho''(\psi(x)) \varrho'(\psi(x)) + 6 \varrho'(\psi(x))}{\psi'(x)} - 15 \varrho'(\psi(x))^3\right).
	\end{align*}
	We note that each derivative of $\psi$ is bounded, based on the fact that $\varrho$ is by definition in $\mathcal{C}_{0}(\mathbb{R})$, hence $\varrho\circ\psi = 1/\psi' \in\mathcal{C}(\mathbb{T})$ and any power of $1/\psi'$ is bounded, too.
	Additionally we assumed that the first $k$ derivatives of $\varrho$ are in $\mathcal{C}_{0}(\mathbb{R})$, too.
	Therefore, with constants $C_k>0$ and $C>0$, for all $k\in\mathbb{N}$ we can estimate
	\begin{align*}
		\left| \frac{\mathrm{d}^k}{\mathrm{d}x^k} [\psi](x) \right| \leq C_k |\psi'(x)|^{2k-1}
	\end{align*}
	and for the Bell polynomials $B_{k,\ell}$ in \eqref{eq:FdiB_formula} we then estimate
	\begin{align} \label{eq:Bell_estimate1}
		&|B_{k,\ell}(\psi'(x), \psi''(x),\ldots, \psi^{(k-\ell+1)}(x))| \\
		&\leq C\cdot B_{k,\ell}(|\psi'(x)|, |\psi'(x)|^3,\ldots, |\psi'(x)|^{2(k-\ell+1)-1}). \nonumber
	\end{align}
	The Bell polynomials were defined according to the rules to partition a number $k\in\mathbb{N}$ into a sum of $\ell\in\{1,2,\ldots,k\}$ natural numbers $j_1,\ldots,j_\ell \in \mathbb{N}$, that are given by
	\begin{align*}
		j_1+j_2+j_3+\ldots+j_{k-\ell+1}=\ell, \\
		j_1+2j_2+3j_3+\ldots+(k-\ell+1)j_{k-\ell+1} = k.
	\end{align*}
	Substracting the first rule from two times the second rule results in the equation
	\begin{align*}
		j_1+3j_2+5j_3+\ldots+(2(k-\ell+1)-1)j_{k-\ell+1} = 2k-\ell
	\end{align*}
	which reveals that in the polynomials $B_{k,\ell}(|\psi'(x)|, |\psi'(x)|^3,\ldots, |\psi'(x)|^{2(k-\ell+1)-1})$ as defined in \eqref{eq:Bell_estimate1}
	the highest appearing power of $|\psi'|$ is ${2k-1}$ for $\ell=1$.
	By extracting $|\psi'(x)|^{2k-1}$ from each $B_{k,\ell}$ the remaining polynomials consist only of powers of $1/\psi'$, that are all bounded.
	Hence, in \eqref{eq:Bell_estimate1} we further estimate
	\begin{align} \label{eq:Bell_estimate2}
		&|B_{k,\ell}(\psi'(x), \psi''(x),\ldots, \psi^{(k-\ell+1)}(x))| \\
		&\leq C \left| \psi'(x)^{2k-1} \frac{B_{k,\ell}(|\psi'(x)|, |\psi'(x)|^3,\ldots, |\psi'(x)|^{2(k-\ell+1)-1})}{\psi'(x)^{2k-1}} \right| \nonumber\\
		&\leq C'  \left| \psi'(x)^{2k-1} \right| \nonumber
	\end{align}
	with constants $C,C'>0$.
	
	We go back to the derivatives of $h\circ\psi$ in \eqref{eq:FdiB_formula} and estimate them individually.	
	For $k=0$ we simply estimate 
	\begin{align} \label{eq:Leibniz_term_estimate1}
		\left\|h\circ \psi\right\|_{L_2(\mathbb{T})}
		&= \left(
		\int_{-\frac{1}{2}}^{\frac{1}{2}} \left| h(\psi(x)) \, \psi'(x)^{-\frac{1}{2}} \right|^{2} \psi'(x) \,\mathrm{d}x
		\right)^{\frac{1}{2}} \\
		&\leq \| \psi'(\circ)^{-\frac{1}{2}} \|_{L_{\infty}(\mathbb{T})}
		\left( \int_{-\infty}^{\infty} \left| h(y) \right|^{2} \,\mathrm{d}y \right)^{\frac{1}{2}}, \nonumber
	\end{align}
	which exists if $\psi'(\circ)^{-\frac{1}{2}}\in L_{\infty}(\mathbb{T})$.
	With the {Fa\'{a} di Bruno} formula \eqref{eq:FdiB_formula} and the upper bound \eqref{eq:Bell_estimate2} for $\frac{\mathrm{d}^{k}}{\mathrm{d}x^{k}}\left[ h\circ\psi \right](x)$ we estimate
	\begin{align} \label{eq:Leibniz_term_estimate2}
		\left\|\frac{\mathrm{d}^{k}}{\mathrm{d}x^{k}}\left[ h\circ\psi \right](x)\right\|_{L_2(\mathbb{T})}
		&= \left(
		\int_{-\frac{1}{2}}^{\frac{1}{2}} \left|   
		\frac{\mathrm{d}^{k}}{\mathrm{d}x^{k}}\left[ h\circ\psi \right](x) \, \psi'(x)^{-\frac{1}{2}} 
		\right|^{2} \psi'(x) \,\mathrm{d}x
		\right)^{\frac{1}{2}} \nonumber\\
		&\leq \| \psi'(x)^{2k-\frac{3}{2}} \|_{L_{\infty}(\mathbb{T})}
		\left(	\int_{-\frac{1}{2}}^{\frac{1}{2}} \left|    \sum_{j=1}^{k} h^{(j)}(\psi(x))	\right|^{2} \,\psi'(x) \,\mathrm{d}x	\right)^{\frac{1}{2}}  \nonumber\\
		&\leq C\cdot \| \psi'(x)^{2k-\frac{3}{2}} \|_{L_{\infty}(\mathbb{T})}
		\sum_{j=1}^{k} \left\| \frac{\mathrm{d}^{j}}{\mathrm{d}y^{j}} \left[h(y)\right] \right\|_{L_{2}(\mathbb{R})}.
	\end{align}
	By inserting \eqref{eq:Leibniz_term_estimate1} and \eqref{eq:Leibniz_term_estimate2} into \eqref{eq:Leibniz_single_term} we in total have for all $\ell\in\mathbb{N}_{0}$ the estimate
	\begin{align}\label{eq:proof_univar_summary}
		&\left\| \frac{\mathrm{d}^{\ell}}{\mathrm{d}x^{\ell}}[f](x) \right\|_{L_2(\mathbb{T})} \\
		&\leq \sum_{k=0}^{\ell} \binom{\ell}{k}
		\left\|
		\frac{\mathrm{d}^{k}}{\mathrm{d}x^{k}}\left[ h\circ\psi \right](x) 
		\frac{\mathrm{d}^{\ell-k}}{\mathrm{d}x^{\ell-k}}\left[ \sqrt{ (\omega\circ\psi) \, \psi' } \right](x) 
		\right\|_{L_{2}(\mathbb{T})} \nonumber\\ 
		&\leq \left\| \frac{\mathrm{d}^{\ell}}{\mathrm{d}x^{\ell}}\left[ \sqrt{ (\omega\circ\psi) \, \psi' } \right](x) \,\psi'(x)^{-\frac{1}{2}} \right\|_{L_{\infty}(\mathbb{T})} \left\| h \right\|_{L_{2}(\mathbb{R})} \nonumber\\
		&\quad + C \cdot 
		\sum_{k=1}^{\ell} \binom{\ell}{k}
		\left\| \frac{\mathrm{d}^{\ell-k}}{\mathrm{d}x^{\ell-k}}\left[ \sqrt{ (\omega\circ\psi) \, \psi' } \right](x) \,\psi'(x)^{2k-\frac{3}{2}} \right\|_{L_{\infty}(\mathbb{T})}
		\sum_{j=1}^{k} \left\| \frac{\mathrm{d}^{j}}{\mathrm{d}y^{j}}\left[ h(y) \right] \right\|_{L_{2}(\mathbb{R})} \nonumber\\
		&\leq C' \max_{k=0,\ldots,\ell} \left\| \frac{\mathrm{d}^{\ell-k}}{\mathrm{d}x^{\ell-k}}\left[ \sqrt{ (\omega\circ\psi) \, \psi' } \right](x)\, \psi'(x)^{\max(-\frac{1}{2},2k-\frac{3}{2})} \right\|_{L_{\infty}(\mathbb{T})} \times \nonumber\\
			&\quad\times\left( \left\| h \right\|_{L_{2}(\mathbb{R})} + \sum_{k=1}^{\ell} \|h\|_{H^{k}(\mathbb{R})} \right) \nonumber\\
		&\leq C' \max_{k=0,\ldots,\ell} \left\| \frac{\mathrm{d}^{\ell-k}}{\mathrm{d}x^{\ell-k}}\left[ \sqrt{ (\omega\circ\psi) \, \psi' } \right](x)\, \psi'(x)^{\max(-\frac{1}{2},2k-\frac{3}{2})} \right\|_{L_{\infty}(\mathbb{T})} (\ell+1)\|h\|_{H^{\ell}(\mathbb{R})} \nonumber
	\end{align}
	with constants $C, C' \geq 1$.
	This upper bound exists as long as the $L_{\infty}$-norms are finite and $h\in H^{m}(\mathbb{R})$.
	In total we finally estimate
	\begin{align*}
		&\|f\|_{H^{m}(\mathbb{T})} 
		= \left( \sum_{\ell=0}^{m} \left\| \frac{\mathrm{d}^{\ell}}{\mathrm{d}x^{\ell}}[f](x) \right\|_{L_2(\mathbb{T})}^2 \right)^{\frac{1}{2}} \\
		&\leq C \max_{\ell=0,\ldots,m} \left( \max_{k=0,\ldots,\ell} \left\| \frac{\mathrm{d}^{\ell-k}}{\mathrm{d}x^{\ell-k}}\left[ \sqrt{ (\omega\circ\psi) \, \psi' } \right](x)\, \psi'(x)^{\max(-\frac{1}{2},2k-\frac{3}{2})} \right\|_{L_{\infty}(\mathbb{T})} \right) \times\\
		&\quad\times \left( \sum_{\ell = 0}^{m} (\ell+1)^2\|h\|_{H^{\ell}(\mathbb{R})}^2 \right)^{\frac{1}{2}} \\
		&\leq C \max_{\ell=0,\ldots,m} \left( 
		\max_{k=0,\ldots,\ell} \left\| \frac{\mathrm{d}^{\ell-k}}{\mathrm{d}x^{\ell-k}}\left[ \sqrt{ (\omega\circ\psi) \, \psi' } \right](x)\, \psi'(x)^{\max(-\frac{1}{2},2k-\frac{3}{2})} \right\|_{L_{\infty}(\mathbb{T})} 
		\right) \times\nonumber\\
		&\quad\times(m+1)^{\frac{3}{2}} \|h\|_{H^{m}(\mathbb{R})}. 
	\end{align*}
\end{proof}

Next, we generalize the previous Theorem by proving its multivariate version.
Again, to simplify the notation in \eqref{def:differential_def_mult} of the $d$-variate differential operator $D^{\mathbf m}[f](\mathbf x)$ with both ${\mathbf m = (m_1,\ldots,m_d) \in\mathbb{N}_{0}^{d}}$ and ${\mathbf x = (x_1,\ldots,x_d)^{\top}\in\mathbb{R}^d}$
we use equivalent expressions for certain (partial) derivatives and state explicit arguments in various norms.
When differentiating a multivariate function $f$ with respect to the $j$-th coordinate $m_j$-times we write
\begin{align*}
	\partial^{m_j}[f](\mathbf x) := \frac{\partial^{m_j}}{\partial x_{j}^{m_j}}[f](\mathbf x) = D^{(0,\ldots,0,m_j,0,\ldots,0)}[f](\mathbf x).
\end{align*}
For the first and $\ell$-th derivatives of univariate functions with $\ell\in\mathbb{N}$ we use the notation 
\begin{align*}
	\psi_j'(x_j) := \frac{\mathrm{d}}{\mathrm{d}x_j}[\psi_j](x_j)
	\quad \text{and} \quad
	\psi_j^{(\ell)}(x_j) := \frac{\mathrm{d}^{\ell}}{\mathrm{d}x_{j}^{\ell}}[\psi_j](x_j).
\end{align*}
Similar to \eqref{eq:f_is_transformed_h} we consider multivariate transformed functions $f\in L_{2}(\mathbb{T}^d)$ of the form 
\begin{align}\label{eq:f_is_transformed_h_mult}
	f(\mathbf x) 
	&= (h\circ\psi)(\mathbf x) \sqrt{(\omega\circ\psi)(\mathbf x) D^{\mathbf 1}[\psi](\mathbf x)} \nonumber\\
	&= h(\psi_{1}(x_1),\ldots,\psi_{d}(x_d)) \prod_{k=1}^{d}\sqrt{\omega_k(\psi_k(x_k)) \psi'_{k}(x_k)},
	\quad \mathbf x\in\mathbb{T}^d,
\end{align}
that are the result of applying the multivariate change of variables 
\begin{align*}
	\mathbf y = (y_1,\ldots,y_d)^{\top} = (\psi_{1}(x_{1}),\ldots,\psi_{d}(x_{d}))^{\top} = \psi(\mathbf x)
\end{align*}
as defined in \eqref{eq:psi_weighted_mult} to a function $h\in L_2(\mathbb{R}^d,\omega)$ with a product weight $\omega$ as in \eqref{eq:omega_weighted_mult} and for which we have the identity
\begin{align} \label{eq:transformed_function_mult}
	\|h\|_{L_{2}(\mathbb{R}^d,\omega)}^2
	&= \int_{\mathbb{R}^d} |h(\mathbf y)|^2 \, \omega(\mathbf y) \,\mathrm{d}\mathbf y \\
	&= \int_{\mathbb{T}^d} |(h\circ\psi)(\mathbf x)|^2 \, (\omega\circ\psi)(\mathbf x) \, D^{\mathbf 1}[\psi](\mathbf x) \, \mathrm{d}\mathbf x
	= \|f\|_{L_{2}(\mathbb{T}^d)}^2. \nonumber
\end{align}
Again, we derive a set of sufficient $L_{\infty}$-conditions on the multivariate transformation $\psi$ and the product weight $\omega$, that determine when a function $h\in L_{2}(\mathbb{R}^d,\omega) \cap H_{\mathrm{mix}}^{m}(\mathbb{R}^d)$ can be transformed by $\psi$ into an $f \in H_{\mathrm{mix}}^{m}(\mathbb{T}^d)$ of form \eqref{eq:f_is_transformed_h_mult}.

\begin{theorem} \label{thm:Hm_composition_criteria_mult}
	Let the dimension $d\in\mathbb{N}$, $m\in\mathbb{N}_{0}$, a $d$-variate transformation ${\psi:(-\frac{1}{2},\frac{1}{2})^d\to\mathbb{R}^d}$ as defined in \eqref{eq:psi_weighted_mult} with the $d$-variate density function $\varrho(\mathbf y) = \prod_{j=1}^{d}\varrho_{j}(y_j)$ of $\psi$ as in \eqref{def:varrho_mult}, 
	a non-negative product weight function $\omega:\mathbb{R}^d\to[0,\infty)$ as in \eqref{eq:omega_weighted_mult},
	a multivariate function ${h\in L_2(\mathbb{R}^d,\omega)\cap H_{\mathrm{mix}}^{m}(\mathbb{R}^d)}$
	and the corresponding transformed functions $f$ of the form \eqref{eq:f_is_transformed_h_mult} be given.
	
	If for all multi-indices $\mathbf m = (m_1,\ldots,m_d)^{\top} \in\mathbb{N}_{0}^{d}, \|\mathbf m\|_{\ell_{\infty}} \leq m$ and all ${j_{\ell}=0,\ldots,m}$, ${\ell=1,\ldots,d}$ we have 
	\begin{align} \label{eq:Hm_composition_criteria_mult}
		&{\partial^{j_{\ell}}\left[ \varrho \right](y_{\ell})\in\mathcal{C}_{0}(\mathbb{R})}, \qquad
		{\partial^{j_{\ell}}\left[ \psi \right](x_{\ell})\in\mathcal{C}((-1/2,1/2))} \\
		\text{and} \quad 
		&\max_{j_\ell=0,\ldots,m_{\ell}} \left\| \partial^{m_{\ell}-j_{\ell}}\left[\sqrt{(\omega_{\ell}\circ\psi_{\ell})\psi'_{\ell}}\right](x_{\ell}) \, \psi_{\ell}'(x_{\ell})^{\max(-\frac{1}{2},2j_{\ell}-\frac{3}{2})} \right\|_{L_{\infty}(\mathbb{T})}
		<\infty, \nonumber
	\end{align}
	then we have ${f \in H_{\mathrm{mix}}^{m}(\mathbb{T}^d)}$.
\end{theorem}
\begin{proof}
	For $h\in L_{2}(\mathbb{R}^d,\omega) \cap H_{\mathrm{mix}}^{m}(\mathbb{R}^d)$ with $m\in\mathbb{N}_{0}$ and a transformation $\psi$ as defined in $\eqref{eq:psi_weighted_mult}$
	we consider the function $f$ as given in \eqref{eq:f_is_transformed_h_mult}.
	In order to prove that $f\in H_{\mathrm{mix}}^{m}(\mathbb{T}^d)$ we have to show that for all multi-indices $\mathbf m = (m_1,\ldots,m_d)^{\top}\in\mathbb{N}_{0}^{d}$ with $\|\mathbf m\|_{\ell_{\infty}} \leq m$ we have $\left\| D^{\mathbf m}[f](\mathbf x) \right\|_{L_2(\mathbb{T}^d)} < \infty$. 
	
	Let $\mathbf m = (m_1,\ldots,m_d)^{\top}\in\mathbb{N}_{0}^{d}$ be any multi-index with $\|\mathbf m\|_{\ell_{\infty}} \leq m$. 
	For a multivariate transformed function of the form \eqref{eq:f_is_transformed_h_mult} we have
	\begin{align}\label{eq:mult_proof_L2T}
		\left\| D^{\mathbf m}[f](\mathbf x) \right\|_{L_2(\mathbb{T}^d)} = \left( \int_{\mathbb T^d} \left| D^{\mathbf m}\left[(h\circ\psi)\prod_{k=1}^{d}\sqrt{(\omega_k\circ\psi_k) \psi'_{k}}\right](\mathbf x) 
		\right|^2 \mathrm{d}\mathbf x \right)^{\frac{1}{2}}. 
	\end{align}
	Based on the product weight function in the transformed function $f$ in \eqref{eq:f_is_transformed_h_mult} we have
	\begin{align}\label{eq:mult_proof_tmp2}
		&D^{\mathbf m}\left[(h\circ\psi)\prod_{k=1}^{d}\sqrt{(\omega_k\circ\psi_k)\psi'_{k}}\right](\mathbf x) \nonumber\\
		&= \partial^{m_d}\left[ \ldots \partial^{m_2}\left[ \partial^{m_1}\left[
			(h\circ\psi)\prod_{k=1}^{d}\sqrt{(\omega_k\circ\psi_k)\psi'_{k}}
		\right](x_1) \right](x_2) \ldots \right](x_d).
	\end{align}
	By applying the Leibniz formula as in \eqref{eq:Leibniz_estimate} we obtain for all $\ell=1,\ldots,d$
	\begin{align}\label{eq:mult_proof1}
		&\partial^{m_{\ell}}\left[ (h\circ\psi)\prod_{k=1}^{d}\sqrt{(\omega_k\circ\psi_k)\psi'_{k}}	\right](x_{\ell}) \nonumber\\
		&= \sum_{j_{\ell}=0}^{m_{\ell}} \binom{m_{\ell}}{j_{\ell}} \partial^{j_{\ell}}[h\circ \psi](x_{\ell}) \,\partial^{m_{\ell}-j_{\ell}}\left[\prod_{k=1}^{d}\sqrt{(\omega_k\circ\psi_k)\psi'_{k}}\right](x_{\ell})
	\end{align}
	and
	in total rewrite the expression in \eqref{eq:mult_proof_tmp2} as 
	\begin{align} \label{eq:mult_proof_leibniz_applied}
		D^{\mathbf m}\left[ (h\circ\psi)\prod_{k=1}^{d}\sqrt{(\omega_k\circ\psi_k)\psi'_{k}}\right](\mathbf x) 
		&= \sum_{j_1=0}^{m_1} \binom{m_1}{j_1} \ldots \sum_{j_d=0}^{m_d} \binom{m_d}{j_d}
		D^{(j_1,\ldots,j_d)}[h\circ \psi](\mathbf x) \,  \times\\
		&\quad\times D^{(m_1-j_1,\ldots,m_d-j_d)} \left[\prod_{k=1}^{d}\sqrt{(\omega_k\circ\psi_k)\psi'_{k}}\right](\mathbf x). \nonumber
	\end{align}
	Next, we apply the {Fa\'{a} di Bruno} formula \eqref{eq:FdiB_formula} to each univariate $j_k$-th derivative of $h\circ\psi$ in \eqref{eq:mult_proof1} so that for ${\ell}=1,\ldots,d$ we have
	\begin{align}
		 \partial^{j_{\ell}}[h\circ \psi](\mathbf x) =
		 \begin{cases} \label{eq:mult_proof_FdB}
		 	\displaystyle \sum_{i_{\ell}=1}^{j_{\ell}} \partial^{i_{\ell}}[h](\psi(\mathbf x)) B_{j_{\ell},i_{\ell}}(\psi'_{\ell}(x_{\ell}),\ldots, \psi^{(j_{\ell}-i_{\ell}+1)}_{\ell}(x_{\ell})) &\quad: j_{\ell}\in\mathbb{N}, \\
		 	h(\psi(\mathbf x)) &\quad: j_{\ell}=0,
		 \end{cases}
	\end{align}
	i.e., $B_{0,i_{\ell}}(\psi'_{\ell}(x_{\ell}),\psi''_{\ell}(x_{\ell}),\ldots, \psi^{(j_{\ell}-i_{\ell}+1)}_{\ell}(x_{\ell})) = 1$.
	
	We combine the norm $\left\| D^{\mathbf m}[f](\mathbf x) \right\|_{L_2(\mathbb{T}^d)}$ in $\eqref{eq:mult_proof_L2T}$ 
	with the expression resulting from applying the Leibniz formula to $D^{\mathbf m}[f]$ in $\eqref{eq:mult_proof_leibniz_applied}$ and the subsequent application of the {Fa\'{a} di Bruno} formula in \eqref{eq:mult_proof_FdB}.
	Then we estimate
	\begin{align} \label{eq:mult_proof_combi1}
		\left\| D^{\mathbf m}[f](\mathbf x) \right\|_{L_2(\mathbb{T}^d)}
		&\leq \sum_{j_1=0,\ldots,j_d=0}^{m_1,\ldots,m_d} \prod_{{\ell}=1}^{d} \binom{m_{\ell}}{j_{\ell}}
		\sum_{i_1=1,\ldots,i_d=1}^{j_1,\ldots,j_d}
		\left( \int_{\mathbb T^d}
			| D^{(i_1,\ldots,i_d)}[h](\psi(\mathbf x))|^2 \right.\times \\
		&\quad\times	
			\left. 
			 \prod_{{\ell}=1}^{d} |B_{j_{\ell},i_{\ell}}(\psi'_{\ell}(x_{\ell}),\ldots, \psi^{(j_{\ell}-i_{\ell}+1)}_{\ell}(x_{\ell}))|^2 
			\right.\times\nonumber\\
		&\quad\times
			\left. 
			\left| D^{(m_1-j_1,\ldots,m_d-j_d)}\left[\prod_{k=1}^{d}\sqrt{(\omega_k\circ\psi_k)\psi'_{k}}\right](\mathbf x) \right|^2
		\,\mathrm{d}\mathbf x
		\right)^{\frac{1}{2}}.\nonumber
	\end{align}
	In the multivariate integral appearing in \eqref{eq:mult_proof_combi1} we estimate each coordinate separately with the univariate arguments of the previous proof by fixing all but one coordinate one after another.
	Recalling the arguments in \eqref{eq:Bell_estimate2}, if all appearing derivatives of $\psi_{\ell}$ are in $\mathcal{C}(-1/2,1/2)$ and the corresponding derivatives of the density $\varrho_{\ell}$ are in $\mathcal{C}_{0}(\mathbb{R})$ then for all Bell polynomials $B_{j_{\ell},i_{\ell}}$ with $j_{\ell}\geq 1$ appearing in \eqref{eq:mult_proof_FdB} and \eqref{eq:mult_proof_combi1} there is some constant $C>0$ so that we can estimate
	\begin{align*}
		|B_{j_{\ell},i_{\ell}}(\psi'_{\ell}(x_{\ell}),\psi''_{\ell}(x_{\ell}),\ldots, \psi^{(j_{\ell}-i_{\ell}+1)}_{\ell}(x_{\ell}))| \leq C |\psi'_{\ell}(x_{\ell})|^{2j_{\ell}-1}.
	\end{align*}
	Analogously, to \eqref{eq:Leibniz_term_estimate1} and \eqref{eq:Leibniz_term_estimate2} for each $\ell=1,\ldots,d$ we have to separate the summand for $j_{\ell}=0$ from the summands corresponding to $j_{\ell} = 1,\ldots,d$.
	Starting with $\ell=1$ we estimate \eqref{eq:mult_proof_combi1} as in \eqref{eq:proof_univar_summary} after inserting the productive one $1 = \psi'_1(x_1)\frac{1}{\psi'_1(x_1)}$, so that
	\begin{align*}
		&\left\| D^{\mathbf m}[f](\mathbf x) \right\|_{L_2(\mathbb{T}^d)} \\
		&\leq C_1 \binom{m_{1}}{\lceil\frac{m_1}{2}\rceil} \max_{j_1=0,\ldots,m_{1}} \left\| \partial^{m_1-j_1}\left[\sqrt{(\omega_1\circ\psi_1)\psi'_{1}}\right](x_1) \, \psi_1'(x_1)^{\max(-\frac{1}{2},2j_{1}-\frac{3}{2})} \right\|_{L_{\infty}(\mathbb{T})} \times \\
		&\quad\times \sum_{j_2=0,\ldots,j_d=0}^{m_2,\ldots,m_d} \prod_{{\ell}=2}^{d} \binom{m_{\ell}}{j_{\ell}}
		\sum_{i_2=1,\ldots,i_d=1}^{j_2,\ldots,j_d}
		\left( \int_{\mathbb T^{d-1}} \right.\times \nonumber\\
		&\left.\quad\times
		\int_{\mathbb T} | D^{(i_1,\ldots,i_d)}[h](\psi_1(x_1),\ldots,\psi_d(x_d))|^2 \psi'_1(x_1) \,\mathrm{d}x_1 \right.\times \nonumber\\
		&\quad\times	
		\left. 
		\prod_{{\ell}=2}^{d} |B_{j_{\ell},i_{\ell}}(\psi'_{\ell}(x_{\ell}),\ldots, \psi^{(j_{\ell}-i_{\ell}+1)}_{\ell}(x_{\ell}))|^2 \right.\times\\
		&\quad\times \left.\left| D^{(m_2-j_2,\ldots,m_d-j_d)}\left[\prod_{k=2}^{d}\sqrt{(\omega_k\circ\psi_k)\psi'_{k}}\right](x_2,\ldots,x_d) \right|^2
		\,\mathrm{d}(x_2,\ldots,x_d)\right)^{\frac{1}{2}}.
	\end{align*}
	After repeating this process for $\ell=2,\ldots,d$ and inserting the inverse transformations ${x_{\ell} = \psi_{\ell}^{-1}(y_{\ell})}$ for all $\ell = 1,\ldots,d$ we end up with the estimate
	\begin{align*}
	&\left\| D^{\mathbf m}[f](\mathbf x) \right\|_{L_2(\mathbb{T}^d)} \\
	&\leq \prod_{\ell=1}^{d} C_{\ell} \binom{m_{\ell}}{\lceil\frac{m_\ell}{2}\rceil} \max_{j_\ell=0,\ldots,m_{\ell}} \left\| \partial^{m_{\ell}-j_{\ell}}\left[\sqrt{(\omega_{\ell}\circ\psi_{\ell})\psi'_{\ell}}\right](x_{\ell}) \, \psi_{\ell}'(x_{\ell})^{\max(-\frac{1}{2},2j_{\ell}-\frac{3}{2})} \right\|_{L_{\infty}(\mathbb{T})} \times \\
		&\quad\times \sum_{j_1=0,\ldots,j_d=0}^{m_1,\ldots,m_d} \sum_{i_1=1,\ldots,i_d=1}^{j_1,\ldots,j_d}
		\left( \int_{\mathbb T^d} 
		| D^{(i_1,\ldots,i_d)}[h](\psi_1(x_1),\ldots,\psi_d(x_d))|^2 \prod_{\ell=1}^{d} \psi'_{\ell}(x_{\ell}) \,\mathrm{d}\mathbf x \right)^{\frac{1}{2}} \\
	&\leq C \prod_{\ell=1}^{d} 
		\max_{j_\ell=0,\ldots,m_{\ell}} \left\| \partial^{m_{\ell}-j_{\ell}}\left[\sqrt{(\omega_{\ell}\circ\psi_{\ell})\psi'_{\ell}}\right](x_{\ell}) \, \psi_{\ell}'(x_{\ell})^{\max(-\frac{1}{2},2j_{\ell}-\frac{3}{2})} \right\|_{L_{\infty}(\mathbb{T})} \times \\
		&\quad\times \sum_{j_1=0,\ldots,j_d=0}^{m_1,\ldots,m_d} \sum_{i_1=1,\ldots,i_d=1}^{j_1,\ldots,j_d}
		\left( \int_{\mathbb R^d}
		| D^{(i_1,\ldots,i_d)}[h](y_1,\ldots,y_d)|^2 \,\mathrm{d}\mathbf y \right)^{\frac{1}{2}} \\
	&\leq C \prod_{\ell=1}^{d} 
		\max_{j_\ell=0,\ldots,m_{\ell}} \left\| \partial^{m_{\ell}-j_{\ell}}\left[\sqrt{(\omega_{\ell}\circ\psi_{\ell})\psi'_{\ell}}\right](x_{\ell}) \, \psi_{\ell}'(x_{\ell})^{\max(-\frac{1}{2},2j_{\ell}-\frac{3}{2})} \right\|_{L_{\infty}(\mathbb{T})} \times\nonumber\\
	&\quad\times
		\sum_{j_1=0,\ldots,j_d=0}^{m_1,\ldots,m_d} \|h\|_{H_{\mathrm{mix}}^{j}(\mathbb{R}^d)},
	\end{align*}
	with $j = \max(j_1,\ldots,j_d)$ in the last estimate.
	
	Since the previous estimate is valid for all multi-indices $\mathbf m = (m_1,\ldots,m_d)^{\top}\in\mathbb{N}_{0}^{d}$ with ${\|\mathbf m\|_{\ell_{\infty}} \leq m}$ we finally estimate
	\begin{align*}
		&\|f\|_{H_{\mathrm{mix}}^m(\mathbb{T}^d)}
		= \left( \sum_{\|\mathbf m\|_{\ell_{\infty}}\leq m} \left\| D^{\mathbf m}[f](\mathbf x) \right\|_{L_2(\mathbb{T}^d)}^2 \right)^{\frac{1}{2}}
		= \left( \sum_{m_1=0,\ldots,m_d=0}^{m,\ldots,m} \left\| D^{\mathbf m}[f](\mathbf x) \right\|_{L_2(\mathbb{T}^d)}^2 \right)^{\frac{1}{2}} \\
		&\leq C \prod_{\ell=1}^{d} \max_{m_{\ell}=0,\ldots,m} \times\\
		&\quad \times	
			\left(
			\max_{j_\ell=0,\ldots,m_{\ell}} \left\| \partial^{m_{\ell}-j_{\ell}}\left[\sqrt{(\omega_{\ell}\circ\psi_{\ell})\psi'_{\ell}}\right](x_{\ell}) \, \psi_{\ell}'(x_{\ell})^{\max(-\frac{1}{2},2j_{\ell}-\frac{3}{2})} \right\|_{L_{\infty}(\mathbb{T})}
		\right) \times\\
		&\quad \times (m+1)^{d} \|h\|_{H_{\mathrm{mix}}^{m}(\mathbb{R}^d)}. 
	\end{align*}
\end{proof}

\subsection{Approximation of transformed functions}
We establish two specific approximation error bounds for functions defined on $\mathbb{R}^d$ based on the approximation error bounds on the torus $\mathbb{T}^d$ that we recalled in Theorems~\ref{thm:L_infty_approx_error_torus} and \ref{eq:L_2_approximation_error_bound}.
The corresponding proofs rely heavily on the previously introduced sufficient conditions in Theorem~\ref{thm:Hm_composition_criteria_mult} that describe when Sobolev functions ${h\in L_2(\mathbb{R}^d,\omega)\cap H_{\mathrm{mix}}^{m}(\mathbb{R}^d)}$ with a multivariate weight function $\omega:\mathbb{R}^d\to[0,\infty)$ as given in \eqref{eq:omega_weighted_mult} can be transformed into Sobolev functions of dominated mixed smoothness on $\mathbb{T}^d$ of the form \eqref{eq:f_is_transformed_h_mult} by multivariate transformations $\psi:(-\frac{1}{2},\frac{1}{2})^d\to\mathbb{R}^d$ as given in \eqref{eq:psi_weighted_mult}.

At first, we fix some notation of certain multivariate objects.
Based on the definition of a {rank-$1$} lattice $\Lambda(\mathbf z,M)$ in \eqref{def:rank_one_lattice} we define a \textsl{transformed {rank-$1$} lattice} as 
\begin{align}\label{eq:Def_trafo_lattice}
	\Lambda_{\psi}(\mathbf z,M) := \left\{ \mathbf y_j := \psi(\mathbf x_j) : \mathbf x_j\in\Lambda(\mathbf z,M), j = 0,\ldots,M-1 \right\}.
\end{align}
Accordingly, we denote the \emph{transformed reconstructing {rank-$1$} lattice} by $\Lambda_{\psi}(\mathbf z,M,I)$.

Besides the weight function $\omega$, also the density $\varrho$ of the transformation $\psi$ is of product form as defined in \eqref{def:varrho_mult}, i.e., it is the product of univariate densities $\varrho_j(y_j), j=1,\ldots,d$.
Hence, based on the functions $\varphi_{k}$ given in \eqref{eq:transformed_basis_functions} this product form extends to 
\begin{align}\label{eq:transformed_basis_functions_mult}
	{\varphi_{\mathbf k}(\mathbf y) := \prod_{j=1}^{d}\varphi_{k_j}(y_j)}.
\end{align}
Similar to \eqref{def:weighted_scalar_product}, the multivariate weighted $L_{2}(\mathbb{R}^d,\omega)$ scalar product reads as
\begin{align} 
	(h_1, h_2)_{L_2\left(\mathbb{R}^d, \omega \right)}
	:= \int_{\mathbb{R}^d} \prod_{j=1}^{d}\omega_j(y_j) \, h_1(\mathbf y) \, \overline{h_2(\mathbf y)} \, \mathrm d\mathbf y
\end{align}
and similar to \eqref{def:FouCoeff_of_h} the multivariate Fourier coefficients are naturally given with respect to this scalar product as
\begin{align} \label{def:FC_trafo_multivar}
	\hat h_{\mathbf k} 
	= (h, \varphi_{\mathbf k})_{L_2\left(\mathbb{R}^d, \omega \right)}.
\end{align}
As before in \eqref{def:Fourier_part_sum_of_h} we define the multivariate Fourier partial sum as
\begin{align*}
	S_{I}h(\mathbf y)
	:= \sum_{\mathbf k\in I} \hat h_{\mathbf k} \, \varphi_{\mathbf k}(\mathbf y).
\end{align*}
Suppose $f\in L_{2}(\mathbb{T}^d)$, then for each $I\subset\mathbb{Z}^d$ the system $\left\{\varphi_{\mathbf k}\right\}_{\mathbf k\in I}$ spans the space of transformed trigonometric polynomials
\begin{align}
	\label{def:trig_poly_trafo_mult}
	\Pi_{I,\psi} := \mathrm{span}\left\{ \sqrt{\frac{\varrho(\circ)}{\omega(\circ)}} \, \mathrm{e}^{2\pi\mathrm i \mathbf k\cdot\psi^{-1}(\circ)} : \mathbf k \in I \right\}.
\end{align}
Similar to \eqref{eq:exact_integration_formula}, for transformed trigonometric polynomials $h \in \Pi_{I,\psi}$, transformed lattice nodes 
${\mathbf y_j\in\Lambda_{\psi}(\mathbf z,M,I)}$ 
and $\mathbf k\in I$ we have the exact integration property of the form 
\begin{align} 
	\label{eq:exact_integration_prop_2}
	\hat h_{\mathbf k}
	&= \int_{\mathbb{R}^d} h(\mathbf y) \, \sqrt{\varrho(\mathbf y)\,\omega(\mathbf y)} \, \mathrm{e}^{-2\pi\mathrm i \mathbf k\cdot\psi^{-1}(\mathbf y)} \, \mathrm dy 
	= \int_{\mathbb{T}^d} f(\mathbf x) \,\mathrm{e}^{-2\pi\mathrm{i}\mathbf k\cdot\mathbf x} \,\mathrm{d}\mathbf x \nonumber \\
	&= \frac{1}{M} \sum_{j =0}^{M-1} f(\mathbf x_j) \,\mathrm{e}^{-2\pi\mathrm{i}\mathbf k\cdot\mathbf x_j} 
	= \frac{1}{M} \sum_{j=0}^{M-1} h(\mathbf y_j)\, \sqrt{\frac{\varrho(\mathbf y_j)}{\omega(\mathbf y_j)}} \, \mathrm{e}^{-2\pi\mathrm i \mathbf k\cdot\psi^{-1}(\mathbf y_j)}
	=\hat h_{\mathbf k}^{\Lambda}.
\end{align}
Generally, the multivariate approximated Fourier coefficients of the form
\begin{align*}
	\hat h_{\mathbf k}^{\Lambda} 
	= \frac{1}{M} \sum_{j=0}^{M-1} h(\mathbf y_j)\, \sqrt{\frac{\varrho(\mathbf y_j)}{\omega(\mathbf y_j)}} \, \mathrm{e}^{-2\pi\mathrm i \mathbf k\cdot\psi^{-1}(\mathbf y_j)} 
	= \frac{1}{M} \sum_{j=0}^{M-1} h(\mathbf y_j)\, \overline{\varphi_{\mathbf k}(\mathbf y_j)}
\end{align*}
approximate the multivariate Fourier coefficients $\hat h_{\mathbf k}$.
Finally, the multivariate version of the approximated Fourier partial sum is given by
\begin{align}
	\label{def:mult_approximated_FPS_FC}
	S_{I}^{\Lambda}h(\mathbf x) := \sum_{\mathbf k\in I} \hat h_{\mathbf k}^{\Lambda} \, \varphi_{\mathbf k}(\mathbf y).
\end{align}
Similarly to the $\mathcal{H}^{\beta}(\mathbb{T}^d)$-norm in \eqref{def:HbetaRaum} we define a norm of weighted Fourier coefficients $\hat h_{\mathbf k}$ of the form
\begin{align*}
	\|h\|_{\mathcal{H}^{m}(\mathbb{R}^d)}^2
	:= \sum_{\mathbf k\in\mathbb{Z}^d} \omega_{\mathrm{hc}}(\mathbf k)^{2m} |\hat h_{\mathbf k}|^2.
\end{align*}

With these rewritten objects we transfer the approximation error bounds in Theorems~\ref{thm:L_infty_approx_error_torus} and \ref{eq:L_2_approximation_error_bound} for functions defined on the torus to $\mathbb{R}^d$.

\subsubsection{$L_{\infty}$-approximation error}
Based on the $L_{\infty}(\mathbb{T}^d)$-approximation error bound \eqref{eq:torus_infty_approx_bound} and the conditions proposed in Theorem~\ref{thm:Hm_composition_criteria_mult} we prove a similar upper bound for the approximation error $\|h - S_{I_N^d}^{\Lambda} h\|$ in terms of a weighted $L_{\infty}$-norm on $\mathbb{R}^d$.
\begin{theorem} \label{thm:L_infty_approx_error_multivar}
	Let $d\in\mathbb{N}$, $m\in\mathbb{N}_{0}$, a hyperbolic cross $I_{N}^{d}$ with $N\geq 2^{d+1}$ and a reconstructing {rank-$1$} lattice $\Lambda(\mathbf z, M, I_{N}^{d})$ be given.	
	Let $\psi$ be a multivariate transformation as defined in \eqref{eq:psi_weighted_mult} with its corresponding density function $\varrho$ in product form \eqref{def:varrho_mult}.
	Let $\omega$ be a weight function as in \eqref{eq:omega_weighted_mult}
	and we consider a multivariate function ${h\in L_2(\mathbb{R}^d,\omega)\cap H_{\mathrm{mix}}^{m}(\mathbb{R}^d)}$.
	Let $\lambda > \frac{1}{2}$.
	Furthermore, 
	for all multi-indices $\mathbf m = (m_1,\ldots,m_d)^{\top} \in\mathbb{N}_{0}^{d}$ with $\|\mathbf m\|_{\ell_{\infty}} \leq m$ and all $j_{\ell}=0,\ldots,m, \ell=1,\ldots,d$ we assume 
	\begin{align*}
		{\partial^{j_{\ell}}\left[ \varrho \right](y_{\ell})\in\mathcal{C}_{0}(\mathbb{R})}, \qquad
		{\partial^{j_{\ell}}\left[ \psi \right](x_{\ell})\in\mathcal{C}((-1/2,1/2))}
	\end{align*}
	and
	\begin{align*}
		\max_{j_\ell=0,\ldots,m_{\ell}} \left\| \partial^{m_{\ell}-j_{\ell}}\left[\sqrt{(\omega_{\ell}\circ\psi_{\ell})\psi'_{\ell}}\right](x_{\ell}) \, \psi_{\ell}'(x_{\ell})^{\max(-\frac{1}{2},2j_{\ell}-\frac{3}{2})} \right\|_{L_{\infty}(\mathbb{T})}
		<\infty.
	\end{align*}
	
	Then there is an approximation error estimate of the form
	\begin{align*}
		\left\|h - S_{I_N^d}^{\Lambda} h \right\|_{L_{\infty}\big(\mathbb{R}^d, \sqrt{\frac{\omega}{\varrho}}\big)} 
		\lesssim N^{-m+\lambda} \|h\|_{\mathcal{H}^{m}(\mathbb{R}^d)}.\end{align*}\end{theorem}
\begin{proof}
	Let $m\in\mathbb{N}, d\in\mathbb{N}$ and let ${h\in L_2(\mathbb{R}^d,\omega)\cap H_{\mathrm{mix}}^{m}(\mathbb{R}^d)}$.
	By assumption the criteria in Theorem~\ref{thm:Hm_composition_criteria_mult} are fulfilled and thus the transformed function $f$ of the form \eqref{eq:f_is_transformed_h_mult} is in $H_{\mathrm{mix}}^{m}(\mathbb{T}^d)$.
	This $f$ is also in $\mathcal{H}^{m}(\mathbb{T}^d)$, due to the norm equivalence \eqref{eq:Hs_norm_equivalence} and furthermore have a continuous representative, because of the inclusion ${\mathcal{H}^{m}(\mathbb{T}^d) \hookrightarrow \mathcal{A}^{m-\lambda}(\mathbb T^d) \hookrightarrow \mathcal{C}(\mathbb{T}^d)}$ with $\lambda>\frac{1}{2}$ as in \eqref{eq:Wiener_algebra_inclusion}.
	Hence, for $f\in\mathcal{A}^{m-\lambda}(\mathbb{T}^d)\cap\mathcal{C}(\mathbb{T}^d)$ we have the approximation error bound 
	\begin{align}
		\label{eq:Linf_proof_1}
		\|f - S_{I_N^d}^{\Lambda} f\|_{L_{\infty}(\mathbb{T}^d)}
		\leq 2 N^{-m+\lambda} \|f\|_{\mathcal{A}^{m-\lambda}(\mathbb{T}^d)}
	\end{align}
	as stated in Theorem~\ref{thm:L_infty_approx_error_torus}.
	
	With the inverse transformation $\mathbf x=\psi^{-1}(\mathbf y)$ we have
	\begin{align*}
		\hat h_{\mathbf k}
= \left(h, \varphi_{\mathbf k} \right)_{L_2\left(\mathbb{R}^d, \omega \right)}
		= (f,\mathrm{e}^{2\pi\mathrm i \mathbf k\cdot\circ})_{L_2(\mathbb{T}^d)}
		= \hat f_{\mathbf k}
	\end{align*}
	and 
	\begin{align}
		\label{eq:Linf_proof_2}
		\|h\|_{\mathcal{H}^{m}(\mathbb{R}^d)}^2
		= \sum_{\mathbf k\in\mathbb{Z}^d} \omega_{\mathrm{hc}}(\mathbf k)^{2m} |\hat h_{\mathbf k}|^2
		= \sum_{\mathbf k\in\mathbb{Z}^d} \omega_{\mathrm{hc}}(\mathbf k)^{2m} |\hat f_{\mathbf k}|^2		
		= \|f\|_{\mathcal{H}^{m}(\mathbb{T}^d)}^2,
	\end{align}
	as well as
	\begin{align*}
		\|h - S_{I_{N}^{d}}h\|_{L_{\infty}\big(\mathbb{R}^d, \sqrt{\frac{\omega}{\varrho}}\big)}
		&= \esssup_{\mathbf y\in\mathbb{R}^{d}} \left| \sqrt{\frac{\omega(\mathbf y)}{\varrho(\mathbf y)}} \left( h(\mathbf y) - \sum_{\mathbf k\in I_{N}^{d}} \hat h_{\mathbf k} \, \varphi_{\mathbf k}(\mathbf y)\right) \right| \\
		&= \esssup_{\mathbf y\in\mathbb{R}^{d}} \left| h(\mathbf y)\,\sqrt{\frac{\omega(\mathbf y)}{\varrho(\mathbf y)}} - \sum_{\mathbf k\in I_{N}^{d}} \hat h_{\mathbf k} \, \mathrm{e}^{2\pi\mathrm i \mathbf k\cdot\psi^{-1}(\mathbf{y})} \right| \\
		&= \esssup_{\mathbf x\in\mathbb{T}^{d}} \left| h(\psi(\mathbf x))\,\sqrt{\omega(\psi(\mathbf x))\psi'(\mathbf x)} - \sum_{\mathbf k\in I_{N}^{d}} \hat h_{\mathbf k} \, \mathrm{e}^{2\pi\mathrm i \mathbf k\cdot\mathbf{x}} \right| \\
		&= \|f - S_{I_{N}^{d}}f\|_{L_{\infty}(\mathbb{T}^d)}
	\end{align*}
	and 
	\begin{align}
		\label{eq:Linf_proof_3}
		\|h - S_{I_{N}^{d}}^{\Lambda}h\|_{L_{\infty}\big(\mathbb{R}^d, \sqrt{\frac{\omega}{\varrho}}\big)} 
		= \|f - S_{I_{N}^{d}}^{\Lambda}f\|_{L_{\infty}(\mathbb{T}^d)}.
	\end{align} 
	In total, by combining \eqref{eq:Linf_proof_3}, \eqref{eq:Linf_proof_1}, \eqref{eq:Wiener_algebra_inclusion2}, and \eqref{eq:Linf_proof_2} we estimated for $f\in\mathcal{H}^{m}(\mathbb{T}^d)\cap\mathcal{C}(\mathbb{T}^d)$ that the approximation error can be bounded by
	\begin{align*}
		\|h - S_{I_N^d}^{\Lambda} h\|_{L_{\infty}\big(\mathbb{R}^d, \sqrt{\frac{\omega}{\varrho}}\big)}
		&= \|f - S_{I_N^d}^{\Lambda} f\|_{L_{\infty}(\mathbb{T}^d)}
		\leq 2 N^{-m+\lambda} \|f\|_{\mathcal{A}^{m-\lambda}(\mathbb{T}^d)} \\
		&\leq 2 C_{d,\lambda} N^{-m+\lambda} \|f\|_{\mathcal{H}^{m}(\mathbb{T}^d)}
		= 2 C_{d,\lambda} N^{-m+\lambda} \|h\|_{\mathcal{H}^{m}(\mathbb{R}^d)}
		< \infty
	\end{align*}
	with $\lambda > \frac{1}{2}$ and some constant $C_{d,\lambda} > 1$.
\end{proof}

\subsubsection{$L_2$-approximation error}
Similarly, based on the $L_{2}(\mathbb{T}^d)$-approximation error bound \eqref{eq:H_beta_error_bound} and the conditions proposed in Theorem~\ref{thm:Hm_composition_criteria_mult} we prove an upper bound for the approximation error $\|h - S_{I_N^d}^{\Lambda} h\|$ in terms of a weighted $L_{2}$-norm on $\mathbb{R}^d$.
\begin{theorem}\label{thm:Hm_approx_error_decay_multivar}
	Let $d\in\mathbb{N}$, $m\in\mathbb{N}_0$, a hyperbolic cross $I_{N}^{d}$ with $N\geq 2^{d+1}$ and a reconstructing {rank-$1$} lattice $\Lambda(\mathbf z, M, I_{N}^{d})$ be given.
	Let $\psi$ be a multivariate transformation as in \eqref{eq:psi_weighted_mult}
	and $\omega$ be a multivariate weight function as in \eqref{eq:omega_weighted_mult}.
	We consider a multivariate function ${h\in L_2(\mathbb{R}^d,\omega)\cap H_{\mathrm{mix}}^{m}(\mathbb{R}^d)}$.
	Furthermore, for all multi-indices $\mathbf m = (m_1,\ldots,m_d)^{\top} \in\mathbb{N}_{0}^{d}$ with $\|\mathbf m\|_{\ell_{\infty}} \leq m$ and all $j_{\ell}=0,\ldots,m, \ell=1,\ldots,d$ we assume
	\begin{align*}
		{\partial^{j_{\ell}}\left[ \varrho \right](y_{\ell})\in\mathcal{C}_{0}(\mathbb{R})}, \qquad
		{\partial^{j_{\ell}}\left[ \psi \right](x_{\ell})\in\mathcal{C}((-1/2,1/2))}
	\end{align*}
	and
	\begin{align*}
		\max_{j_\ell=0,\ldots,m_j} \left\| \partial^{m_{\ell}-j_{\ell}}\left[\sqrt{(\omega_{\ell}\circ\psi_{\ell})\psi'_{\ell}}\right](x_{\ell}) \, \psi_{\ell}'(x_{\ell})^{\max(-\frac{1}{2},2j_{\ell}-\frac{3}{2})} \right\|_{L_{\infty}(\mathbb{T})}
		<\infty.
	\end{align*}
	
	Then there is an approximation error estimate of the form
	\begin{align*}
		\|h - S_{I_{N}^{d}}^{\Lambda}h\|_{L_{2}\left(\mathbb{R}^d, \omega \right)} 
		\lesssim N^{-m} (\log N)^{(d-1)/2} \|h\|_{\mathcal{H}^{m}(\mathbb{R}^d)}.
	\end{align*}
\end{theorem}
\begin{proof}
	Let $m\in\mathbb{N}, d\in\mathbb{N}$ and let ${h\in L_2(\mathbb{R}^d,\omega)\cap H_{\mathrm{mix}}^{m}(\mathbb{R}^d)}$.
	By assumption the criteria in Theorem~\ref{thm:Hm_composition_criteria_mult} are fulfilled and thus the transformed function $f$ of the form \eqref{eq:f_is_transformed_h_mult} is in $H^{m}(\mathbb{T}^d)$.
	This $f$ is also in $\mathcal{H}^{m}(\mathbb{T}^d)$, due to the norm equivalence \eqref{eq:Hs_norm_equivalence} and they furthermore have a continuous representative, because of the inclusion $\mathcal{H}^{m}(\mathbb{T}^d) \hookrightarrow \mathcal{C}(\mathbb{T}^d)$ as in \eqref{eq:Wiener_algebra_inclusion}.
	For $f\in\mathcal{H}^{m}(\mathbb{T}^d) \cap \mathcal{C}(\mathbb{T}^d)$ Theorem~\ref{eq:L_2_approximation_error_bound} yields the approximation error bound of the form
	\begin{align}\label{eq:L2_proof_1}
	\|f - S_{I_{N}^{d}}^{\Lambda}f\|_{L_{2}(\mathbb{T}^d)} 
	\leq C_{d,\beta} N^{-\beta} (\log N)^{(d-1)/2} \|f\|_{\mathcal{H}^{\beta}(\mathbb{T}^d)}
	\end{align}
	with some constant $C_{d,\beta} := C(d,\beta)>0$.
	With the inverse transformation $\bm x=\psi^{-1}(\bm y)$ we have
	\begin{align*} 
		\hat h_{\mathbf k}
		= \left(h, \varphi_{\mathbf k} \right)_{L_2\left(\mathbb{R}^d, \omega\right)}
		= (f,\mathrm{e}^{2\pi\mathrm i \mathbf k\cdot\circ})_{L_2(\mathbb{T}^d)}
		= \hat f_{\mathbf k},
	\end{align*}
	and
	\begin{align*}
		\|h\|_{\mathcal{H}^{m}(\mathbb{R}^d)}^2
		= \sum_{\mathbf k\in\mathbb{Z}^d} \omega_{\mathrm{hc}}(\mathbf k)^{2m} |\hat h_{\mathbf k}|^2
		= \sum_{\mathbf k\in\mathbb{Z}^d} \omega_{\mathrm{hc}}(\mathbf k)^{2m} |\hat f_{\mathbf k}|^2		
		= \|f\|_{\mathcal{H}^{m}(\mathbb{T}^d)}^2
	\end{align*}
	as in \eqref{eq:Linf_proof_2}, as well as 
	\begin{align} \label{eq:L2_proof_3}
		\|h - S_{I_{N}^{d}}h\|_{L_{2}\left(\mathbb{R}^d, \omega\right)}^2
		&= \int_{\mathbb{R}^d} \left| h(\mathbf y) - \sum_{\mathbf k\in I_{N}^{d}} \hat h_{\mathbf k} \, \varphi_{\mathbf k}(\mathbf y) \right|^2 \omega(\bm y) \,\mathrm{d}\mathbf y
		= \|f - S_{I_{N}^{d}}f\|_{L_{2}(\mathbb{T}^d)}^2
	\intertext{and}
		\|h - S_{I_{N}^{d}}^{\Lambda}h\|_{L_{2}\left(\mathbb{R}^d, \omega \right)}
		&= \|f - S_{I_{N}^{d}}^{\Lambda}f\|_{L_{2}(\mathbb{T}^d)}. \nonumber
	\end{align}
	In total, by combining \eqref{eq:L2_proof_3}, \eqref{eq:L2_proof_1}, and \eqref{eq:Linf_proof_2} we estimated for $f\in\mathcal{H}^{m}(\mathbb{T}^d)\cap\mathcal{C}(\mathbb{T}^d)$ that the approximation error can be bounded by
	\begin{align*}
		\|h - S_{I_N^d}^{\Lambda} h\|_{L_{2}\big(\mathbb{R}^d, \omega\big)}
		= \|f - S_{I_N^d}^{\Lambda} f\|_{L_{2}(\mathbb{T}^d)}
		&\lesssim C_{d,\beta} N^{-\beta} (\log N)^{(d-1)/2} \|f\|_{\mathcal{H}^{\beta}(\mathbb{T}^d)} \\
		&= C_{d,\beta} N^{-\beta} (\log N)^{(d-1)/2} \|h\|_{\mathcal{H}^{m}(\mathbb{R}^d)}
		< \infty
	\end{align*}
	with some constant $C_{d,\beta} > 0$.
\end{proof}

\section{Algorithms}
In this chapter we start denoting the parameters $\bm\eta,\bm\mu\in\mathbb{R}^d$ 
in families of multivariate parameterized transformations $\psi(\circ,\bm\eta)$ as in \eqref{eq:psi_weighted_mult}, 
in families of multivariate parameterized weight functions $\omega(\circ, \bm \mu)$ as in \eqref{eq:omega_weighted_mult} 
and in all related functions and objects.

We adapt the algorithms described in \cite[Algorithm~3.1 and 3.2]{kaemmererdiss} that are based on one-dimensional fast Fourier transforms (FFTs).
They are used for the fast reconstruction of approximated Fourier coefficients $\hat h_{\mathbf{k}}^{\Lambda}$ and the evaluation of a transformed multivariate trigonometric polynomials, in particular the approximated Fourier series $S_{I}^{\Lambda} h$, both given in \eqref{def:mult_approximated_FPS_FC}.
This is denoted as matrix-vector-products of the form
\begin{align} \label{eq: matrix vector products}
	\mathbf h = \mathbf A\mathbf{\hat h}
	\quad \text{and} \quad 
	\mathbf{\hat{h}} = M^{-1}\mathbf A^*\mathbf{h}
\end{align}
with $\bm\eta,\bm\mu\in\mathbb{R}^d, \mathbf{h} := \left(h(\mathbf y_j)\, \sqrt{\frac{\omega(\mathbf y_j, \bm\mu)}{\varrho(\mathbf y_j, \bm\eta)}}\right)_{j=0,\ldots,M-1}$ for $\mathbf y_j \in \Lambda_{\psi(\circ,\bm\eta)}(\mathbf z,M)$, $\hat{\mathbf{h}} := (\hat h_{\mathbf k})_{\mathbf k \in I_N}$ and the transformed Fourier matrices $\mathbf A$ and $\mathbf A^{*}$ given by
\begin{align*}
	\mathbf A 
	&:= \left( \mathrm{e}^{2 \pi \mathrm{i} \mathbf k \cdot\psi^{-1}(\mathbf y_j,\bm\eta)} \right)_{\mathbf y_j \in \Lambda_{\psi(\circ,\bm\eta)}(\mathbf z,M),\, \mathbf k\in I} \in \mathbb{C}^{M \times |I|}, \nonumber\\
	\mathbf A^{*} 
	&:= \left( \mathrm{e}^{-2 \pi \mathrm{i} \mathbf k \cdot\psi^{-1}(\mathbf y_j,\bm\eta)} \right)_{\mathbf k\in I,\, \mathbf y_j \in \Lambda_{\psi(\circ,\bm\eta)}(\mathbf z,M)}
	\in \mathbb{C}^{|I| \times M}.
\end{align*}
We incorporate the previously described idea that the functions ${h\in L_2(\mathbb{R}^d,\omega)\cap H_{\mathrm{mix}}^{m}(\mathbb{R}^d)}$ are transformed into functions $f$ on the torus $\mathbb T^d$ that are of the form \eqref{eq:f_is_transformed_h_mult} via transformations $\mathbf x_j = \psi(\mathbf y_j,\bm\eta)$, so that we have samples of the form
\begin{align*}
	h(\mathbf y_j)\, \sqrt{\frac{\omega(\mathbf y_j, \bm\mu)}{\varrho(\mathbf y_j, \bm\eta)}} 
	= h(\psi(\mathbf x_j,\bm\eta)) \, \sqrt{ \omega(\psi(\mathbf x_j,\bm\eta), \bm\mu) \, \psi'(\mathbf x_j, \bm\eta) } 
	= f(\mathbf x_j,\bm\eta,\bm\mu) = f(\mathbf x_j),
\end{align*}
depending on the particular choices for $\bm\eta,\bm\mu\in\mathbb{R}^d$.
\begin{remark}
	We identify $\mathbb T^d$ with different cubes.
	On the one hand, when defining {rank-$1$} lattices $\Lambda(\mathbf z, M)$ in \eqref{def:rank_one_lattice} we identify it with $[0,1)^d$.
	On the other hand, in order to apply the transformations $\psi$ we need to consider $\mathbb T^d \simeq [-\frac{1}{2},\frac{1}{2})^d$, which we achieve by reassigning all lattice points $\mathbf x_j\in\Lambda(\mathbf z, M)$ via
	\begin{align*}
		{\mathbf{x}}_j \mapsto \left( \left( \mathbf{x}_j + \frac{\mathbf 1}{2} \right) \bmod \mathbf 1 \right) - \frac{\mathbf 1}{2}
	\end{align*}
	for all $j=0,\ldots,M-1$.
	
	We already showcased in Figure~\ref{fig:Plot of trafo examples} that the definition of $\psi$ in \eqref{def:Trafo_def} allows a range of functions with different slopes, which manifested in algebraic or exponential density functions $\varrho$.
	In Figure~\ref{fig:Lattice_and_transformed_lattice} we highlight these differences once more with transformed {rank-$1$} lattices $\Lambda_{\psi(\circ,\bm\eta)}(\mathbf z, M)$ as defined in \eqref{eq:Def_trafo_lattice}.
	We consider the two-dimensional {rank-$1$} lattice $\Lambda(\mathbf z, M)$ generated by $\mathbf z = (1,3)^{\top}$ and $M = 31$.
	We compare the transformed lattices for the algebraic transformation and the error function transformation of the form \eqref{eq:param_trafo_explicit_example} in their two-dimensional versions given by
	\begin{align}
		\label{eq:2D_example_trafos}
		\psi(\mathbf x, \bm\eta) = \left(\frac{2\eta_1 x_1}{\sqrt{1-4x_1^2}}, \frac{2\eta_2 x_2}{\sqrt{1-4x_2^2}}\right)^{\top},
		\psi(\mathbf x, \bm\eta) = \left(\eta_1\mathrm{erf}^{-1}(2x_1), \eta_2\mathrm{erf}^{-1}(2x_2)\right)^{\top}.
	\end{align}
	For $\bm\eta = (\eta_1, \eta_2)^{\top} = (1, 1)^{\top}$ the graphs in the center and on the right of Figure~\ref{fig:Lattice_and_transformed_lattice} reveal that the algebraic transformation causes a wider spread of the lattice nodes close to the center, 
	whereas the slope of the error function transformation increases hugely towards the boundary points which we only notice for larger values $M$ and much finer lattices with more nodes closer to the boundary of the cube $(-\frac{1}{2},\frac{1}{2})^2$.
\end{remark}
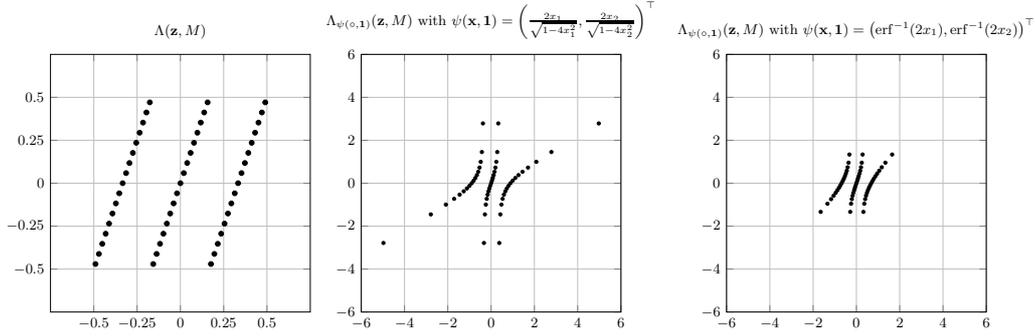
\begin{figure}[tb]
	\centering
\begin{tikzpicture}[scale=0.6]
	\begin{axis}[scatter/classes = { a = {mark=o, draw=black} },
	font=\footnotesize,	grid = both,
	xmax = 0.75,	xmin = -0.75,	ymax = 0.75,	ymin = -0.75,
xtick={-0.5,-0.25,0,0.25,0.5}, ytick={-0.5,-0.25,0,0.25,0.5},
	title = {$\Lambda(\mathbf z, M)$},
	unit vector ratio*=1 1 1
	]
	\addplot[scatter ,only marks, mark size=1.5, scatter src = explicit symbolic]coordinates{
		(   0,    0)  (0.01961, 0.05882)  (0.03922, 0.1176)  (0.05882, 0.1765)  (0.07843, 0.2353)  (0.09804, 0.2941)  (0.1176, 0.3529)  (0.1373, 0.4118)  (0.1569, 0.4706)  (0.1765, -0.4706)  (0.1961, -0.4118) (0.2157, -0.3529)  (0.2353, -0.2941)  (0.2549, -0.2353)  (0.2745, -0.1765)  (0.2941, -0.1176)  (0.3137, -0.05882) (0.3333,    0)  (0.3529, 0.05882)  (0.3725, 0.1176)  (0.3922, 0.1765)  (0.4118, 0.2353)  (0.4314, 0.2941)  (0.451, 0.3529)  (0.4706, 0.4118)  (0.4902, 0.4706)  (-0.4902, -0.4706)  (-0.4706, -0.4118)  (-0.451, -0.3529)  (-0.4314, -0.2941)  (-0.4118, -0.2353)  (-0.3922, -0.1765)  (-0.3725, -0.1176)  (-0.3529, -0.05882)  (-0.3333,    0)  (-0.3137, 0.05882)  (-0.2941, 0.1176)  (-0.2745, 0.1765)  (-0.2549, 0.2353)  (-0.2353, 0.2941)  (-0.2157, 0.3529)  (-0.1961, 0.4118)  (-0.1765, 0.4706)  (-0.1569, -0.4706)  (-0.1373, -0.4118)  (-0.1176, -0.3529)  (-0.09804, -0.2941)  (-0.07843, -0.2353)  (-0.05882, -0.1765)  (-0.03922, -0.1176)  (-0.01961, -0.05882)   };
	\end{axis}
	\end{tikzpicture}
	\begin{tikzpicture}[scale=0.6]
	\begin{axis}[scatter/classes = { a = {mark=o, draw=black} },
	font=\footnotesize,	grid = both,
	xmax = 6,	xmin = -6,	ymax = 6,	ymin = -6,
title = {$\Lambda_{\psi(\circ,\bm 1)}(\mathbf z, M)$ with $\psi(\mathbf x, \mathbf 1) = \left(\frac{2x_1}{\sqrt{1-4x_1^2}}, \frac{2x_2}{\sqrt{1-4x_2^2}}\right)^{\top}$ },
	unit vector ratio*=1 1 1
	]
	\addplot[scatter ,only marks, mark size=1.0, scatter src = explicit symbolic]coordinates{
		(   0,    0)  (0.03925, 0.1185)  (0.07867, 0.2421)  (0.1185, 0.3772)  (0.1588, 0.5333)  ( 0.2, 0.7274)  (0.2421, 0.9965)  (0.2855, 1.452)  (0.3304, 2.785)  (0.3772, -2.785)  (0.4263, -1.452)  (0.4781, -0.9965)  (0.5333, -0.7274)  (0.5926, -0.5333)  (0.6569, -0.3772)  (0.7274, -0.2421)  (0.8058, -0.1185)  (0.8944,    0)  (0.9965, 0.1185)  (1.117, 0.2421)  (1.264, 0.3772)  (1.452, 0.5333)  (1.706, 0.7274)  (2.089, 0.9965)  (2.785, 1.452)  (4.975, 2.785)  (-4.975, -2.785)  (-2.785, -1.452)  (-2.089, -0.9965)  (-1.706, -0.7274)  (-1.452, -0.5333)  (-1.264, -0.3772)  (-1.117, -0.2421)  (-0.9965, -0.1185)  (-0.8944,    0)  (-0.8058, 0.1185)  (-0.7274, 0.2421)  (-0.6569, 0.3772)  (-0.5926, 0.5333)  (-0.5333, 0.7274)  (-0.4781, 0.9965)  (-0.4263, 1.452)  (-0.3772, 2.785)  (-0.3304, -2.785)  (-0.2855, -1.452)  (-0.2421, -0.9965)  (-0.2, -0.7274)  (-0.1588, -0.5333)  (-0.1185, -0.3772)  (-0.07867, -0.2421)  (-0.03925, -0.1185)};
	\end{axis}
	\end{tikzpicture}
	\begin{tikzpicture}[scale=0.6]
	\begin{axis}[scatter/classes = { a = {mark=o, draw=black} },
	font=\footnotesize,	grid = both,
	xmax = 6,	xmin = -6,	ymax = 6,	ymin = -6,
title = {$\Lambda_{\psi(\circ,\bm 1)}(\mathbf z, M)$ with $\psi(\mathbf x, \mathbf 1) = \left(\mathrm{erf}^{-1}(2x_1), \mathrm{erf}^{-1}(2x_2)\right)^{\top}$},
	unit vector ratio*=1 1 1
	]
	\addplot[scatter ,only marks, mark size=1.0, scatter src = explicit symbolic]coordinates{
		(   0,    0)  (0.03477, 0.1046)  (0.06962, 0.2116)  (0.1046, 0.3238) (0.1399, 0.4447)  (0.1756, 0.5804)  (0.2116, 0.7418)  (0.2483, 0.9558)  (0.2856, 1.336)  (0.3238, -1.336)  (0.3629, -0.9558)  (0.4031, -0.7418)  (0.4447, -0.5804)  (0.4879, -0.4447)  (0.533, -0.3238)  (0.5804, -0.2116)  (0.6305, -0.1046)  (0.6841,    0)  (0.7418, 0.1046)  (0.8051, 0.2116)  (0.8755, 0.3238)  (0.9558, 0.4447)  (1.051, 0.5804)  (1.17, 0.7418)  (1.336, 0.9558)  (1.65, 1.336)  (-1.65, -1.336)  (-1.336, -0.9558)  (-1.17, -0.7418)  (-1.051, -0.5804)  (-0.9558, -0.4447)  (-0.8755, -0.3238) (-0.8051, -0.2116)  (-0.7418, -0.1046)  (-0.6841,    0)  (-0.6305, 0.1046)  (-0.5804, 0.2116)  (-0.533, 0.3238)  (-0.4879, 0.4447)  (-0.4447, 0.5804)  (-0.4031, 0.7418)  (-0.3629, 0.9558)  (-0.3238, 1.336)  (-0.2856, -1.336)  (-0.2483, -0.9558)  (-0.2116, -0.7418)  (-0.1756, -0.5804)  (-0.1399, -0.4447) 	(-0.1046, -0.3238)  (-0.06962, -0.2116)  (-0.03477, -0.1046) 
	};
	\end{axis}
	\end{tikzpicture}
	\caption{A two-dimensional lattice $\Lambda(\mathbf z,M)$ with $\mathbf z = (1,3)^{\top}, M = 31$ on the left and the resulting transformed lattice $\Lambda_{\psi(\circ,\bm\eta)}(\mathbf z, M)$ for the algebraic transformation in the center and for the error function transformation on the right, as given in \eqref{eq:2D_example_trafos} and both used with $\bm\eta=\mathbf 1$.
	}
	\label{fig:Lattice_and_transformed_lattice}
\end{figure}

\subsection{Evaluation of transformed multivariate trigonometric polynomials}
\begin{algorithm}[tb]
\caption{Evaluation at {rank-$1$} lattice}
\label{alg:LFFT_eval}
\begin{tabular}{p{2.0cm}p{4.5cm}p{6.0cm}}
	Input: 	
	& $M\in\mathbb{N}$ & lattice size of $\Lambda_{\psi(\circ,\bm\eta)}(\mathbf z,M)$\\
	& $\mathbf z\in\mathbb{Z}^d$ & generating vector of $\Lambda_{\psi(\circ,\bm\eta)}(\mathbf z,M)$\\
	& $I\subset\mathbb{Z}^d$ & frequency set of finite cardinality\\
	& $\mathbf{\hat h} = \left(\hat{h}_{\mathbf k}\right)_{\mathbf k\in I}$ & Fourier coefficients of $h\in\Pi_{I,\psi(\circ,\bm\eta)}$
\end{tabular}
\begin{algorithmic}
	\STATE $\mathbf {\hat{g}}=\left(0\right)_{l=0}^{M-1}$
	\FOR{\textbf{each} $\mathbf k\in I$}
	\STATE  $\hat{g}_{\mathbf k\cdot\mathbf z\bmod{M}} = \hat{g}_{\mathbf k\cdot\mathbf z\bmod{M}} + \hat{h}_{\mathbf k}$
	\ENDFOR
	\STATE $\mathbf h=\mathrm{iFFT\_1D}(\mathbf{\hat{g}})$
	\STATE $\mathbf h=M\mathbf h$
\end{algorithmic}
\begin{tabular}{p{2.0cm}p{5.2cm}p{7.0cm}}
	Output: & $\mathbf h = \mathbf A\mathbf{\hat h} = \left(h(\mathbf y_j)\, \sqrt{\frac{\omega(\mathbf y_j, \bm\mu)}{\varrho(\mathbf y_j, \bm\eta)}}\right)_{j=0}^{M-1}$ & function values of $h\in\Pi_{I,\psi(\circ,\bm\eta)}$
\end{tabular}
\end{algorithm}
\noindent
Given a frequency set $I\subset\mathbb{Z}^d$ of finite cardinality $|I|<\infty$ we consider the multivariate trigonometric polynomial $h\in\Pi_{I,\psi(\circ,\bm\eta)}$ as in \eqref{def:trig_poly_trafo_mult}
with Fourier coefficients $\hat{h}_{\mathbf k} $. The evaluation of $h$ at lattice points $\mathbf y_j \in \Lambda_{\psi(\circ,\bm\eta)}(\mathbf z, M)$ simplifies to
\begin{align*}
	h(\mathbf y_j)\, \sqrt{\frac{\omega(\mathbf y_j, \bm\mu)}{\varrho(\mathbf y_j, \bm\eta)}}
	&= \sum_{\mathbf k\in I} \hat{h}_{\mathbf k} \,\mathrm{e}^{2\pi\mathrm i \mathbf{k} \cdot \psi^{-1}(\mathbf y_j,\bm\eta)} \\
	&= \sum_{\ell = 0}^{M-1} 
		\left( \sum_{ \substack{\mathbf k\in I,\\ \mathbf k \cdot \mathbf z \equiv \ell \,(\bmod{M})} } \hat{h}_{\mathbf k}
		\right)	\,\mathrm{e}^{2\pi\mathrm i \ell \frac{j}{M}} 
	= \sum_{\ell = 0}^{M-1} \hat g_\ell \,\mathrm{e}^{2\pi\mathrm i \ell \frac{j}{M}},
\end{align*}
with
\begin{align*}
	\hat g_\ell
	= \sum_{ \substack{\mathbf k\in I,\\ \mathbf k\cdot\mathbf z \equiv \ell \,(\bmod{M})} } \hat{h}_{\mathbf k}. 
\end{align*}
In total, the evaluation of such a function is realized by simply pre-computing $(\hat g_\ell)_{\ell=0}^{M-1}$
and applying a one-dimensional inverse fast Fourier transform, see Algorithm~\ref{alg:LFFT_eval}.

\subsection{Reconstruction of transformed multivariate trigonometric polynomials}
\begin{algorithm}[tb]
\caption{Reconstruction from sampling values along a transformed reconstructing {rank-$1$} lattice} \label{alg:LFFT_recon}
\begin{tabular}{p{2.0cm}p{4.5cm}p{7.0cm}}
	Input: 	
	& $I\subset\mathbb{Z}^d$ & frequency set of finite cardinality\\
	& $M\in\mathbb{N}$ & lattice size of $\Lambda_{\psi(\circ,\bm\eta)}(\mathbf z,M,I)$\\
	& $\mathbf z\in\mathbb{Z}^d$ & generating vector of $\Lambda_{\psi(\circ,\bm\eta)}(\mathbf z,M,I)$\\
	& $\mathbf{h} = \left(h(\mathbf y_j)\, \sqrt{\frac{\omega(\mathbf y_j, \bm\mu)}{\varrho(\mathbf y_j, \bm\eta)}}\right)_{j=0}^{M-1}$ & function values of $h\in\Pi_{I,\psi(\circ,\bm\eta)}$
\end{tabular}
\begin{algorithmic}
	\STATE $\mathbf {\hat{g}}=\mathrm{FFT\_1D}(\mathbf h)$
	\FOR{\textbf{each} $\mathbf k\in I$}
	\STATE  $\hat{h}_{\mathbf k}=\frac{1}{M}\hat{g}_{\mathbf k\cdot\mathbf z\bmod{M}}$
	\ENDFOR
\end{algorithmic}
\begin{tabular}{p{2.0cm}p{4.5cm}p{7.0cm}}
	Output: & $\mathbf{\hat{h}} = M^{-1}\mathbf A^*\mathbf{h} = \left(\hat{h}_{\mathbf k}\right)_{\mathbf k\in I}$ & Fourier coefficients supported on $I$
\end{tabular}
\end{algorithm}
\noindent
For the reconstruction of a multivariate trigonometric polynomial $h\in\Pi_{I,\psi(\circ,\bm\eta)}$ as in \eqref{def:trig_poly_trafo_mult} from lattice points $\mathbf y_j\in\Lambda_{\psi(\circ,\bm\eta)}(\mathbf z, M, I)$ we utilize the exact integration property \eqref{eq:exact_integration_prop_2} and the fact that we have
\begin{align}
	\label{eq:reconR1L_matrix_prop}
	\sum_{j=0}^{M-1} \left( \mathrm{e}^{2\pi\mathrm i\frac{(\mathbf k-\mathbf h)\cdot\mathbf z}{M} } \right)^{j} = 
	\begin{cases}
		M & \text{for } \mathbf{k}\cdot\mathbf{z} \equiv \mathbf{k}\cdot\mathbf{h} \,(\bmod{M}), \\
		0 & \text{otherwise},
	\end{cases}
\end{align}
and thus $\mathbf A^* \mathbf A = M \textbf{I}$ with $\textbf{I}\in\mathbb{C}^{|I|\times|I|}$ being the identity matrix.
For fixed parameters $\bm\eta,\bm\mu\in\mathbb{R}^d$ we have input sample points of the form 
\begin{align*}
	h(\mathbf y_j)\, \sqrt{\frac{\omega(\mathbf y_j, \bm\mu)}{\varrho(\mathbf y_j, \bm\eta)}}
	= h(\psi(\mathbf x_j,\bm\eta)) \, \sqrt{ \omega(\psi(\mathbf x_j,\bm\eta), \bm\mu) \, \psi'(\mathbf x_j, \bm\eta) } 
	= f(\mathbf x_j,\bm\eta,\bm\mu) 
	= f(\mathbf x_j) .
\end{align*}
For the reconstruction of the Fourier coefficients $\hat{h}_{\mathbf k}$ we use a single one-dimensional fast Fourier transform. The entries of the resulting vector $\left( \hat g_\ell \right)_{\ell=0}^{M-1}$ are renumbered by means of the unique inverse mapping $\mathbf k \mapsto\mathbf k \cdot\mathbf z\bmod{M}$, see Algorithm~\ref{alg:LFFT_recon}.

\subsection{Discrete approximation error}
In order to use Algorithms~\ref{alg:LFFT_eval} and \ref{alg:LFFT_recon} to illustrate the proposed error bounds of Theorems~\ref{thm:L_infty_approx_error_multivar} and \ref{thm:Hm_approx_error_decay_multivar}
we sample the approximated Fourier partial sum $S_{I}^{\Lambda}h$ in order to discretize and thus approximate the error $\|h - S_{I}^{\Lambda}h\|_{L_{\infty}\big(\mathbb{R}^d, \sqrt{\frac{\omega(\circ,\bm\mu)}{\varrho(\circ,\bm\eta)}}\big)}$ that is equal to $\|f - S_{I}^{\Lambda}f\|_{L_{\infty}(\mathbb{T}^d)}$ as shown in the proof of Theorem~\ref{thm:L_infty_approx_error_multivar}.
Based on the given sample data in the vector $\mathbf{h}=\left(h(\mathbf y_j)\, \sqrt{\frac{\omega(\mathbf y_j, \bm\mu)}{\varrho(\mathbf y_j, \bm\eta)}}\right)_{j=0}^{M-1}$ with lattice points ${\mathbf y_j\in\Lambda_{\psi(\circ,\bm\eta)}(\mathbf z, M, I)}$ we apply Algorithm~\ref{alg:LFFT_recon} yielding a vector of approximated Fourier coefficients via ${\mathbf{\hat{h}} = M^{-1}\mathbf A^*\mathbf{h}}$, which we immediately put into Algorithm~\ref{alg:LFFT_eval}.
After applying both algorithms we have computed the vector $\mathbf{h}_{\mathrm{approx}} := M^{-1} \mathbf A \mathbf A^*\mathbf{h} = \left( \sqrt{\frac{\omega(\mathbf y_j, \bm\mu)}{\varrho(\mathbf y_j, \bm\eta)}}\,S_{I_N^d}^{\Lambda} h(\mathbf y_j)\right)_{j=0}^{M-1}$.

In \cite[Corollary~1]{Kae2013} it was shown under mild assumptions that for each frequency set ${I\subset\mathbb{Z}^d}$ that induces a reconstructing {rank-$1$} lattice, there is an $M\in\mathbb{N}$ such that ${|I| \leq M \lesssim |I|^2}$.
Furthermore, in \eqref{eq:reconR1L_matrix_prop} we already observed that for a reconstruction {rank-$1$} lattice $\Lambda_{\psi(\circ,\bm\eta)}(\mathbf z,M, I)$ we have $\mathbf A^* \mathbf A = M \textbf{I}$ with $\textbf{I}\in\mathbb{C}^{|I|\times|I|}$ being the identity matrix.
However, $\mathbf A \mathbf A^* \in\mathbb{C}^{M\times M}$ is generally not an identity matrix.
Hence, there is a gap between the initially given values $\mathbf h$ and the resulting vector $\mathbf h_{\mathrm{approx}}$ that we quantify with the \emph{discrete approximation error} 
\begin{align}
	\label{eq:Discrete_ellinfty_error}
	\|\mathbf h - \mathbf h_{\mathrm{approx}}\|_{\ell_{\infty}} := \max_{j=0,\ldots,M-1} \left| \sqrt{\frac{\omega(\mathbf y_j, \bm\mu)}{\varrho(\mathbf y_j, \bm\eta)}}\, \left( h(\mathbf y_j)- S_{I_N^d}^{\Lambda} h(\mathbf y_j) \right) \right|.
\end{align}
But it's important to note, that we only discuss this particular discretization approach that is exclusively sampling on the {rank-$1$} lattice nodes and doesn't measure the quality of the approximation at any point outside the {rank-$1$} lattice.
Nevertheless, for hyperbolic crosses $I_{N}^{d}$ we still have the upper bound
\begin{align}
	\label{eq:discrete_upper_bound}
	\|\mathbf h - \mathbf h_{\mathrm{approx}}\|_{\ell_{\infty}} 
	&\leq \|h - S_{I_{N}^{d}}^{\Lambda}h\|_{L_{\infty}\left(\mathbb{R}^d,\sqrt{\frac{\omega(\circ,\bm\mu)}{\varrho(\circ,\bm\eta)}}\right)}\\
	&= \|f - S_{I_{N}^{d}}^{\Lambda}f\|_{L_{\infty}(\mathbb{T}^d)}
	\leq 2 N^{-m} \|h\|_{\mathcal{H}^{m}(\mathbb{R}^d)} \nonumber
\end{align} 
for appropriately chosen parameters $\bm\eta,\bm\mu\in\mathbb{R}^d$ as shown in Theorem~\ref{thm:L_infty_approx_error_multivar}.
Hence, the theoretical results predict a certain decay rate of the discretized approximation error for increasing ${N\in\mathbb{N}}$ with fixed $m\in\mathbb{N}$ and suitably chosen parameter $\bm\eta$ and $\bm\mu$.

On the other hand, for the $L_{2}$-approximation error we lack a similar discretization approach.
However, by Theorem~\ref{thm:Hm_approx_error_decay_multivar} we know that for fixed $m\in\mathbb{N}$ and suitably chosen parameters $\bm\eta$ and $\bm\mu$ the error $\|h - S_{I_{N}^{d}}^{\Lambda}h\|_{L_{2}(\mathbb{R}^d,\omega)} = \|f - S_{I_{N}^{d}}^{\Lambda}f\|_{L_{2}(\mathbb{T}^d)}$
is bounded above by ${N^{-m}(\log N)^{(d-1)/2}\|f\|_{\mathcal{H}^{m}(\mathbb{T}^d)}}$.
By Parseval's equation we have 
\begin{align}
	\label{eq:Discrete_ell2_error}
	\|f - S_{I_{N}^{d}}^{\Lambda}f\|_{L_{2}(\mathbb{T}^d)}^2
	= \sum_{\mathbf k\in\mathbb{Z}^d} |\hat{f_{\mathbf k}} - \hat{f}_{\mathbf k}^{\Lambda}|^2 
	&= \sum_{\mathbf k\in\mathbb{Z}^d\setminus I_{N}^{d}} |\hat{f_{\mathbf k}}|^2 + \sum_{\mathbf k\in I_{N}^{d}} |\hat{f_{\mathbf k}} - \hat{f}_{\mathbf k}^{\Lambda}|^2 \nonumber\\
	&= \|f\|_{L_{2}(\mathbb{T}^d)}^2 + \sum_{\mathbf k\in I_{N}^{d}} \left( |\hat{f_{\mathbf k}} - \hat{f}_{\mathbf k}^{\Lambda}|^2 - |\hat{f_{\mathbf k}}|^2 \right).
\end{align}
Hence, we can evaluate the $L_2$-approximation error if we used Algorithm~\ref{alg:LFFT_recon} to reconstruct the approximated Fourier coefficients $\hat{f}_{\mathbf k}^{\Lambda}$ and if it is possible to calculate the Fourier coefficients $\hat{f_{\mathbf k}}$ for all ${\mathbf k\in I_{N}^{d}}$.
Later on we present an example where the Fourier coefficients $\hat{f_{\mathbf k}}$ can be computed for all ${\mathbf k \in \mathbb{Z}^d}$. 
Generally this isn't possible, so that we have to resort to the theoretical approach based on norm equivalences presented earlier in this paper in order to obtain the information if the Fourier coefficients $\hat{f_{\mathbf k}}$ are square summable.

\section{Examples}
Based on the algebraic transformation \eqref{eq:algebraic_trafo} and the error function transformations \eqref{eq:error_function_trafo}
we discuss certain choices for test functions $h$ and weight functions $\omega$ for which the proposed smoothness conditions~\eqref{eq:Hm_composition_criteria_mult} in Theorem~\ref{thm:Hm_composition_criteria_mult} are fulfilled.
In both cases we proceed similarly: 
We fix a family of multivariate weight functions $\omega(\circ,\bm\mu),\bm\mu\in\mathbb{R}^d$ as well as the test function $h$ in $L_2(\mathbb{R}^d,\omega(\circ,\bm\mu))\cap H_{\mathrm{mix}}^{m}(\mathbb{R}^d)$.
Then we fix a family of multivariate transformations $\psi(\circ,\bm\eta),\bm\eta\in\mathbb{R}^d$ of the form \eqref{eq:param_trafo_explicit_example}.
Afterwards we calculate lower bounds for $\bm\mu$ and $\bm\eta$ such that ${f(\mathbf x,\bm\eta,\bm\mu) := h(\psi(\mathbf x,\bm\eta)) \, \sqrt{ \omega(\psi(\mathbf x,\bm\eta),\bm\mu) \, \psi'(\mathbf x,\bm\eta) }}$ is in $H_{\mathrm{mix}}^m(\mathbb{T}^d)$ for Sobolev-smoothness orders $m=0,1,2,3$.
Finally we switch to dimension $d=2$ and based on the calculated parameter bounds, we use Algorithms~\ref{alg:LFFT_eval} and \ref{alg:LFFT_recon} for numerical tests of the $L_{\infty}$-approximation error bound proposed in Theorems~\ref{thm:L_infty_approx_error_multivar} and discuss the possibility to evaluate Fourier coefficients $\hat h_{\mathbf k}$.

Throughout this section we repeatedly specify parameter vectors that have the same number in each entry, for which we recall the short notation of just using a single bold number, e.g., $\mathbf 1 = (1,\ldots,1)^{\top}$ that appeared earlier in the definition of rank-$1$ lattices $\Lambda(\mathbf z,M)$ in \eqref{def:rank_one_lattice}.

\subsection{Algebraic transformation}
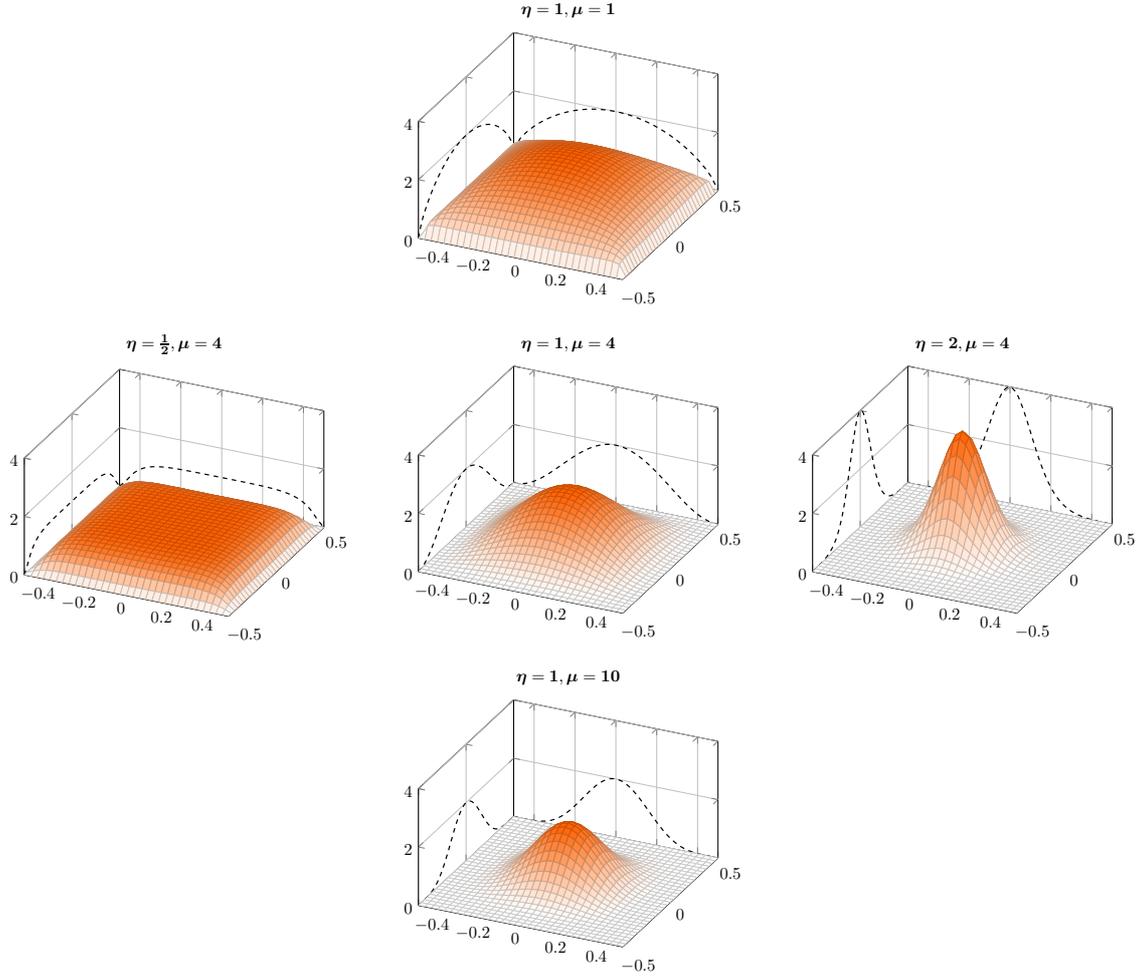
\begin{figure}[t]\begin{minipage}{.325\linewidth}
	\mbox{}
\end{minipage}
\begin{minipage}{.325\linewidth}
\centering
	\begin{tikzpicture}[
		scale=0.575,
		declare function = {eta1 = 1;},	declare function = {mu1 = 1;},
		declare function = {eta2 = 1;}, declare function = {mu2 = 1;},
		declare function = {funiv(\ETA,\MU)=(1/( 1+(4*\ETA^2*x^2)/(1-4*x^2) )) * sqrt( ((1-4*x^2)/(1-4*x^2*(1-\ETA^2)))^(\MU) * 2*\ETA*(1-4*x^2)^(-3/2) * 2*\ETA);},
		declare function = {fmult(\ETA,\ETAa,\MU,\MUa)=(1/( 1+(4*\ETA^2*x^2)/(1-4*x^2)+(4*\ETAa^2*y^2)/(1-4*y^2) )) * sqrt( ((1-4*x^2)/(1-4*x^2*(1-\ETA^2)))^(\MU) * 2*\ETA*(1-4*x^2)^(-3/2)) * sqrt( ((1-4*y^2)/(1-4*y^2*(1-\ETAa^2)))^\MUa * 2*\ETAa*(1-4*y^2)^(-3/2));}
	]
	\begin{axis}[view={25}{40},
		colormap name  = whitered,
		grid = both,
		title = {$\bm\eta=\mathbf 1, \bm\mu = \mathbf{1}$},
		samples = 33,
		zmin = 0, zmax = 4,
		legend style={at={(0.5,1.05)}, anchor=south,legend columns=2,legend cell align=left, font=\small}
	]
		\addplot3 [domain=-1/2:1/2,samples=31, samples y=0, dashed, thick, smooth]
		      (x,1/2,{funiv(eta1,mu1)});
	    \addplot3 [domain=-1/2:1/2,samples=31, samples y=0, dashed, thick, smooth]
		      (-1/2,x,{funiv(eta1,mu1)});
		\addplot3[surf,domain=-0.5:0.5, domain y=-0.5:0.5] {fmult(eta1,eta2,mu1,mu2)};
	\end{axis}
	\end{tikzpicture}
\end{minipage}
\begin{minipage}{.325\linewidth}
	\mbox{}
\end{minipage}
\par\medskip
\begin{minipage}{.325\linewidth}
	\centering
	\begin{tikzpicture}[
		scale=0.575,
		declare function = {eta1 = 0.5;}, declare function = {mu1 = 4;}, 
		declare function = {eta2 = 0.5;}, declare function = {mu2 = 4;},
		declare function = {funiv(\ETA,\MU)=(1/( 1+(4*\ETA^2*x^2)/(1-4*x^2) )) * sqrt( ((1-4*x^2)/(1-4*x^2*(1-\ETA^2)))^(\MU) * 2*\ETA*(1-4*x^2)^(-3/2)* 2*\ETA);},
		declare function = {fmult(\ETA,\ETAa,\MU,\MUa)=(1/( 1+(4*\ETA^2*x^2)/(1-4*x^2)+(4*\ETAa^2*y^2)/(1-4*y^2) )) * sqrt( ((1-4*x^2)/(1-4*x^2*(1-\ETA^2)))^(\MU) * 2*\ETA*(1-4*x^2)^(-3/2)) * sqrt( ((1-4*y^2)/(1-4*y^2*(1-\ETAa^2)))^\MUa * 2*\ETAa*(1-4*y^2)^(-3/2));}
	]
	\begin{axis}[view={25}{40},
		colormap name  = whitered,
		grid = both,
		title = {$\bm\eta=\mathbf{\frac{1}{2}}, \bm\mu = \mathbf 4$},
		samples = 33,
		zmin = 0, zmax = 4,
		legend style={at={(0.5,1.05)}, anchor=south,legend columns=2,legend cell align=left, font=\small}
	]
		\addplot3 [domain=-1/2:1/2,samples=31, samples y=0, dashed, thick, smooth]
  			(x,1/2,{funiv(eta1,mu1)});
		\addplot3 [domain=-1/2:1/2,samples=31, samples y=0, dashed, thick, smooth]
			(-1/2,x,{funiv(eta1,mu1)});
		\addplot3[surf,domain=-0.5:0.5, domain y=-0.5:0.5] {fmult(eta1,eta2,mu1,mu2)};
	\end{axis}
	\end{tikzpicture}
\end{minipage}
\begin{minipage}{.325\linewidth}
	\centering
	\begin{tikzpicture}[
		scale=0.575,
		declare function = {eta1 = 1;}, declare function = {mu1 = 4;}, 
		declare function = {eta2 = 1;}, declare function = {mu2 = 4;},
		declare function = {funiv(\ETA,\MU)=(1/( 1+(4*\ETA^2*x^2)/(1-4*x^2) )) * sqrt( ((1-4*x^2)/(1-4*x^2*(1-\ETA^2)))^(\MU) * 2*\ETA*(1-4*x^2)^(-3/2)* 2*\ETA);},
		declare function = {fmult(\ETA,\ETAa,\MU,\MUa)=(1/( 1+(4*\ETA^2*x^2)/(1-4*x^2)+(4*\ETAa^2*y^2)/(1-4*y^2) )) * sqrt( ((1-4*x^2)/(1-4*x^2*(1-\ETA^2)))^(\MU) * 2*\ETA*(1-4*x^2)^(-3/2)) * sqrt( ((1-4*y^2)/(1-4*y^2*(1-\ETAa^2)))^\MUa * 2*\ETAa*(1-4*y^2)^(-3/2));}
	]
	\begin{axis}[view={25}{40},
		colormap name  = whitered,
		grid = both,
		title = {$\bm\eta=\mathbf{1}, \bm\mu = \mathbf 4$},
		samples = 33,
		zmin = 0, zmax = 4,
		legend style={at={(0.5,1.05)}, anchor=south,legend columns=2,legend cell align=left, font=\small}
	]
		\addplot3 [domain=-1/2:1/2,samples=31, samples y=0, dashed, thick, smooth]
  			(x,1/2,{funiv(eta1,mu1)});
		\addplot3 [domain=-1/2:1/2,samples=31, samples y=0, dashed, thick, smooth]
			(-1/2,x,{funiv(eta1,mu1)});
		\addplot3[surf,domain=-0.5:0.5, domain y=-0.5:0.5] {fmult(eta1,eta2,mu1,mu2)};
	\end{axis}
	\end{tikzpicture}
\end{minipage}
\begin{minipage}{.325\linewidth}
	\centering
	\begin{tikzpicture}[
		scale=0.575,
		declare function = {eta1 = 2;}, declare function = {mu1 = 4;}, 
		declare function = {eta2 = 2;}, declare function = {mu2 = 4;},
		declare function = {funiv(\ETA,\MU)=(1/( 1+(4*\ETA^2*x^2)/(1-4*x^2) )) * sqrt( ((1-4*x^2)/(1-4*x^2*(1-\ETA^2)))^(\MU) * 2*\ETA*(1-4*x^2)^(-3/2)* 2*\ETA);},
		declare function = {fmult(\ETA,\ETAa,\MU,\MUa)=(1/( 1+(4*\ETA^2*x^2)/(1-4*x^2)+(4*\ETAa^2*y^2)/(1-4*y^2) )) * sqrt( ((1-4*x^2)/(1-4*x^2*(1-\ETA^2)))^(\MU) * 2*\ETA*(1-4*x^2)^(-3/2)) * sqrt( ((1-4*y^2)/(1-4*y^2*(1-\ETAa^2)))^\MUa * 2*\ETAa*(1-4*y^2)^(-3/2));}
	]
	\begin{axis}[view={25}{40},
		colormap name  = whitered,
		title = {$\bm\eta=\mathbf{2}, \bm\mu = \mathbf 4$},
		grid = both,
		samples = 33,
		zmin = 0, zmax = 4,
		legend style={at={(0.5,1.05)}, anchor=south,legend columns=2,legend cell align=left, font=\small}
	]
		\addplot3 [domain=-1/2:1/2,samples=31, samples y=0, dashed, thick, smooth]
  			(x,1/2,{funiv(eta1,mu1)});
		\addplot3 [domain=-1/2:1/2,samples=31, samples y=0, dashed, thick, smooth]
			(-1/2,x,{funiv(eta1,mu1)});
		\addplot3[surf,domain=-0.5:0.5, domain y=-0.5:0.5] {fmult(eta1,eta2,mu1,mu2)};
	\end{axis}
	\end{tikzpicture}
\end{minipage}
\par\medskip
\begin{minipage}{.325\linewidth}
	\mbox{}
\end{minipage}
\begin{minipage}{.325\linewidth}
\centering
	\begin{tikzpicture}[
		scale=0.575,
		declare function = {eta1 = 1;}, declare function = {mu1 = 10;}, 
		declare function = {eta2 = 1;}, declare function = {mu2 = 10;},
		declare function = {funiv(\ETA,\MU)=(1/( 1+(4*\ETA^2*x^2)/(1-4*x^2) )) * sqrt( ((1-4*x^2)/(1-4*x^2*(1-\ETA^2)))^(\MU) * 2*\ETA*(1-4*x^2)^(-3/2)* 2*\ETA);},
		declare function = {fmult(\ETA,\ETAa,\MU,\MUa)=(1/( 1+(4*\ETA^2*x^2)/(1-4*x^2)+(4*\ETAa^2*y^2)/(1-4*y^2) )) * sqrt( ((1-4*x^2)/(1-4*x^2*(1-\ETA^2)))^(\MU) * 2*\ETA*(1-4*x^2)^(-3/2)) * sqrt( ((1-4*y^2)/(1-4*y^2*(1-\ETAa^2)))^\MUa * 2*\ETAa*(1-4*y^2)^(-3/2));}
	]
	\begin{axis}[title = {$\bm\eta=\mathbf{1}, \bm\mu = \mathbf{10}$},
		view={25}{40},
		colormap name  = whitered,
		grid = both,
		samples = 33,
		zmin = 0, zmax = 4,
		legend style={at={(0.5,1.05)}, anchor=south,legend columns=2,legend cell align=left, font=\small}
	]
		\addplot3 [domain=-1/2:1/2,samples=31, samples y=0, dashed, thick, smooth]
  			(x,1/2,{funiv(eta1,mu1)});
		\addplot3 [domain=-1/2:1/2,samples=31, samples y=0, dashed, thick, smooth]
			(-1/2,x,{funiv(eta1,mu1)});
		\addplot3[surf,domain=-0.5:0.5, domain y=-0.5:0.5] {fmult(eta1,eta2,mu1,mu2)};
	\end{axis}
	\end{tikzpicture}
\end{minipage}
\begin{minipage}{.325\linewidth}
	\mbox{}
\end{minipage}
\caption{Plots of the two-dimensional transformed function $f(\circ,\bm\eta,\bm\mu)$ for various combinations of the parameters $\bm\mu$ and $\bm\eta$ with an algebraic weight function $\omega(\circ,\bm\mu)$ in \eqref{eq:alebraic_weight_function_parametrized} and the algebraic transformation $\psi(\circ,\bm\eta)$ in \eqref{eq:algebraic_trafo_parametrized}. Horizontally $\bm\mu=(4,4)^{\top}$ is fixed, vertically $\bm\eta = (1,1)^{\top}$ is fixed. The individual univariate functions $f((x_1,0),\bm\eta,\bm\mu),f((0,x_2),\bm\eta,\bm\mu)$ are shown with dashed lines.}
\label{fig:Plots_f_with_algTrafo2}
\end{figure}

The test function is of the form
\begin{align}\label{eq:exemplary_h}
	h(\mathbf y) &= \frac{1}{1 + \|\mathbf y\|_{\ell_2}^{2}}
\end{align}
with $\|\mathbf y\|_{\ell_2} := \sqrt{y_1^2+\ldots+y_d^2}$ for $\mathbf y\in\mathbb{R}^d$.
According to \cite[pp.~363-364]{boyd00} in $d=1$ this function is rather difficult to approximate by classical approximation methods.
We fix the algebraic weight function \eqref{eq:alebraic_weight_function_parametrized_univar} in its multivariate version of the form
\begin{align} \label{eq:alebraic_weight_function_parametrized}
	\omega(\mathbf y,\bm\mu) := \prod_{j=1}^{d}\left( \frac{1}{1+y_j^2} \right)^{\mu_j}
\end{align}
with ${\bm\mu = (\mu_1,\ldots,\mu_d)^{\top}\in\mathbb{R}^d}$,
and the algebraic transformation 
${\psi(\mathbf x,\bm\eta) = ((\psi_j(x_j,\eta_j))_{j=1}^{d})^{\top}}$ 
in the form \eqref{eq:param_trafo_explicit_example} with $\mathbf x\in(-\frac{1}{2},\frac{1}{2})^d$, the parameter ${\bm\eta = (\eta_1,\ldots,\eta_d)^{\top}\in\mathbb{R}^d}$ and its univariate components given by
\begin{align} \label{eq:algebraic_trafo_parametrized}
	\psi_j(x_j,\eta_j) &= \frac{2 \eta_j x_j}{(1-4x_j^2)^{\frac{1}{2}}},
	\quad \psi_j'(x_j,\eta_j) = \frac{2 \eta_j}{(1-4x_j^2)^{\frac{3}{2}}},\\
	\psi_j^{-1}(y_j,\eta_j) &= \frac{y_j}{2 (\eta_j^2+y_j^2)^{\frac{1}{2}}},
	\quad \varrho_j(y_j,\eta_j) = \frac{1}{2 (\eta_j^2+y_j^2)^{\frac{3}{2}}}. \nonumber
\end{align}
For $\eta_j = 1$ we stated the definition of $\psi_j(\circ,1)$ earlier in \eqref{eq:algebraic_trafo}.
For the resulting weighted Hilbert space $L_{2}(\mathbb{R}^d,\omega(\circ,\bm\mu))$ we have a system $\left\{\varphi_{\mathbf k}\right\}_{\mathbf k\in\mathbb{Z}^d}$ of product functions given in \eqref{eq:transformed_basis_functions_mult} with univariate components $(\varphi_{k_j})_{j=1}^{d}$ as in \eqref{eq:transformed_basis_functions} of the form
\begin{align*}
	\varphi_{k_j}(y_j,\eta_j,\mu_j)
	:= \frac{ 1 }{\sqrt{2}} (1+y_j^2)^{\frac{\mu_j}{2}} (\eta_j^2+y_j^2)^{-\frac{3}{4}} \, \mathrm{e}^{\pi\mathrm i k_j y_j (\eta_j^2+y_j^2)^{-\frac{1}{2}}},
\end{align*}
that are orthogonal with respect to the weighted scalar product
\begin{align*}
	(h_1, h_2)_{L_2\left(\mathbb{R}^d, \omega(\circ,\bm\mu) \right)}
	= \pi^{-\frac{d}{2}}\,\int_{\mathbb{R}^d} \prod_{j=1}^{d}(1+y_j^2)^{-\mu_j} \, h_1(\mathbf y) \, \overline{h_2(\mathbf y)} \, \mathrm d\mathbf y
\end{align*}
and the Fourier coefficients $\hat h_{\mathbf k}$ of an arbitrary function $h\in L_{2}(\mathbb{R}^d,\omega(\circ,\bm\mu))$ are of the form
\begin{align*}
	\hat h_{\mathbf k}
	&:= \left(h, \varphi_{\mathbf k} \right)_{L_2\left(\mathbb{R}^d,\omega(\circ,\bm \mu) \right)} \\
	&= \int_{\mathbb{R}^d} h(\mathbf y) \, \overline{\varphi_{\mathbf k}(\mathbf y,\bm\eta,\bm\mu)} \, \omega(\mathbf y,\bm\mu) \, \mathrm d\mathbf y \\
	&= 2^{-\frac{d}{2}} \int_{\mathbb{R}^d} h(\mathbf y) \, \prod_{j=1}^{d}(1+y_j^2)^{-\frac{\mu_j}{2}} (\eta_j^2+y_j^2)^{-\frac{3}{4}} \, \mathrm{e}^{-\pi\mathrm i k_j y_j (\eta_j^2+y_j^2)^{-\frac{1}{2}}} \, \mathrm d\mathbf y.
\end{align*}
The test function $h$ in \eqref{eq:exemplary_h} combined with the weight function \eqref{eq:alebraic_weight_function_parametrized} and the transformations \eqref{eq:algebraic_trafo_parametrized} lead to transformed functions $f$ in the sense of \eqref{eq:f_is_transformed_h_mult} of the form
\begin{align}
	\label{eq:Algebraic_trafo_function}
	f(\mathbf x,\bm\eta,\bm\mu)
	&= h(\psi_1(x_1,\eta_1),\ldots,\psi_d(x_d,\eta_d)) \, \prod_{j=1}^{d}\sqrt{ \omega_j(\psi_j(x_j,\eta_j),\mu_j) \, \psi_j'(x_j,\eta_j) } \nonumber\\
	&= \left( 1+ \sum_{j=1}^{d} \frac{4\eta_j^2 x_j^2}{1-4x_j^2}\right)^{-1} \, \prod_{j=1}^{d}\sqrt{\left( \frac{1-4x_j^2}{1-4(1-\eta_j^2)x_j^2} \right)^{\mu_j} 2\eta_j \left( 1-4x_j^2 \right)^{-\frac{3}{2}}}
\end{align}
In Figure~\ref{fig:Plots_f_with_algTrafo2} we have a side-by-side comparison of the graphs of these transformed functions with $d=2$ for fixed ${\bm\mu=(4,4)^{\top}}$ with varied ${\bm\eta = (\eta_1,\eta_2)\in\mathbb{R}^2}$, ${1/2\leq \eta_1,\eta_2\leq 2}$ and for fixed $\bm\eta=(1,1)^{\top}$ with varied ${\bm\mu=(\mu_1,\mu_2)^{\top}}$, ${0\leq \mu_1,\mu_2\leq 10}$.

We proceed to determine the values $\bm\eta,\bm\mu\in\mathbb{R}$ for which $f(\circ, \bm\eta, \bm\mu)$ as in \eqref{eq:Algebraic_trafo_function} is element of $H_{\mathrm{mix}}^{m}(\mathbb{T}^d)$ by investigating conditions~\eqref{eq:Hm_composition_criteria_mult} in Theorem~\ref{thm:Hm_composition_criteria_mult}.
First of all, we observe that for $\eta_1,\ldots,\eta_d > 0$ the components $\psi_1,\ldots,\psi_d$ of the function $\psi(\circ,\bm\eta)$ in \eqref{eq:algebraic_trafo_parametrized} are transformations in the sense of \eqref{def:Trafo_def} by being increasing, continuously differentiable, and invertible functions.
Furthermore, for all $\ell=1,\ldots,d$ it's easy to check that its first three derivatives of all $\psi_j(\circ,\eta_j)$ are in fact continuous on $(-\frac{1}{2},\frac{1}{2})$ for $\eta_j > 0$ and that the first three derivatives of $\varrho_j(\circ,\eta_j)$ are in $\mathcal{C}_{0}(\mathbb{R})$ for all non-zero $\eta_j\in\mathbb{R}$.
Finally, we check the $L_{\infty}$-conditions~\eqref{eq:Hm_composition_criteria_mult} in Theorem~\ref{thm:Hm_composition_criteria_mult} for $m=0,1,2,3$. 
We suppose that for $\ell=1,\ldots,d$ we have $m=m_{\ell}$ and need to check that the appearing $L_{\infty}(\mathbb{T})$-norms are finite for all $j_{\ell}=0,\ldots,m$:
\begin{itemize}
\item 
	Let $m=0$, then we only have the condition
	\begin{align*}
		&\left\| \sqrt{\omega_{\ell}(\psi_{\ell}(x_{\ell},\eta_{\ell}),\mu_{\ell}) \, \psi_{\ell}'(x_{\ell},\eta_{\ell})} \, (\psi_{\ell}'(x_{\ell},\eta_{\ell}))^{-\frac{1}{2}} \right\|_{L_{\infty}(\mathbb{T})} \\
		&= \left\| \left( \frac{1-4x_{\ell}^2}{1-4(1-\eta_{\ell}^2)x_{\ell}^2} \right)^{\mu_{\ell}} \right\|_{L_{\infty}(\mathbb{T})}
	\end{align*}
	which is finite for $\mu_{\ell} \geq 0$.
\item 
	Let $m=1$. We have to check two conditions.
	For $j_{\ell}=0$ we have
	\begin{align*}
		&\left\| \frac{\partial}{\partial x_{\ell}}\left[\sqrt{\omega_{\ell}(\psi_{\ell}(x_{\ell},\eta_{\ell}),\mu_{\ell}) \, \psi_{\ell}'(x_{\ell},\eta_{\ell})} \right] \psi_{\ell}'(x_{\ell},\eta_{\ell})^{-\frac{1}{2}} \right\|_{L_{\infty}(\mathbb{T})} \\
		&= \left\| x_{\ell}\left( \frac{(1-4x_{\ell}^2)}{1+4(\eta_{\ell}^2-1)x_{\ell}^2} \right)^{\frac{\mu_{\ell}}{2}-1} 
		\frac{(-(\mu_{\ell}^2+3)\mu_{\ell}+6(1+(\eta_{\ell}^2-1)x_{\ell}^2))}{(1+4(\eta_{\ell}^2-1)x_{\ell}^2)^{2}} 
		\right\|_{L_{\infty}(\mathbb{T})}
	\end{align*}
	and this is finite if $\mu_{\ell} > 2$.
	
	For $j_{\ell}=1$ we have
	\begin{align*}
		&\left\| \sqrt{\omega_{\ell}(\psi_{\ell}(x_{\ell},\eta_{\ell}),\mu_{\ell}) \, \psi_{\ell}'(x_{\ell},\eta_{\ell})} \, (\psi_{\ell}'(x_{\ell},\eta_{\ell}))^{\frac{1}{2}} \right\|_{L_{\infty}(\mathbb{T})} \\
		&= \left\| 2\eta_{\ell} \left(\frac{1-4x_{\ell}^2}{1+4(\eta_{\ell}^2-1)x_{\ell}^2}\right)^{\frac{\mu_{\ell}}{2}} (1-4x_{\ell}^2)^{- \frac{3}{2}}  \right\|_{L_{\infty}(\mathbb{T})}
	\end{align*}
	and this is finite for $\mu_{\ell} > 3$.
\item 
	Likewise, after checking the individual conditions we conclude that for $m=2$ we have an lower bound of $\mu_{\ell} > 9$ and for $m=3$ it is $\mu_{\ell} > 15$.
\end{itemize}
In total, $f$ is at least an $L_{2}(\mathbb{T}^d)$-function for all $\mu_{1},\ldots,\mu_{d} \geq 0$, it is at least in $H_{\mathrm{mix}}^{1}(\mathbb{T}^d)$ for ${\mu_1,\ldots,\mu_d> 3}$, at least in $H_{\mathrm{mix}}^{2}(\mathbb{T}^d)$ for ${\mu_1,\ldots,\mu_d> 9}$ and is at least an $H_{\mathrm{mix}}^{2}(\mathbb{T}^d)$-function for ${\mu_1,\ldots,\mu_d > 15}$.
Apparently, the parameters $\eta_1,\ldots,\eta_d$ in the transformation $\psi(\circ,\bm\eta)$ don't have an impact on the Sobolev-smoothness of $f(\circ, \bm\eta, \bm\mu)$ as in \eqref{eq:Algebraic_trafo_function}, according to this specific set of conditions.
In other words, if $\bm\eta$ is able to control the smoothness of $f$ then we can't recognize it with these conditions -- at least for this particular combination of transformation $\psi$ and weight function $\omega$.
	
\subsubsection{$L_{\infty}$-approximation error discussion}
	Next we discuss the application of the multivariate $L_{\infty}(\mathbb{R}^d)$-approximation error bound in Theorem~\ref{thm:L_infty_approx_error_multivar} for $d=2$ with the two-dimensional test function $h$ in \eqref{eq:exemplary_h}, the weight function \eqref{eq:alebraic_weight_function_parametrized}, the transformations \eqref{eq:algebraic_trafo_parametrized} and the resulting transformed functions $f$ given in \eqref{eq:Algebraic_trafo_function}.
	
\begin{figure}[t]
	\centering
		\begin{tikzpicture}[baseline,scale=0.75]
		\begin{axis}[
		ymode = log,
		enlargelimits=false,
		xmin=0, xmax=90, ymin=1e-12, ymax=1e-0,
		ytick={1e-1,1e-2,1e-3,1e-4,1e-5,1e-6,1e-7,1e-8,1e-9,1e-10,1e-11},
		grid=both, 
		xlabel={$N$}, 
		ylabel={$\|\mathbf h - \mathbf h_{\mathrm{approx}}\|_{\ell_\infty}/\|\mathbf h\|_{\ell_\infty}$},
		legend style={at={(0.5,1.05)}, anchor=south,legend columns=4,legend cell align=left, font=\small,  
		},
		xminorticks=false,
		yminorticks=false
		]
		\addplot[mark options={solid}, blue, mark=o, mark size=1.0] coordinates {
			(8, 1.1118e-02)  (9, 1.5324e-02)  (10, 1.0303e-02)  (11, 1.3566e-02)  (12, 1.3373e-02)  (13, 1.1843e-02)  (14, 1.2930e-02)  (15, 1.2922e-02)  (16, 1.2707e-02)  (17, 1.4225e-02)  (18, 1.4076e-02)  (19, 1.5402e-02)  (20, 1.4482e-02)  (21, 1.4720e-02)  (22, 1.5578e-02)  (23, 1.6645e-02)  (24, 1.5438e-02)  (25, 1.5651e-02)  (26, 1.6390e-02)  (27, 1.6741e-02)  (28, 1.6374e-02)  (29, 1.7149e-02)  (30, 1.6310e-02)  (31, 1.7005e-02)  (32, 1.6784e-02)  (33, 1.7100e-02)  (34, 1.7651e-02)  (35, 1.7381e-02)  (36, 1.6611e-02)  (37, 1.7146e-02)  (38, 1.7623e-02)  (39, 1.7918e-02)  (40, 1.7173e-02)  (41, 1.7637e-02)  (42, 1.7204e-02)  (43, 1.7637e-02)  (44, 1.7615e-02)  (45, 1.7345e-02)  (46, 1.7709e-02)  (47, 1.8090e-02)  (48, 1.7438e-02)  (49, 1.7506e-02)  (50, 1.7387e-02)  (51, 1.7619e-02)  (52, 1.7646e-02)  (53, 1.7959e-02)  (54, 1.7730e-02)  (55, 1.7661e-02)  (56, 1.7233e-02)  (57, 1.7438e-02)  (58, 1.7700e-02)  (59, 1.7965e-02)  (60, 1.7287e-02)  (61, 1.7536e-02)  (62, 1.7772e-02)  (63, 1.7577e-02)  (64, 1.7436e-02)  (65, 1.7407e-02)  (66, 1.7269e-02)  (67, 1.7483e-02)  (68, 1.7541e-02)  (69, 1.7705e-02)  (70, 1.7363e-02)  (71, 1.7559e-02)  (72, 1.7009e-02)  (73, 1.7194e-02)  (74, 1.7373e-02)  (75, 1.7337e-02)  (76, 1.7397e-02)  (77, 1.7297e-02)  (78, 1.7217e-02)  (79, 1.7380e-02)  (80, 1.7007e-02)  
		};
		\addlegendentry{$\bm\mu=\mathbf{0}$}
		\addplot[mark options={solid},red!75!yellow,mark=x,mark size=1.0] coordinates {
			(8, 3.6615e-04)  (9, 1.9107e-04)  (10, 2.2361e-04)  (11, 2.3331e-04)  (12, 1.3640e-04)  (13, 1.4013e-04)  (14, 1.3856e-04)  (15, 1.0883e-04)  (16, 8.6951e-05)  (17, 8.4319e-05)  (18, 7.1837e-05)  (19, 6.8143e-05)  (20, 5.2474e-05)  (21, 4.7983e-05)  (22, 4.3797e-05)  (23, 4.1023e-05)  (24, 3.3485e-05)  (25, 2.9990e-05)  (26, 2.7120e-05)  (27, 2.5143e-05)  (28, 2.1663e-05)  (29, 2.1916e-05)  (30, 1.7363e-05)  (31, 1.7377e-05)  (32, 1.6245e-05)  (33, 1.6825e-05)  (34, 1.6806e-05)  (35, 1.4571e-05)  (36, 1.2332e-05)  (37, 1.2686e-05)  (38, 1.2570e-05)  (39, 1.2616e-05)  (40, 1.1019e-05)  (41, 1.1127e-05)  (42, 9.3268e-06)  (43, 9.4248e-06)  (44, 9.1777e-06)  (45, 8.5019e-06)  (46, 8.3216e-06)  (47, 8.3058e-06)  (48, 7.2598e-06)  (49, 6.8762e-06)  (50, 6.4420e-06)  (51, 6.3180e-06)  (52, 6.1550e-06)  (53, 6.1038e-06)  (54, 5.4997e-06)  (55, 5.3112e-06)  (56, 4.8055e-06)  (57, 4.6888e-06)  (58, 4.5613e-06)  (59, 4.5069e-06)  (60, 4.0671e-06)  (61, 4.0666e-06)  (62, 4.0981e-06)  (63, 3.6779e-06)  (64, 3.4619e-06)  (65, 3.4909e-06)  (66, 3.2849e-06)  (67, 3.3274e-06)  (68, 3.3236e-06)  (69, 3.3151e-06)  (70, 3.0600e-06)  (71, 3.0868e-06)  (72, 2.6636e-06)  (73, 2.6882e-06)  (74, 2.6888e-06)  (75, 2.6691e-06)  (76, 2.6484e-06)  (77, 2.5184e-06)  (78, 2.4269e-06)  (79, 2.4348e-06)  (80, 2.2474e-06)
		};
		\addlegendentry{$\bm\mu=\mathbf{4}$}
		\addplot[mark options={solid},red!25!yellow,mark=triangle,mark size=1.5] coordinates {
			(8, 3.6670e-04)  (9, 1.1495e-04)  (10, 5.2275e-05)  (11, 5.9904e-05)  (12, 1.5702e-05)  (13, 1.7951e-05)  (14, 1.0685e-05)  (15, 9.5958e-06)  (16, 3.4796e-06)  (17, 3.3378e-06)  (18, 3.0850e-06)  (19, 3.2320e-06)  (20, 1.3955e-06)  (21, 1.1829e-06)  (22, 1.1057e-06)  (23, 1.0970e-06)  (24, 7.5623e-07)  (25, 4.2469e-07)  (26, 4.1598e-07)  (27, 4.2516e-07)  (28, 3.7053e-07)  (29, 3.7531e-07)  (30, 1.7156e-07)  (31, 1.7585e-07)  (32, 1.6591e-07)  (33, 1.6041e-07)  (34, 1.5594e-07)  (35, 1.0893e-07)  (36, 7.2137e-08)  (37, 7.1777e-08)  (38, 6.9894e-08)  (39, 6.6598e-08)  (40, 5.4009e-08)  (41, 5.2922e-08)  (42, 3.1381e-08)  (43, 3.0908e-08)  (44, 3.0337e-08)  (45, 2.5819e-08)  (46, 2.5377e-08)  (47, 2.6231e-08)  (48, 1.8193e-08)  (49, 1.3945e-08)  (50, 1.3368e-08)  (51, 1.3428e-08)  (52, 1.3492e-08)  (53, 1.3744e-08)  (54, 1.1045e-08)  (55, 1.0997e-08)  (56, 7.1413e-09)  (57, 7.1160e-09)  (58, 7.1700e-09)  (59, 7.2363e-09)  (60, 6.1380e-09)  (61, 6.1642e-09)  (62, 6.1555e-09)  (63, 4.6698e-09)  (64, 3.8012e-09)  (65, 3.8125e-09)  (66, 3.4448e-09)  (67, 3.4423e-09)  (68, 3.3800e-09)  (69, 3.3386e-09)  (70, 2.7721e-09)  (71, 2.7631e-09)  (72, 1.9144e-09)  (73, 1.9103e-09)  (74, 1.8977e-09)  (75, 1.8733e-09)  (76, 1.8340e-09)  (77, 1.6077e-09)  (78, 1.5371e-09)  (79, 1.5228e-09)  (80, 1.2184e-09)  
		};
		\addlegendentry{$\bm\mu=\mathbf{10}$}
		\addplot[mark options={solid},red!50!blue,mark=square,mark size=1.5] coordinates {
			(8, 2.9163e-03)  (9, 4.0847e-04)  (10, 4.0898e-04)  (11, 4.0916e-04)  (12, 5.6782e-05)  (13, 5.9886e-05)  (14, 5.9966e-05)  (15, 2.9371e-05)  (16, 5.5110e-06)  (17, 5.6882e-06)  (18, 1.9923e-06)  (19, 2.0049e-06)  (20, 9.7161e-07)  (21, 4.7891e-07)  (22, 4.8724e-07)  (23, 4.9044e-07)  (24, 4.4241e-07)  (25, 1.1259e-07)  (26, 1.1880e-07)  (27, 1.1707e-07)  (28, 8.9803e-08)  (29, 8.9517e-08)  (30, 3.4030e-08)  (31, 3.2195e-08)  (32, 2.1114e-08)  (33, 2.1408e-08)  (34, 2.0986e-08)  (35, 1.7745e-08)  (36, 5.7926e-09)  (37, 5.7837e-09)  (38, 5.6973e-09)  (39, 5.7283e-09)  (40, 4.8147e-09)  (41, 4.8650e-09)  (42, 1.9257e-09)  (43, 1.8528e-09)  (44, 1.8420e-09)  (45, 1.3845e-09)  (46, 1.4127e-09)  (47, 1.4005e-09)  (48, 1.1646e-09)  (49, 7.4199e-10)  (50, 5.1250e-10)  (51, 4.6856e-10)  (52, 4.9814e-10)  (53, 5.0357e-10)  (54, 4.5854e-10)  (55, 3.6542e-10)  (56, 2.2103e-10)  (57, 1.8707e-10)  (58, 1.9305e-10)  (59, 1.9417e-10)  (60, 1.3137e-10)  (61, 1.3267e-10)  (62, 1.3180e-10)  (63, 1.0243e-10)  (64, 7.8147e-11)  (65, 6.3028e-11)  (66, 5.2974e-11)  (67, 5.3022e-11)  (68, 5.2408e-11)  (69, 5.7202e-11)  (70, 4.6650e-11)  (71, 4.6873e-11)  (72, 2.4918e-11)  (73, 2.4666e-11)  (74, 2.4597e-11)  (75, 2.1724e-11)  (76, 2.1952e-11)  (77, 2.1898e-11)  (78, 1.9757e-11)  (79, 1.9943e-11)  (80, 1.4026e-11)       
		};
		\addlegendentry{$\bm\mu=\mathbf{16}$}
		\end{axis}
		\end{tikzpicture}
	\caption{Comparison of discrete $\ell_\infty$-approximation error $\| \mathbf h - \mathbf h_{\text{approx}} \|_{\ell_\infty}/\|\mathbf h\|_{\ell_\infty}$ of two-dimensional test function \eqref{eq:exemplary_h} in combination with 
		the algebraic transformation $\psi(\circ,\bm\eta)$ \eqref{eq:algebraic_trafo_parametrized} 
		and the algebraic weight function $\omega(\circ,\bm\mu)$ \eqref{eq:alebraic_weight_function_parametrized}
		in their two-dimensional versions with fixed $\bm\eta = \mathbf 1$ and $\bm\mu \in \{ \mathbf 0, \mathbf 4, \mathbf{10}, \mathbf{16}\}$ .}
	\label{fig:numeric_ell2error_algebraic}
\end{figure}
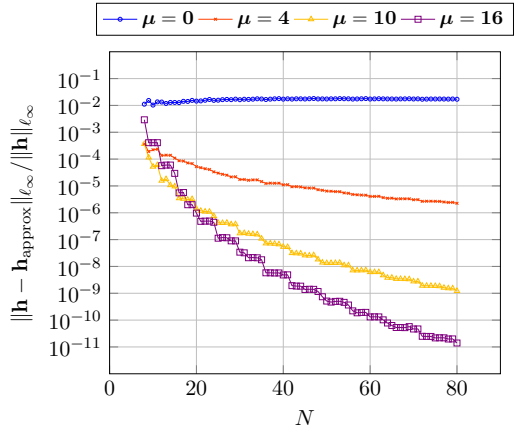

	Let a reconstructing {rank-$1$} lattice $\Lambda(\mathbf z, M, I_{N}^{2})$ with $N\geq 8$ be given. 
	We already evaluated the sufficient conditions proposed in Theorem~\ref{thm:Hm_composition_criteria_mult}, yielding lower bounds for $\bm\mu\geq \mathbf 0$ such that $f$ is at least of Sobolev-smoothness order $m=0,1,2,3$, i.e., $f\in H_{\mathrm{mix}}^{m}(\mathbb{T}^2)$ and thus $f\in \mathcal{H}^{m}(\mathbb{T}^2)$.	
	We fix $\lambda=1$ and for ${m\in\mathbb{N}_{0}}$ we choose $\bm\mu,\bm\eta\in\mathbb{R}^2$ such that $f\in \mathcal{H}^{m+1}(\mathbb{T}^2)\hookrightarrow \mathcal{A}^{m}(\mathbb{T}^2)$.
	As outlined in \eqref{eq:discrete_upper_bound} we expect the discrete approximation error \eqref{eq:Discrete_ellinfty_error} to be bounded by
	\begin{align} \label{eq:discrete_bounds_exact}
		\|\mathbf h - \mathbf h_{\mathrm{approx}}\|_{\ell_{\infty}} 
		\leq \|f - S_{I_N^2}^{\Lambda} f\|_{L_{\infty}(\mathbb{T}^2)}
		\lesssim 
		\begin{cases}
			N^{0} \quad &\text{ for } \quad \mu_j\geq 0,\\
			N^{-1} \quad &\text{ for } \quad \mu_j > 3, \\
			N^{-2} \quad &\text{ for } \quad \mu_j > 9, \\
			N^{-3} \quad &\text{ for } \quad \mu_j > 15.
		\end{cases}
	\end{align}
	For $N=8,\ldots,80$, $\bm\eta=\mathbf 1$ and $\bm\mu\in\{\mathbf 0,\mathbf{4},\mathbf{10},\mathbf{16}\}$ we actually observe this behavior for the relative discrete approximation error $\|\mathbf h - \mathbf h_{\mathrm{approx}}\|_{\ell_{\infty}}/\|\mathbf h\|_{\ell_{\infty}}$ as seen in the left plot of Figure~\ref{fig:numeric_ell2error_algebraic}.
	
\subsubsection{$L_{2}$-approximation error discussion}
	We switch to dimension $d=1$. 
	In Theorem~\ref{thm:Hm_approx_error_decay_multivar} we proved that when $f$ of the form \eqref{eq:f_is_transformed_h} is in $\mathcal{H}^{m}(\mathbb{T})\cap\mathcal{C}(\mathbb{T})$ we have
	\begin{align*}
		\|h - S_{I_{N}^{1}}^{\Lambda}h\|_{L_{2}\left(\mathbb{R}, \omega \right)} 
		=\|f - S_{I_{N}^{1}}^{\Lambda}f\|_{L_{2}\left(\mathbb{T}\right)} 
		\lesssim N^{-m}.
	\end{align*}
	For one particular special case with explicitly computable Fourier coefficients $\hat f_k$ we observe that their rate of decay is consistent with the theoretical propositions.
	The conditions of Theorem~\ref{thm:Hm_composition_criteria} yielded for $m=1,2,3$
	\begin{align}
		\label{eq:Param_bounds_alg_general}
		f \in
		\begin{cases}
		 	H^{1}(\mathbb{T}) \quad &\text{ for } \quad \mu > 3,\\
			H^{2}(\mathbb{T}) \quad &\text{ for } \quad \mu > 9,\\
			H^{3}(\mathbb{T}) \quad &\text{ for } \quad \mu > 15.
		\end{cases}
	\end{align}
	We compare these lower bounds with the specific lower bounds for the chosen $h$ in \eqref{eq:exemplary_h}.
	Fixing $\eta = 1$ the transformed function $f$ in \eqref{eq:Algebraic_trafo_function} simplifies to
	\begin{align*}
		f(x,1,\mu) = f(x,\mu) := (1-4x^2)^{\frac{1}{2}(\mu+\frac{1}{2})}.
	\end{align*}
	We then explicitly calculate that
	\begin{align*}
		\int_{\mathbb{T}} \left|\frac{\mathrm{d}}{\mathrm{d}x}[f](x,\mu)\right|^2 \,\mathrm{d}x
		= \left(\mu+\frac{1}{2}\right)^2 \,\int_{\mathbb{T}} 16x^2 \left|1-4x^2\right|^{\mu-\frac{3}{2}} \,\mathrm{d}x
		< \infty
	\end{align*}
	for $\mu\geq \frac{3}{2}$, as well as
	\begin{align*}
		\int_{\mathbb{T}} \left|\frac{\mathrm{d}^2}{\mathrm{d}x^2}[f](x,\mu)\right|^2 \,\mathrm{d}x
		= 4\left(2\mu+1\right)^2 \, \int_{\mathbb{T}} \left(1-(4\mu - 2)x^2\right)^{2} \, \left|1-4x^2\right|^{\mu-\frac{7}{2}} \,\mathrm{d}x
		< \infty
	\end{align*}
	for $\mu\geq \frac{7}{2}$ and so forth, which is summarized for $m=1,2,3$ as
	\begin{align}\label{eq:Param_bounds_alg_example}
		f \in
		\begin{cases}
			H^{1}(\mathbb{T}) \quad &\text{ for } \quad \mu \geq \frac{3}{2},\\
			H^{2}(\mathbb{T}) \quad &\text{ for } \quad \mu \geq \frac{7}{2},\\
			H^{3}(\mathbb{T}) \quad &\text{ for } \quad \mu \geq \frac{11}{2}.
		\end{cases}
	\end{align}
	Due to the norm equivalence \eqref{eq:Hs_norm_equivalence} we know that the absolute Fourier coefficients $|\hat{f}_{k}|$ of a function $f\in H_{\mathrm{mix}}^{m}(\mathbb{T})$ decay at least as fast as $|k|^{-m}$.
	In our particular example with the above $f(x,\mu)$ we have a decay twice as fast, which is observed by considering $k\in\mathbb{Z}\setminus\{0\}$ and calculating that
	\begin{align*}
		|\hat{f}_{k}|
		= 
		\sqrt{2}\cdot
		\begin{cases}
			\displaystyle \int_{\mathbb{T}} (1-4x^2) \, \mathrm{e}^{-2\pi\mathrm{i}k x} \, \mathrm dx = \frac{2 \sqrt{2}}{\pi^2 |k|^2} \quad &\text{ for } \quad \mu = \frac{3}{2},\\
			\displaystyle \int_{\mathbb{T}} (1-4x^2)^2 \, \mathrm{e}^{-2\pi\mathrm{i}k x} \, \mathrm dx = \frac{24 \sqrt{2}}{\pi^4 |k|^4} \quad &\text{ for } \quad \mu = \frac{7}{2},\\
			\displaystyle \int_{\mathbb{T}} (1-4x^2)^3 \, \mathrm{e}^{-2\pi\mathrm{i}k x} \, \mathrm dx = \frac{48 \sqrt{2}|\pi^2 |k|^2-15|}{\pi^6 |k|^6} \quad &\text{ for } \quad \mu = \frac{11}{2}.
		\end{cases}
	\end{align*}
	
	The general $L_{\infty}$-parameter bounds in \eqref{eq:Param_bounds_alg_general} look relatively coarse in comparison to the exact bounds in \eqref{eq:Param_bounds_alg_example}.
	However, generally we can't compute the Fourier coefficients $\hat f_{k}$ of a transformed function $f$, which makes the conditions proposed in Theorem~\ref{thm:Hm_composition_criteria} so powerful, as they work independent of the particular choice of ${h\in L_{2}(\mathbb{R},\omega)\cap H^{m}(\mathbb{R})}$ for the cost of yielding not the most precise lower parameter bounds.

\subsection{Error~function transformation}
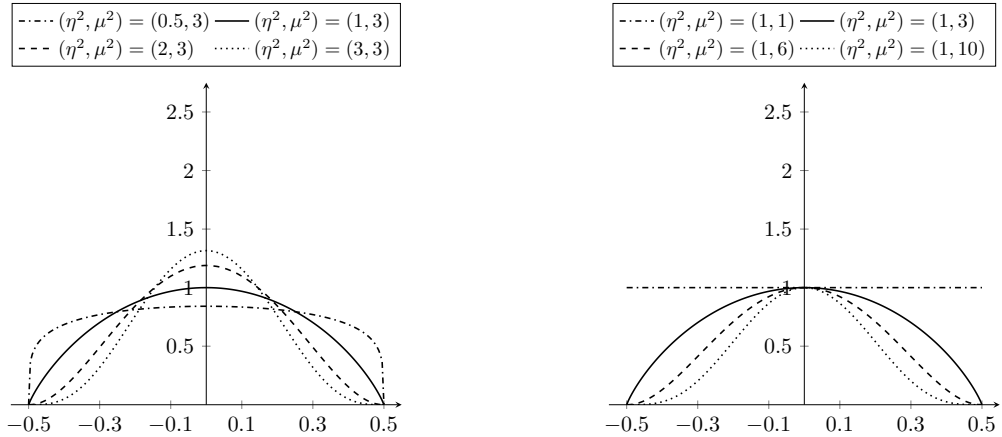
\begin{figure}[t]
\begin{minipage}{.5\linewidth}
	\centering
\begin{tikzpicture}[scale=0.75]
	\pgfmathsetmacro\MU{3} \pgfmathsetmacro\ETA{1/2}
	\pgfmathsetmacro\MUa{3} \pgfmathsetmacro\ETAa{1}
	\pgfmathsetmacro\MUb{3} \pgfmathsetmacro\ETAb{2}
	\pgfmathsetmacro\MUc{3} \pgfmathsetmacro\ETAc{3}
	\begin{axis}[samples=100,  
	xtick={-0.5, -0.3, -0.1, 0, 0.1, 0.3, 0.5},
xmin=-0.55, xmax=0.55, ymin=0, ymax=2.75,
	axis x line=center, axis y line=center,
	every axis plot/.append style={thick},
	legend style={at={(0.5,1.05)}, anchor=south,legend columns=2,legend cell align=left, font=\small}
	]
\addplot[dashdotted, mark size = 2] coordinates { 
		(-0.5, 0.000e+00) (-0.495, 3.669e-01) (-0.49, 4.275e-01) (-0.485, 4.668e-01) (-0.48, 4.963e-01) (-0.475, 5.202e-01) (-0.47, 5.404e-01) (-0.465, 5.579e-01) (-0.46, 5.733e-01) (-0.455, 5.871e-01) (-0.45, 5.996e-01) (-0.445, 6.111e-01) (-0.44, 6.216e-01) (-0.435, 6.314e-01) (-0.43, 6.405e-01) (-0.425, 6.490e-01) (-0.42, 6.570e-01) (-0.415, 6.645e-01) (-0.41, 6.717e-01) (-0.405, 6.784e-01) (-0.4, 6.848e-01) (-0.395, 6.909e-01) (-0.39, 6.967e-01) (-0.385, 7.023e-01) (-0.38, 7.076e-01) (-0.375, 7.127e-01) (-0.37, 7.176e-01) (-0.365, 7.223e-01) (-0.36, 7.267e-01) (-0.355, 7.311e-01) (-0.35, 7.352e-01) (-0.345, 7.392e-01) (-0.34, 7.431e-01) (-0.335, 7.468e-01) (-0.33, 7.504e-01) (-0.325, 7.539e-01) (-0.32, 7.573e-01) (-0.315, 7.605e-01) (-0.31, 7.637e-01) (-0.305, 7.667e-01) (-0.3, 7.696e-01) (-0.295, 7.725e-01) (-0.29, 7.752e-01) (-0.285, 7.779e-01) (-0.28, 7.805e-01) (-0.275, 7.830e-01) (-0.27, 7.854e-01) (-0.265, 7.878e-01) (-0.26, 7.901e-01) (-0.255, 7.923e-01) (-0.25, 7.944e-01) (-0.245, 7.965e-01) (-0.24, 7.985e-01) (-0.235, 8.004e-01) (-0.23, 8.023e-01) (-0.225, 8.042e-01) (-0.22, 8.059e-01) (-0.215, 8.077e-01) (-0.21, 8.093e-01) (-0.205, 8.109e-01) (-0.2, 8.125e-01) (-0.195, 8.140e-01) (-0.19, 8.154e-01) (-0.185, 8.169e-01) (-0.18, 8.182e-01) (-0.175, 8.195e-01) (-0.17, 8.208e-01) (-0.165, 8.220e-01) (-0.16, 8.232e-01) (-0.155, 8.243e-01) (-0.15, 8.254e-01) (-0.145, 8.265e-01) (-0.14, 8.275e-01) (-0.135, 8.285e-01) (-0.13, 8.294e-01) (-0.125, 8.303e-01) (-0.12, 8.311e-01) (-0.115, 8.320e-01) (-0.11, 8.327e-01) (-0.105, 8.335e-01) (-0.1, 8.342e-01) (-0.095, 8.348e-01) (-0.09, 8.355e-01) (-0.085, 8.361e-01) (-0.08, 8.366e-01) (-0.075, 8.371e-01) (-0.07, 8.376e-01) (-0.065, 8.381e-01) (-0.06, 8.385e-01) (-0.055, 8.389e-01) (-0.05, 8.392e-01) (-0.045, 8.396e-01) (-0.04, 8.398e-01) (-0.035, 8.401e-01) (-0.03, 8.403e-01) (-0.025, 8.405e-01) (-0.02, 8.406e-01) (-0.015, 8.407e-01) (-0.01, 8.408e-01) (-0.005, 8.409e-01) (0, 8.409e-01) (0.005, 8.409e-01) (0.01, 8.408e-01) (0.015, 8.407e-01) (0.02, 8.406e-01) (0.025, 8.405e-01) (0.03, 8.403e-01) (0.035, 8.401e-01) (0.04, 8.398e-01) (0.045, 8.396e-01) (0.05, 8.392e-01) (0.055, 8.389e-01) (0.06, 8.385e-01) (0.065, 8.381e-01) (0.07, 8.376e-01) (0.075, 8.371e-01) (0.08, 8.366e-01) (0.085, 8.361e-01) (0.09, 8.355e-01) (0.095, 8.348e-01) (0.1, 8.342e-01) (0.105, 8.335e-01) (0.11, 8.327e-01) (0.115, 8.320e-01) (0.12, 8.311e-01) (0.125, 8.303e-01) (0.13, 8.294e-01) (0.135, 8.285e-01) (0.14, 8.275e-01) (0.145, 8.265e-01) (0.15, 8.254e-01) (0.155, 8.243e-01) (0.16, 8.232e-01) (0.165, 8.220e-01) (0.17, 8.208e-01) (0.175, 8.195e-01) (0.18, 8.182e-01) (0.185, 8.169e-01) (0.19, 8.154e-01) (0.195, 8.140e-01) (0.2, 8.125e-01) (0.205, 8.109e-01) (0.21, 8.093e-01) (0.215, 8.077e-01) (0.22, 8.059e-01) (0.225, 8.042e-01) (0.23, 8.023e-01) (0.235, 8.004e-01) (0.24, 7.985e-01) (0.245, 7.965e-01) (0.25, 7.944e-01) (0.255, 7.923e-01) (0.26, 7.901e-01) (0.265, 7.878e-01) (0.27, 7.854e-01) (0.275, 7.830e-01) (0.28, 7.805e-01) (0.285, 7.779e-01) (0.29, 7.752e-01) (0.295, 7.725e-01) (0.3, 7.696e-01) (0.305, 7.667e-01) (0.31, 7.637e-01) (0.315, 7.605e-01) (0.32, 7.573e-01) (0.325, 7.539e-01) (0.33, 7.504e-01) (0.335, 7.468e-01) (0.34, 7.431e-01) (0.345, 7.392e-01) (0.35, 7.352e-01) (0.355, 7.311e-01) (0.36, 7.267e-01) (0.365, 7.223e-01) (0.37, 7.176e-01) (0.375, 7.127e-01) (0.38, 7.076e-01) (0.385, 7.023e-01) (0.39, 6.967e-01) (0.395, 6.909e-01) (0.4, 6.848e-01) (0.405, 6.784e-01) (0.41, 6.717e-01) (0.415, 6.645e-01) (0.42, 6.570e-01) (0.425, 6.490e-01) (0.43, 6.405e-01) (0.435, 6.314e-01) (0.44, 6.216e-01) (0.445, 6.111e-01) (0.45, 5.996e-01) (0.455, 5.871e-01) (0.46, 5.733e-01) (0.465, 5.579e-01) (0.47, 5.404e-01) (0.475, 5.202e-01) (0.48, 4.963e-01) (0.485, 4.668e-01) (0.49, 4.275e-01) (0.495, 3.669e-01) (0.5, 0.000e+00)
	};
	\addlegendentry{$(\eta^2,\mu^2)=(\ETA,\MU)$};
	\addplot[mark size = 2] coordinates { 
		(-0.5, 0.000e+00) (-0.495, 3.625e-02) (-0.49, 6.681e-02) (-0.485, 9.493e-02) (-0.48, 1.214e-01) (-0.475, 1.465e-01) (-0.47, 1.706e-01) (-0.465, 1.937e-01) (-0.46, 2.160e-01) (-0.455, 2.376e-01) (-0.45, 2.585e-01) (-0.445, 2.788e-01) (-0.44, 2.986e-01) (-0.435, 3.178e-01) (-0.43, 3.366e-01) (-0.425, 3.548e-01) (-0.42, 3.727e-01) (-0.415, 3.901e-01) (-0.41, 4.071e-01) (-0.405, 4.237e-01) (-0.4, 4.399e-01) (-0.395, 4.558e-01) (-0.39, 4.713e-01) (-0.385, 4.865e-01) (-0.38, 5.014e-01) (-0.375, 5.160e-01) (-0.37, 5.303e-01) (-0.365, 5.442e-01) (-0.36, 5.579e-01) (-0.355, 5.713e-01) (-0.35, 5.844e-01) (-0.345, 5.973e-01) (-0.34, 6.099e-01) (-0.335, 6.222e-01) (-0.33, 6.343e-01) (-0.325, 6.461e-01) (-0.32, 6.577e-01) (-0.315, 6.691e-01) (-0.31, 6.802e-01) (-0.305, 6.911e-01) (-0.3, 7.018e-01) (-0.295, 7.122e-01) (-0.29, 7.224e-01) (-0.285, 7.324e-01) (-0.28, 7.422e-01) (-0.275, 7.518e-01) (-0.27, 7.611e-01) (-0.265, 7.703e-01) (-0.26, 7.792e-01) (-0.255, 7.880e-01) (-0.25, 7.965e-01) (-0.245, 8.049e-01) (-0.24, 8.131e-01) (-0.235, 8.210e-01) (-0.23, 8.288e-01) (-0.225, 8.364e-01) (-0.22, 8.438e-01) (-0.215, 8.510e-01) (-0.21, 8.580e-01) (-0.205, 8.649e-01) (-0.2, 8.715e-01) (-0.195, 8.780e-01) (-0.19, 8.843e-01) (-0.185, 8.904e-01) (-0.18, 8.964e-01) (-0.175, 9.022e-01) (-0.17, 9.078e-01) (-0.165, 9.132e-01) (-0.16, 9.185e-01) (-0.155, 9.235e-01) (-0.15, 9.285e-01) (-0.145, 9.332e-01) (-0.14, 9.378e-01) (-0.135, 9.422e-01) (-0.13, 9.464e-01) (-0.125, 9.505e-01) (-0.12, 9.544e-01) (-0.115, 9.582e-01) (-0.11, 9.617e-01) (-0.105, 9.652e-01) (-0.1, 9.684e-01) (-0.095, 9.715e-01) (-0.09, 9.744e-01) (-0.085, 9.772e-01) (-0.08, 9.798e-01) (-0.075, 9.823e-01) (-0.07, 9.846e-01) (-0.065, 9.867e-01) (-0.06, 9.887e-01) (-0.055, 9.905e-01) (-0.05, 9.921e-01) (-0.045, 9.936e-01) (-0.04, 9.950e-01) (-0.035, 9.961e-01) (-0.03, 9.972e-01) (-0.025, 9.980e-01) (-0.02, 9.987e-01) (-0.015, 9.993e-01) (-0.01, 9.997e-01) (-0.005, 9.999e-01) (0, 1.000e+00) (0.005, 9.999e-01) (0.01, 9.997e-01) (0.015, 9.993e-01) (0.02, 9.987e-01) (0.025, 9.980e-01) (0.03, 9.972e-01) (0.035, 9.961e-01) (0.04, 9.950e-01) (0.045, 9.936e-01) (0.05, 9.921e-01) (0.055, 9.905e-01) (0.06, 9.887e-01) (0.065, 9.867e-01) (0.07, 9.846e-01) (0.075, 9.823e-01) (0.08, 9.798e-01) (0.085, 9.772e-01) (0.09, 9.744e-01) (0.095, 9.715e-01) (0.1, 9.684e-01) (0.105, 9.652e-01) (0.11, 9.617e-01) (0.115, 9.582e-01) (0.12, 9.544e-01) (0.125, 9.505e-01) (0.13, 9.464e-01) (0.135, 9.422e-01) (0.14, 9.378e-01) (0.145, 9.332e-01) (0.15, 9.285e-01) (0.155, 9.235e-01) (0.16, 9.185e-01) (0.165, 9.132e-01) (0.17, 9.078e-01) (0.175, 9.022e-01) (0.18, 8.964e-01) (0.185, 8.904e-01) (0.19, 8.843e-01) (0.195, 8.780e-01) (0.2, 8.715e-01) (0.205, 8.649e-01) (0.21, 8.580e-01) (0.215, 8.510e-01) (0.22, 8.438e-01) (0.225, 8.364e-01) (0.23, 8.288e-01) (0.235, 8.210e-01) (0.24, 8.131e-01) (0.245, 8.049e-01) (0.25, 7.965e-01) (0.255, 7.880e-01) (0.26, 7.792e-01) (0.265, 7.703e-01) (0.27, 7.611e-01) (0.275, 7.518e-01) (0.28, 7.422e-01) (0.285, 7.324e-01) (0.29, 7.224e-01) (0.295, 7.122e-01) (0.3, 7.018e-01) (0.305, 6.911e-01) (0.31, 6.802e-01) (0.315, 6.691e-01) (0.32, 6.577e-01) (0.325, 6.461e-01) (0.33, 6.343e-01) (0.335, 6.222e-01) (0.34, 6.099e-01) (0.345, 5.973e-01) (0.35, 5.844e-01) (0.355, 5.713e-01) (0.36, 5.579e-01) (0.365, 5.442e-01) (0.37, 5.303e-01) (0.375, 5.160e-01) (0.38, 5.014e-01) (0.385, 4.865e-01) (0.39, 4.713e-01) (0.395, 4.558e-01) (0.4, 4.399e-01) (0.405, 4.237e-01) (0.41, 4.071e-01) (0.415, 3.901e-01) (0.42, 3.727e-01) (0.425, 3.548e-01) (0.43, 3.366e-01) (0.435, 3.178e-01) (0.44, 2.986e-01) (0.445, 2.788e-01) (0.45, 2.585e-01) (0.455, 2.376e-01) (0.46, 2.160e-01) (0.465, 1.937e-01) (0.47, 1.706e-01) (0.475, 1.465e-01) (0.48, 1.214e-01) (0.485, 9.493e-02) (0.49, 6.681e-02) (0.495, 3.625e-02) (0.5, 0.000e+00)  
	};
	\addlegendentry{$(\eta^2,\mu^2)=(\ETAa,\MUa)$};
	\addplot[dashed, mark size = 2] coordinates { 
		(-0.5, 0.000e+00) (-0.495, 2.974e-04) (-0.49, 1.372e-03) (-0.485, 3.302e-03) (-0.48, 6.102e-03) (-0.475, 9.769e-03) (-0.47, 1.429e-02) (-0.465, 1.963e-02) (-0.46, 2.579e-02) (-0.455, 3.272e-02) (-0.45, 4.041e-02) (-0.445, 4.883e-02) (-0.44, 5.794e-02) (-0.435, 6.772e-02) (-0.43, 7.815e-02) (-0.425, 8.919e-02) (-0.42, 1.008e-01) (-0.415, 1.130e-01) (-0.41, 1.257e-01) (-0.405, 1.389e-01) (-0.4, 1.526e-01) (-0.395, 1.668e-01) (-0.39, 1.814e-01) (-0.385, 1.964e-01) (-0.38, 2.117e-01) (-0.375, 2.274e-01) (-0.37, 2.435e-01) (-0.365, 2.599e-01) (-0.36, 2.765e-01) (-0.355, 2.934e-01) (-0.35, 3.105e-01) (-0.345, 3.279e-01) (-0.34, 3.455e-01) (-0.335, 3.632e-01) (-0.33, 3.811e-01) (-0.325, 3.991e-01) (-0.32, 4.172e-01) (-0.315, 4.355e-01) (-0.31, 4.538e-01) (-0.305, 4.722e-01) (-0.3, 4.906e-01) (-0.295, 5.090e-01) (-0.29, 5.275e-01) (-0.285, 5.459e-01) (-0.28, 5.644e-01) (-0.275, 5.827e-01) (-0.27, 6.010e-01) (-0.265, 6.193e-01) (-0.26, 6.374e-01) (-0.255, 6.555e-01) (-0.25, 6.734e-01) (-0.245, 6.912e-01) (-0.24, 7.089e-01) (-0.235, 7.264e-01) (-0.23, 7.437e-01) (-0.225, 7.608e-01) (-0.22, 7.778e-01) (-0.215, 7.945e-01) (-0.21, 8.110e-01) (-0.205, 8.273e-01) (-0.2, 8.433e-01) (-0.195, 8.590e-01) (-0.19, 8.745e-01) (-0.185, 8.898e-01) (-0.18, 9.047e-01) (-0.175, 9.194e-01) (-0.17, 9.337e-01) (-0.165, 9.477e-01) (-0.16, 9.614e-01) (-0.155, 9.748e-01) (-0.15, 9.878e-01) (-0.145, 1.000e+00) (-0.14, 1.013e+00) (-0.135, 1.025e+00) (-0.13, 1.036e+00) (-0.125, 1.047e+00) (-0.12, 1.058e+00) (-0.115, 1.069e+00) (-0.11, 1.079e+00) (-0.105, 1.088e+00) (-0.1, 1.098e+00) (-0.095, 1.106e+00) (-0.09, 1.115e+00) (-0.085, 1.123e+00) (-0.08, 1.130e+00) (-0.075, 1.137e+00) (-0.07, 1.144e+00) (-0.065, 1.150e+00) (-0.06, 1.156e+00) (-0.055, 1.161e+00) (-0.05, 1.166e+00) (-0.045, 1.170e+00) (-0.04, 1.174e+00) (-0.035, 1.178e+00) (-0.03, 1.181e+00) (-0.025, 1.183e+00) (-0.02, 1.185e+00) (-0.015, 1.187e+00) (-0.01, 1.188e+00) (-0.005, 1.189e+00) (0, 1.189e+00) (0.005, 1.189e+00) (0.01, 1.188e+00) (0.015, 1.187e+00) (0.02, 1.185e+00) (0.025, 1.183e+00) (0.03, 1.181e+00) (0.035, 1.178e+00) (0.04, 1.174e+00) (0.045, 1.170e+00) (0.05, 1.166e+00) (0.055, 1.161e+00) (0.06, 1.156e+00) (0.065, 1.150e+00) (0.07, 1.144e+00) (0.075, 1.137e+00) (0.08, 1.130e+00) (0.085, 1.123e+00) (0.09, 1.115e+00) (0.095, 1.106e+00) (0.1, 1.098e+00) (0.105, 1.088e+00) (0.11, 1.079e+00) (0.115, 1.069e+00) (0.12, 1.058e+00) (0.125, 1.047e+00) (0.13, 1.036e+00) (0.135, 1.025e+00) (0.14, 1.013e+00) (0.145, 1.000e+00) (0.15, 9.878e-01) (0.155, 9.748e-01) (0.16, 9.614e-01) (0.165, 9.477e-01) (0.17, 9.337e-01) (0.175, 9.194e-01) (0.18, 9.047e-01) (0.185, 8.898e-01) (0.19, 8.745e-01) (0.195, 8.590e-01) (0.2, 8.433e-01) (0.205, 8.273e-01) (0.21, 8.110e-01) (0.215, 7.945e-01) (0.22, 7.778e-01) (0.225, 7.608e-01) (0.23, 7.437e-01) (0.235, 7.264e-01) (0.24, 7.089e-01) (0.245, 6.912e-01) (0.25, 6.734e-01) (0.255, 6.555e-01) (0.26, 6.374e-01) (0.265, 6.193e-01) (0.27, 6.010e-01) (0.275, 5.827e-01) (0.28, 5.644e-01) (0.285, 5.459e-01) (0.29, 5.275e-01) (0.295, 5.090e-01) (0.3, 4.906e-01) (0.305, 4.722e-01) (0.31, 4.538e-01) (0.315, 4.355e-01) (0.32, 4.172e-01) (0.325, 3.991e-01) (0.33, 3.811e-01) (0.335, 3.632e-01) (0.34, 3.455e-01) (0.345, 3.279e-01) (0.35, 3.105e-01) (0.355, 2.934e-01) (0.36, 2.765e-01) (0.365, 2.599e-01) (0.37, 2.435e-01) (0.375, 2.274e-01) (0.38, 2.117e-01) (0.385, 1.964e-01) (0.39, 1.814e-01) (0.395, 1.668e-01) (0.4, 1.526e-01) (0.405, 1.389e-01) (0.41, 1.257e-01) (0.415, 1.130e-01) (0.42, 1.008e-01) (0.425, 8.919e-02) (0.43, 7.815e-02) (0.435, 6.772e-02) (0.44, 5.794e-02) (0.445, 4.883e-02) (0.45, 4.041e-02) (0.455, 3.272e-02) (0.46, 2.579e-02) (0.465, 1.963e-02) (0.47, 1.429e-02) (0.475, 9.769e-03) (0.48, 6.102e-03) (0.485, 3.302e-03) (0.49, 1.372e-03) (0.495, 2.974e-04) (0.5, 0.000e+00)   
	};
	\addlegendentry{$(\eta^2,\mu^2)=(\ETAb,\MUb)$};
	\addplot[dotted, mark size = 2] coordinates { 
		(-0.5, 0.000e+00) (-0.495, 2.271e-06) (-0.49, 2.622e-05) (-0.485, 1.069e-04) (-0.48, 2.855e-04) (-0.475, 6.062e-04) (-0.47, 1.114e-03) (-0.465, 1.852e-03) (-0.46, 2.865e-03) (-0.455, 4.194e-03) (-0.45, 5.879e-03) (-0.445, 7.956e-03) (-0.44, 1.046e-02) (-0.435, 1.343e-02) (-0.43, 1.689e-02) (-0.425, 2.086e-02) (-0.42, 2.538e-02) (-0.415, 3.046e-02) (-0.41, 3.613e-02) (-0.405, 4.240e-02) (-0.4, 4.929e-02) (-0.395, 5.680e-02) (-0.39, 6.495e-02) (-0.385, 7.375e-02) (-0.38, 8.320e-02) (-0.375, 9.330e-02) (-0.37, 1.041e-01) (-0.365, 1.155e-01) (-0.36, 1.275e-01) (-0.355, 1.402e-01) (-0.35, 1.535e-01) (-0.345, 1.675e-01) (-0.34, 1.821e-01) (-0.335, 1.973e-01) (-0.33, 2.131e-01) (-0.325, 2.294e-01) (-0.32, 2.463e-01) (-0.315, 2.638e-01) (-0.31, 2.817e-01) (-0.305, 3.002e-01) (-0.3, 3.192e-01) (-0.295, 3.386e-01) (-0.29, 3.584e-01) (-0.285, 3.787e-01) (-0.28, 3.994e-01) (-0.275, 4.204e-01) (-0.27, 4.417e-01) (-0.265, 4.633e-01) (-0.26, 4.853e-01) (-0.255, 5.074e-01) (-0.25, 5.298e-01) (-0.245, 5.524e-01) (-0.24, 5.751e-01) (-0.235, 5.980e-01) (-0.23, 6.210e-01) (-0.225, 6.440e-01) (-0.22, 6.671e-01) (-0.215, 6.902e-01) (-0.21, 7.133e-01) (-0.205, 7.364e-01) (-0.2, 7.593e-01) (-0.195, 7.822e-01) (-0.19, 8.049e-01) (-0.185, 8.274e-01) (-0.18, 8.497e-01) (-0.175, 8.718e-01) (-0.17, 8.937e-01) (-0.165, 9.153e-01) (-0.16, 9.365e-01) (-0.155, 9.574e-01) (-0.15, 9.780e-01) (-0.145, 9.981e-01) (-0.14, 1.018e+00) (-0.135, 1.037e+00) (-0.13, 1.056e+00) (-0.125, 1.074e+00) (-0.12, 1.092e+00) (-0.115, 1.109e+00) (-0.11, 1.126e+00) (-0.105, 1.142e+00) (-0.1, 1.158e+00) (-0.095, 1.172e+00) (-0.09, 1.187e+00) (-0.085, 1.200e+00) (-0.08, 1.213e+00) (-0.075, 1.225e+00) (-0.07, 1.237e+00) (-0.065, 1.247e+00) (-0.06, 1.257e+00) (-0.055, 1.267e+00) (-0.05, 1.275e+00) (-0.045, 1.283e+00) (-0.04, 1.290e+00) (-0.035, 1.296e+00) (-0.03, 1.301e+00) (-0.025, 1.306e+00) (-0.02, 1.309e+00) (-0.015, 1.312e+00) (-0.01, 1.314e+00) (-0.005, 1.316e+00) (0, 1.316e+00) (0.005, 1.316e+00) (0.01, 1.314e+00) (0.015, 1.312e+00) (0.02, 1.309e+00) (0.025, 1.306e+00) (0.03, 1.301e+00) (0.035, 1.296e+00) (0.04, 1.290e+00) (0.045, 1.283e+00) (0.05, 1.275e+00) (0.055, 1.267e+00) (0.06, 1.257e+00) (0.065, 1.247e+00) (0.07, 1.237e+00) (0.075, 1.225e+00) (0.08, 1.213e+00) (0.085, 1.200e+00) (0.09, 1.187e+00) (0.095, 1.172e+00) (0.1, 1.158e+00) (0.105, 1.142e+00) (0.11, 1.126e+00) (0.115, 1.109e+00) (0.12, 1.092e+00) (0.125, 1.074e+00) (0.13, 1.056e+00) (0.135, 1.037e+00) (0.14, 1.018e+00) (0.145, 9.981e-01) (0.15, 9.780e-01) (0.155, 9.574e-01) (0.16, 9.365e-01) (0.165, 9.153e-01) (0.17, 8.937e-01) (0.175, 8.718e-01) (0.18, 8.497e-01) (0.185, 8.274e-01) (0.19, 8.049e-01) (0.195, 7.822e-01) (0.2, 7.593e-01) (0.205, 7.364e-01) (0.21, 7.133e-01) (0.215, 6.902e-01) (0.22, 6.671e-01) (0.225, 6.440e-01) (0.23, 6.210e-01) (0.235, 5.980e-01) (0.24, 5.751e-01) (0.245, 5.524e-01) (0.25, 5.298e-01) (0.255, 5.074e-01) (0.26, 4.853e-01) (0.265, 4.633e-01) (0.27, 4.417e-01) (0.275, 4.204e-01) (0.28, 3.994e-01) (0.285, 3.787e-01) (0.29, 3.584e-01) (0.295, 3.386e-01) (0.3, 3.192e-01) (0.305, 3.002e-01) (0.31, 2.817e-01) (0.315, 2.638e-01) (0.32, 2.463e-01) (0.325, 2.294e-01) (0.33, 2.131e-01) (0.335, 1.973e-01) (0.34, 1.821e-01) (0.345, 1.675e-01) (0.35, 1.535e-01) (0.355, 1.402e-01) (0.36, 1.275e-01) (0.365, 1.155e-01) (0.37, 1.041e-01) (0.375, 9.330e-02) (0.38, 8.320e-02) (0.385, 7.375e-02) (0.39, 6.495e-02) (0.395, 5.680e-02) (0.4, 4.929e-02) (0.405, 4.240e-02) (0.41, 3.613e-02) (0.415, 3.046e-02) (0.42, 2.538e-02) (0.425, 2.086e-02) (0.43, 1.689e-02) (0.435, 1.343e-02) (0.44, 1.046e-02) (0.445, 7.956e-03) (0.45, 5.879e-03) (0.455, 4.194e-03) (0.46, 2.865e-03) (0.465, 1.852e-03) (0.47, 1.114e-03) (0.475, 6.062e-04) (0.48, 2.855e-04) (0.485, 1.069e-04) (0.49, 2.622e-05) (0.495, 2.271e-06) (0.5, 0.000e+00) 
	};
	\addlegendentry{$(\eta^2,\mu^2)=(\ETAc,\MUc)$};
	\end{axis}
	\end{tikzpicture}
\end{minipage}
\begin{minipage}{.5\linewidth}
	\centering
\begin{tikzpicture}[scale=0.75]
	\pgfmathsetmacro\MU{1} \pgfmathsetmacro\ETA{1}
	\pgfmathsetmacro\MUa{3} \pgfmathsetmacro\ETAa{1}
	\pgfmathsetmacro\MUb{6} \pgfmathsetmacro\ETAb{1}
	\pgfmathsetmacro\MUc{10} \pgfmathsetmacro\ETAc{1}
	\begin{axis}[samples=100,  
	xtick={-0.5, -0.3, -0.1, 0, 0.1, 0.3, 0.5},
xmin=-0.55, xmax=0.55, ymin=0, ymax=2.75,
	axis x line=center, axis y line=center,
	every axis plot/.append style={thick},
	legend style={at={(0.5,1.05)}, anchor=south,legend columns=2,legend cell align=left, font=\small}
	]
\addplot[dashdotted, mark size = 2] coordinates{ (-0.5, 1) (0.5, 1) };
	\addlegendentry{$(\eta^2,\mu^2)=(\ETA,\MU)$};
	\addplot[mark size = 2] coordinates { 
		(-0.5, 0.000e+00) (-0.495, 3.625e-02) (-0.49, 6.681e-02) (-0.485, 9.493e-02) (-0.48, 1.214e-01) (-0.475, 1.465e-01) (-0.47, 1.706e-01) (-0.465, 1.937e-01) (-0.46, 2.160e-01) (-0.455, 2.376e-01) (-0.45, 2.585e-01) (-0.445, 2.788e-01) (-0.44, 2.986e-01) (-0.435, 3.178e-01) (-0.43, 3.366e-01) (-0.425, 3.548e-01) (-0.42, 3.727e-01) (-0.415, 3.901e-01) (-0.41, 4.071e-01) (-0.405, 4.237e-01) (-0.4, 4.399e-01) (-0.395, 4.558e-01) (-0.39, 4.713e-01) (-0.385, 4.865e-01) (-0.38, 5.014e-01) (-0.375, 5.160e-01) (-0.37, 5.303e-01) (-0.365, 5.442e-01) (-0.36, 5.579e-01) (-0.355, 5.713e-01) (-0.35, 5.844e-01) (-0.345, 5.973e-01) (-0.34, 6.099e-01) (-0.335, 6.222e-01) (-0.33, 6.343e-01) (-0.325, 6.461e-01) (-0.32, 6.577e-01) (-0.315, 6.691e-01) (-0.31, 6.802e-01) (-0.305, 6.911e-01) (-0.3, 7.018e-01) (-0.295, 7.122e-01) (-0.29, 7.224e-01) (-0.285, 7.324e-01) (-0.28, 7.422e-01) (-0.275, 7.518e-01) (-0.27, 7.611e-01) (-0.265, 7.703e-01) (-0.26, 7.792e-01) (-0.255, 7.880e-01) (-0.25, 7.965e-01) (-0.245, 8.049e-01) (-0.24, 8.131e-01) (-0.235, 8.210e-01) (-0.23, 8.288e-01) (-0.225, 8.364e-01) (-0.22, 8.438e-01) (-0.215, 8.510e-01) (-0.21, 8.580e-01) (-0.205, 8.649e-01) (-0.2, 8.715e-01) (-0.195, 8.780e-01) (-0.19, 8.843e-01) (-0.185, 8.904e-01) (-0.18, 8.964e-01) (-0.175, 9.022e-01) (-0.17, 9.078e-01) (-0.165, 9.132e-01) (-0.16, 9.185e-01) (-0.155, 9.235e-01) (-0.15, 9.285e-01) (-0.145, 9.332e-01) (-0.14, 9.378e-01) (-0.135, 9.422e-01) (-0.13, 9.464e-01) (-0.125, 9.505e-01) (-0.12, 9.544e-01) (-0.115, 9.582e-01) (-0.11, 9.617e-01) (-0.105, 9.652e-01) (-0.1, 9.684e-01) (-0.095, 9.715e-01) (-0.09, 9.744e-01) (-0.085, 9.772e-01) (-0.08, 9.798e-01) (-0.075, 9.823e-01) (-0.07, 9.846e-01) (-0.065, 9.867e-01) (-0.06, 9.887e-01) (-0.055, 9.905e-01) (-0.05, 9.921e-01) (-0.045, 9.936e-01) (-0.04, 9.950e-01) (-0.035, 9.961e-01) (-0.03, 9.972e-01) (-0.025, 9.980e-01) (-0.02, 9.987e-01) (-0.015, 9.993e-01) (-0.01, 9.997e-01) (-0.005, 9.999e-01) (0, 1.000e+00) (0.005, 9.999e-01) (0.01, 9.997e-01) (0.015, 9.993e-01) (0.02, 9.987e-01) (0.025, 9.980e-01) (0.03, 9.972e-01) (0.035, 9.961e-01) (0.04, 9.950e-01) (0.045, 9.936e-01) (0.05, 9.921e-01) (0.055, 9.905e-01) (0.06, 9.887e-01) (0.065, 9.867e-01) (0.07, 9.846e-01) (0.075, 9.823e-01) (0.08, 9.798e-01) (0.085, 9.772e-01) (0.09, 9.744e-01) (0.095, 9.715e-01) (0.1, 9.684e-01) (0.105, 9.652e-01) (0.11, 9.617e-01) (0.115, 9.582e-01) (0.12, 9.544e-01) (0.125, 9.505e-01) (0.13, 9.464e-01) (0.135, 9.422e-01) (0.14, 9.378e-01) (0.145, 9.332e-01) (0.15, 9.285e-01) (0.155, 9.235e-01) (0.16, 9.185e-01) (0.165, 9.132e-01) (0.17, 9.078e-01) (0.175, 9.022e-01) (0.18, 8.964e-01) (0.185, 8.904e-01) (0.19, 8.843e-01) (0.195, 8.780e-01) (0.2, 8.715e-01) (0.205, 8.649e-01) (0.21, 8.580e-01) (0.215, 8.510e-01) (0.22, 8.438e-01) (0.225, 8.364e-01) (0.23, 8.288e-01) (0.235, 8.210e-01) (0.24, 8.131e-01) (0.245, 8.049e-01) (0.25, 7.965e-01) (0.255, 7.880e-01) (0.26, 7.792e-01) (0.265, 7.703e-01) (0.27, 7.611e-01) (0.275, 7.518e-01) (0.28, 7.422e-01) (0.285, 7.324e-01) (0.29, 7.224e-01) (0.295, 7.122e-01) (0.3, 7.018e-01) (0.305, 6.911e-01) (0.31, 6.802e-01) (0.315, 6.691e-01) (0.32, 6.577e-01) (0.325, 6.461e-01) (0.33, 6.343e-01) (0.335, 6.222e-01) (0.34, 6.099e-01) (0.345, 5.973e-01) (0.35, 5.844e-01) (0.355, 5.713e-01) (0.36, 5.579e-01) (0.365, 5.442e-01) (0.37, 5.303e-01) (0.375, 5.160e-01) (0.38, 5.014e-01) (0.385, 4.865e-01) (0.39, 4.713e-01) (0.395, 4.558e-01) (0.4, 4.399e-01) (0.405, 4.237e-01) (0.41, 4.071e-01) (0.415, 3.901e-01) (0.42, 3.727e-01) (0.425, 3.548e-01) (0.43, 3.366e-01) (0.435, 3.178e-01) (0.44, 2.986e-01) (0.445, 2.788e-01) (0.45, 2.585e-01) (0.455, 2.376e-01) (0.46, 2.160e-01) (0.465, 1.937e-01) (0.47, 1.706e-01) (0.475, 1.465e-01) (0.48, 1.214e-01) (0.485, 9.493e-02) (0.49, 6.681e-02) (0.495, 3.625e-02) (0.5, 0.000e+00) 
	};
	\addlegendentry{$(\eta^2,\mu^2)=(\ETAa,\MUa)$};
	\addplot[dashed, mark size = 2] coordinates { 
		(-0.5, 0.000e+00) (-0.495, 2.501e-04) (-0.49, 1.154e-03) (-0.485, 2.776e-03) (-0.48, 5.132e-03) (-0.475, 8.215e-03) (-0.47, 1.201e-02) (-0.465, 1.651e-02) (-0.46, 2.169e-02) (-0.455, 2.752e-02) (-0.45, 3.398e-02) (-0.445, 4.106e-02) (-0.44, 4.872e-02) (-0.435, 5.695e-02) (-0.43, 6.571e-02) (-0.425, 7.500e-02) (-0.42, 8.477e-02) (-0.415, 9.502e-02) (-0.41, 1.057e-01) (-0.405, 1.168e-01) (-0.4, 1.284e-01) (-0.395, 1.403e-01) (-0.39, 1.525e-01) (-0.385, 1.651e-01) (-0.38, 1.780e-01) (-0.375, 1.913e-01) (-0.37, 2.048e-01) (-0.365, 2.185e-01) (-0.36, 2.325e-01) (-0.355, 2.467e-01) (-0.35, 2.611e-01) (-0.345, 2.757e-01) (-0.34, 2.905e-01) (-0.335, 3.054e-01) (-0.33, 3.204e-01) (-0.325, 3.356e-01) (-0.32, 3.509e-01) (-0.315, 3.662e-01) (-0.31, 3.816e-01) (-0.305, 3.971e-01) (-0.3, 4.125e-01) (-0.295, 4.281e-01) (-0.29, 4.436e-01) (-0.285, 4.591e-01) (-0.28, 4.746e-01) (-0.275, 4.900e-01) (-0.27, 5.054e-01) (-0.265, 5.208e-01) (-0.26, 5.360e-01) (-0.255, 5.512e-01) (-0.25, 5.663e-01) (-0.245, 5.812e-01) (-0.24, 5.961e-01) (-0.235, 6.108e-01) (-0.23, 6.254e-01) (-0.225, 6.398e-01) (-0.22, 6.540e-01) (-0.215, 6.681e-01) (-0.21, 6.820e-01) (-0.205, 6.956e-01) (-0.2, 7.091e-01) (-0.195, 7.224e-01) (-0.19, 7.354e-01) (-0.185, 7.482e-01) (-0.18, 7.608e-01) (-0.175, 7.731e-01) (-0.17, 7.851e-01) (-0.165, 7.969e-01) (-0.16, 8.084e-01) (-0.155, 8.197e-01) (-0.15, 8.306e-01) (-0.145, 8.413e-01) (-0.14, 8.516e-01) (-0.135, 8.617e-01) (-0.13, 8.714e-01) (-0.125, 8.808e-01) (-0.12, 8.899e-01) (-0.115, 8.987e-01) (-0.11, 9.071e-01) (-0.105, 9.152e-01) (-0.1, 9.229e-01) (-0.095, 9.303e-01) (-0.09, 9.373e-01) (-0.085, 9.440e-01) (-0.08, 9.503e-01) (-0.075, 9.563e-01) (-0.07, 9.619e-01) (-0.065, 9.671e-01) (-0.06, 9.719e-01) (-0.055, 9.764e-01) (-0.05, 9.805e-01) (-0.045, 9.842e-01) (-0.04, 9.875e-01) (-0.035, 9.904e-01) (-0.03, 9.929e-01) (-0.025, 9.951e-01) (-0.02, 9.969e-01) (-0.015, 9.982e-01) (-0.01, 9.992e-01) (-0.005, 9.998e-01) (0, 1.000e+00) (0.005, 9.998e-01) (0.01, 9.992e-01) (0.015, 9.982e-01) (0.02, 9.969e-01) (0.025, 9.951e-01) (0.03, 9.929e-01) (0.035, 9.904e-01) (0.04, 9.875e-01) (0.045, 9.842e-01) (0.05, 9.805e-01) (0.055, 9.764e-01) (0.06, 9.719e-01) (0.065, 9.671e-01) (0.07, 9.619e-01) (0.075, 9.563e-01) (0.08, 9.503e-01) (0.085, 9.440e-01) (0.09, 9.373e-01) (0.095, 9.303e-01) (0.1, 9.229e-01) (0.105, 9.152e-01) (0.11, 9.071e-01) (0.115, 8.987e-01) (0.12, 8.899e-01) (0.125, 8.808e-01) (0.13, 8.714e-01) (0.135, 8.617e-01) (0.14, 8.516e-01) (0.145, 8.413e-01) (0.15, 8.306e-01) (0.155, 8.197e-01) (0.16, 8.084e-01) (0.165, 7.969e-01) (0.17, 7.851e-01) (0.175, 7.731e-01) (0.18, 7.608e-01) (0.185, 7.482e-01) (0.19, 7.354e-01) (0.195, 7.224e-01) (0.2, 7.091e-01) (0.205, 6.956e-01) (0.21, 6.820e-01) (0.215, 6.681e-01) (0.22, 6.540e-01) (0.225, 6.398e-01) (0.23, 6.254e-01) (0.235, 6.108e-01) (0.24, 5.961e-01) (0.245, 5.812e-01) (0.25, 5.663e-01) (0.255, 5.512e-01) (0.26, 5.360e-01) (0.265, 5.208e-01) (0.27, 5.054e-01) (0.275, 4.900e-01) (0.28, 4.746e-01) (0.285, 4.591e-01) (0.29, 4.436e-01) (0.295, 4.281e-01) (0.3, 4.125e-01) (0.305, 3.971e-01) (0.31, 3.816e-01) (0.315, 3.662e-01) (0.32, 3.509e-01) (0.325, 3.356e-01) (0.33, 3.204e-01) (0.335, 3.054e-01) (0.34, 2.905e-01) (0.345, 2.757e-01) (0.35, 2.611e-01) (0.355, 2.467e-01) (0.36, 2.325e-01) (0.365, 2.185e-01) (0.37, 2.048e-01) (0.375, 1.913e-01) (0.38, 1.780e-01) (0.385, 1.651e-01) (0.39, 1.525e-01) (0.395, 1.403e-01) (0.4, 1.284e-01) (0.405, 1.168e-01) (0.41, 1.057e-01) (0.415, 9.502e-02) (0.42, 8.477e-02) (0.425, 7.500e-02) (0.43, 6.571e-02) (0.435, 5.695e-02) (0.44, 4.872e-02) (0.445, 4.106e-02) (0.45, 3.398e-02) (0.455, 2.752e-02) (0.46, 2.169e-02) (0.465, 1.651e-02) (0.47, 1.201e-02) (0.475, 8.215e-03) (0.48, 5.132e-03) (0.485, 2.776e-03) (0.49, 1.154e-03) (0.495, 2.501e-04) (0.5, 0.000e+00)    
	};
	\addlegendentry{$(\eta^2,\mu^2)=(\ETAb,\MUb)$};
	\addplot[dotted, mark size = 2] coordinates { 
		(-0.5, 0.000e+00) (-0.495, 3.286e-07) (-0.49, 5.149e-06) (-0.485, 2.502e-05) (-0.48, 7.559e-05) (-0.475, 1.763e-04) (-0.47, 3.495e-04) (-0.465, 6.194e-04) (-0.46, 1.012e-03) (-0.455, 1.553e-03) (-0.45, 2.271e-03) (-0.445, 3.192e-03) (-0.44, 4.344e-03) (-0.435, 5.752e-03) (-0.43, 7.444e-03) (-0.425, 9.442e-03) (-0.42, 1.177e-02) (-0.415, 1.446e-02) (-0.41, 1.752e-02) (-0.405, 2.097e-02) (-0.4, 2.484e-02) (-0.395, 2.914e-02) (-0.39, 3.388e-02) (-0.385, 3.909e-02) (-0.38, 4.476e-02) (-0.375, 5.092e-02) (-0.37, 5.757e-02) (-0.365, 6.472e-02) (-0.36, 7.237e-02) (-0.355, 8.053e-02) (-0.35, 8.919e-02) (-0.345, 9.837e-02) (-0.34, 1.081e-01) (-0.335, 1.182e-01) (-0.33, 1.289e-01) (-0.325, 1.401e-01) (-0.32, 1.518e-01) (-0.315, 1.639e-01) (-0.31, 1.766e-01) (-0.305, 1.896e-01) (-0.3, 2.032e-01) (-0.295, 2.171e-01) (-0.29, 2.315e-01) (-0.285, 2.463e-01) (-0.28, 2.614e-01) (-0.275, 2.769e-01) (-0.27, 2.928e-01) (-0.265, 3.090e-01) (-0.26, 3.255e-01) (-0.255, 3.423e-01) (-0.25, 3.593e-01) (-0.245, 3.766e-01) (-0.24, 3.941e-01) (-0.235, 4.117e-01) (-0.23, 4.296e-01) (-0.225, 4.476e-01) (-0.22, 4.656e-01) (-0.215, 4.838e-01) (-0.21, 5.021e-01) (-0.205, 5.203e-01) (-0.2, 5.386e-01) (-0.195, 5.569e-01) (-0.19, 5.751e-01) (-0.185, 5.933e-01) (-0.18, 6.113e-01) (-0.175, 6.292e-01) (-0.17, 6.470e-01) (-0.165, 6.646e-01) (-0.16, 6.820e-01) (-0.155, 6.991e-01) (-0.15, 7.160e-01) (-0.145, 7.326e-01) (-0.14, 7.489e-01) (-0.135, 7.649e-01) (-0.13, 7.805e-01) (-0.125, 7.958e-01) (-0.12, 8.106e-01) (-0.115, 8.250e-01) (-0.11, 8.390e-01) (-0.105, 8.525e-01) (-0.1, 8.655e-01) (-0.095, 8.780e-01) (-0.09, 8.900e-01) (-0.085, 9.015e-01) (-0.08, 9.124e-01) (-0.075, 9.227e-01) (-0.07, 9.324e-01) (-0.065, 9.415e-01) (-0.06, 9.500e-01) (-0.055, 9.579e-01) (-0.05, 9.651e-01) (-0.045, 9.717e-01) (-0.04, 9.776e-01) (-0.035, 9.828e-01) (-0.03, 9.873e-01) (-0.025, 9.912e-01) (-0.02, 9.944e-01) (-0.015, 9.968e-01) (-0.01, 9.986e-01) (-0.005, 9.996e-01) (0, 1.000e+00) (0.005, 9.996e-01) (0.01, 9.986e-01) (0.015, 9.968e-01) (0.02, 9.944e-01) (0.025, 9.912e-01) (0.03, 9.873e-01) (0.035, 9.828e-01) (0.04, 9.776e-01) (0.045, 9.717e-01) (0.05, 9.651e-01) (0.055, 9.579e-01) (0.06, 9.500e-01) (0.065, 9.415e-01) (0.07, 9.324e-01) (0.075, 9.227e-01) (0.08, 9.124e-01) (0.085, 9.015e-01) (0.09, 8.900e-01) (0.095, 8.780e-01) (0.1, 8.655e-01) (0.105, 8.525e-01) (0.11, 8.390e-01) (0.115, 8.250e-01) (0.12, 8.106e-01) (0.125, 7.958e-01) (0.13, 7.805e-01) (0.135, 7.649e-01) (0.14, 7.489e-01) (0.145, 7.326e-01) (0.15, 7.160e-01) (0.155, 6.991e-01) (0.16, 6.820e-01) (0.165, 6.646e-01) (0.17, 6.470e-01) (0.175, 6.292e-01) (0.18, 6.113e-01) (0.185, 5.933e-01) (0.19, 5.751e-01) (0.195, 5.569e-01) (0.2, 5.386e-01) (0.205, 5.203e-01) (0.21, 5.021e-01) (0.215, 4.838e-01) (0.22, 4.656e-01) (0.225, 4.476e-01) (0.23, 4.296e-01) (0.235, 4.117e-01) (0.24, 3.941e-01) (0.245, 3.766e-01) (0.25, 3.593e-01) (0.255, 3.423e-01) (0.26, 3.255e-01) (0.265, 3.090e-01) (0.27, 2.928e-01) (0.275, 2.769e-01) (0.28, 2.614e-01) (0.285, 2.463e-01) (0.29, 2.315e-01) (0.295, 2.171e-01) (0.3, 2.032e-01) (0.305, 1.896e-01) (0.31, 1.766e-01) (0.315, 1.639e-01) (0.32, 1.518e-01) (0.325, 1.401e-01) (0.33, 1.289e-01) (0.335, 1.182e-01) (0.34, 1.081e-01) (0.345, 9.837e-02) (0.35, 8.919e-02) (0.355, 8.053e-02) (0.36, 7.237e-02) (0.365, 6.472e-02) (0.37, 5.757e-02) (0.375, 5.092e-02) (0.38, 4.476e-02) (0.385, 3.909e-02) (0.39, 3.388e-02) (0.395, 2.914e-02) (0.4, 2.484e-02) (0.405, 2.097e-02) (0.41, 1.752e-02) (0.415, 1.446e-02) (0.42, 1.177e-02) (0.425, 9.442e-03) (0.43, 7.444e-03) (0.435, 5.752e-03) (0.44, 4.344e-03) (0.445, 3.192e-03) (0.45, 2.271e-03) (0.455, 1.553e-03) (0.46, 1.012e-03) (0.465, 6.194e-04) (0.47, 3.495e-04) (0.475, 1.763e-04) (0.48, 7.559e-05) (0.485, 2.502e-05) (0.49, 5.149e-06) (0.495, 3.286e-07) (0.5, 0.000e+00)  
	};
	\addlegendentry{$(\eta^2,\mu^2)=(\ETAc,\MUc)$};
	\end{axis}
	\end{tikzpicture}
\end{minipage}
	\caption{Plots of the univariate transformed function $f$ for various combinations of the parameters $\mu$ and $\eta$ with a Gaussian weight function $\omega$ \eqref{eq:gaussian_weight_param} and the error function transformation \eqref{eq:erf_trafo_param}. On the left hand side with fixed $\mu^2=3$, on the right hand side with fixed $\eta^2 = 1$.}
	\label{fig:Plots_f_with_erfTrafo}
\end{figure}

In this section we settle for the constant function $h(\mathbf y) \equiv 1$. 
We could choose $h(\mathbf y)=\mathrm{e}^{-\|\mathbf y\|_{\ell_2}^2}$ or even the algebraic function $h(\mathbf y)=\frac{1}{1+\|\mathbf y\|_{\ell_2}^2}$ in \eqref{eq:exemplary_h}, but they all impose the same problem that we will not able to compute their Fourier coefficients $\hat{h}_{\mathbf k}$.
We proceed in the same way as in the previous section with the algebraic transformation.
We fix the multivariate version of the Gaussian weight function \eqref{eq:gaussian_weight_param_univar}, reading as
\begin{align}\label{eq:gaussian_weight_param}
	\omega(\mathbf y,\bm\mu) = \frac{1}{\pi^{\frac{d}{2}}}\,\prod_{j=1}^{d}\mathrm{e}^{-\mu_j^2 y_j^2}
\end{align}
with $\bm\mu\in\mathbb{R}^d$, as well as the error function transformation $\psi(\mathbf x,\bm\eta) = ((\psi_j(x_j,\eta_j))_{j=1}^{d})^{\top}$
in the form \eqref{eq:param_trafo_explicit_example} with $\mathbf x\in(-\frac{1}{2},\frac{1}{2})^d$, the parameter ${\bm\eta = (\eta_1,\ldots,\eta_d)^{\top}\in\mathbb{R}^d}$ and its univariate components given by
\begin{align}\label{eq:erf_trafo_param}
	\psi_j(x_j,\eta_j) &= \eta_j \, \mathrm{erf}^{-1}(2x_j) ,
	\quad \psi_j'(x_j,\eta_j) = \eta_j \, \sqrt{\pi} \, \mathrm{e}^{(\mathrm{erf}^{-1}(2x_j))^2} \\
	\psi^{-1}(y_j,\eta_j) &= \frac{1}{2}\,\mathrm{erf}\left( \frac{y_j}{\eta_j} \right),
	\quad \varrho(y_j,\eta_j) = \frac{1}{\sqrt{\pi \eta_j^2}}\,\mathrm{e}^{-\left(\frac{y_j}{\eta_j}\right)^2}. \nonumber
\end{align}
For $\eta_j=1$ we stated the definition of $\psi_j(\circ,1)$ already in $\eqref{eq:error_function_trafo}$.
For the resulting weighted Hilbert space $L_{2}(\mathbb{R}^d,\omega(\circ,\bm\mu))$ we have a system $\left\{\varphi_{\mathbf k}\right\}_{\mathbf k\in\mathbb{Z}^d}$ of product functions given in \eqref{eq:transformed_basis_functions_mult} with univariate components $(\varphi_{k_j})_{j=1}^{d}$ as in \eqref{eq:transformed_basis_functions} of the form
\begin{align*}
	\varphi_{k_j}(y_j,\eta_j,\mu_j)
	= \frac{1}{\eta_j} \, \mathrm{e}^{\frac{1}{2}(\mu_j^2-\frac{1}{\eta_j^2})y_j^2 + \pi\mathrm i k_j\,\mathrm{erf}\left(\frac{y_j}{\eta_j}\right)},
\end{align*}
that are orthogonal with respect to the weighted scalar product
\begin{align*}
	(h_1, h_2)_{L_2\left(\mathbb{R}^d, \omega(\circ, \bm\mu) \right)}
	= \frac{1}{\pi^{\frac{d}{2}}}\,\int_{\mathbb{R}^d} \prod_{j=1}^{d}\mathrm{e}^{-\mu_j^2 y_j^2} \, h_1(\mathbf y) \, \overline{h_2(\mathbf y)} \, \mathrm d\mathbf y
\end{align*}
and the Fourier coefficients $\hat f_{\mathbf k}$ of an arbitrary function $h\in L_{2}(\mathbb{R}^d,\omega(\circ,\bm\mu))$ are of the form
\begin{align*}
	\hat h_{\mathbf k}
	:= \left(h, \varphi_{\mathbf k} \right)_{L_2\left(\mathbb{R}^d, \omega(\circ, \bm\mu) \right)}
	&= \int_{\mathbb{R}^d} h(\mathbf y) \, \overline{\varphi_{\mathbf k}(\mathbf y,\bm\eta,\bm\mu)} \, \omega(\mathbf y,\bm\mu) \, \mathrm d\mathbf y \\
	&= \int_{\mathbb{R}^d} h(\mathbf y) \, \prod_{j=1}^{d}\frac{1}{\eta_j} \, \mathrm{e}^{\frac{1}{2}(\mu_j^2-\frac{1}{\eta_j^2})y_j^2 - \pi\mathrm i k_j\,\mathrm{erf}\left(\frac{y_j}{\eta_j}\right)} \, \frac{1}{\sqrt{\pi}}\,\mathrm{e}^{-\mu_j^2 y_j^2} \, \mathrm d\mathbf y \\
	&= \pi^{-\frac{d}{2}}\prod_{j=1}^{d}\frac{1}{\eta_j} \int_{\mathbb{R}^d} h(\mathbf y) \,\prod_{j=1}^{d}\mathrm{e}^{- \pi\mathrm i k_j\,\mathrm{erf}\left(\frac{y_j}{\eta_j}\right)} \, \mathrm{e}^{-\frac{1}{2}(\mu_j^2+\frac{1}{\eta_j^2})y_j^2} \, \mathrm d\mathbf y.
\end{align*}
The constant test function $h(y)\equiv 1$ combined with the weight function \eqref{eq:gaussian_weight_param} and the transformations \eqref{eq:erf_trafo_param} lead to transformed functions $f$ in the sense of \eqref{eq:f_is_transformed_h_mult} of the form
\begin{align}
	\label{eq:Error_trafo_function}
	f(\mathbf x,\bm\eta,\bm\mu)
	&= h(\psi_1(x_1,\eta_1),\ldots,\psi_d(x_d,\eta_d)) \, \prod_{j=1}^{d}\sqrt{ \omega_j(\psi_j(x_j,\eta_j),\mu_j) \, \psi_j'(x_j,\eta_j) } \nonumber\\
	&= \prod_{j=1}^{d}\eta_j^{\frac{1}{2}} \, \mathrm{e}^{\frac{1}{2}\left(1- \mu_j^2 \, \eta_j^2 \right) \mathrm{erf}^{-1}(2x_j)^2}.
\end{align}
In Figure~\ref{fig:Plots_f_with_erfTrafo} we have a side-by-side comparison of the graphs of these transformed functions with $d=1$ for fixed ${\mu^2 = 3}$ with varied ${1/2\leq \eta^2\leq 3}$ and for fixed $\eta=1$ with varied ${1\leq \mu^2\leq 10}$.

We proceed to determine the values $\bm\eta,\bm\mu\in\mathbb{R}^d$ for which $f(\circ, \bm\eta, \bm\mu)$ as in \eqref{eq:Error_trafo_function} is element of $H_{\mathrm{mix}}^{m}(\mathbb{T}^d)$ by investigating conditions~\eqref{eq:Hm_composition_criteria_mult} in Theorem~\ref{thm:Hm_composition_criteria_mult}.
First of all, we observe that for $\eta_1,\ldots,\eta_d > 0$ the components $\psi_1,\ldots,\psi_d$ of the function $\psi(\circ,\bm\eta)$ in \eqref{eq:erf_trafo_param} are transformations in the sense of \eqref{def:Trafo_def} by being increasing, continuously differentiable and invertible functions.
Furthermore, for all $\ell=1,\ldots,d$ it's easy to check that its first three derivatives of all $\psi_j(\circ,\eta_j)$ are in fact continuous on $(-\frac{1}{2},\frac{1}{2})$ for $\eta_j > 0$ and that the first three derivatives of $\varrho_j(\circ,\eta_j)$ are in $\mathcal{C}_{0}(\mathbb{R})$ for all non-zero $\eta_j\in\mathbb{R}$.
Finally, we check the $L_{\infty}$-conditions~\eqref{eq:Hm_composition_criteria_mult} in Theorem~\ref{thm:Hm_composition_criteria_mult} for $m=0,1,2,3$. We suppose that for $\ell=1,\ldots,d$ we have $m=m_{\ell}$ and need to check that the appearing $L_{\infty}(\mathbb{T})$-norms are finite for all $j_{\ell}=0,\ldots,m$:
\begin{itemize}
\item 
	Let $m=0$. We have
	\begin{align*}
		\left\| \sqrt{\omega(\psi(x,\eta),\mu)} \right\|_{L_{\infty}(\mathbb{T})}
		&= \pi^{-\frac{1}{4}} \left\| \mathrm{e}^{-\frac{1}{2}\eta^2\mu^2\,\mathrm{erf}^{-1}(2x)^2} \right\|_{L_{\infty}(\mathbb{T})}
		< \infty
	\end{align*}
	for $\eta^2\mu^2 \geq 0$.
\item 
	Let $m=1$. We have to check two conditions.
	For $j_{\ell}=0$ we have
	\begin{align*}
		&\left\| \frac{\partial}{\partial x_{\ell}}\left[\sqrt{\omega_{\ell}(\psi_{\ell}(x_{\ell},\eta_{\ell}),\mu_{\ell}) \, \psi_{\ell}'(x_{\ell},\eta_{\ell})} \right] \psi_{\ell}'(x_{\ell},\eta_{\ell})^{-\frac{1}{2}} \right\|_{L_{\infty}(\mathbb{T})} \\
		&= \pi^{\frac{1}{4}} \left|\frac{\eta_{\ell}^2 - \mu_{\ell}^2}{\eta_{\ell}^2}\right| \left\| \mathrm{erf}^{-1}(2x_{\ell}) \, \mathrm{e}^{-\frac{1}{2}(\eta_{\ell}^2\mu_{\ell}^2-2)\,\mathrm{erf}^{-1}(2x_{\ell})^2} \right\|_{L_{\infty}(\mathbb{T})}
\end{align*}
	being finite for $\eta^2\mu^2 > 2$.
	For $j_{\ell}=1$ we have
	\begin{align*}
		&\left\| \sqrt{\omega_{\ell}(\psi_{\ell}(x_{\ell},\eta_{\ell}),\mu_{\ell}) \, \psi_{\ell}'(x_{\ell},\eta_{\ell})} \, (\psi_{\ell}'(x_{\ell},\eta_{\ell}))^{\frac{1}{2}} \right\|_{L_{\infty}(\mathbb{T})} \\
		&= \pi^{\frac{1}{4}}
		\left\| \mathrm{e}^{ -\frac{1}{2}(\mu_{\ell}^2 \eta_{\ell}^2-2)(\mathrm{erf}^{-1}(2x_{\ell})^2 } \right\|_{L_{\infty}(\mathbb{T})}
	\end{align*}
	and this is finite if the exponent is negative or zero, which is the case for $\eta^2\mu^2 \geq 2$.
\item 
	Let $m=2$. We check three conditions.
	For $j_{\ell}=0$
	\begin{align*}
		\left\| \frac{\partial^2}{\partial x_{\ell}^2}\left[\sqrt{\omega_{\ell}(\psi_{\ell}(x_{\ell},\eta_{\ell}),\mu_{\ell}) \, \psi_{\ell}'(x_{\ell},\eta_{\ell})} \right] \psi_{\ell}'(x_{\ell},\eta_{\ell})^{-\frac{1}{2}} \right\|_{L_{\infty}(\mathbb{T})} < \infty 
	\end{align*}
	for all $\eta_{\ell}^2\mu_{\ell}^2 > 4$.
	For $j_{\ell}=1$
	\begin{align*}
		\left\| \frac{\partial}{\partial x_{\ell}}\left[\sqrt{\omega_{\ell}(\psi_{\ell}(x_{\ell},\eta_{\ell}),\mu_{\ell}) \, \psi_{\ell}'(x_{\ell},\eta_{\ell})} \right] \psi_{\ell}'(x_{\ell},\eta_{\ell})^{\frac{1}{2}} \right\|_{L_{\infty}(\mathbb{T})} < \infty
	\end{align*}
	for all $\eta_{\ell}^2\mu_{\ell}^2 > 4$.
	For $j_{\ell}=2$
	\begin{align*}
		\left\| \sqrt{\omega_{\ell}(\psi_{\ell}(x_{\ell},\eta_{\ell}),\mu_{\ell}) \, \psi_{\ell}'(x_{\ell},\eta_{\ell})} \, (\psi_{\ell}'(x_{\ell},\eta_{\ell}))^{\frac{5}{2}} \right\|_{L_{\infty}(\mathbb{T})} < \infty
	\end{align*} 
	for all $\eta_{\ell}^2\mu_{\ell}^2 \geq 6$.
\item 
	For $m=3$ the individual conditions for $k=0,1,2,3$ are finite in case of ${\eta_{\ell}^2\mu_{\ell}^2 > 6}$, ${\eta_{\ell}^2\mu_{\ell}^2 > 6}$, ${\eta_{\ell}^2\mu_{\ell}^2 > 8}$ and ${\eta_{\ell}^2\mu_{\ell}^2 \geq 10}$ respectively.
	Hence, we need ${\eta_{\ell}^2\mu_{\ell}^2 \geq 10}$ in order to have $f\in H_{\mathrm{mix}}^{3}(\mathbb{T}^d)$.
\end{itemize}
In total we calculated that
\begin{align*}
	f \in
	\begin{cases}
H_{\mathrm{mix}}^{1}(\mathbb{T}^d) \quad &\text{ for } \quad \eta_{\ell}^2\mu_{\ell}^2 \geq 2,\\
		H_{\mathrm{mix}}^{2}(\mathbb{T}^d) \quad &\text{ for } \quad \eta_{\ell}^2\mu_{\ell}^2 \geq 6,\\
		H_{\mathrm{mix}}^{3}(\mathbb{T}^d) \quad &\text{ for } \quad \eta_{\ell}^2\mu_{\ell}^2 \geq 10.
	\end{cases}
\end{align*}

Contrary to the previous section concerned with the algebraic transformation \eqref{eq:algebraic_trafo_parametrized}, the $L_2$-approximation error can't be discussed this time as we can't compute the Fourier coefficients
\begin{align*}
	\hat{f}_{\mathbf k}
	= \int_{\mathbb{T}^d} f(\mathbf x,\bm\eta,\bm\mu) \, \mathrm{e}^{- 2 \pi\mathrm i \mathbf k\mathbf x} \, \mathrm d\mathbf x
	= \int_{\mathbb{T}} h(\psi(\mathbf x,\bm\eta)) \, \prod_{j=1}^{d}\eta_j^{\frac{1}{2}} \, \mathrm{e}^{\frac{1}{2}\left(1- \mu_j^2 \, \eta_j^2 \right) \mathrm{erf}^{-1}(2x_j)^2} \, \mathrm{e}^{- 2 \pi\mathrm i k_j x_j} \, \mathrm d\mathbf x
\end{align*}
regardless of the chosen $h$.
Even for trivial choices of $h$ we're not able to integrate the transformed weight function.

Hence, we only discuss the application of the weighted $L_{\infty}(\mathbb{R}^d)$-approximation error bound from Theorem~\ref{thm:L_infty_approx_error_multivar} for dimension $d=2$.
With the constant test function for $d=2$ given by $h(\mathbf y) = h(y_1,y_2) \equiv 1$, the weight function \eqref{eq:gaussian_weight_param} and the transformations \eqref{eq:erf_trafo_param},
the corresponding transformed functions $f$ in \eqref{eq:Error_trafo_function} read as
\begin{align*}
	{f(\mathbf x) = \prod_{j=1}^{2}\sqrt{\omega_j(\psi_j(x_j,\eta_j),\mu_j) \, \psi_j'(x_j,\eta_j)}}.
\end{align*}
Let $N\geq 8$, the two-dimensional hyperbolic cross $I_{N}^{2}$ as in \eqref{def:hyperbolic_cross} and a reconstructing {rank-$1$} lattice $\Lambda(\mathbf z, M, I_{N}^{2})$ be given. 
We already evaluated the sufficient conditions proposed in Theorem~\ref{thm:Hm_composition_criteria}, yielding lower bounds for $\bm\eta,\bm\mu\geq \mathbf 0$ such that $f$ is at least of Sobolev-smoothness order $m=0,1,2,3$, i.e., $f\in H_{\mathrm{mix}}^{m}(\mathbb{T}^2)$ and thus $f\in \mathcal{H}^{m}(\mathbb{T}^2)$.
We fix $\lambda=1$ and for $m\in\mathbb{N}_{0}$ we choose ${\bm\eta =(\eta_1,\eta_2)^{\top},\bm\mu =(\mu_1,\mu_2)^{\top}\in\mathbb{R}^2}$ such that ${f\in \mathcal{H}^{m+1}(\mathbb{T}^2)\hookrightarrow \mathcal{A}^{m}(\mathbb{T}^2)}$.
As outlined in \eqref{eq:discrete_upper_bound} we expect the discrete approximation error \eqref{eq:Discrete_ellinfty_error} to be bounded by
\begin{align*}
	\|\mathbf h - \mathbf h_{\mathrm{approx}}\|_{\ell_{\infty}} \leq
	\|f - S_{I_N^d}^{\Lambda} f\|_{L_{\infty}(\mathbb{T}^2)}
	\lesssim 
	\begin{cases}
		N^{0} \quad &\text{ for } \quad \eta_j^2\mu_j^2>0,\\
		N^{-1} \quad &\text{ for } \quad \eta_j^2\mu_j^2 > 2, \\
		N^{-2} \quad &\text{ for } \quad \eta_j^2\mu_j^2 \geq 6, \\
		N^{-3} \quad &\text{ for } \quad \eta_j^2\mu_j^2 \geq 10.
	\end{cases}
\end{align*}
\begin{figure}[t]
\begin{minipage}{0.5\linewidth}
\begin{tikzpicture}[scale=0.75]
		\begin{axis}[
		ymode = log,
		enlargelimits=false,
		xmin=0, xmax=90, ymin=1e-16, ymax=1e-0,
		ytick={1e-1,1e-2,1e-3,1e-4,1e-5,1e-6,1e-7,1e-8,1e-9,1e-10,1e-11,1e-12,1e-13,1e-14,1e-15,1e-16},
grid=both, 
		xlabel={$N$}, 
		ylabel={$\|\mathbf h - \mathbf h_{\mathrm{approx}}\|_{\ell_\infty}/\|\mathbf h\|_{\ell_\infty}$},
		title = {$h(\mathbf y) = 1$},
legend style={at={(0.5,1.20)}, anchor=south,legend columns=2,legend cell align=left, font=\small,  
		},
		xminorticks=false,
		yminorticks=false
		]
		\addplot[mark options={solid}, blue, mark=*, mark size=1.0] coordinates {
			(8, 4.4409e-16)  (9, 2.2204e-16)  (10, 4.4409e-16)  (11, 5.5511e-16)  (12, 5.5511e-16)  (13, 3.3307e-16)  (14, 4.4409e-16)  (15, 4.4409e-16)  (16, 2.3315e-15)  (17, 3.3307e-16)  (18, 3.7748e-15)  (19, 5.5511e-16)  (20, 3.3307e-15)  (21, 4.4409e-16)  (22, 3.4417e-15)  (23, 4.4409e-16)  (24, 5.5511e-16)  (25, 3.3307e-16)  (26, 5.5511e-16)  (27, 4.7740e-15)  (28, 4.4409e-16)  (29, 4.4409e-16)  (30, 1.8874e-15)  (31, 3.4417e-15)  (32, 7.7716e-16)  (33, 1.4433e-15)  (34, 3.3307e-16)  (35, 2.4425e-15)  (36, 6.6613e-16)  (37, 5.5511e-16)  (38, 6.6613e-16)  (39, 5.5511e-16)  (40, 6.6613e-16)  (41, 2.4425e-15)  (42, 5.5511e-16)  (43, 3.3307e-15)  (44, 5.5511e-16)  (45, 6.6613e-16)  (46, 5.3291e-15)  (47, 4.7740e-15)  (48, 4.2188e-15)  (49, 7.1054e-15)  (50, 4.4409e-16)  (51, 8.8818e-16)  (52, 5.5511e-16)  (53, 8.8818e-16)  (54, 2.6645e-15)  (55, 3.3307e-16)  (56, 7.7716e-16)  (57, 1.9984e-15)  (58, 3.9968e-15)  (59, 1.2212e-15)  (60, 3.2196e-15)  (61, 6.6613e-16)  (62, 9.9920e-16)  (63, 3.6637e-15)  (64, 2.9976e-15)  (65, 5.5511e-16)  (66, 4.4409e-16)  (67, 3.7748e-15)  (68, 5.5511e-16)  (69, 4.4409e-16)  (70, 7.7716e-16)  (71, 5.5511e-16)  (72, 4.4409e-16)  (73, 8.4377e-15)  (74, 5.5511e-16)  (75, 6.5503e-15)  (76, 3.3307e-16)  (77, 7.7716e-16)  (78, 5.5511e-16)  (79, 1.1102e-15)  (80, 7.7716e-16) 
		};
		\addlegendentry{$(\bm\eta,\bm\mu)=(\mathbf 1,\mathbf 1)$}
		\addplot[mark options={solid},red!75!yellow,mark=x,mark size=1.5] coordinates {
			(8, 1.8789e-02)  (9, 2.2904e-02)  (10, 1.2211e-02)  (11, 2.2071e-02)  (12, 2.0997e-02)  (13, 1.3103e-02)  (14, 1.9999e-02)  (15, 2.0084e-02)  (16, 1.2841e-02)  (17, 1.8599e-02)  (18, 1.8969e-02)  (19, 1.4077e-02)  (20, 1.7532e-02)  (21, 1.7442e-02)  (22, 1.3623e-02)  (23, 1.6380e-02)  (24, 1.6039e-02)  (25, 1.3633e-02)  (26, 1.5261e-02)  (27, 1.5462e-02)  (28, 1.3567e-02)  (29, 1.4439e-02)  (30, 1.4256e-02)  (31, 1.3910e-02)  (32, 1.3494e-02)  (33, 1.3650e-02)  (34, 1.3807e-02)  (35, 1.3728e-02)  (36, 1.3004e-02)  (37, 1.3378e-02)  (38, 1.3509e-02)  (39, 1.3628e-02)  (40, 1.2969e-02)  (41, 1.3284e-02)  (42, 1.2876e-02)  (43, 1.3169e-02)  (44, 1.2945e-02)  (45, 1.2751e-02)  (46, 1.2844e-02)  (47, 1.3094e-02)  (48, 1.2511e-02)  (49, 1.2601e-02)  (50, 1.2423e-02)  (51, 1.2498e-02)  (52, 1.2338e-02)  (53, 1.2537e-02)  (54, 1.2301e-02)  (55, 1.2269e-02)  (56, 1.1914e-02)  (57, 1.1975e-02)  (58, 1.2038e-02)  (59, 1.2202e-02)  (60, 1.1638e-02)  (61, 1.1792e-02)  (62, 1.1847e-02)  (63, 1.1737e-02)  (64, 1.1544e-02)  (65, 1.1523e-02)  (66, 1.1348e-02)  (67, 1.1477e-02)  (68, 1.1381e-02)  (69, 1.1424e-02)  (70, 1.1188e-02)  (71, 1.1305e-02)  (72, 1.0908e-02)  (73, 1.1018e-02)  (74, 1.1059e-02)  (75, 1.0977e-02)  (76, 1.0899e-02)  (77, 1.0879e-02)  (78, 1.0751e-02)  (79, 1.0847e-02)  (80, 1.0556e-02)
		};
		\addlegendentry{$(\bm\eta,\bm\mu)=(\mathbf 1,\sqrt{\mathbf 3})$}
		\addplot[mark options={solid},red!25!yellow,mark=triangle*,mark size=1.5] coordinates {
			(8, 1.6611e-04)  (9, 1.0213e-04)  (10, 9.0042e-05)  (11, 1.2709e-04)  (12, 9.8451e-05)  (13, 8.5408e-05)  (14, 9.6622e-05)  (15, 8.2238e-05)  (16, 7.2218e-05)  (17, 7.7151e-05)  (18, 6.5724e-05)  (19, 6.0206e-05)  (20, 6.0471e-05)  (21, 5.4699e-05)  (22, 4.9740e-05)  (23, 4.9705e-05)  (24, 4.4228e-05)  (25, 4.1853e-05)  (26, 4.0571e-05)  (27, 3.7867e-05)  (28, 3.5208e-05)  (29, 3.4447e-05)  (30, 3.1584e-05)  (31, 3.0087e-05)  (32, 2.8770e-05)  (33, 2.7072e-05)  (34, 2.5832e-05)  (35, 2.4901e-05)  (36, 2.3338e-05)  (37, 2.2618e-05)  (38, 2.1652e-05)  (39, 2.0569e-05)  (40, 1.9418e-05)  (41, 1.8833e-05)  (42, 1.7725e-05)  (43, 1.7338e-05)  (44, 1.6427e-05)  (45, 1.5869e-05)  (46, 1.5362e-05)  (47, 1.4829e-05)  (48, 1.4037e-05)  (49, 1.3664e-05)  (50, 1.3088e-05)  (51, 1.2588e-05)  (52, 1.2168e-05)  (53, 1.1741e-05)  (54, 1.1321e-05)  (55, 1.1104e-05)  (56, 1.0591e-05)  (57, 1.0228e-05)  (58, 9.9253e-06)  (59, 9.6117e-06)  (60, 9.2057e-06)  (61, 9.1099e-06)  (62, 8.7231e-06)  (63, 8.5187e-06)  (64, 8.2739e-06)  (65, 8.0137e-06)  (66, 7.7202e-06)  (67, 7.5649e-06)  (68, 7.2656e-06)  (69, 7.0642e-06)  (70, 6.9413e-06)  (71, 6.7047e-06)  (72, 6.5199e-06)  (73, 6.4557e-06)  (74, 6.2235e-06)  (75, 6.0623e-06)  (76, 5.8883e-06)  (77, 5.7339e-06)  (78, 5.5545e-06)  (79, 5.5068e-06)  (80, 5.2759e-06) 
		};
		\addlegendentry{$(\bm\eta,\bm\mu)=(\mathbf 1,\sqrt{\mathbf 6})$}
		\addplot[mark options={solid},red!50!blue,mark=square,mark size=1.0] coordinates {
			(8, 2.1708e-04)  (9, 1.5934e-04)  (10, 9.3054e-05)  (11, 1.1412e-04)  (12, 4.7945e-05)  (13, 5.6347e-05)  (14, 4.9568e-05)  (15, 3.6278e-05)  (16, 2.5462e-05)  (17, 2.7806e-05)  (18, 1.8725e-05)  (19, 2.0069e-05)  (20, 1.3972e-05)  (21, 1.1404e-05)  (22, 1.0786e-05)  (23, 1.1312e-05)  (24, 7.1762e-06)  (25, 6.5926e-06)  (26, 6.3576e-06)  (27, 5.6663e-06)  (28, 4.5493e-06)  (29, 4.7011e-06)  (30, 3.3400e-06)  (31, 3.4433e-06)  (32, 2.8872e-06)  (33, 2.6460e-06)  (34, 2.5960e-06)  (35, 2.3152e-06)  (36, 1.6673e-06)  (37, 1.7064e-06)  (38, 1.6805e-06)  (39, 1.5834e-06)  (40, 1.2424e-06)  (41, 1.2658e-06)  (42, 1.0268e-06)  (43, 1.0447e-06)  (44, 9.4279e-07)  (45, 8.0539e-07)  (46, 7.9672e-07)  (47, 8.0812e-07)  (48, 6.3334e-07)  (49, 6.1011e-07)  (50, 5.5187e-07)  (51, 5.2932e-07)  (52, 4.8808e-07)  (53, 4.9384e-07)  (54, 4.3132e-07)  (55, 4.0396e-07)  (56, 3.4437e-07)  (57, 3.3215e-07)  (58, 3.3004e-07)  (59, 3.3326e-07)  (60, 2.6075e-07)  (61, 2.6328e-07)  (62, 2.6192e-07)  (63, 2.3759e-07)  (64, 2.1639e-07)  (65, 2.0514e-07)  (66, 1.8410e-07)  (67, 1.8561e-07)  (68, 1.7629e-07)  (69, 1.7216e-07)  (70, 1.5277e-07)  (71, 1.5387e-07)  (72, 1.2581e-07)  (73, 1.2672e-07)  (74, 1.2628e-07)  (75, 1.1772e-07)  (76, 1.1272e-07)  (77, 1.0737e-07)  (78, 9.9009e-08)  (79, 9.9601e-08)  (80, 8.6612e-08)     
		};
		\addlegendentry{$(\bm\eta,\bm\mu)=(\mathbf 1,\sqrt{\mathbf{10}})$}
		\end{axis}
		\end{tikzpicture}
	\end{minipage}
	\begin{minipage}{0.5\linewidth}
\begin{tikzpicture}[scale=0.75]
		\begin{axis}[
		ymode = log,
		enlargelimits=false,
		xmin=0, xmax=90, ymin=1e-9, ymax=1e-0,
		ytick={1e-1,1e-2,1e-3,1e-4,1e-5,1e-6,1e-7,1e-8,1e-9,1e-10,1e-11},
grid=both, 
		xlabel={$N$}, 
		ylabel={$\|\mathbf h - \mathbf h_{\mathrm{approx}}\|_{\ell_\infty}/\|\mathbf h\|_{\ell_\infty}$},
		title = {$h(\mathbf y) = \mathrm{e}^{-y_1^2-y_2^2}$},
legend style={at={(0.5,1.20)}, anchor=south,legend columns=2,legend cell align=left, font=\small},
		xminorticks=false,
		yminorticks=false
		]
		\addplot[mark options={solid}, blue, mark=*, mark size=1.0] coordinates {
			(8, 1.8789e-02)  (9, 2.2904e-02)  (10, 1.2211e-02)  (11, 2.2071e-02)  (12, 2.0997e-02)  (13, 1.3103e-02)  (14, 1.9999e-02)  (15, 2.0084e-02)  (16, 1.2841e-02)  (17, 1.8599e-02)  (18, 1.8969e-02)  (19, 1.4077e-02)  (20, 1.7532e-02)  (21, 1.7442e-02)  (22, 1.3623e-02)  (23, 1.6380e-02)  (24, 1.6039e-02)  (25, 1.3633e-02)  (26, 1.5261e-02)  (27, 1.5462e-02)  (28, 1.3567e-02)  (29, 1.4439e-02)  (30, 1.4256e-02)  (31, 1.3910e-02)  (32, 1.3494e-02)  (33, 1.3650e-02)  (34, 1.3807e-02)  (35, 1.3728e-02)  (36, 1.3004e-02)  (37, 1.3378e-02)  (38, 1.3509e-02)  (39, 1.3628e-02)  (40, 1.2969e-02)  (41, 1.3284e-02)  (42, 1.2876e-02)  (43, 1.3169e-02)  (44, 1.2945e-02)  (45, 1.2751e-02)  (46, 1.2844e-02)  (47, 1.3094e-02)  (48, 1.2511e-02)  (49, 1.2601e-02)  (50, 1.2423e-02)  (51, 1.2498e-02)  (52, 1.2338e-02)  (53, 1.2537e-02)  (54, 1.2301e-02)  (55, 1.2269e-02)  (56, 1.1914e-02)  (57, 1.1975e-02)  (58, 1.2038e-02)  (59, 1.2202e-02)  (60, 1.1638e-02)  (61, 1.1792e-02)  (62, 1.1847e-02)  (63, 1.1737e-02)  (64, 1.1544e-02)  (65, 1.1523e-02)  (66, 1.1348e-02)  (67, 1.1477e-02)  (68, 1.1381e-02)  (69, 1.1424e-02)  (70, 1.1188e-02)  (71, 1.1305e-02)  (72, 1.0908e-02)  (73, 1.1018e-02)  (74, 1.1059e-02)  (75, 1.0977e-02)  (76, 1.0899e-02)  (77, 1.0879e-02)  (78, 1.0751e-02)  (79, 1.0847e-02)  (80, 1.0556e-02)
		};
		\addlegendentry{$(\bm\eta,\bm\mu)=(\mathbf 1,\mathbf 1)$}
		\addplot[mark options={solid},red!75!yellow,mark=x,mark size=1.5] coordinates {
			(8, 1.0912e-03)  (9, 1.1690e-03)  (10, 7.0072e-04)  (11, 8.0753e-04)  (12, 7.0518e-04)  (13, 6.5537e-04)  (14, 6.0715e-04)  (15, 5.6625e-04)  (16, 4.8187e-04)  (17, 5.2696e-04)  (18, 4.4600e-04)  (19, 4.7778e-04)  (20, 4.0339e-04)  (21, 3.8553e-04)  (22, 3.7438e-04)  (23, 3.9263e-04)  (24, 3.2085e-04)  (25, 3.1860e-04)  (26, 3.1177e-04)  (27, 3.0306e-04)  (28, 2.7415e-04)  (29, 2.8304e-04)  (30, 2.4526e-04)  (31, 2.5233e-04)  (32, 2.3275e-04)  (33, 2.2759e-04)  (34, 2.2437e-04)  (35, 2.1713e-04)  (36, 1.8930e-04)  (37, 1.9338e-04)  (38, 1.9105e-04)  (39, 1.8784e-04)  (40, 1.6926e-04)  (41, 1.7229e-04)  (42, 1.5794e-04)  (43, 1.6053e-04)  (44, 1.5277e-04)  (45, 1.4477e-04)  (46, 1.4339e-04)  (47, 1.4540e-04)  (48, 1.3088e-04)  (49, 1.3025e-04)  (50, 1.2484e-04)  (51, 1.2333e-04)  (52, 1.1855e-04)  (53, 1.1994e-04)  (54, 1.1308e-04)  (55, 1.1103e-04)  (56, 1.0399e-04)  (57, 1.0289e-04)  (58, 1.0218e-04)  (59, 1.0318e-04)  (60, 9.3103e-05)  (61, 9.3986e-05)  (62, 9.3410e-05)  (63, 9.0336e-05)  (64, 8.6658e-05)  (65, 8.5389e-05)  (66, 8.1566e-05)  (67, 8.2231e-05)  (68, 8.0066e-05)  (69, 7.9433e-05)  (70, 7.5744e-05)  (71, 7.6308e-05)  (72, 7.0373e-05)  (73, 7.0882e-05)  (74, 7.0537e-05)  (75, 6.8727e-05)  (76, 6.7164e-05)  (77, 6.6340e-05)  (78, 6.4006e-05)  (79, 6.4410e-05)  (80, 6.0850e-05)
		};
		\addlegendentry{$(\bm\eta,\bm\mu)=(\mathbf 1,\sqrt{\mathbf 3})$}
		\addplot[mark options={solid},red!25!yellow,mark=triangle*,mark size=1.5] coordinates {
			(8, 2.3962e-04)  (9, 1.3450e-04)  (10, 1.1239e-04)  (11, 1.1145e-04)  (12, 6.9902e-05)  (13, 7.6372e-05)  (14, 7.6966e-05)  (15, 6.2660e-05)  (16, 5.3530e-05)  (17, 5.8402e-05)  (18, 4.6651e-05)  (19, 5.0001e-05)  (20, 4.0722e-05)  (21, 3.5966e-05)  (22, 3.6070e-05)  (23, 3.7758e-05)  (24, 2.8433e-05)  (25, 2.7175e-05)  (26, 2.7245e-05)  (27, 2.5361e-05)  (28, 2.2440e-05)  (29, 2.3116e-05)  (30, 1.8689e-05)  (31, 1.9204e-05)  (32, 1.7404e-05)  (33, 1.6501e-05)  (34, 1.6533e-05)  (35, 1.5364e-05)  (36, 1.2559e-05)  (37, 1.2812e-05)  (38, 1.2831e-05)  (39, 1.2361e-05)  (40, 1.0652e-05)  (41, 1.0825e-05)  (42, 9.5124e-06)  (43, 9.6563e-06)  (44, 9.1108e-06)  (45, 8.2384e-06)  (46, 8.2464e-06)  (47, 8.3499e-06)  (48, 7.1553e-06)  (49, 7.0062e-06)  (50, 6.6103e-06)  (51, 6.4425e-06)  (52, 6.1544e-06)  (53, 6.2175e-06)  (54, 5.7092e-06)  (55, 5.4837e-06)  (56, 4.9800e-06)  (57, 4.8709e-06)  (58, 4.8745e-06)  (59, 4.9162e-06)  (60, 4.2079e-06)  (61, 4.2430e-06)  (62, 4.2460e-06)  (63, 3.9879e-06)  (64, 3.7713e-06)  (65, 3.6516e-06)  (66, 3.4165e-06)  (67, 3.4410e-06)  (68, 3.3402e-06)  (69, 3.2896e-06)  (70, 3.0561e-06)  (71, 3.0757e-06)  (72, 2.7119e-06)  (73, 2.7289e-06)  (74, 2.7304e-06)  (75, 2.6099e-06)  (76, 2.5448e-06)  (77, 2.4719e-06)  (78, 2.3496e-06)  (79, 2.3623e-06)  (80, 2.1670e-06) 
		};
		\addlegendentry{$(\bm\eta,\bm\mu)=(\mathbf 1,\sqrt{\mathbf 6})$}
		\addplot[mark options={solid},red!50!blue,mark=square,mark size=1.0] coordinates {
			(8, 7.7922e-05)  (9, 3.0210e-05)  (10, 1.7557e-05)  (11, 1.6967e-05)  (12, 7.4736e-06)  (13, 8.0272e-06)  (14, 5.8982e-06)  (15, 6.0976e-06)  (16, 4.2921e-06)  (17, 4.2925e-06)  (18, 3.0561e-06)  (19, 3.0327e-06)  (20, 2.2836e-06)  (21, 2.1180e-06)  (22, 1.5862e-06)  (23, 1.5352e-06)  (24, 1.1233e-06)  (25, 1.1148e-06)  (26, 8.7204e-07)  (27, 8.1325e-07)  (28, 6.2655e-07)  (29, 6.5431e-07)  (30, 4.7030e-07)  (31, 4.9390e-07)  (32, 4.0042e-07)  (33, 3.6354e-07)  (34, 3.1041e-07)  (35, 3.2181e-07)  (36, 2.4250e-07)  (37, 2.3785e-07)  (38, 2.1056e-07)  (39, 1.9081e-07)  (40, 1.6023e-07)  (41, 1.6170e-07)  (42, 1.3277e-07)  (43, 1.3133e-07)  (44, 1.1109e-07)  (45, 1.0689e-07)  (46, 9.1255e-08)  (47, 9.1888e-08)  (48, 7.6995e-08)  (49, 7.4628e-08)  (50, 6.4493e-08)  (51, 6.2983e-08)  (52, 5.4577e-08)  (53, 5.2860e-08)  (54, 4.6621e-08)  (55, 4.5067e-08)  (56, 3.9569e-08)  (57, 3.8709e-08)  (58, 3.4001e-08)  (59, 3.3014e-08)  (60, 2.9180e-08)  (61, 2.8318e-08)  (62, 2.6449e-08)  (63, 2.4656e-08)  (64, 2.2282e-08)  (65, 2.2624e-08)  (66, 2.0030e-08)  (67, 2.0180e-08)  (68, 1.8808e-08)  (69, 1.7746e-08)  (70, 1.6713e-08)  (71, 1.6835e-08)  (72, 1.4344e-08)  (73, 1.4431e-08)  (74, 1.3841e-08)  (75, 1.3134e-08)  (76, 1.2261e-08)  (77, 1.2158e-08)  (78, 1.0979e-08)  (79, 1.1016e-08)  (80, 1.0143e-08)     
		};
		\addlegendentry{$(\bm\eta,\bm\mu)=(\mathbf 1,\sqrt{\mathbf{10}})$}
		\end{axis}
		\end{tikzpicture}
	\end{minipage}
	\caption{Comparison of discrete $\ell_\infty$-approximation error $\| \mathbf h - \mathbf h_{\text{approx}} \|_{\ell_\infty}/\|\mathbf h\|_{\ell_\infty}$ when using 
		the Gaussian weight function $\omega(\circ,\bm\mu)$ as in \eqref{eq:gaussian_weight_param} 
		and the error function transformation $\psi(\circ,\bm\eta)$ as in \eqref{eq:erf_trafo_param}
		with $\bm\mu \in \{\mathbf 1, \sqrt{\mathbf 3}, \sqrt{\mathbf 6}, \sqrt{\mathbf{10}} \}$ and fixed $\bm\eta = \mathbf 1$.}
	\label{fig:numeric_ell2error_errorfct}
\end{figure}
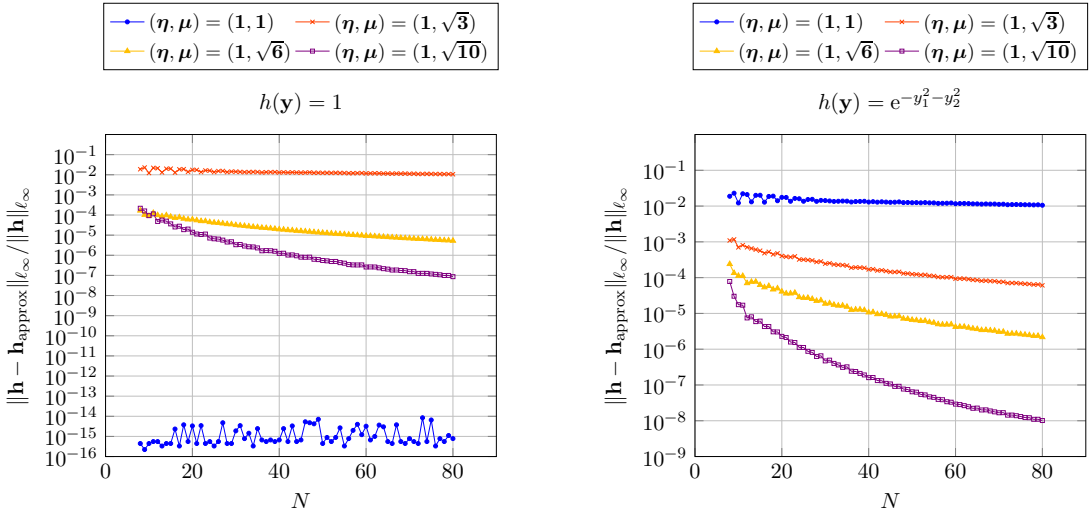
We actually observe this behavior numerically, as showcased in Figure~\ref{fig:numeric_ell2error_errorfct}, where we show in the left graph the approximation error decay of the constant test function $h(\mathbf y) \equiv 1$ for ${N=8,\ldots,80}$, fixed ${\bm\eta = \mathbf 1}$ and $\bm\mu\in\{\mathbf 1, \sqrt{\mathbf 3}, \sqrt{\mathbf 6}, \sqrt{\mathbf{10}}\}$.
The outlier with $\bm\eta = \bm\mu = \mathbf 1$ is explained by the fact the corresponding Fourier coefficients are trivial as these parameters lead to a constant weight function $\omega(\mathbf y) \equiv 1$.
We repeat this numerical test with the non-constant test function ${h(\mathbf y) = \mathrm{e}^{-y_1^2-y_2^2}}$ that is in $L_{2}(\mathbb{R}^2,\omega(\circ,\bm\mu))$ for all $\bm\mu\in\mathbb{R}^2$ with $\mu_1,\mu_2 >-2$.
Then we have a similar decay of the discrete approximation error as shown in the right graph of Figure~\ref{fig:numeric_ell2error_errorfct}.

\section{Remarks on multiple {rank-1} lattices and sparse frequency sets}
Now that we are able to construct functions on the torus $\mathbb{T}^d$ with a guaranteed minimal Sobolev-smoothness degree $m\in\mathbb{N}_{0}$, we adapt the techniques of both multiple {rank-$1$} lattices \cite{Kae16} and sparse FFT algorithms \cite{PoVo14}.
Usually we consider the algebraic test function in \eqref{eq:exemplary_h} that was given by
\begin{align*} h(\mathbf y) = \frac{1}{1 + \|\mathbf y\|_{\ell_2}^2}, \quad \mathbf y\in\mathbb{R}^d.
\end{align*}
\subsection{Multiple {rank-1} lattices}
In Lemma~\ref{lem:existence_of_M} we recalled that under mild assumptions it's possible to generate a reconstructing {rank-$1$} lattice $\Lambda(\mathbf z,M,I)$ with some frequency set $I\subset\mathbb{Z}^d$ of finite cardinality $|I|<\infty$ such that 
\begin{align*}
	|I| \leq M \leq |I|^2.
\end{align*}
Even though this upper bound is independent of the dimension $d$, if the lattice size $M$ is usually close to $|I|^2$ and is therefore still pretty large.
In order to overcome this limitation of the single rank\mbox{-}1 lattice approach
L. K\"ammerer suggested the use of multiple {rank-$1$} lattices which are obtained by taking a union of $s$ {rank-$1$} lattices 
\begin{align*}
	\Lambda(\mathbf z_1, M_1,\ldots, \mathbf z_s, M_s) := \bigcup_{j=1,\ldots,s} \Lambda(\mathbf z_j, M_j),
\end{align*}
see \cite{Kae16, Kae17}. 
Then it's possible to determine a reconstructing sampling set for multivariate trigonometric polynomials in $\Pi_{I}$ supported on the given frequency set $I$ with a probability of at least $1-\delta_s$, where
\begin{align*}
	\delta_s=C_1\,{\mathrm e}^{-C_2s}
\end{align*}
is an upper bound on the probability that the approach fails and $C_1,C_2>0$ are constants.
In \cite{Kae17} it was proven that the upper bound on the lattice size improves with high probability to 
\begin{align*}
	M\leq C |I|\log|I|
\end{align*}
for these particular reconstructing lattices.
For the adaptation of this approach in the context of families of transformations $\psi(\circ,\bm\eta)$ with $\bm\eta\in\mathbb{R}^d$ we analogously consider unions of $s$ transformed {rank-$1$} lattices 
\begin{align*}
	\Lambda_{\psi(\circ,\bm\eta)}(\mathbf z_1, M_1,\ldots, \mathbf z_s, M_s) := \bigcup_{j=1,\ldots,s} \Lambda_{\psi(\circ,\bm\eta)}(\mathbf z_j, M_j)
\end{align*}
in order to sample the test function $h\in L_2(\mathbb{R}^d, \omega)$.

\begin{figure}[t]
	\centering
	\begin{minipage}{0.5\linewidth}
\begin{tikzpicture}[baseline,scale=0.75]
		\begin{axis}[
		ymode = log,
		enlargelimits=false,
		xmin=0, xmax=90, ymin=1e-12, ymax=1e-0,
		ytick={1e-1,1e-2,1e-3,1e-4,1e-5,1e-6,1e-7,1e-8,1e-9,1e-10,1e-11},
grid=both, 
		xlabel={$N$}, 
		ylabel={$\|\mathbf h - \mathbf h_{\mathrm{approx}}\|_{\ell_\infty}/\|\mathbf h\|_{\ell_\infty}$},
legend style={at={(0.5,1.05)}, anchor=south,legend columns=2,legend cell align=left, font=\small,  
		},
		xminorticks=false,
		yminorticks=false
		]
		\addplot[mark options={solid}, blue, mark=*, mark size=1.0] coordinates {
			(8, 6.8638e-02)  (9, 6.2104e-02)  (10, 4.9204e-02)  (11, 4.8943e-02)  (12, 3.9966e-02)  (13, 3.9167e-02)  (14, 3.9482e-02)  (15, 4.0028e-02)  (16, 3.7429e-02)  (17, 3.9248e-02)  (18, 3.8473e-02)  (19, 3.9902e-02)  (20, 3.7869e-02)  (21, 3.7468e-02)  (22, 3.7934e-02)  (23, 3.6672e-02)  (24, 3.5485e-02)  (25, 3.6071e-02)  (26, 3.6225e-02)  (27, 3.5314e-02)  (28, 3.5413e-02)  (29, 3.3265e-02)  (30, 3.5018e-02)  (31, 3.3603e-02)  (32, 3.2760e-02)  (33, 3.0635e-02)  (34, 3.2597e-02)  (35, 3.4373e-02)  (36, 3.2853e-02)  (37, 3.2947e-02)  (38, 3.3813e-02)  (39, 3.0932e-02)  (40, 3.1714e-02)  (41, 3.3167e-02)  (42, 3.2625e-02)  (43, 2.9557e-02)  (44, 3.2522e-02)  (45, 3.0400e-02)  (46, 2.9113e-02)  (47, 2.9907e-02)  (48, 2.2369e-02)  (49, 2.9641e-02)  (50, 3.1149e-02)  (51, 3.0573e-02)  (52, 2.9462e-02)  (53, 3.1092e-02)  (54, 2.8222e-02)  (55, 2.2376e-02)  (56, 2.8425e-02)  (57, 2.7029e-02)  (58, 2.8941e-02)  (59, 2.9837e-02)  (60, 2.9157e-02)  (61, 2.9242e-02)  (62, 2.7088e-02)  (63, 2.7674e-02)  (64, 2.8837e-02)  (65, 2.7044e-02)  (66, 2.8032e-02)  (67, 2.7177e-02)  (68, 2.8175e-02)  (69, 2.8623e-02)  (70, 2.6734e-02)  (71, 2.7783e-02)  (72, 2.5566e-02)  (73, 2.6235e-02)  (74, 2.6036e-02)  (75, 2.6244e-02)  (76, 2.6482e-02)  (77, 2.6214e-02)  (78, 2.5272e-02)  (79, 2.5087e-02)  (80, 2.6007e-02)   
		};
		\addlegendentry{$(\bm\eta,\bm\mu)=(\mathbf 1,\mathbf 0)$}
		\addplot[mark options={solid},red!75!yellow,mark=square*,mark size=1.0] coordinates {
			(8, 7.6414e-04)  (9, 4.5162e-04)  (10, 3.7231e-04)  (11, 3.1595e-04)  (12, 1.8792e-04)  (13, 1.6654e-04)  (14, 1.3497e-04)  (15, 1.2233e-04)  (16, 7.9869e-05)  (17, 7.7522e-05)  (18, 7.1416e-05)  (19, 6.8717e-05)  (20, 4.8129e-05)  (21, 4.6326e-05)  (22, 4.1942e-05)  (23, 4.3866e-05)  (24, 3.4242e-05)  (25, 3.0869e-05)  (26, 2.7832e-05)  (27, 2.9101e-05)  (28, 2.4590e-05)  (29, 2.4778e-05)  (30, 1.9620e-05)  (31, 1.9877e-05)  (32, 1.7417e-05)  (33, 1.7706e-05)  (34, 1.6477e-05)  (35, 1.4964e-05)  (36, 1.2602e-05)  (37, 1.2691e-05)  (38, 1.2074e-05)  (39, 1.2140e-05)  (40, 1.0470e-05)  (41, 1.0216e-05)  (42, 8.9328e-06)  (43, 8.9531e-06)  (44, 8.5466e-06)  (45, 7.7470e-06)  (46, 7.7190e-06)  (47, 7.7293e-06)  (48, 6.7958e-06)  (49, 6.3088e-06)  (50, 5.9554e-06)  (51, 5.9105e-06)  (52, 5.9482e-06)  (53, 5.8629e-06)  (54, 5.2076e-06)  (55, 5.1428e-06)  (56, 4.5385e-06)  (57, 4.6663e-06)  (58, 4.5924e-06)  (59, 4.6401e-06)  (60, 4.0922e-06)  (61, 4.1335e-06)  (62, 4.2071e-06)  (63, 3.8141e-06)  (64, 3.6049e-06)  (65, 3.5659e-06)  (66, 3.2965e-06)  (67, 3.2743e-06)  (68, 3.3810e-06)  (69, 3.3356e-06)  (70, 3.0655e-06)  (71, 3.0632e-06)  (72, 2.5410e-06)  (73, 2.6645e-06)  (74, 2.6503e-06)  (75, 2.6515e-06)  (76, 2.6690e-06)  (77, 2.5114e-06)  (78, 2.4180e-06)  (79, 2.3257e-06)  (80, 2.2069e-06) 
		};
		\addlegendentry{$(\bm\eta,\bm\mu)=(\mathbf 1,\mathbf{4})$}
		\addplot[mark options={solid},red!25!yellow,mark=triangle*,mark size=1.5] coordinates {
			(8, 3.7913e-04)  (9, 1.6100e-04)  (10, 6.2517e-05)  (11, 6.3171e-05)  (12, 1.9441e-05)  (13, 1.9461e-05)  (14, 1.1859e-05)  (15, 1.1752e-05)  (16, 3.9492e-06)  (17, 3.9309e-06)  (18, 3.7268e-06)  (19, 3.7176e-06)  (20, 1.4858e-06)  (21, 1.4275e-06)  (22, 1.2286e-06)  (23, 1.2318e-06)  (24, 8.3472e-07)  (25, 5.0166e-07)  (26, 4.4341e-07)  (27, 4.5511e-07)  (28, 3.7564e-07)  (29, 3.7584e-07)  (30, 1.7983e-07)  (31, 1.7933e-07)  (32, 1.6252e-07)  (33, 1.6067e-07)  (34, 1.5639e-07)  (35, 1.0832e-07)  (36, 7.2172e-08)  (37, 7.2020e-08)  (38, 7.1254e-08)  (39, 6.8502e-08)  (40, 5.5994e-08)  (41, 5.6480e-08)  (42, 3.2138e-08)  (43, 3.1980e-08)  (44, 3.2698e-08)  (45, 2.8210e-08)  (46, 2.8259e-08)  (47, 2.8153e-08)  (48, 1.9998e-08)  (49, 1.5314e-08)  (50, 1.4673e-08)  (51, 1.4141e-08)  (52, 1.4140e-08)  (53, 1.3999e-08)  (54, 1.1322e-08)  (55, 1.1123e-08)  (56, 7.4852e-09)  (57, 7.2882e-09)  (58, 7.2013e-09)  (59, 7.2482e-09)  (60, 6.1521e-09)  (61, 6.1208e-09)  (62, 6.1252e-09)  (63, 4.6569e-09)  (64, 3.8241e-09)  (65, 3.8137e-09)  (66, 3.4416e-09)  (67, 3.4449e-09)  (68, 3.3616e-09)  (69, 3.3416e-09)  (70, 2.7555e-09)  (71, 2.7730e-09)  (72, 1.9049e-09)  (73, 1.8960e-09)  (74, 1.8968e-09)  (75, 1.8871e-09)  (76, 1.8655e-09)  (77, 1.6078e-09)  (78, 1.5807e-09)  (79, 1.5778e-09)  (80, 1.2253e-09)   
		};
		\addlegendentry{$(\bm\eta,\bm\mu)=(\mathbf 1,\mathbf{10})$}
		\addplot[mark options={solid},red!50!blue,mark=square,mark size=1.0] coordinates {
			(8, 2.9745e-03)  (9, 4.2933e-04)  (10, 4.0915e-04)  (11, 4.0912e-04)  (12, 7.1303e-05)  (13, 7.0960e-05)  (14, 7.0602e-05)  (15, 2.9369e-05)  (16, 6.2199e-06)  (17, 6.0169e-06)  (18, 2.0824e-06)  (19, 2.1056e-06)  (20, 9.7079e-07)  (21, 5.4965e-07)  (22, 5.4431e-07)  (23, 5.4983e-07)  (24, 4.4990e-07)  (25, 1.1470e-07)  (26, 1.1910e-07)  (27, 1.1838e-07)  (28, 9.1304e-08)  (29, 9.1271e-08)  (30, 3.4315e-08)  (31, 3.4067e-08)  (32, 2.1835e-08)  (33, 2.1512e-08)  (34, 2.1277e-08)  (35, 1.8228e-08)  (36, 5.7801e-09)  (37, 5.7943e-09)  (38, 5.7839e-09)  (39, 5.9061e-09)  (40, 4.9992e-09)  (41, 4.9619e-09)  (42, 1.9620e-09)  (43, 1.8996e-09)  (44, 1.8849e-09)  (45, 1.4190e-09)  (46, 1.4336e-09)  (47, 1.4369e-09)  (48, 1.2150e-09)  (49, 8.0128e-10)  (50, 5.3793e-10)  (51, 4.8833e-10)  (52, 5.1980e-10)  (53, 5.1722e-10)  (54, 4.6225e-10)  (55, 3.8549e-10)  (56, 2.2724e-10)  (57, 1.9614e-10)  (58, 1.9521e-10)  (59, 1.9485e-10)  (60, 1.3581e-10)  (61, 1.3559e-10)  (62, 1.3368e-10)  (63, 1.0241e-10)  (64, 7.8715e-11)  (65, 6.6308e-11)  (66, 5.2845e-11)  (67, 5.2870e-11)  (68, 5.2974e-11)  (69, 5.8122e-11)  (70, 4.7129e-11)  (71, 4.6641e-11)  (72, 2.5101e-11)  (73, 2.6192e-11)  (74, 2.8689e-11)  (75, 2.2198e-11)  (76, 2.2448e-11)  (77, 2.7823e-11)  (78, 2.0070e-11)  (79, 2.0376e-11)  (80, 1.4087e-11)     
		};
		\addlegendentry{$(\bm\eta,\bm\mu)=(\mathbf 1,\mathbf{16})$}
		\end{axis}
		\end{tikzpicture}
	\end{minipage}
	\caption{Comparison of discrete $\ell_\infty$-approximation error $\| \mathbf h - \mathbf h_{\text{approx}} \|_{\ell_\infty}/\|\mathbf h\|_{\ell_\infty}$ of test function 		\eqref{eq:exemplary_h} 
		for multiple rank-$1$ lattices $\Lambda_{\psi(\circ,\bm\eta)}(\mathbf z_1, M_1,\ldots, \mathbf z_s, M_s)$ with 
		the algebraic transformation $\psi(\circ,\bm\eta)$ \eqref{eq:algebraic_trafo_parametrized} 
		and the algebraic weight function $\omega(\circ,\bm\mu)$ \eqref{eq:alebraic_weight_function_parametrized}
		in their two-dimensional versions with fixed $\bm\eta = \mathbf 1$ and $\bm\mu \in \{ \mathbf 0, \mathbf 4, \mathbf{10}, \mathbf{16}\}$ .}
	\label{fig:numeric_algebraic_multr1l}
\end{figure}
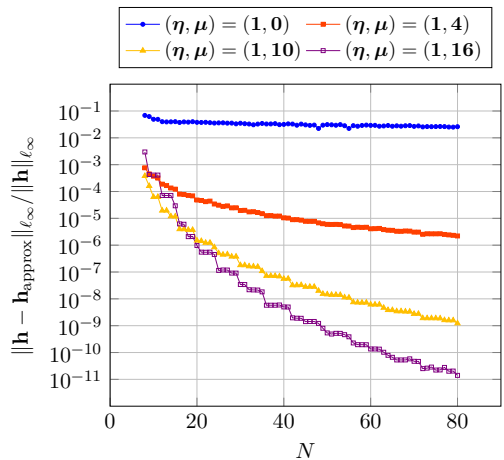

For an example in dimension $d=2$ we consider the test function $h$ in \eqref{eq:exemplary_h}, the algebraic weight function 
\begin{align*}
	\omega(\mathbf y,\bm\mu) := \left( \frac{1}{1+y_1^2} \right)^{\mu_1} \left( \frac{1}{1+y_2^2} \right)^{\mu_2}
\end{align*}
and algebraic transformation 
\begin{align} \label{eq:alg_trafo_two_dim}
	\psi(\mathbf x, \bm\eta) = \left(\frac{2\eta_1 x_1}{\sqrt{1-4x_1^2}}, \frac{2\eta_2 x_2}{\sqrt{1-4x_2^2}}\right)^{\top}
\end{align}
based on their univariate versions in \eqref{eq:alebraic_weight_function_parametrized} and \eqref{eq:algebraic_trafo_parametrized}.
We consider the sample data vector ${\mathbf{h}=\left(h(\mathbf y_j)\, \sqrt{\frac{\omega(\mathbf y_j, \bm\mu)}{\varrho(\mathbf y_j, \bm\eta)}}\right)_{j=0}^{M-1}}$ 
and the corresponding approximated data vector of the form ${\mathbf{h}_{\mathrm{approx}} = \left(\sqrt{\frac{\omega(\mathbf y_j, \bm\mu)}{\varrho(\mathbf y_j, \bm\eta)}}\,S_{I_{N}^{d}}^{\Lambda}h(\mathbf y_j)\right)_{j=0}^{M-1}}$
with lattice points $\mathbf y_j$ from the multiple rank-$1$ lattice ${\Lambda_{\psi(\circ,\bm\eta)}(\mathbf z_1, M_1,\ldots, \mathbf z_s, M_s,I)}$ transformed by the algebraic transformation $\psi(\circ, \bm\eta)$ in \eqref{eq:alg_trafo_two_dim}.
In \eqref{eq:discrete_bounds_exact} we already discussed that the discretized approximation error defined in \eqref{eq:Discrete_ellinfty_error} is bounded above by
\begin{align*}
	\|\mathbf h - \mathbf h_{\mathrm{approx}}\|_{\ell_{\infty}} 
	\leq \|f - S_{I_N^2}^{\Lambda} f\|_{L_{\infty}(\mathbb{T}^2)}
	\lesssim N^{-m}
\end{align*}
for $\mu_j \geq 0$ if $m=0$ and for $\mu_j > 3m$ if $m=1,2,3$.
Similarly to the results of the numerical test with single rank-$1$ lattices shown in Figure~\ref{fig:numeric_ell2error_algebraic},
we achieve this behavior of the relative discrete approximation error $\|\mathbf h - \mathbf h_{\mathrm{approx}}\|_{\ell_{\infty}}/\|\mathbf h\|_{\ell_{\infty}}$
when applying the multiple rank-1 algorithms described in \cite{Kae16,Kae17}.
In particular we adapted \cite[Algorithm~6]{Kae17}.
For $N=8,\ldots,80$, ${\bm\mu\in\{\mathbf 0,\mathbf{4},\mathbf{10},\mathbf{16}\}}$ and $\bm\eta=\mathbf 1$ we initialize this Algorithm with the parameters ${c=30}, {n=30}$ and $\delta=0.5$ and still have the proposed decay rates of the discrete approximation errors as seen in Figure~\ref{fig:numeric_algebraic_multr1l}.
A major advantage of this approach is that we don't need to construct the generating vector $\mathbf z$ via component-by-component construction methods, that generally takes quite some time.
\subsection{The construction of sparse frequency sets}
Likewise, once we set up the transformed function $f$ on the torus of the form \eqref{eq:f_is_transformed_h}, we can make use of dimension incremental algorithms -- the sparse fast Fourier transforms (sparse FFT), see \cite{PoVo14} and \cite{Vo_FFTr1l} -- that reconstruct sparse multivariate trigonometric polynomials with an unknown support in a frequency domain $I\subset\mathbb{Z}^d$. 
Based on component-by-component construction of {rank-$1$} lattices the approach of \cite[Algorithm~1 and Algorithm~2]{PoVo14} describes a dimension incremental construction of a frequency set $I\subset\mathbb{Z}^d$  belonging to the non-zero or approximately largest Fourier coefficients. 
This is achieved by restricting the search space to a full grid $[-N,N]^d\cap\mathbb{Z}^d$ of refinement ${N\in\mathbb{N}}$ and by assuming that the cardinality of the support of the multivariate trigonometric polynomial is bounded by a sparsity $s\in\mathbb{N}$. 
Then we end up with up to $s$ non-zero Fourier coefficients $\hat f_{\mathbf k}$ of the respective test function $f$.

We adapt these algorithms for transformed reconstructing {rank-$1$} lattices $\Lambda_{\psi(\circ,\bm\eta)}(\mathbf z,M, I)$ by again calculating the relative discretized approximation error $\|\mathbf h - \mathbf h_{\mathrm{approx}}\|_{\ell_{\infty}}/\|\mathbf h\|_{\ell_{\infty}}$ as in \eqref{eq:Discrete_ellinfty_error} with samples $\mathbf{h}=\left(h(\mathbf y_j)\, \sqrt{\frac{\omega(\mathbf y_j, \bm\mu)}{\varrho(\mathbf y_j, \bm\eta)}}\right)_{j=0}^{M-1}$ and $\mathbf{h}_{\mathrm{approx}} = \left( \sqrt{\frac{\omega(\mathbf y_j, \bm\mu)}{\varrho(\mathbf y_j, \bm\eta)}}\,S_{I}^{\Lambda} h(\mathbf y_j)\right)_{j=0}^{M-1}$ but use an unknown frequency set $I$ 
with cardinality $|I|=s$ that was constructed via a dimensional incremental construction method as outlined above.

\subsubsection{Example for the algebraic transformation}
We use the algebraic test function \eqref{eq:exemplary_h}
in combination with the multivariate version of algebraic weight function \eqref{eq:alebraic_weight_function_parametrized}
and the multivariate algebraic transformation based on \eqref{eq:algebraic_trafo_parametrized}, reading as
\begin{align*} 
	\omega(\mathbf y,\bm\mu) = \prod_{j=1}^{d} \left( \frac{1}{1+y_j^2} \right)^{\mu_j}, \quad
	\psi(\mathbf x,\bm\eta) = \left( \frac{2 \eta_1 x_1}{\sqrt{1-4x_{1}^2}}, \ldots, \frac{2 \eta_d x_d}{\sqrt{1-4x_{d}^2}} \right)^{\top}
\end{align*}
with $\bm\mu = \mathbf 4$ and $\bm\eta = \mathbf 1$.
Earlier we used a similar setup for $d=2$ where we choose a hyperbolic cross $I_{N}^{d}$ as the frequency set. 

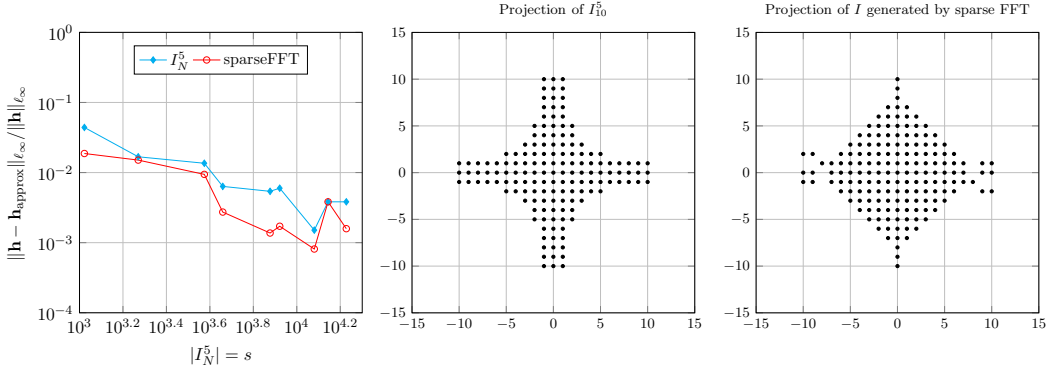
\begin{figure}[]
	\centering
\begin{tikzpicture}[baseline,scale=0.6]
	\begin{loglogaxis}[
	enlargelimits=false,
	ymax=1e-0, ymin=1e-4, xmax=2e4,	xmin=1e3,
	ytick={1e0,1e-1,1e-2,1e-3,1e-4,1e-5,1e-6},
	height=0.5\textwidth, width=0.5\textwidth, 
	grid=both, 
xlabel={$|I_N^5| = s$}, 
	ylabel={$\| \mathbf h - \mathbf h_{\text{approx}} \|_{\ell_{\infty}}/\|\mathbf h\|_{\ell_{\infty}}$},
	legend style={at={(0.5,0.85)}, anchor=south,legend columns=4,legend cell align=left, font=\small,  
	},
	xminorticks=false, yminorticks=false
	]	
	\addplot[mark options={solid},cyan,mark=diamond*,mark size=2.0] coordinates {
		(1053, 4.4028e-02)  (1863, 1.6795e-02)  (3753, 1.3580e-02)  (4563, 6.3476e-03)  (7533, 5.4000e-03) (8343, 5.9863e-03)  (12033, 1.5080e-03)  (13923, 3.8265e-03)  (16893, 3.8265e-03) 
	};
	\addlegendentry{$I_{N}^{5}$}
	\addplot[mark options={solid},red,mark=o,mark size=2.0] coordinates {
		(1053, 1.8706e-02)  (1863, 1.5066e-02)  (3753, 9.4225e-03)  (4563, 2.7314e-03)  (7533, 1.3724e-03) (8343, 1.7157e-03)  (12033, 8.1170e-04)  (13923, 3.8265e-03)  (16893, 1.5830e-03)
	};
	\addlegendentry{sparseFFT}
	\end{loglogaxis}
	\end{tikzpicture}
\begin{tikzpicture}[baseline,scale=0.6]
	\begin{axis}[
	scatter/classes = { a = {mark=o, draw=black} },
	font=\footnotesize,
	grid = both,
	xmax = 15, xmin = -15, ymax = 15, ymin = -15,
	height=0.5\textwidth, width=0.5\textwidth,
	title = {Projection of $I_{10}^{5}$},
	unit vector ratio*=1 1 1
	]
	\addplot[scatter ,only marks, mark size=1, scatter src = explicit symbolic]
	coordinates{
		(-10,-1)(-10,0)(-10,1)(-9,-1)(-9,0)(-9,1)(-8,-1)(-8,0)(-8,1)(-7,-1)(-7,0)(-7,1)(-6,-1)(-6,0)(-6,1)(-5,-2)(-5,-1)(-5,0)(-5,1)(-5,2)(-4,-2)(-4,-1)(-4,0)(-4,1)(-4,2)(-3,-3)(-3,-2)(-3,-1)(-3,0)(-3,1)(-3,2)(-3,3)(-2,-5)(-2,-4)(-2,-3)(-2,-2)(-2,-1)(-2,0)(-2,1)(-2,2)(-2,3)(-2,4)(-2,5)(-1,-10)(-1,-9)(-1,-8)(-1,-7)(-1,-6)(-1,-5)(-1,-4)(-1,-3)(-1,-2)(-1,-1)(-1,0)(-1,1)(-1,2)(-1,3)(-1,4)(-1,5)(-1,6)(-1,7)(-1,8)(-1,9)(-1,10)(0,-10)(0,-9)(0,-8)(0,-7)(0,-6)(0,-5)(0,-4)(0,-3)(0,-2)(0,-1)
		(0,0)(0,1)(0,2)(0,3)(0,4)(0,5)(0,6)(0,7)(0,8)(0,9)(0,10)(1,-10)(1,-9)(1,-8)(1,-7)(1,-6)(1,-5)(1,-4)(1,-3)(1,-2)(1,-1)(1,0)(1,1)(1,2)(1,3)(1,4)(1,5)(1,6)(1,7)(1,8)(1,9)(1,10)(2,-5)(2,-4)(2,-3)(2,-2)(2,-1)(2,0)(2,1)(2,2)(2,3)(2,4)(2,5)(3,-3)(3,-2)(3,-1)(3,0)(3,1)(3,2)(3,3)(4,-2)(4,-1)(4,0)(4,1)(4,2)(5,-2)(5,-1)(5,0)(5,1)(5,2)(6,-1)(6,0)(6,1)(7,-1)(7,0)(7,1)(8,-1)(8,0)(8,1)(9,-1)(9,0)(9,1)(10,-1)(10,0)(10,1)
	};
	\end{axis}
	\end{tikzpicture}
\begin{tikzpicture}[baseline,scale=0.6]
	\begin{axis}[
	scatter/classes = { a = {mark=o, draw=black} },
	font=\footnotesize,
	grid = both,
	xmax = 15, xmin = -15, ymax = 15, ymin = -15,
	height=0.5\textwidth, width=0.5\textwidth,
	title = {Projection of $I$ generated by sparse FFT},
	unit vector ratio*=1 1 1
	]
	\addplot[scatter ,only marks, mark size=1, scatter src = explicit symbolic]
	coordinates{
		(-10,-1) (-10,0) (-10,2) (-9,-1) (-9,0) (-9,2) (-8,1) (-7,-1) (-7,0) (-7,1) (-6,-2) (-6,-1) (-6,0) (-6,1) (-6,2) (-5,-3) (-5,-2) (-5,-1) (-5,0) (-5,1) (-5,2) (-5,3) (-4,-4) (-4,-3) (-4,-2) (-4,-1) (-4,0) (-4,1) (-4,2) (-4,3) (-4,4) (-3,-5) (-3,-4) (-3,-3) (-3,-2) (-3,-1) (-3,0) (-3,1) (-3,2) (-3,3) (-3,4) (-3,5) (-2,-6) (-2,-5) (-2,-4) (-2,-3) (-2,-2) (-2,-1) (-2,0) (-2,1) (-2,2) (-2,3) (-2,4) (-2,5) (-2,6) (-1,-7) (-1,-6) (-1,-5) (-1,-4) (-1,-3) (-1,-2) (-1,-1) (-1,0) (-1,1) (-1,2) (-1,3) (-1,4) (-1,5) (-1,6) (-1,7) (0,-10) (0,-9) (0,-8) (0,-7) (0,-6) (0,-5) (0,-4) (0,-3) (0,-2) (0,-1) (0,0) (0,1) (0,2) (0,3) (0,4) (0,5) (0,6) (0,7) (0,8) (0,9) (0,10) (1,-7) (1,-6) (1,-5) (1,-4) (1,-3) (1,-2) (1,-1) (1,0) (1,1) (1,2) (1,3) (1,4) (1,5) (1,6) (1,7) (2,-6) (2,-5) (2,-4) (2,-3) (2,-2) (2,-1) (2,0) (2,1) (2,2) (2,3) (2,4) (2,5) (2,6) (3,-5) (3,-4) (3,-3) (3,-2) (3,-1) (3,0) (3,1) (3,2) (3,3) (3,4) (3,5) (4,-4) (4,-3) (4,-2) (4,-1) (4,0) (4,1) (4,2) (4,3) (4,4) (5,-3) (5,-2) (5,-1) (5,0) (5,1) (5,2) (5,3) (6,-2) (6,-1) (6,0) (6,1) (6,2) (7,-1) (7,0) (7,1) (8,-1) (9,-2) (9,0) (9,1) (10,-2) (10,0) (10,1)
	};
	\end{axis}
	\end{tikzpicture}
	\caption{Relative discrete approximation error $\| \mathbf h - \mathbf h_{\text{approx}} \|_{\ell_\infty}/\| \mathbf h \|_{\ell_\infty}$ in dimension $d=5$ for the algebraic transformation with the hyperbolic cross $I_{N}^{d}$ with $N=2,\ldots,10$ compared to the frequency set generated by the sparse FFT algorithm (left). 
	In the center and on the right are the two-dimensional projections of $I_{10}^{5}$ and of the frequency set generated by the sparse FFT algorithm.} 
	\label{fig:Algebraic_error_compared_with_sFFT}
\end{figure}

Now we let the sparse FFT algorithm \cite[Algorithm~2]{PoVo14} determine a suitable frequency set $I$.
For dimension ${d=5}$ and for each ${N=2,3,\ldots,10}$ we choose algorithm called 'a2r1l' in \cite{Vo_FFTr1l} and use the cardinality of the hyperbolic crosses $I_{N}^{5}$ as the sparsity parameter {'sparsity\_s'}~$= s = |I_N^5|$.
As expected, the resulting discretized relative approximation errors $\|\mathbf h - \mathbf h_{\mathrm{approx}}\|_{\ell_{\infty}}/\|\mathbf h\|_{\ell_{\infty}}$ are just as good as the ones where we fixed the hyperbolic cross $I_{N}^{d}$, but the two-dimensional projections of both frequency sets to their first two coordinates differ substantially in size and shape, even though they have the same cardinality, as seen in Figure~\ref{fig:Algebraic_error_compared_with_sFFT}.

\subsubsection{Example for the error function transformation and logarithmic transformation}
The sparse FFT algorithm is especially interesting for the error function transformation $\eqref{eq:error_function_trafo}$ and the logarithmic transformation \eqref{eq:logarithmic_trafo}, because we can't calculate the transformed Fourier coefficients $\hat h_{\mathbf k}$ given in \eqref{def:FC_trafo_multivar}. 
Again we simply let the sparse FFT algorithm \cite[Algorithm~2]{PoVo14} construct a suitable frequency set $I$ depending on the sparsity $s\in\mathbb{N}$.

We return to dimension $d = 2$, use
\begin{align*}
	h(\mathbf y) = \mathrm{e}^{-y_1^2-y_2^2}
\end{align*}
as the test function and consider the constant weight function $\omega(\bm y) \equiv 1$.
We again apply two different transformations.
The two-dimensional error function transformation 
\begin{align*}
	\psi(\bm x,\bm\eta) = (\eta_1\mathrm{erf}^{-1}(2x_1), \eta_2\mathrm{erf}^{-1}(2x_2))^{\top},
\end{align*}
which we consider for $\bm\eta = \mathbf 1$, is based on its univariate version given in \eqref{eq:erf_trafo_param}.
The two-dimensional logarithmic transformation is also based on its univariate version given in \eqref{eq:logarithmic_trafo} and reads as
\begin{align*}
	\psi(\bm x, \bm\eta) := \left(\eta_1\log\left(\frac{1+2x_1}{1-2x_1}\right), \eta_2\log\left(\frac{1+2x_2}{1-2x_2}\right)\right)^{\top},
\end{align*}
which we consider only for $\bm\eta = \mathbf 1$, too.

\begin{figure}[t]
	\centering
	\begin{tikzpicture}[scale=0.75]
	\begin{axis}[
	scatter/classes = { a = {mark=o, draw=black} },
	font=\footnotesize,
	grid = both,
	xmax = 20,	xmin = -20,	ymax = 20,	ymin = -20,
title = {$\psi(\bm x,\bm 1) = \mathrm{erf}^{-1}(2x_1) \, \mathrm{erf}^{-1}(2x_2)$, $s=100$},
	unit vector ratio*=1 1 1
	]
	\addplot[scatter ,only marks, mark size=1, scatter src = explicit symbolic]
	coordinates{
		(0,0) (-1,0) (1,0) (0,-1) (0,1) (-2,0) (2,0) (0,-2) (0,2) (-1,-1) (1,1) (-1,1) (1,-1) (-3,0) (3,0) (0,-3) (0,3) (-4,0) (4,0) (0,-4) (0,4) (-1,2) (1,-2) (-1,-2) (1,2) (-2,-1) (2,1) (-2,1) (2,-1) (-5,0) (5,0) (0,-5) (0,5) (-6,0) (6,0) (0,-6) (0,6) (-1,3) (1,-3) (-1,-3) (1,3) (-3,-1) (3,1) (-3,1) (3,-1) (-7,0) (7,0) (0,-7) (0,7) (-8,0) (8,0) (-2,2) (2,-2) (-2,-2) (2,2) (0,-8) (0,8) (-1,4) (1,-4) (-1,-4) (1,4) (-4,-1) (4,1) (-4,1) (4,-1) (-9,0) (9,0) (0,-9) (0,9) (-10,0) (10,0) (-1,5) (1,-5) (-1,-5) (1,5) (-5,-1) (5,1) (0,-10) (0,10) (-5,1) (5,-1) (-11,0) (11,0) (-3,2) (3,-2) (-3,-2) (3,2) (-2,3) (2,-3) (-2,-3) (2,3) (0,-11) (0,11) (-12,0) (12,0) (-1,6) (1,-6) (-1,-6) (1,6) (-6,-1)
	};
	\end{axis}
	\end{tikzpicture}
	\begin{tikzpicture}[scale=0.75]
	\begin{axis}[
	scatter/classes = { a = {mark=o, draw=black} },
	font=\footnotesize,
	grid = both,
	xmax = 20,	xmin = -20,	ymax = 20,	ymin = -20,
title = {$\psi(\bm x, \bm 1) = \log\left(\frac{1+2x_1}{1-2x_1}\right) \log\left(\frac{1+2x_2}{1-2x_2}\right)$, $s=100$},
	unit vector ratio*=1 1 1
	]
	\addplot[scatter ,only marks, mark size=1, scatter src = explicit symbolic]
	coordinates{
		(0,0) (0,-1) (0,1) (-1,0) (1,0) (-1,1) (1,-1) (-1,-1) (1,1) (-2,0) (2,0) (0,-2) (0,2) (-2,-1) (2,1) (-1,2) (1,-2) (-1,-2) (1,2) (-2,1) (2,-1) (-3,0) (3,0) (0,-3) (0,3) (-2,-2) (2,2) (-2,2) (2,-2) (-3,-1) (3,1) (-1,-3) (1,3) (-1,3) (1,-3) (-3,1) (3,-1) (-2,3) (2,-3) (-2,-3) (2,3) (-3,-2) (3,2) (-3,2) (3,-2) (-4,0) (4,0) (0,-4) (0,4) (-4,-1) (4,1) (-1,-4) (1,4) (-1,4) (1,-4) (-4,1) (4,-1) (0,-6) (0,6) (-6,0) (6,0) (-3,3) (3,-3) (-3,-3) (3,3) (0,-5) (0,5) (-5,0) (5,0) (-4,-2) (4,2) (-2,4) (2,-4) (-2,-4) (2,4) (-4,2) (4,-2) (0,-7) (0,7) (-7,0) (7,0) (-6,1) (6,-1) (-1,-6) (1,6) (-1,6) (1,-6) (-6,-1) (6,1) (-5,1) (5,-1) (-1,5) (1,-5) (-1,-5) (1,5) (-5,-1) (5,1) (-7,1) (7,-1) (-1,-7)   		
	};
	\end{axis}
	\end{tikzpicture}
	
	\begin{tikzpicture}[scale=0.75]
	\begin{axis}[
	scatter/classes = { a = {mark=o, draw=black} },
	font=\footnotesize,
	grid = both,
	xmax = 20,	xmin = -20,	ymax = 20,	ymin = -20,
title = {$\psi(\bm x,\bm 1) = \mathrm{erf}^{-1}(2x_1) \, \mathrm{erf}^{-1}(2x_2)$, $s=500$},
	unit vector ratio*=1 1 1
	]
	\addplot[scatter ,only marks, mark size=1, scatter src = explicit symbolic]
	coordinates{
		(0,0) (-1,0) (1,0) (0,-1) (0,1) (-2,0) (2,0) (0,-2) (0,2) (-1,-1) (1,1) (-1,1) (1,-1) (-3,0) (3,0) (0,-3) (0,3) (-4,0) (4,0) (0,-4) (0,4) (-1,2) (1,-2) (-1,-2) (1,2) (-2,-1) (2,1) (-2,1) (2,-1) (-5,0) (5,0) (0,-5) (0,5) (-6,0) (6,0) (0,-6) (0,6) (-1,3) (1,-3) (-1,-3) (1,3) (-3,-1) (3,1) (-3,1) (3,-1) (-7,0) (7,0) (0,-7) (0,7) (-8,0) (8,0) (-2,2) (2,-2) (-2,-2) (2,2) (0,-8) (0,8) (-1,4) (1,-4) (-1,-4) (1,4) (-4,-1) (4,1) (-4,1) (4,-1) (-9,0) (9,0) (0,-9) (0,9) (-10,0) (10,0) (-1,5) (1,-5) (-1,-5) (1,5) (-5,-1) (5,1) (0,-10) (0,10) (-5,1) (5,-1) (-11,0) (11,0) (-3,2) (3,-2) (-3,-2) (3,2) (-2,3) (2,-3) (-2,-3) (2,3) (0,-11) (0,11) (-12,0) (12,0) (-1,6) (1,-6) (-1,-6) (1,6) (-6,-1) (6,1) (-6,1) (6,-1) (-13,0) (13,0) (0,-12) (0,12) (-14,0) (14,0) (-1,7) (1,-7) (-1,-7) (1,7) (-7,-1) (7,1) (0,-13) (0,13) (-4,2) (4,-2) (-4,-2) (4,2) (-2,4) (2,-4) (-2,-4) (2,4) (-7,1) (7,-1) (-15,0) (15,0) (0,-14) (0,14) (-1,8) (1,-8) (-1,-8) (1,8) (-8,-1) (8,1) (-16,0) (16,0) (-3,3) (3,-3) (-3,-3) (3,3) (-17,0) (17,0) (-8,1) (8,-1) (0,-15) (0,15) (-5,2) (5,-2) (-5,-2) (5,2) (-1,9) (1,-9) (-1,-9) (1,9) (-9,-1) (9,1) (-2,5) (2,-5) (-2,-5) (2,5) (-18,0) (18,0) (-10,-1) (10,1) (-19,0) (19,0) (0,-16) (0,16) (-1,10) (1,-10) (-1,-10) (1,10) (-9,1) (9,-1) (-6,2) (6,-2) (-4,3) (4,-3) (-4,-3) (4,3) (-6,-2) (6,2) (-3,4) (3,-4) (-3,-4) (3,4) (-20,0) (20,0) (-2,6) (2,-6) (-2,-6) (2,6) (-11,-1) (11,1) (-1,11) (1,-11) (-1,-11) (1,11) (0,-17) (0,17) (-20,1) (20,-1) (-7,2) (7,-2) (-12,-1) (12,1) (-7,-2) (7,2) (-10,1) (10,-1) (-1,12) (1,-12) (-1,-12) (1,12) (-2,7) (2,-7) (-2,-7) (2,7) (-5,3) (5,-3) (-5,-3) (5,3) (-19,1) (19,-1) (-3,5) (3,-5) (-13,-1) (13,1) (-3,-5) (3,5) (-8,2) (8,-2) (-4,4) (4,-4) (-4,-4) (4,4) (-8,-2) (8,2) (-1,13) (1,-13) (-1,-13) (1,13) (0,-18) (0,18) (-14,-1) (14,1) (-18,1) (18,-1) (-2,8) (2,-8) (-2,-8) (2,8) (-9,2) (9,-2) (-6,3) (6,-3) (-15,-1) (15,1) (-11,1) (11,-1) (-6,-3) (6,3) (-9,-2) (9,2) (-1,14) (1,-14) (-1,-14) (1,14) (-3,6) (3,-6) (-3,-6) (3,6) (-10,2) (10,-2) (-16,-1) (16,1) (-2,9) (2,-9) (-2,-9) (2,9) (-17,1) (17,-1) (-5,4) (5,-4) (-5,-4) (5,4) (-10,-2) (10,2) (-4,5) (4,-5) (-4,-5) (4,5) (-17,-1) (17,1) (-7,3) (7,-3) (-11,2) (11,-2) (-7,-3) (7,3) (-1,15) (1,-15) (-1,-15) (1,15) (-18,-1) (18,1) (-12,2) (12,-2) (-3,7) (3,-7) (-3,-7) (3,7) (-11,-2) (11,2) (0,-19) (0,19) (-2,10) (2,-10) (-2,-10) (2,10) (-19,-1) (19,1) (-13,2) (13,-2) (-8,3) (8,-3) (-20,-1) (20,1) (-8,-3) (8,3) (-6,4) (6,-4) (-12,1) (12,-1) (-14,2) (14,-2) (-12,-2) (12,2) (-6,-4) (6,4) (-20,2) (20,-2) (-1,16) (1,-16) (-1,-16) (1,16) (-16,1) (16,-1) (-15,2) (15,-2) (-4,6) (4,-6) (-19,2) (19,-2) (-4,-6) (4,6) (-16,2) (16,-2) (-18,2) (18,-2) (-17,2) (17,-2) (-5,5) (5,-5) (-5,-5) (5,5) (-2,11) (2,-11) (-2,-11) (2,11) (-9,3) (9,-3) (-3,8) (3,-8) (-3,-8) (3,8) (-13,-2) (13,2) (-9,-3) (9,3) (-7,4) (7,-4) (-14,-2) (14,2) (-7,-4) (7,4) (-10,3) (10,-3) (-1,17) (1,-17) (-1,-17) (1,17) (-2,12) (2,-12) (-2,-12) (2,12) (-4,7) (4,-7) (-15,-2) (15,2) (-10,-3) (10,3) (-4,-7) (4,7) (-3,9) (3,-9) (-3,-9) (3,9) (-11,3) (11,-3) (-6,5) (6,-5) (-6,-5) (6,5) (-5,6) (5,-6) (-5,-6) (5,6) (-16,-2) (16,2) (-8,4) (8,-4) (-8,-4) (8,4) (-12,3) (12,-3) (-11,-3) (11,3) (-2,13) (2,-13) (-2,-13) (2,13) (-17,-2) (17,2) (-13,3) (13,-3) (-3,10) (3,-10) (-3,-10) (3,10) (-18,-2) (18,2) (-4,8) (4,-8) (-4,-8) (4,8) (-9,4) (9,-4) (-15,1) (15,-1) (-12,-3) (12,3) (-7,5) (7,-5) (-19,-2) (19,2) (-14,3) (14,-3) (-7,-5) (7,5) (-1,18) (1,-18) (-1,-18) (1,18) (-9,-4) (9,4) (-20,-2) (20,2) (-15,3) (15,-3) (-20,3) (20,-3) (-16,3) (16,-3) (-5,7) (5,-7) (-2,14) (2,-14) (-2,-14) (2,14) (-6,6) (6,-6) (-5,-7) (5,7) (-19,3) (19,-3) (-17,3) (17,-3) (-6,-6) (6,6) (-18,3) (18,-3) (-13,-3) (13,3) (-10,4) (10,-4) (-3,11) 	
	};
	\end{axis}
	\end{tikzpicture}
	\begin{tikzpicture}[scale=0.75]
	\begin{axis}[
	scatter/classes = { a = {mark=o, draw=black} },
	font=\footnotesize,
	grid = both,
	xmax = 20,	xmin = -20,	ymax = 20,	ymin = -20,
title = {$\psi(\bm x, \bm 1) = \log\left(\frac{1+2x_1}{1-2x_1}\right) \log\left(\frac{1+2x_2}{1-2x_2}\right)$, $s=500$},
	unit vector ratio*=1 1 1
	]
	\addplot[scatter ,only marks, mark size=1, scatter src = explicit symbolic]
	coordinates{
		(0,0) (0,-1) (0,1) (-1,0) (1,0) (-1,1) (1,-1) (-1,-1) (1,1) (-2,0) (2,0) (0,-2) (0,2) (-2,-1) (2,1) (-1,2) (1,-2) (-1,-2) (1,2) (-2,1) (2,-1) (-3,0) (3,0) (0,-3) (0,3) (-2,-2) (2,2) (-2,2) (2,-2) (-3,-1) (3,1) (-1,-3) (1,3) (-1,3) (1,-3) (-3,1) (3,-1) (-2,3) (2,-3) (-2,-3) (2,3) (-3,-2) (3,2) (-3,2) (3,-2) (-4,0) (4,0) (0,-4) (0,4) (-4,-1) (4,1) (-1,-4) (1,4) (-1,4) (1,-4) (-4,1) (4,-1) (0,-6) (0,6) (-6,0) (6,0) (-3,3) (3,-3) (-3,-3) (3,3) (0,-5) (0,5) (-5,0) (5,0) (-4,-2) (4,2) (-2,4) (2,-4) (-2,-4) (2,4) (-4,2) (4,-2) (0,-7) (0,7) (-7,0) (7,0) (-6,1) (6,-1) (-1,-6) (1,6) (-1,6) (1,-6) (-6,-1) (6,1) (-5,1) (5,-1) (-1,5) (1,-5) (-1,-5) (1,5) (-5,-1) (5,1) (-7,1) (7,-1) (-1,-7) (1,7) (-1,7) (1,-7) (-7,-1) (7,1) (0,-8) (0,8) (-8,0) (8,0) (-8,1) (8,-1) (-1,-8) (1,8) (-1,8) (1,-8) (-8,-1) (8,1) (0,-9) (0,9) (-9,0) (9,0) (-9,1) (9,-1) (-1,-9) (1,9) (-1,9) (1,-9) (-9,-1) (9,1) (-6,2) (6,-2) (-2,6) (2,-6) (-2,-6) (2,6) (-6,-2) (6,2) (0,-10) (0,10) (-10,0) (10,0) (-4,3) (4,-3) (-4,-3) (4,3) (-3,4) (3,-4) (-3,-4) (3,4) (-5,2) (5,-2) (-2,-5) (2,5) (-2,5) (2,-5) (-5,-2) (5,2) (-7,2) (7,-2) (-2,7) (2,-7) (-2,-7) (2,7) (-7,-2) (7,2) (-10,1) (10,-1) (-1,-10) (1,10) (-1,10) (1,-10) (-10,-1) (10,1) (0,-11) (0,11) (-11,0) (11,0) (-8,2) (8,-2) (-2,8) (2,-8) (-2,-8) (2,8) (-8,-2) (8,2) (-11,1) (11,-1) (-1,-11) (1,11) (-1,11) (1,-11) (-11,-1) (11,1) (0,-12) (0,12) (-12,0) (12,0) (-9,2) (9,-2) (-2,9) (2,-9) (-2,-9) (2,9) (-9,-2) (9,2) (-12,1) (12,-1) (-1,-12) (1,12) (-1,12) (1,-12) (-12,-1) (12,1) (-3,6) (3,-6) (-3,-6) (3,6) (-6,-3) (6,3) (-6,3) (6,-3) (0,-13) (0,13) (-13,0) (13,0) (-10,2) (10,-2) (-2,10) (2,-10) (-2,-10) (2,10) (-10,-2) (10,2) (-3,-5) (3,5) (-3,5) (3,-5) (-5,-3) (5,3) (-5,3) (5,-3) (-3,7) (3,-7) (-3,-7) (3,7) (-7,-3) (7,3) (-7,3) (7,-3) (-13,1) (13,-1) (-1,-13) (1,13) (-1,13) (1,-13) (-13,-1) (13,1) (-4,4) (4,-4) (-4,-4) (4,4) (0,-14) (0,14) (-14,0) (14,0) (-11,2) (11,-2) (-2,11) (2,-11) (-2,-11) (2,11) (-11,-2) (11,2) (-3,8) (3,-8) (-3,-8) (3,8) (-8,-3) (8,3) (-8,3) (8,-3) (-14,1) (14,-1) (-1,-14) (1,14) (-1,14) (1,-14) (-14,-1) (14,1) (0,-15) (0,15) (-15,0) (15,0) (-12,2) (12,-2) (-2,12) (2,-12) (-2,-12) (2,12) (-12,-2) (12,2) (-3,9) (3,-9) (-3,-9) (3,9) (-9,-3) (9,3) (-9,3) (9,-3) (-15,1) (15,-1) (-1,-15) (1,15) (-1,15) (1,-15) (-15,-1) (15,1) (0,-16) (0,16) (-16,0) (16,0) (-13,2) (13,-2) (-2,13) (2,-13) (-2,-13) (2,13) (-13,-2) (13,2) (-3,10) (3,-10) (-3,-10) (3,10) (-10,-3) (10,3) (-10,3) (10,-3) (-4,6) (4,-6) (-4,-6) (4,6) (-6,-4) (6,4) (-6,4) (6,-4) (-16,1) (16,-1) (-1,-16) (1,16) (-1,16) (1,-16) (-16,-1) (16,1) (-4,-5) (4,5) (-4,5) (4,-5) (-5,-4) (5,4) (-5,4) (5,-4) (-4,7) (4,-7) (-4,-7) (4,7) (-7,-4) (7,4) (-7,4) (7,-4) (0,-17) (0,17) (-14,2) (14,-2) (-17,0) (17,0) (-2,14) (2,-14) (-2,-14) (2,14) (-14,-2) (14,2) (-3,11) (3,-11) (-3,-11) (3,11) (-11,-3) (11,3) (-11,3) (11,-3) (-17,1) (17,-1) (-1,-17) (1,17) (-1,17) (1,-17) (-4,8) (4,-8) (-4,-8) (4,8) (-8,-4) (8,4) (-8,4) (8,-4) (-17,-1) (17,1) (-15,2) (15,-2) (-2,15) (2,-15) (-2,-15) (2,15) (-15,-2) (15,2) (0,-18) (0,18) (-3,12) (3,-12) (-3,-12) (3,12) (-12,-3) (12,3) (-12,3) (12,-3) (-18,0) (18,0) (-18,1) (18,-1) (-4,9) (4,-9) (-4,-9) (4,9) (-9,-4) (9,4) (-9,4) (9,-4) (-16,2) (16,-2) (-1,-18) (1,18) (-1,18) (1,-18) (-6,-6) (6,6) (-6,6) (6,-6) (-2,16) (2,-16) (-2,-16) (2,16) (-3,13) (3,-13) (-3,-13) (3,13) (-13,-3) (13,3) (-13,3) (13,-3) (-16,-2) (16,2) (-18,-1) (18,1) (-6,5) (6,-5) (-6,-5) (6,5) (-5,-6) (5,6) (-5,6) (5,-6) (0,-19) (0,19) (-6,-7) (6,7) (-6,7) (6,-7) (-7,-6) (7,6) (-7,6) (7,-6) (-17,2) (17,-2) (-19,1) (19,-1) (-4,10) (4,-10) (-4,-10) (4,10) (-10,-4) (10,4) (-10,4) (10,-4) (-19,0) (19,0) (-5,5) (5,-5) (-5,-5) (5,5) (-7,5) (7,-5) (-5,-7)
	};
	\end{axis}
	\end{tikzpicture}
	\caption{Two-dimensional frequency sets $I_{N}$ with $N=20$ and $s\in\{100,500\}$ for the error function transformation (left column) and the logarithmic transformation (right column).}
	\label{fig:spare index sets log erf}
\end{figure}
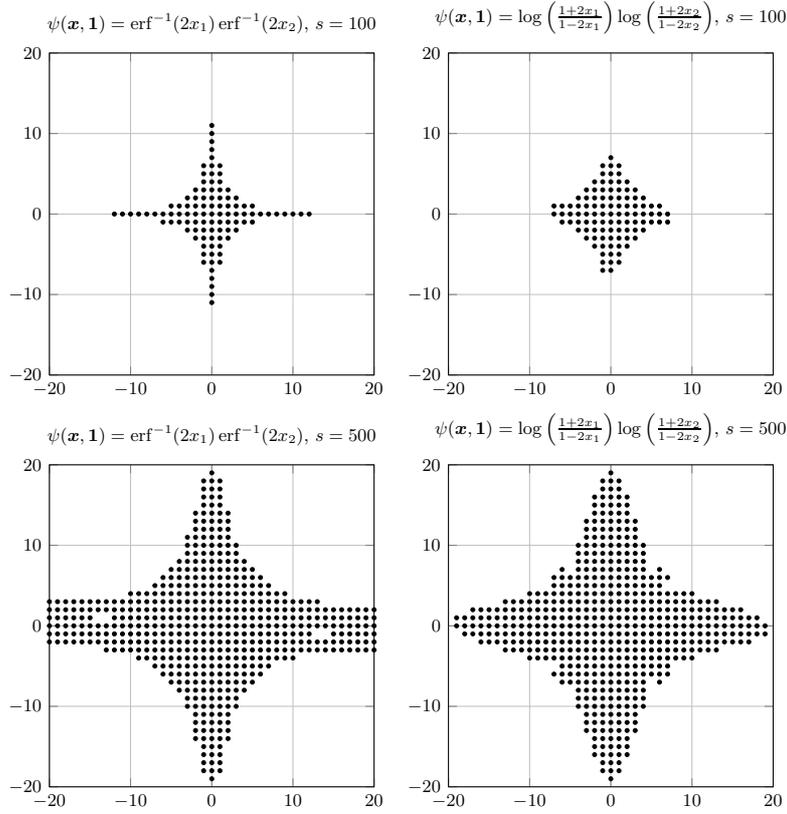

At first we fix the refinement $N=20$. Then the full $41\times 41$-integer grid contains ${(2\cdot 20 +1)^2 = 1681}$ elements. 
Again, we initialize the algorithm 'a2r1l' in \cite{Vo_FFTr1l} with the default threshold parameter 'threshold\_theta' of 1e-12 and denote the sparsity parameter 'sparsity\_s' as $s\in\mathbb{N}$.
For the sparsities $s = 100$ and $s = 500$ the error function transformation leads to a frequency set $I_N$ that reminds us of a hyperbolic cross, whereas the logarithmic transformation leads to a frequency set that could also resemble an appropriately scaled unit ball $\left\{ \mathbf x\in\mathbb{Z}^2 : \left(|x_{1}|^p + |x_{2}|^p\right)^{\frac{1}{p}} \leq N \right\}$
of a two-dimensional sequence space $\ell_p$ with $0<p<1$, see Figure~\ref{fig:spare index sets log erf}.

\begin{figure}[t]
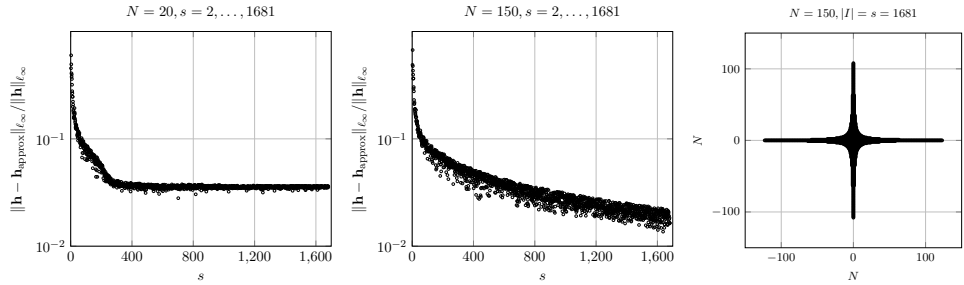

\centering

	\caption{Relative approximation errors $\|\mathbf h - \mathbf h_{\mathrm{approx}}\|_{\ell_{\infty}}/\|\mathbf h\|_{\ell_{\infty}}$ calculated by the spare FFT algorithm for the error function transformation with $\bm \eta = \mathbf 1, s=2,\ldots,1681$ and refinement $N=20$ on the left and $N=150$ in the center. The automatically constructed frequency set $I$ for $N=150,s=1681$ is shown on the right.}
	\label{fig:spare_index_sets_erf}
\end{figure}

Finally, we focus on the two-dimensional error function transformation with $\bm \eta = \mathbf 1$ and compare the corresponding relative approximation errors $\|\mathbf h - \mathbf h_{\mathrm{approx}}\|_{\ell_{\infty}}/\|\mathbf h\|_{\ell_{\infty}}$ calculated by the spare FFT algorithm in two different setups.
At first we keep the refinement $N=20$ and consider increasing sparsity parameters $s=2,\ldots,1681$. Hence, for small sparsities we have frequency sets that look like hyperbolic crosses as shown in the left column of Figure~\ref{fig:spare index sets log erf}, whose branches along the central axes become thicker as $s$ increases and eventually we end up with the full $41\times 41$-grid.
Based on these frequency sets, the relative approximation errors stagnate at a certain point, shown on the left in Figure~\ref{fig:spare_index_sets_erf}, because the relatively small refinement value forces the algorithm to consider frequencies within the $41 \times 41$ grid, that don't significantly improve the approximation of $h$.
In comparison we raise the refinement to $N=150$ and let the sparsity parameter $s$ run from $2$ to $1681$ again, so that the resulting frequency sets have the same cardinality as before but keep their hyperbolic cross like shape, which is shown on the right of Figure~\ref{fig:spare_index_sets_erf}, where we have the frequency set constructed by the sparse FFT algorithm for $N=150$ and $s=1681$.
With these frequency sets we now have steadily decreasing relative approximation errors for increasing sparsity values, as plotted in the center of Figure~\ref{fig:spare_index_sets_erf}.

\subsubsection{Example with the tangens transformation}
Finally, we consider a different algebraic test function that is given in product form by
\begin{align*}
	h(\mathbf y) = \prod_{j=1}^{d} \frac{1}{1+y_j^2}.
\end{align*}
Additionally we consider the constant weight function $\omega(\bm y) \equiv \bm 1$ and the multivariate tangens transformation 
\begin{align*}
	\psi(\bm x, \bm \eta) := \left( \eta_1 \tan(\pi x_1), \ldots, \eta_d \tan(\pi x_d) \right)^{\top}
\end{align*}
with $\bm \eta = \mathbf 1$ based on the univariate version defined in \eqref{eq:tan_trafo}. 
The resulting transformed function is of the form
\begin{align}
	\label{eq:trafo_fct_tan}
	f(\bm x, \bm 1, \bm \mu) := \prod_{j=1}^{d} \frac{1}{1 + \tan(\pi x_j)^2}.
\end{align}
This product form extends to the corresponding Fourier coefficients, i.e., $\hat h_{\mathbf k} = \prod_{j=1}^{d} \hat h_{k_j}$
and the one-dimensional Fourier coefficients $\hat h_{k_j}$ are of the form
\begin{align} \label{eq:exact_tan_FC}
	\hat h_{k_j}
	= \int_{-1/2}^{1/2} \frac{\mathrm{e}^{-2 \pi \mathrm{i} k_j x_j}}{1 + \tan(\pi x_j)^2} \mathrm d x_j
	= \int_{-1/2}^{1/2} \cos(\pi x_j)^2 \,\mathrm{e}^{-2 \pi \mathrm{i} k_j x_j} \mathrm d x_j 
	=
	\begin{cases}
	\frac{1}{2} & \text{for } k_j = 0, \\
	\frac{1}{4} & \text{for } |k_j| = 1,\\
	0 & \text{otherwise}.
	\end{cases}
\end{align}
Hence, over a full grid  $[-N,N]^d\cap\mathbb{Z}^d$ with $(2 N+1)^d$ points there are just $3^d$ non-zero multivariate Fourier coefficients $\hat h_{\mathbf k}$.
Again we check this for dimension $d = 12$ and $N = 4$ with the sparse FFT algorithm, see \cite[Algorithm~2]{PoVo14} and \cite{Vo_FFTr1l}.
We initialize this algorithm with the test function $f$ in \eqref{eq:trafo_fct_tan}, choose the algorithm name 'a2r1l', set the sparsity parameter 'sparsity\_s' to $10^6$ and the threshold parameter 'threshold\_theta' to 1e-12.
This results in an exact reconstruction as the algorithm indeed only detected the $3^{12} = 531441$ out of $(2\cdot 4 + 1)^{12} \approx 2.8\cdot 10^{11}$ possible frequencies, corresponding to the $12$-dimensional integer unit cube of radius $1$ for which the transformed Fourier coefficients $\hat h_{k_j}$ are non-zero, as calculated in \eqref{eq:exact_tan_FC}.

\section{Conclusion}
In this paper we considered functions $h\in L_{2}(\mathbb{R}^d,\omega(\circ,\bm\mu)) \cap H_{\mathrm{mix}}^{m}(\mathbb{R}^d)$ with a parameterized weight function $\omega(\circ,\bm\mu):\mathbb{R}^d\to[0,\infty), \bm\mu\in\mathbb{R}^d$ and discussed strategies for transforming them into functions $f$ on the torus $\mathbb T^d$.
A parameterized transformation $\psi(\circ,\bm\eta): (-\frac{1}{2},\frac{1}{2})^d\to\mathbb{R}^d$ with $\bm\eta\in\mathbb{R}^d$ in combination with the weight function $\omega(\circ,\bm\mu)$ let us control the degree of smoothness $m\in\mathbb{N}$ of a function $h$ defined on $\mathbb{R}^d$, that is preserved under the change of variables $\psi(\circ,\bm\eta)$. 
Hence, the parameters $\bm\eta$ and $\bm\mu$ control which Sobolev space $H_{\mathrm{mix}}^{m}(\mathbb{T}^d)$ the transformed functions ${f(\circ,\bm\eta,\bm\mu) = h(\psi(\circ,\bm\eta))\,\sqrt{ \omega(\psi(\circ,\bm\eta),\bm\mu)\,\psi'(\circ,\bm\eta) }}$ belong to.
Due to the embedding of the Sobolev space $H_{\mathrm{mix}}^{m}(\mathbb{T}^d)$ into the Wiener algebra $\mathcal{A}(\mathbb{T}^d)$ of functions with absolutely summable Fourier coefficients, we have information on the rate of decay of the Fourier coefficients $\hat f_{\mathbf k}$ and $\hat h_{\mathbf k}$ without having to calculate them -- which in a lot of cases is not possible in the first place.
Thus, the essential theoretical $L_2$- and $L_{\infty}$-approximation error bounds on the torus $\mathbb{T}^d$, proposed in \cite[Theorem~2.30]{volkmerdiss} and \cite[Theorem~3.3]{KaPoVo13}, can be transferred to $\mathbb{R}^d$ by means of the inverse transformation ${\psi^{-1}(\circ,\bm\eta): \mathbb{R}^d \to (-\frac{1}{2},\frac{1}{2})^d}$.
Furthermore, only slight modifications are necessary to incorporate such transformations into the algorithms based on single reconstructing {rank-$1$} lattices for the evaluation and the reconstruction of transformed multivariate trigonometric polynomials presented in \cite[Algorithm~3.1 and 3.2]{kaemmererdiss}.
Algorithms based on multiple reconstructing {rank-$1$} lattices \cite{Kae16} and sparse fast Fourier transformations \cite{PoVo14} can be adjusted, too.

Our numerical tests show that these algorithms are still working within the proposed upper bounds for the approximation error. 
Additionally, special cases in which we can actually calculate the Fourier coefficients confirm the theoretical parameter bounds for $\bm\mu$ and $\bm\eta$, that are sufficient to achieve a certain degree of Sobolev-smoothness under a change of variables.
In several examples we apply the adapted multiple {rank-$1$} lattice methods and adjusted dimension incremental construction methods for sparse frequency sets.

\section*{Acknowledgements}
The authors thank Lutz K\"ammerer, Tino Ullrich and Toni Volkmer for fruitful discussions during the preparation of this paper.
Furthermore we thank Lutz K\"ammerer  for the multiple lattice software,
see \cite{Ka_LFFT, Kae17}, and Toni Volkmer for the software 'sparseFFTr1l' for 
computing the sparse fast Fourier transform based on reconstructing {rank-$1$} lattices in a dimension incremental way, see \cite{Vo_FFTr1l}.
The first named author gratefully acknowledges the support by the funding of the European Union and the
Free State of Saxony (ESF).

\bibliographystyle{abbrv}

\end{document}